\input amstex
\documentstyle{amsppt}
\document

\magnification 1000

\def\gen{\frak{g}}
\def\slen{\frak{s}\frak{l}}

\def\ten{\frak{t}}
\def\aen{\frak{a}}
\def\ben{\frak{b}}

\def\ken{\frak{k}}
\def\een{\frak{e}}

\def\hen{\frak{h}}

\def\len{\frak{l}}

\def\nen{\frak{n}}

\def\qen{\frak{q}}
\def\cen{\frak{c}}

\def\sen{\frak{s}}
\def\uen{\frak{u}}

\def\zen{\frak{z}}

\def\Aen{\frak{A}}

\def\Cen{\frak{C}}
\def\Den{\frak{D}}

\def\a{{\alpha}}
\def\g{{\gamma}}
\def\o{{\omega}}
\def\l{{\lambda}}
\def\b{{\beta}}
\def\eps{{\varepsilon}}

\def\1b{{\bold 1}}

\def\Ab{{\bold A}}

\def\Gb{{\bold G}}

\def\Kb{{\bold K}}

\def\Hb{{\bold H}}

\def\Pb{{\bold P}}

\def\Rb{{\bold R}}

\def\Sb{{\bold S}}

\def\Lb{{\bold L}}

\def\Extb{{\bold{Ext}}}

\def\rhov{{\check \rho}}
\def\thetav{{\check\theta}}
\def\av{{\check a}}

\def\ov{{\check o}}
\def\alphav{{\check\alpha}}
\def\omegav{{\check\omega}}
\def\lambdav{{\check\lambda}}
\def\muv{{\check\mu}}
\def\nuv{{\check\nu}}

\def\ad{{\roman{ad}}}

\def\Spec{{\roman{Spec\,}}}

\def\sgnb{{\bold{sgn}}}
\def\trivb{{\bold{triv}}}

\def\Aut{\text{Aut}}

\def\Ext{\text{Ext}}

\def\Hom{{\roman{Hom}}}
\def\Gal{{\roman{Gal}}}
\def\cox{{\text{\sl cox}}}

\def\Res{\text{Res}}

\def\dim{{\roman{dim}}}

\def\Ind{\text{Ind}}

\def\End{{\roman{End}}}

\def\Id{{\roman{id}}}
\def\Irr{{\roman{Irr}}}

\def\tr{\roman{tr}}

\def\Ker{\roman{Ker}\,}

\def\Gr{{\roman{Gr}}}

\def\Lie{\roman{Lie}}

\def\SL{{\roman{SL}}}

\def\re{{\roman{re}}}
\def\EN{{{\text{EN}}}}
\def\ERS{{{\text{ERS}}}}
\def\HRS{{{\text{HRS}}}}
\def\RN{{{\text{RN}}}}
\def\reg{{{\text{RS}}}}
\def\re{{{\text{re}}}}

\def\AA{{\Bbb A}}

\def\CC{{\Bbb C}}

\def\FF{{\Bbb F}}

\def\GG{{\Bbb G}}

\def\PP{{\Bbb P}}
\def\QQ{{\Bbb Q}}
\def\RR{{\Bbb R}}

\def\ZZ{{\Bbb Z}}

\def\Ac{{\Cal A}}
\def\Bc{{\Cal B}}
\def\Cc{{\Cal C}}

\def\Ec{{\Cal E}}
\def\Fc{{\Cal F}}

\def\Hc{{\Cal H}}
\def\HHc{{\Hc\!\!\!\Hc}}
\def\Ic{{\Cal I}}

\def\Lc{{\Cal L}}
\def\Mc{{\Cal M}}
\def\Nc{{\Cal N}}
\def\Oc{{\Cal O}}
\def\Pc{{\Cal P}}

\def\Sc{{\Cal S}}

\def\Uc{{\Cal U}}
\def\Vc{{\Cal V}}
\def\Wc{{\Cal W}}
\def\Xc{{\Cal X}}
\def\Yc{{\Cal Y}}

\def\mod{{\text{mod}}}

\def\and{{\text{and}}}

\def\adm{{\roman{adm}}}
\def\rat{{\roman{rat}}}

\def\ss{\scriptstyle}

\def\qed{\hfill $\sqcap \hskip-6.5pt \sqcup$}        % White box
\overfullrule=0pt                                    % No black boxes
\def\newpage{{\vfill\break}}

\def\7dag{{{\!\!\!\!\!\!\!\dag}}}
\def\6dag{{{\!\!\!\!\!\!\dag}}}
\def\5dag{{{\!\!\!\!\!\dag}}}
\def\4dag{{{\!\!\!\!\dag}}}
\def\3dag{{{\!\!\!\dag}}}
\def\2dag{{{\!\!\dag}}}
\def\1dag{{{\!\dag}}}

\def\pro{{\lim\limits_{\lla}}}

\def\la{{\langle}}
\def\ra{{\rangle}}
\def\lla{{\longleftarrow}}
\def\lra{{{\longrightarrow}}}

\newdimen\Squaresize\Squaresize=14pt
\newdimen\Thickness\Thickness=0.5pt
\def\Square#1{\hbox{\vrule width\Thickness
          \alphaox to \Squaresize{\hrule height \Thickness\vss
          \hbox to \Squaresize{\hss#1\hss}
          \vss\hrule height\Thickness}
          \unskip\vrule width \Thickness}
          \kern-\Thickness}
\def\Vsquare#1{\alphaox{\Square{$#1$}}\kern-\Thickness}

\nologo

%%%%%%%%%%%%%%%%%%%%%%%%%%%%%%%%%%%%%%%%%%%%%%%%%%%%%%%%%%%%%%%%%%%%%%%%%%%%%%%%
\topmatter
\title Finite dimensional representations
of DAHA and affine Springer fibers :
the spherical case
\endtitle
\rightheadtext{DAHA and affine Springer fibers}
\abstract
We classify finite dimensional simple spherical representations of
rational double affine Hecke algebras,
and we study a remarkable family of
finite dimensional simple spherical representations of
double affine Hecke algebras.
\endabstract
\author M. Varagnolo, E. Vasserot\endauthor
\address D\'epartement de Math\'ematiques,
Universit\'e de Cergy-Pontoise, 2 av. A. Chauvin,
BP 222, 95302 Cergy-Pontoise Cedex, France,
Fax : 01 34 25 66 45\endaddress
\email michela.varagnolo\@math.u-cergy.fr\endemail
\address D\'epartement de Math\'ematiques,
Universit\'e Paris 7,
175 rue du Chevaleret,
75013 Paris, France,
Fax : 01 44 27 78 18
\endaddress
\email vasserot\@math.jussieu.fr\endemail
\thanks
2000{\it Mathematics Subject Classification.}
Primary ??; Secondary ??.
\endthanks
\endtopmatter
\document

\head Introduction \endhead

Double affine Hecke algebras, DAHA for short, have been introduced by Cherednik
about 15 years ago to prove MacDonald conjectures.
The understanding of their representation theory has progressed 
very much recently, in particular by the classification of the simple
modules in the category $O$ in \cite{V1} 
(when the parameters are not roots of unity).
The latter is very similar to Kazhdan-Lusztig classification of simple modules
of affine Hecke algebras. One can show that any simple module
in the category $O$  is the top of a module induced from an 
affine Hecke subalgebra. See A.3.6 below.
However, the representation theory of DAHA has some specific feature
which has no analogue for affine Hecke algebras.
For instance, it is very difficult to classify the finite dimensional
simple modules.

This can be approached in several ways.
The DAHA, denoted by $\Hb$, admits two remarkable degenerated forms.
The first one, the degenerated DAHA, denoted by $\Hb'$, is an analogue 
of the degenerate Hecke algebras introduced by Drinfeld and Lusztig.
Its representation theory is more or less the same as that of $\Hb$.
See \cite{VV}, \cite{L1}.
The second one has been introduced by Etingof and Ginzburg in \cite{EG} and
is called the rational DAHA (or rational Cherednik algebra).
We'll denote it by $\Hb''$.

It has been believed, see \cite{BEG} for instance,
that finite dimensional representations of $\Hb''$ and $\Hb$ are the same.
A little attention shows that indeed
the category of finite dimensional $\Hb''$-modules
embeds strictly into the category of finite dimensional $\Hb$-modules.
While the parametrization of isomorphism classes of spherical
finite dimensional
simple $\Hb''$-modules involves only one rational number, see
2.8.2(c), the parametrization
of isomorphism classes of spherical finite dimensional simple $\Hb$-modules
involves also some irreducible local system, coming from
the Langlands parameters in \cite{V1}. 
It is not difficult to find examples of
finite dimensional $\Hb$-modules which are not
$\Hb''$-modules, see 2.3.5 below.
After this paper was writen we have had a discussion with P. Etingof
who has been able to prove how to recover the missing representations. 
See \cite{E} for details.

In this paper we concentrate on the spherical finite dimensional modules.
The case of non spherical modules can probably be done with similar techniques.
We'll come back to this elsewhere.
The paper contains two main results.

First we classify all spherical finite dimensional simple
$\Hb''$-modules in Theorem 2.8.1.
Since the finite dimensional simple
$\Hb''$-modules belong to the category O, each of them is 
the top of a standard module.
The spherical one are the top of a polynomial representation 
(= a standard module induced from the trivial representation
of the Weyl group).
So they are labelled by the value of the parameter of
$\Hb''$, which is a rational number $c=k/m$ with $(k,m)=1$ and $m>0$.
Surprisingly, the classification we get is extremely simple and nice.
The spherical finite dimensional simple modules correspond to the integers
$k,m$ such that $k<0$ and $m$ is an ellipic number, i.e.,
the integer $m$ is the order of an elliptic element of the Weyl group
which is regular in Springer's sense.
For instance, in type $E_8$ there are 12 elliptic numbers.
The only known cases before were the case where $m$ is the
Coxeter number in arbitrary type,
and the dihedral types (in particular all rank 2 types).
Notice that in this paper we assume 
that $\Hb''$ is crystallographic with equal parameters.
The proof is as follows.
Any simple spherical finite dimensional $\Hb''$-module 
$M''$ has also the structure
of a simple spherical $\Hb$-module, denoted by $M$.
The algebra $\Hb''$ has two remarkable polynomial subalgebras,
yielding, under induction, 
two representations.
They are called the polynomial representations.
A spherical finite dimensional $\Hb''$-module
is a quotient of both polynomial representations.
Using this, one can identify $M$ with the top of 
a standard $\Hb$-module with explicit Langlands parameters.
See \cite{V1} for the terminology.
A case by case analysis using the Fourier-Sato transform of perverse sheaves 
shows that this explicit module is finite dimensional precisely when
$m$ is elliptic.

The classification of 
all spherical finite dimensional simple $\Hb$-modules
may be deduced from the knowledge of
the spherical finite dimensional simple $\Hb''$-modules.
Another question is to understand the representation of $\Hb$
in the homology of the (elliptic homogeneous) affine Springer fibers.
They enter in the geometric construction of $\Hb$ in \cite{V1}.
The homology of affine Springer fibers is difficult to compute in general.
We describe explicitely all the spherical Jordan-H\"older
factors (modulo a technical hypothesis).
This classification involves
interesting combinatorial objects which already appear
in local Langlands correspondence for $p$-adic groups. 
See \cite{R2}, \cite{R3} for the combinatorial aspects, and
\cite{R4}, \cite{GR}, \cite{DR} for the representation theoretic aspects.
Affine Hecke algebras are related to unramified Langlands 
correspondence via Bernstein's functor.
DAHA's seem to be related to the tamely ramified
correspondence. 

More precisely, let $k$ be a non-Archimedean local field of characteristic 
zero with residue field $\ken=\FF_q$.
Fix an algebraic closure $\bar k$ of $k$.
Let $K$ be the maximal unramified extension of $k$ in $\bar k$, 
and $K_t$ be the maximal tame extension of $K$ in $\bar k$. 
The reduction map yields an exact sequence
$$1\to\Ic_t\to\Gal(K_t,k)\to\Gal(\bar\ken,\ken)\to 1,$$
where $\Ic_t$ is the tame inertia subgroup. 
Let $F\in\Gal(K_t,k)$ be a lifting of the Frobenius element in 
$\Gal(\bar\ken,\ken)$.
The tame Weil group $\Wc_t\subset\Gal(K_t,k)$ is isomorphic 
to the semi-direct product $\la F\ra\ltimes\Ic_t$, 
where $\la F\ra$ is the subgroup generated by $F$
and
$$F\gamma F^{-1}=\gamma^q,\quad\forall\gamma\in\Ic_t.$$
The tame Weil-Deligne group is the semi-direct
product $\Wc'_t=\Wc_t\ltimes\CC$, where $\Wc_t$ acts on $z\in\CC$ as follows
$$FzF^{-1}=qz,\quad \gamma z\gamma^{-1}=z,\quad\forall\gamma\in\Ic_t.$$

Let $\Gb$ be a split connected simple group-$k$-scheme of adjoint type.
Each continuous cocycle $c:\Gal(\bar k/k)\to\Gb$ yields a new 
$k$-structure $\Gb_c$ on $\Gb$.
If $c,c'$ define the same class in the Galois cohomology group
$H^1(k,\Gb)$, then the group-$k$-schemes $\Gb_c$, $\Gb_{c'}$
are isomorphic.
The group-$k$-schemes $\Gb_\omega$, with $\omega\in H^1(k,\Gb)$, 
are called the pure inner forms of $\Gb$.

The Langlands dual group of $\Gb$ is a complex Chevalley group $G_\CC$.
Let $G_0$, $T_0$ be the sets of $\CC$-points of $G_\CC$ and of a maximal
torus $T_\CC\subset G_\CC$.
Let $X_0$ be the group of characters of the torus $T_0$,
$\check Y_0$ be the root lattice,
and $W_0$ be the Weyl group.

A tamely ramified Langlands parameter (or TRLP) is a continuous homomorphism 
$$\varphi:\Wc'_t\to G_0.$$ 
Assume that $\varphi$ is a regular elliptic TRLP (or ERTRLP), i.e., that
$Z_{G_0}(\varphi(\Ic_t))$ is a maximal torus
and that $Z_{G_0}(\varphi(\Wc_t))$ is finite.
Since $\varphi$ is continuous, its restriction to $\Ic_t$ factors through
the group of units of a finite extension of $\ken$.
Thus $\varphi(\Ic_t)$ is the cyclic group 
generated by a regular semisimple element $s_\varphi$ 
of order prime to $q$. 
Set $n_\varphi=\varphi(F)$. 
Since $\varphi$ vanishes on the subgroup $\CC\subset\Wc'_t$, 
it is uniquely determined by the pair $(s_\varphi, n_\varphi).$
After conjugating by $G_0$, we may assume that
$Z_{G_0}(\varphi(\Ic_t))$ equals $T_0$.
Then $n_\varphi\in N_{G_0}(T_0),$ 
and its image $w_\varphi$ in $W_0$ satisfies the relation 
${}^{w_\varphi}s_\varphi=s_\varphi^{q}.$
Let $T_{0,\reg}\subset T_0$ be the subset consisting of regular semisimple
elements.
The assignement
$\varphi\mapsto(s_\varphi,w_\varphi)$ is a bijection from the
set of $G_0$-conjugacy classes of ERTRLP's to the set of $W_0$-conjugacy classes
of pairs $(s,w)\in T_{0,\reg}\times W_0$ such that
$${}^{w}(s^q)=s,\quad
T_0^w\roman{\ is\ finite}.$$
Notice that if $m=o(w)$ then we have $s^{q^m}=s$.
See \cite{DR, rem.~4.1.1} for details.
To simplify we'll assume that $w$ is regular and elliptic.

%In particular its conjugacy class is uniquely determined by the order
%$m_\varphi=o(w_\varphi)$ of $w_\varphi$.
According to Langlands correspondence, to the 
$G_0$-conjugacy class of $\varphi$ should correspond
a $L$-packet $\pi(\varphi)$
of supercuspidal complex representations 
of the pure inner forms of $\Gb$.
More precisely, let $\Irr(T_0^{w_\varphi})$ be the
group of irreducible characters 
of the finite Abelian group $T_0^{w_\varphi}$.
DeBacker and Reeder construct a $L$-packet 
$$\pi(\varphi)=\{\pi(\varphi,\rho);\rho\in\Irr(T_0^{w_\varphi})\}.$$
A theorem of Kottwitz implies that
$$H^1(k,\Gb)=\Irr(Z(G_0)).$$
Thus, restricting a character to the center $Z(G_0)$ we get a map
$$\Irr(T_0^{w_\varphi})\to H^1(k,\Gb),\quad
\rho\mapsto\omega_\rho,$$
such that $\pi(\varphi,\rho)$ is a representation of the pure inner form
$\Gb_{\omega_\rho}$.
Observe that
$$\Irr(T_0^{w_\varphi})\simeq X_0/(1-w_\varphi)X_0.$$
Thus, if $\rho\in\check Y_0/(1-w_\varphi)X_0$
then $\pi(\varphi,\rho)$ is a representation of the group $\Gb$.

Recall that $G_0$-conjugate triples $(s_\varphi,n_\varphi,\rho)$
should yield isomorphic representations of $\Gb$. 
Indeed, it is proved in
\cite{R4} that we have the following equivariance property
$$\pi(s_\varphi,n_\varphi,\rho)=\pi({}^ws_\varphi,n_\varphi,{}^w\rho),
\quad
\forall w\in Z_{W_0}(w_\varphi).$$
%Since $T^{w_\varphi}$ is finite, the map
%$T_0\to T_0$, $t\mapsto {}^wt\,t^{-1}$ has finite fibers,
%hence is surjective.
%Thus, conjugating $n_\varphi$ by elements of $T_0$ we can change it to 
%any other representative of $w_\varphi$.
Thus the depth-zero supercuspidal complex representations of $\Gb$ we get
are labelled by the $W_0$-conjugacy classes
of triples $(s,w,\rho)$ with $s,w$ as above and
$\rho\in\check Y_0/(1-w_\varphi)X_0$.
In other words, to each $W_0$-conjugacy class of pairs
$(w,[\rho])$ where $w\in W_0$ is regular elliptic 
and $[\rho]$ is a $Z_W(w)$-orbit in 
$\check Y_0/(1-w)X_0$ is attached a set 
of cuspidal representations of $\Gb$.

We prove in 3.3.1, 3.3.6 below that
the set of isomorphism classes of
finite dimensional simple spherical $\Hb$-modules 
which are Jordan-H\"older composition factors of the homology of
an elliptic affine Springer fiber
decomposes into families labelled by conjugacy classes of regular 
elliptic elements in $W_0$.
%Let $w_m$ denote an elliptic element in $W_0$ of order $m$. 
We also prove that the simple modules in the family
corresponding to the $W_0$-orbit of $w$ are labelled 
by a positive integer $k$ which is prime to the order of $w$
and a $Z_W(w)$-orbit in the set $\check Y_0/(1-w)X_0.$ 
This coincidence is certainly not a `hasard'.
%We'll come back to this and
%to the conjecture in \cite{VV2} elsewhere.
In particular, it raises the following interesting problem.

\proclaim{\bf Question}
What is the $p$-adic interpretation of the dimension of the
finite dimensional simple spherical $\Hb$-modules ?
\endproclaim

More precisely, 
given the conjugacy class of an elliptic regular element
$w\in W_0$ of order $m$, a positive integer $k$ prime to $m$,
and a $Z_W(w)$-orbit
$[\rho]$ in $\check Y_0/(1-w)X_0,$
we have a finite dimensional simple $\Hb$-module
$L_{m,k,[\rho]}.$
By 3.4.2 below we have 
$$\dim(L_{m,k,[\rho]})=k^nd_{m,[\rho]},$$
where $n$ is the rank of $G_0$ and 
$d_{m,[\rho]}$ is a positive integer which
depends only on $m$ and $[\rho]$. 
Now, fix a ERTRLP $\varphi$ such that
$w_\varphi=w$ and an element $\rho\in[\rho]$.
%The integer $o(s_\varphi)$ must divide the integer $m_\varphi$.  
%Since the cardinal
%$q$ of the residue field of $k$ is prime to $o(s_\varphi)$,
%It is reasonnable to expect that 
%the integer $k$ should be related to $q$.
We expect that the integer $q^nd_{m,[\rho]}$ 
should play some role for the supercuspidal representation 
$\pi(\varphi,\rho)$ of the group $\Gb$.
 
\head Notation \endhead

A reflection group will mean a group generated by reflections of order 2.
A group generated by complex reflections will be called a 
complex reflection group.

All schemes are assumed to be separated.
By variety we mean a scheme of finite type over an algebraically closed field.
By a coherent sheaf on a not connected scheme we mean a family 
of coherent sheaves on each connected components supported 
on a finite number of components.
By a vector bundle on a not connected scheme we mean a family of vector
bundles on each connected components (without support condition).
By a virtual vector bundle we mean a family of
formal $\CC$-linear combinations of vector bundles on each componenents.

Fix formal variables $\eps_\ell$, with $\ell\neq 0$, such that
$(\eps_{\ell\ell'})^{\ell'}=\eps_{\ell}$.
Set $\eps=\eps_1$, and $\eps^c=\eps^k_m$
for each $c=k/m\in\QQ^\times.$
Put $K=\CC((\eps))$, $A=\CC[[\eps]]$,
and $\bar K=\bigcup_{\ell}\CC((\eps_\ell))$.
Fix formal variables $q,t,\kappa$.
Set $\CC_{q,t}=\CC[q^{\pm 1}, t^{\pm 1}]$,
$\CC_{t}=\CC[t^{\pm 1}]$,
and $\CC_{\kappa}=\CC[\kappa]$.

For any $\CC$-scheme $\Sc_\CC$ 
set 
$\Sc_0=\Sc_\CC(\CC)$, 
$\Sc_+=\Sc_\CC(A)$, 
$\Sc=\Sc_\CC(K)$,
and
$\bar\Sc=\Sc_\CC(\bar K)$.
%and any $\CC$-algebra $\k$
%and let $\Sc_\k$ be the $\k$-scheme $\Sc_\CC\otimes_\CC\k.$
Let $K(\Sc_\CC)$, $H_*(\Sc_\CC,\CC)$
be the complexified Grothendieck group
of coherent sheaves and the Borel-Moore homology with complex coefficients
respectively.
Given an action of a complex algebraic group $G_\CC$ on $\Sc_\CC$
we write also $K^{G_0}(\Sc_\CC)$, $K_{G_0}(\Sc_\CC)$ for
the complexified Grothendieck groups of $G_\CC$-equivariant coherent sheaves
and $G_\CC$-equivariant vector bundles respectively.

For each locally closed subset $S$ we write
$\CC_S$ for the constant sheaf over $S$.
The stalk of a sheaf $\Fc$ at a point $x$ is writen
$\Fc_x$ as usual. 

The cardinal of a finite set $S$ is denoted by $|S|$.
Given an action of a group $G$ on $S$ and subsets 
$S'\subset S$, $G'\subset G$, we write
$$\aligned
&Z_G(S')=\{g\in G;gx=x,\forall x\in S'\},
\\
&N_G(S')=\{g\in G;gx\in S',\forall x\in S'\},
\\
&S^{G'}=\{x\in S;gx=x,\,\forall g\in G'\}.
\endaligned$$ 

We write $\Irr(\Cc)$ for the set of isomorphism classes 
of simple objects in an Abelian category $\Cc$.
If $\Cc$ is the category of complex representations of
an algebra $\Ab$ or a group $G$ we may write
$\Irr(\Ab)$, $\Irr(G)$ for $\Irr(\Cc)$.

\head Contents \endhead

\item{1.} Homogeneous elliptic regular elements in loop Lie algebras
\itemitem{1.1.} Elliptic and regular elements in $W_0$ 
\itemitem{1.2.} Reminder on maximal tori in $G$ 
\itemitem{1.3.} Regular semisimple elements in $\gen$
\item{2.} Finite dimensional representations of rational DAHA's
\itemitem{2.1.} Reminder on DAHA's 
\itemitem{2.2.} Reminder on degenerate DAHA's 
\itemitem{2.3.} Reminder on rational DAHA's 
\itemitem{2.4.} DAHA's and affine Springer fibers
\itemitem{2.5.} Isogenies and simple modules of DAHA's 
\itemitem{2.6.} Fourier-Sato tranform 
\itemitem{2.7.} Geometrization of the polynomial representation 
\itemitem{2.8.} Classification of spherical simple finite dimensional modules
coming from $\Hb''$ 
\item{3.} Finite dimensional representations of DAHA's 
and affine Springer fibers
\itemitem{3.1.} The evaluation map 
\itemitem{3.2.} Homology of affine Springer fibers
\itemitem{3.3.} Spherical simple finite dimensional modules
\itemitem{3.4.} Dimension 
\itemitem{3.5.} The Coxeter case and the sub-Coxeter case
\item{4.} Computations
\itemitem{4.1.} Type $A_n$
\itemitem{4.2.} Type $B_n$
\itemitem{4.3.} Type $C_n$
\itemitem{4.4.} Type $D_n$
\itemitem{4.5.} Type $E_6$
\itemitem{4.6.} Type $E_7$
\itemitem{4.7.} Type $E_8$
\itemitem{4.8.} Type $F_4$
\itemitem{4.9.} Type $G_2$
\item{A.} Appendix 
\itemitem{A.1.} Fourier-Sato transform
\itemitem{A.2.} Induction 
\itemitem{A.3.} Standard modules 
\itemitem{A.4.} Weights of simple modules 
\itemitem{A.5.} Index of notation 
\itemitem{A.6.} Remarks on \cite{V1} 

\head 1. Homogeneous elliptic regular elements in loop Lie algebras\endhead

Let $(X_0,\Delta_0,\Pi_0,Y_0,\check\Delta_0,\check\Pi_0)$ 
be a based root datum.
This means that $X_0,Y_0$ are free Abelian groups in duality via a perfect
pairing 
$$X_0\times Y_0\to\ZZ,
\quad(\l,\lambdav)\mapsto\l\cdot\lambdav,
\leqno(1.0.1)$$
that
$\Delta_0\subset X_0$, $\check\Delta_0\subset Y_0$ 
are root systems with a bijection
$\Delta_0\to\check\Delta_0,$ $a\mapsto\av$ commuting with 
the map $a\mapsto-a$ and such that $a\cdot\check a=2$, and that
$\Pi_0$, $\check\Pi_0$ are bases of $\Delta_0$, $\check\Delta_0$ 
respectively.
We'll always assume that the root systems $\Delta_0$, $\check\Delta_0$ 
are reduced.
Let $n$ be the rank.

Let $W_0$ be the abstract Weyl group,
$\Delta_0^+$ be the set of positive roots,
and $\Delta_0^-$ be the set of negative roots.
Let $\{s_a;a\in\Pi_0\}\subset W_0$ be
the set of simple reflections.
Write $a_i$, $i\in I_0$, for the simple roots, i.e.,
the elements in $\Pi_0$, and
$\av_i$, $i\in I_0$, for the simple coroots.
Let also $\{o_i\}$ be the set of fundamental weights, and
$\{\ov_i\}$ be the set of fundamental coweights.
We have set $I_0=\{1,2,...n\}$.
Set $\rhov=\sum_i\ov_i$,
and $w_0$ be the longest element in $W_0$.

Let $G_\CC$ be the connected reductive group over $\CC$
associated to the root datum. 
Let $T_\CC\subset G_\CC$ be a maximal torus, 
and $\gen_\CC$, $\ten_\CC$ be the Lie algebra $\CC$-schemes
of $G_\CC$, $T_\CC$.
Taking the points over $\CC$,
$A$, $K$, $\bar K$ we get the groups $G_0$, $T_0$, $G_+$, $T_+$, $G$, 
$T$, $\bar G$, 
and the Lie algebras
$\gen_0$, $\ten_0$, $\gen_+$, $\ten_+$, $\gen$, $\ten$, $\bar\gen$.
For any $\CC$-algebra $k$ we abbreviate
$G_k=G_\CC\otimes_\CC k$,
$\gen_k=\gen_\CC\otimes_\CC k$.
The Langlands dual group $\check G_\CC$ is the connected reductive group
associated to the root datum
$(Y_0,\check\Delta_0,X_0,\Delta_0).$

For any group-scheme $H$ let
$X^*(H)$ be the group of characters $H\to\GG_m$ and
$X_*(H)$ the group of cocharacters $\GG_m\to H.$
Notice that $X_0=X^*(T_\CC)$, $Y_0=X_*(T_\CC)$, and
$W_0\simeq N_{G_0}(T_0)/T_0$.

Set $\check V_0=Y_0\otimes\CC$ and $V_0=X_0\otimes\CC$.
The pairing 1.0.1 yields a perfect pairing of $\CC$-vector spaces
$$V_0\times \check V_0\to\CC,
\quad(\l,\lambdav)\mapsto\l\cdot\lambdav.$$
So
$\check V_0$ 
is identified 
with the dual $V^*_0$ of $V_0$.
Fix a $W_0$-invariant positive nondegenerate
symmetric bilinear form $(\ :\ )$ on $V_0$.
It yields a vector space isomorphism 
$$\nu:V_0\to V^*_0=\check V_0,\quad \lambda\mapsto(\lambda:\ ),$$
called the standard isomorphism.
Write $(\ :\ )$ again for the 
bilinear form on $\check V_0$ which coincides with $(\ :\ )$
under $\nu$. 
We have $\nu(a)=2\av/(\av:\av)$ for all roots.
Normalize $(\ :\ )$ in such a way that $\nu(\theta)=\check\theta.$ 
%Thus the pairing 1.0.1 may be viewed as the unique bilinear form
%on $\check V_0$ (or on $V_0$) such that
%$\av\cdot\nu(a)=a(\av)$ for all $a\in X_0$.
Put 
$V_{0,\RR}=X_0\otimes\RR$,
$\check V_{0,\RR}=Y_0\otimes\RR$.

Let $\Delta_\re=\Delta_0\times\ZZ$,
$\Delta=\Delta_\re\cup(\{0\}\times\ZZ)$,
$\Delta^+_\re=(\Delta_0\times\ZZ_{>0})\cup(\Delta_0^+\times\{0\})$,
and $\Pi=\{\a_i\}\cup\{(a,1);a\in\Delta_{0,\min}\}$.
Here $\Delta_{0,\min}$ is the set of minimal roots for the
usual order, and $\a_i=(a_i,0)$.
Elements of $\Pi$ are called simple affine roots.
Put also $\check \Delta_{0,\min}=\{\av;a\in\Delta_{0,\min}\}$,
and $\alphav_i=(\av_i,0).$
The affine Weyl group associated to the group $G_\CC$,
or to the root datum $(X_0,\Delta_0,Y_0,\check\Delta_0)$,
is the group $W$ of affine automorphisms of $\check V_{0,\RR}$
generated by $W_0$
and the translations by elements of $Y_0$.
It is isomorphic to the group $N_G(T)/T_+$.
There is a unique group isomorphism $W\to W_0\ltimes Y_0$ taking
the affine reflection $s_{(a,\ell)}$ to $\xi_{-\ell \av}s_a$,
where we write $w,\xi_{\lambdav}$ for $(w,0)$, $(1,\lambdav)$ respectively.
The group $W$ acts on
$X_0\times\ZZ$ by the same formulas as in 1.0.2 below.
Consider the group $\Omega_W=\{w\in W;w(\Pi)=\Pi\}$.
It is an Abelian group isomorphic to $Y_0/\ZZ\check\Delta_0$.
There is a normal subgroup $W'\subset W$ such that
$W\simeq W'\rtimes\Omega_W$.
It is the Coxeter group generated by the set
$\{s_\a;\a\in\Pi\}$ of simple affine reflections.
See \cite{L1, sect.~1} for details.

We fix a system of root vectors $\{e_a\}\subset\gen_0$.
Put $f_{a}=e_{-a},$ 
and $e_\a=e_a\otimes\eps^\ell\in\gen$ 
for each real affine root $\a=(a,\ell)$.

Unless specified otherwise we'll assume that
$G_\CC$ is a Chevalley group, i.e.,
it is simple and simply connected.
In this case we have $W=W'$.
This group is called
the affine Weyl group of the root system $\Delta_0$ for short.
Let $h$ be the Coxeter number,
$\theta$ be the highest root,
and $\thetav$ be the highest short coroot.
We have $Y_0=\sum_{i\in I_0}\ZZ\av_i$, $X_0=\sum_{i\in I_0}\ZZ o_i$,
and $\Pi=\{\a_i\}_{i\in I}$ with $\a_0=(-\theta,1)$ and $I=I_0\cup\{0\}.$
We'll write $s_i$ for the reflection $s_{a_i}$ or the affine reflection
$s_{\a_i}$.
We hope that this ambiguity in the notation
will not create confusion.
We set $\check Y_0=\sum_i\ZZ a_i$, $\check X_0=\sum_i\ZZ \ov_i$.

If $G_\CC$ is simple of adjoint type then $W$ is 
the extended affine Weyl froup
$W^e=W_0\ltimes\ZZ\check X_0$ of the root system $\Delta_0$, 
and we set $\Omega=\Omega_{W^e}$.
Let $\{\omegav_i\}\subset\check X_0\times\ZZ$
be the set of affine fundamental coweights.
Fix a positive integer $N$ such that $\l\cdot\lambdav\in(1/N)\ZZ$
for each $\l\in X_0$, $\lambdav\in\check X_0$.
The group $W^e$ acts on
$X_0\times(1/N)\ZZ$, $\check X_0\times\ZZ$ by the following formulas 
$$\aligned
{}^w(\mu,\ell)=(\mu,\ell-\mu\cdot\lambdav),
\quad
\hfill
&
{}^w(\muv,\ell)=(\muv+\ell\lambdav,\ell)\quad \roman{if}\ w=\xi_{\lambdav},
\hfill\cr
{}^w(\mu,\ell)=({}^w\mu,\ell),
\quad
\hfill
&
{}^w(\muv,\ell)=({}^w\muv,\ell)\quad \roman{if}\ w\in W_0.
\hfill
\endaligned
\leqno(1.0.2)$$
The pairing such that $(\mu,\ell)\cdot(\muv,k)=\mu\cdot\muv+\ell k$
is $W^e$-invariant.
Let $\ov_r$, $r\in O$, be the minuscule coweights.
For each $r$ the subgroup $Z_{W_0}(\ov_r)$ is generated by $\{s_i\}_{i\neq r}$.
Let $w_r\in W_0$ be such that
$w_0^{-1}w_r$ is the longest element in 
$Z_{W_0}(\ov_r)$.
The action of the element $\pi_r=\xi_{\ov_r}w_r^{-1}$ in $W^e$
on $\check X_0\times\ZZ$ preserves the set of affine fundamental coweights.
The group $\Omega$ is equal to $\{\pi_r\}$.
It is identified the group of automorphisms of 
the affine Dynkin diagram, so that $\pi_r$ is taken to
the unique automorphism such that $0\mapsto r$.
The group isomorphism
$W^e\to W\rtimes\Omega$
takes $w\in W$ to itself,
and $\xi_{\ov_r}$ to $\pi_rw_r$.
Here, we write $w,\pi$ for the elements $(w,1),(1,\pi)$ 
in the semi-direct product.

Set $X=X_0\times\ZZ^2$ and $\check X=\check X_0\times\ZZ^2$.
%The second copy of $\ZZ$ is the set of characters and cocharacters
%of $Z(\hat G)$ respectively.
We'll identify $X_0$, $\check X_0$,
$X_0\times\ZZ$, $\check X_0\times\ZZ$ with the subgroups
$X_0\times\{0\}$, $\check X_0\times\{0\}$,
$X_0\times\ZZ\times\{0\}$, $\check X_0\times\ZZ\times\{0\}$
in the obvious way.
The $W$-action extends to
$X$, $\check X$
so that
$$
\aligned
{}^{s_0}(\l,\ell,z)
=({}^{s_\theta}\l+z\theta,
\ell+\l\cdot\thetav-z, z),
&\quad
{}^{s_i}(\l,\ell,z)
=({}^{s_{i}}\l,\ell,z),
\\
{}^{s_0}(\lambdav,\ell,z)=({}^{s_\theta}\lambdav+\ell\thetav,\ell,
z+\theta\cdot\lambdav-\ell),
&\quad
{}^{s_i}(\lambdav,\ell,z)
=({}^{s_{i}}\lambdav,\ell,z),
\endaligned
$$
for all $i\in I_0.$
See \cite{C3, prop.~17.20} for instance.
%These actions are dual to each other with respect to the pairing
%such that 
%$$(\l,\ell,z)\cdot(\lambdav,k,y)=\l\cdot\lambdav+\ell k+zy.
%\leqno(1.0.3)$$
Let $\{\omega_i\}\subset X$
be the set of affine fundamental weights
and $\delta=(0,1,0)\in X$.
The set of integral dominant weights is
$$X_+=\bigoplus_{i\in I}\ZZ_{\ge 0}\o_i\oplus\ZZ\delta.$$
Notice that $\a_0=(-\theta,1,0)$,
$\a_i=(a_i,0,0)$,
$\omega_0=(0,0,1)$,
$\omega_i=(o_i,0,o_i\cdot\thetav)$,
$\alphav_0=(-\thetav,0,1),$
$\alphav_i=(\av_i,0,0)$,
$\omegav_0=(0,1,0)$,
and
$\omegav_i=(\ov_i,\ov_i\cdot\theta,0)$.
Put $\check\delta=(0,0,1)\in\check X.$ 

Set $V=X\otimes\CC$ and 
$\check V=\check X\otimes\CC$
with the induced $W$-actions.
Let $V_\RR$,
$\check V_\RR$ be the corresponding $\RR$-vector spaces. 
Observe that $V=V_0\times\CC^2$ 
and $\check V=\check V_0\times\CC^2$.
There is an unique $W$-invariant perfect pairing 
$V\times\check V\to\CC$ such that
%$(\lambda,\lambdav)\mapsto\lambdav\cdot\nu(\lambda)$ is $W$-invariant.
$$(v,\ell,z)\cdot(v',\ell',z')=v\cdot v+\ell\ell'+zz'.
\leqno(1.0.3)$$
Identify $\check V$ with the dual $V^*$ of $V$ via 1.0.3.
There are unique $W$-invariant nondegenerate
semidefinite bilinear forms $(\ :\ )$ on 
$V$, $\check V$ which coincide with 
$(\ :\ )$ on $V_0$, $\check V_0$ and satisfy the following relations
$$
(\delta:\omega_0)=(\check\delta:\check\omega_0)=1,
\quad
(\delta:\delta+V_0)=(\omega_0:\omega_0+V_0)=
(\check\delta:\check\delta+\check V_0)=
(\check\omega_0:\check\omega_0+\check V_0)=0.
$$
It yields the standard isomorphism
$$\nu:V\to V^*=\check V,\quad \lambda\mapsto(\lambda:\ ).$$
We have $\nu(\delta)=\check\delta$, 
$\nu(\alpha_0)=\alphav_0$,
$\nu(\omega_0)=\check\omega_0$
and
$\nu(\alpha_i)=(\nu(a_i),0,0)$,
$\nu(\omega_i)=(\nu(o_i),o_i\cdot\check\theta,0)$
for each $i\in I_0$.

\subhead 1.1. Elliptic and regular elements in $W_0$\endsubhead

By an eigenvector of an element of $W_0$ we'll mean an eigenvector in
the reflection representation on $\ten_0$.
A element of $W_0$ is regular iff it has an eigenvector
which is regular (in the sense that no root vanishes on it).
Let $W_0[m]$ be the set of regular elements of order $m$.

An element of $W_0$ is called elliptic (with respect to $\ten_0$) iff
it has no nonzero fixed vector in $\ten_0$.
A conjugacy class in $W_0$ is elliptic,
i.e., it contains elliptic elements, iff it is cuspidal, i.e.,
iff it has an empty intersection with each proper parabolic subgroup of $W_0$.
See \cite{GP, exe.~3.13, 3.16} for instance.

The integers which occur as orders of regular or elliptic regular elements
are called regular or elliptic numbers respectively.
Write EN, RN for the sets of elliptic, regular numbers.
Notice that $1\in\RN$ but $1\notin \EN$.

Let $\varphi_i\in\CC[\ten_\CC]^{W_0}$, $i\in I_0$,
be a set of primitive invariants.
Put $d_i=\deg(\varphi_i).$
For each integer $m\neq 0$ set $I_m=\{i\in I_0;m|d_i\}.$

The following are proved in \cite{S3, thm.~4.2},
\cite{DL, thm.~2.5, cor.~2.9} respectively.

\proclaim{1.1.1. Theorem}
(a)
The sets $W_0[m]$, $m\in\RN$, are conjugacy classes of $W_0$.
The order of the eigenvalue of a regular eigenvector
of any element $w$ in $W_0[m]$ is $m.$
If $V$ is an eigenspace of $w$ containing a regular vector
then $Z_{W_0}(w)$ is a complex reflection group in $V$
whose set of degrees is $\{d_i;i\in I_m\}$.

(b)
Let $V$ be an eigenspace of $w\in W_0[m]$ containing a regular vector.
The complex reflection group $Z_{W_0}(w)$ is irreducible,
its reflection hyperplanes are the traces
of the reflection hyperplanes of $W_0$, and a set of primitive invariants
is given by the restrictions to $V$ of the primitive invariants $\varphi_i$ with
$i\in I_m$.
\endproclaim

\vskip3mm

\noindent{\bf 1.1.2. Example.}
The conjugacy class $\cox=W_0[h]$
consists of the Coxeter elements, i.e.,
of the products of all simple reflections in $W_0$ in any order.
They are elliptic regular.

\vskip3mm

The regular elements have been classified in \cite{S3}
(notice that there is an error in table 2).
Let $W_{0,*}$ be the set of conjugacy classes in $W_0$.
We'll label elements in $W_{0,*}$ as in \cite{C2}.
Write $w_*$ for the conjugacy class of $w$,
and $(w_*)^{[i]}$ for $(w^i)_*$.
The conjugacy classes of elliptic regular elements are the following.

$$\vbox{
\offinterlineskip

\def\tablerule{\noalign{\hrule}}

\def\tableskip{&\omit&height 6pt&\omit&\omit&\omit&\cr}

\halign{
\tabskip= .5em\vrule #
&\it#
&\vrule #
&\vtop{\hsize=7pc\it #\strut}
&#\vrule
&\vtop{\hsize=11pc\it #\strut}
&#\vrule\tabskip=-.1em\cr
\tablerule
\tableskip
& $G_\CC$ && $W_0[m]$&&m&\cr
\tablerule
\tableskip
&$A_n$ && $A_n$&& n+1&\cr
\tablerule
\tableskip
&$C_n$ && $C_n^{[n/r]}$ &&$2r$ with $r|n$&\cr
\tablerule
\tableskip
&$D_n$ && $D_n^{[(n-1)/r]}$

$D_n(a_{n/2-1})^{[n/m]}$ &&
$2r$ with $r|n-1$ and $2r\nmid n-1$

$2r$ with $2r|n$ 
&\cr
\tablerule
\tableskip
&$E_6$ && $E_6^{[12/m]}$

$E_6(a_1)$
&& $3,6,12$

$9$
&\cr
\tablerule
\tableskip
&$E_7$ && $E_7^{[18/m]}$

$E_7(a_1)$&&
$2,6,18$

$14$
&\cr
\tablerule
\tableskip
&$E_8$ && $E_8^{[30/m]}$

$E_8(a_1)^{[24/m]}$

$E_8(a_2)$&&
$2,3,5,6,10,15,30$

$4,8,12,24$

$20$
&\cr
\tablerule
\tableskip
&$F_4$ &&$F_4^{[12/m]}$

$B_4$&&
$2,3,4,6,12$

$8$
&\cr
\tablerule
\tableskip
&$G_2$ &&$G_2^{[6/m]}$ && $2,3,6$&\cr
\tablerule
}}
$$

\noindent{\bf 1.1.3. Remark.}
Write the element $s_1s_2s_3...$ as $123...$ for short.
For the computations in section 3 it is useful to have 
an explicit representative in each class $W_0[m]$.
Recall that $123...n\in W_0[h]$ for all types.
For exceptional types, the classification of cuspidal classes
in \cite{GP, tables B.3-6} yields the following :
$123234\in B_4$,
$12342546\in E_6(a_1)$,
$123425467\in E_7(a_1)$,
$1234254678\in E_8(a_1)$, and
$123425465478\in E_8(a_2).$

\subhead 1.2. Reminder on maximal tori in $G$\endsubhead

Fix complex numbers $0\neq\tau_\ell$, $\ell>0$, such that
$(\tau_{\ell\ell'})^{\ell'}=\tau_{\ell}$
and $\tau=\tau_1$ is not a root of unity.
For each $c$ as above we write $\tau^c$ for $\tau^k_m$.
Let $F_\tau\in\Aut(\bar K)$ be such that
$\eps_\ell\mapsto\tau_\ell\eps_\ell$.
We'll write $F$ for $F_\varpi$,
where $\varpi_\ell=\exp(2i\pi/\ell)$.
Put also $\lambdav_\ell=\lambdav(\varpi_\ell)$
for each $\lambdav\in X_*(G_\CC)$.

%Write $G'_\CC$ for the adjoint form of $G_\CC$.
%Thus $\bar G'=\bar G/Z(\bar G)$ and
%$G'=\{g\in\bar G;g^ {-1}F(g)\in Z(\bar G)\}/Z(\bar G)$.

The set of $G$-conjugacy classes
of maximal tori in $G_K$ (or in $\gen_K$) is in one-to-one
correspondence with $W_{0,*}$ by \cite{KL2, sec.~1}.
See also \cite{D1, sec.~4.2} for more general results.
The conjugacy class associated to a maximal torus is called its type.

Choose a lift of finite order $n_w$ in $N_{G_0}(T_0)$ of $w$.
Set $m'=o(n_w)$.
As $n_w$ is a torsion element of $G_0$,
there is a cocharacter $\lambdav$ such that $n_w=\lambdav_{m'}$.
Set $F_w=(\ad n_w)^{-1}\circ F.$
The automorphism $(\ad\lambdav(\eps_{-m'}))$ takes the $K$-torus
$T_{\bar K}^{F_w}$
into a maximal torus of type $w_*$ in $G_{K}$.
%Let $\ten_{w,K}$ be its Lie algebra.
%So the $K$-scheme $\ten_{K,w}$ represents the functor
%$\{K\text{-algebras}\}\to\{\text{sets}\},$
%$S\mapsto(\ad\lambdav(\eps_{-m'}))(\bar\ten^{F_w}\otimes_KS).$

An element $x\in\bar\gen$ is regular semisimple (or RS) iff its centralizer
$\zen_{\bar\gen}(x)$ is a maximal torus.
The type of a RS element $x\in\gen$ is the type of the maximal torus
$\zen_\gen(x)$.

\subhead 1.3. Regular semisimple elements in $\gen$\endsubhead

An element $x\in\bar\gen$ is topologically nilpotent (or TN) iff
$(\ad x)^\infty=0$ in $\End(\bar\gen)$.
Write TNS for topologically nilpotent semisimple element.
A semisimple element is homogeneous (or HS) of degree $c$
iff it is $\bar G$-conjugate into $\ten_0\otimes\eps^c$.
It is elliptic regular (or ERS) iff it is RS, TN, and
the torus $Z_{G_{K}}(x)^\circ$ is anisotropic
(i.e., its group of cocharacters is trivial).
We will write HRS for HS and RS, and HERS for HS and ERS.

From now on we assume that $m>0$ and $k,m$ are relatively prime.
Given $\Sc\subset\bar\gen$ let $\Sc_\reg,\Sc_\HRS,\Sc_\ERS,\Sc^c\subset\Sc$
be the subsets of RS, HRS, ERS elements and of elements $x$ such that
$(\ad\rhov(\tau^c))F_\tau(x)=\tau^c x$ respectively.

Put $e_R=\sum_ie_{a_i}$ and $f_R=\sum_i(o_i\cdot 2\rhov)f_{a_i}$.
Assume that $\{e_{a_i},\av_i,f_{a_i}\}$ is a
$\sen\len_2$-triple for each $i\in I_0$.
Then $\phi=\{e_R,2\rhov,f_R\}$ is also a $\sen\len_2$-triple.
See \cite{CM, p.~58}.

A theorem of Chevalley states that the restriction map
$\CC[\gen_\CC]\to\CC[\ten_\CC]$ induces an isomorphism
of categorical quotients $\eta:\gen_\CC/G_0\to\ten_\CC/W_0$.
We regard this isomorphism as a $G_0$-invariant map
$\gen_\CC\to\ten_\CC/W_0$.

The $I_0$-uple $\varphi=(\varphi_i)$ of primitive invariants
is an isomorphism $\ten_\CC/W_0\to\AA_\CC^{I_0}.$
Set $f_i=d(\varphi_i\eta)(f_R)$.
View $f_i$ as an element in $\gen_0$ via the Killing form.

The Kostant section in $\gen_\CC$ is the closed subscheme
$\sen_\CC=e_R+\zen(f_R)_\CC$.
See \cite{K5, sec.~ 2.4} for a short review,
and \cite{S2, sect.~7} for further results.

Let $B_0^-,B_0\subset G_0$ be the Borel subgroups associated
to the pairs $(T_0,\Delta_0^-)$, $(T_0,\Delta_0^+)$ respectively.
Let $U_0^-,U_0$ be their unipotent radicals,
and $\ben_0$, $\uen_0$, $\ben_0^-$, $\uen_0^-$ 
be the corresponding Lie subalgebras of $\gen_0$.

\proclaim{1.3.1. Proposition}
(a)
Any $G_0$-orbit in $\gen_{0,\reg}$ meets $\sen_{0,\reg}$ exactly once,
and any $G$-orbit in $\gen_\reg$ meets $\sen_\reg$ exactly once.
%Any element in $\gen$ is $G$-conjugate into $\sen$.

(b)
The composition of the closed embedding $\sen_\CC\subset\gen_\CC$
and $\eta$ is an isomorphism of algebraic varieties
$\sen_\CC\to\ten_\CC/W_0$.

(c)
The $I_0$-uple $(f_i)$ is a basis of $\zen(f_R)_0$
such that $[\rhov,f_i]=(1-d_i)f_i.$
In particular, we have $\sen=e_R+\bigoplus_{i\in I_0}K f_i$.

(d)
The group $U_0^-$ acts on $e_R+\ben_0^-$ by the adjoint action,
and the map $U_0^-\times\sen_0\to e_R+\ben_0^-$,
$(g,x)\mapsto(\ad g)(x)$ is an isomorphism of $U_0^-$-varieties.
\endproclaim

\noindent{\sl Proof :}
Part (a) follows from \cite{K3, thm.~8}.
Observe that
the second claim follows also from the following more precise statement.
Any TN element in $\gen_\reg$ is $G$-conjugate into
$e_R+\eps\gen_+$ by \cite{KL2, sec.~4, cor.~1}.
Thus any element in $e_R+\eps\gen_+$ is $G_+$-conjugate into
$e_R+\eps\zen(f_R)_+$ by \cite{B1, lem.~1}.
Hence, any TN element in $\gen_\reg$ is $G$-conjugate into
$e_R+\eps\zen(f_R)_+$.

Part (b) is proved in \cite{K3, thm.~7},
(c) in \cite{D3, p.~185} and \cite{K3, thm.~9},
(d) in the proof of \cite{K3, prop.~19}.

\qed

\vskip3mm

Therefore the composed map $\varphi\eta$ factors to an isomorphism
$\sen_\CC\to\AA^{I_0}_\CC$, denoted by $\varphi\eta$ again.
Further, each $(\ad G)$-orbits in $\gen_\ERS$ contains a unique 
element in $\sen_\ERS$.

\proclaim{1.3.2. Proposition}
(a)
Any element in $\gen^c_\reg$ is HRS of degree $c$.
Any HRS element of degree $c$ in $\gen$
is $G$-conjugate into $\sen^c_\reg$.

(b)
There are HRS elements of degree $c$ in $\gen$ iff $m\in\RN.$
Any HRS element of degree $c$ in $\gen$ has type $W_0[m]$.
Further, it is ERS iff $m\in\EN$, $k>0$.
\endproclaim

\noindent{\sl Proof :}
(a)
Since $\tau$ is not a root of unity, we have
$\gen^c=\gen\cap(\ad\rhov(\eps^{-c}))(\gen_0\otimes\eps^c).$
Thus any element in $\gen^c_\reg$ is HRS of degree $c$.

Any HRS element of degree $c$ in $\gen$ is $G$-conjugate into
a HRS element of degree $c$ in $\sen$ by 1.3.1(a).
We must prove that the latter belong exactly to $\sen^c_\reg.$
The $\bar K^\times$-action on
$\bar\sen$
such that $$t\bullet x=t^{2}\ad(\rhov(t^{-2}))(x)$$
fixes $e_R$ and has the weight $(2d_i)$ with respect to the decomposition
$\bar\sen=e_R+\bigoplus_{i\in I_0}\bar K f_i.$
For each $x\in\bar\sen$ we have
$\varphi_i\eta(t\bullet x)=t^{2d_i}\varphi_i\eta(x)$,
because $\varphi_i\eta$ is $\bar G$-invariant and homogeneous of degree $d_i$.
%see, for instance, \cite{S2, (7.4)},
Since $x\in\bar\sen^c$ iff $F_\tau(x)=\tau^{c/2}\bullet x$,
and since $\varphi\eta$ yields an $F_\tau$-equivariant bijection 
$\bar\sen\to\bar K^{I_0}$, we have
$$\aligned
x\in\bar\sen^c
&\iff
F_\tau\varphi_i\eta(x)=\tau^{cd_i}\varphi_i\eta(x),\ \forall i
\hfill\cr
&\iff
\varphi\eta(x)\in\bigoplus_{i\in I_0}\CC\eps^{cd_i}.
\hfill
\endaligned$$
Since the map $\varphi\eta$ yields a bijection $\sen\to K^{I_0}$, 
we have also 
$$x\in\sen^c\iff\varphi\eta(x)=\bigoplus_{i\in I_m}\CC\eps^{cd_i}.$$
Now, assume that $x\in\sen$ is HRS of degree $c$.
It is $\bar G$-conjugate into $\ten_{0,\reg}\otimes\eps^c$.
Thus, since $\varphi_i\eta$ is a homogeneous polynomial of degree $d_i$,
we have $\varphi_i\eta(x)\in\CC\eps^{cd_i}$ if $i\in I_m$
and $\varphi_i\eta(x)=0$ else.
Thus $x\in\sen^c$.

(b)
Fix an element $x\in\gen^c_\reg$.
Let $w_*$ be its type.
We have $(\ad\rhov(\eps^c))(x)\in\gen_{0,\reg}\otimes\eps^c.$
Thus there is an element $g'\in\bar G$ such that $x'=(\ad g')(x)$
belongs to $\ten_{0,\reg}\otimes\eps^c$.
On the other hand the torus 
$\zen_\gen(x)$ is $\bar G$-conjugated into $\bar\ten^{F_w}$.
Fix $g''\in \bar G$ such that
$x''=(\ad g'')(x)$ belongs to $\bar\ten^{F_w}_{\reg}$.
Since the group $\bar G$ has a split BN-pair, two
elements in $\bar\ten$ which are $\bar G$-conjugate are indeed $W_0$-conjugate.
Thus $x'$, $x''$ are $W_0$-conjugate.
In particular $x''$ lies in
$(\ten_{0,\reg}\otimes\eps^c)^{F_w}$,
yielding $w(x'')=\varpi^cx''$.
So $w\in W_0[m]$ and $m\in\RN$.

Conversely, if $m\in\RN$
there are $w\in W_0[m]$ and $x\in\ten_{0,\reg}$ such that $w(x)=\varpi^cx$.
Then $x\otimes\eps^c\in(\ten_{0,\reg}\otimes\eps^c)^{F_w}$.
So $x\otimes\eps^c$ is HRS of degree $c$ in $\bar\ten^{F_w}$.
Since $\bar\ten^{F_w}$ embeds in $\gen$,
there is also a degree $c$ element in $\gen_\HRS$.

Finally we prove the last claim.
Notice that HS elements of $\gen$ of degree $c$ are TN iff $k>0$.
Further, a maximal torus of type $w_*$
is anisotropic iff no eigenvalue of $w$ equals 1, because the induced action
of the Frobenius on the group of cocharacters is conjugate to $w$.
See \cite{S4, sec. 13.2.2}.
Finally, we have proved that
each degree $c$ element $x\in\gen_\HRS$ has the type $W_0[m]$.
Thus $x$ is ERS iff $m\in\EN$, $k>0$.

\qed

\vskip3mm

\proclaim{1.3.3. Corollary}
(a) We have $\gen^c_\reg\neq\emptyset\iff m\in\RN.$

(b) We have $\gen^c_\ERS\neq\emptyset\iff m\in\EN, k>0
\Rightarrow\gen^c_\ERS=\gen^c_\reg.$
\endproclaim

\vskip3mm

\noindent
{\bf 1.3.4. Remark.}
If $m\in\RN$ we can construct HRS elements of degree $c$ as follows.
Given an $I_0$-tuple $z=(z_i)\in\CC^{I_0}$ we set
$e_{z}=e_R+\sum_iz_if_i\otimes\eps^{cd_i}.$
Then $$\sen^c=\{e_z;z_i=0,\forall i\notin I_m\}.$$
Further $e_z\in\sen^c_\reg$ if the $I_m$-uple $(z_i)_{i\in I_m}$ is generic.
We'll give elements $e_z\in\sen^c_\reg$ 
for each type and each $c$ in section 4.

\head 2. Finite dimensional representations of rational DAHA's\endhead

\subhead 2.1. Reminder on DAHA's \endsubhead

First, let $G_\CC$ be any connected reductive group.
%Recall that $W=W'\rtimes\Omega_W=W_0\ltimes Y_0$.
The DAHA associated to $G_\CC$,
or to the root datum $(X_0,\Delta_0,Y_0,\check\Delta_0)$,
is the unital associative
$\CC_{q,t}$-algebra $\Hb$ generated by the symbols
$t_w,x_\l$ with $w\in W$, $\l\in X_0\times\ZZ$
such that the $t_w$'s satisfy the braid relations of $W$,
plus 
$$\aligned
&x_{(0,1)}=q,\quad x_\mu x_\l=x_{\l+\mu},
\quad (t_{s_i}-t)(t_{s_i}+t^{-1})=0,
\hfill\\
%&t_{s_i}t_{s_j}t_{s_i}\cdots=t_{s_j}t_{s_i}t_{s_j}\cdots\quad
%(\roman{braid\ group\ relations\ for\ }i\neq j),
%\hfill\\
&x_\lambda t_\pi-t_\pi x_\mu=0
\ \roman{if}\ 
\pi\in\Omega_W
\ \roman{and}\ 
{}^{\pi}\!\mu=\lambda,
\hfill\\
&x_\l t_{s_\a}-t_{s_\a}x_{\l-r\a}=
(t-t^{-1})x_\l(1+x_{-\a}+...+x_{-\a}^{r-1})
\ \roman{if\ }
\a\in\Pi
\ \roman{and}\ 
\alphav\cdot\l=r\ge 0.
\hfill
\endaligned$$
%Let $\star$ be the $\CC$-linear anti-involution of $\Hb$ such that
%$t_i\mapsto t_i^{-1}$, $x_\mu\mapsto x_\mu^{-1}$, $q\mapsto q^{-1}$,
%$t\mapsto t^{-1}$.
For each $\a\in\Pi$ we set 
$$\aligned
\psi_{s_\a}&=(t_{s_\a}+t^{-1})(x_{\a}-1)(t^{-1}x_{\a}-t)^{-1}-1
\hfill\\
&=(x_{\a}-1)(tx_{\a}-t^{-1})^{-1}(t_{s_\a}-t)+1.
\endaligned$$
It is the (normalized) intertwiner operator, 
and it belongs to a ring of quotients of $\Hb$.
The $\psi_{s_\a}$'s satisfy the same relations as the 
$s_\a$'s do. 
In particular the element $\psi_w$ is well-defined for each $w\in W$,
and it is invertible with inverse $\psi_{w^{-1}}$.
Further, we have $\psi_w\,x_\l=x_{{}^w\l}\,\psi_w.$
There is a $\CC_t$-algebra automorphism
$$IM:\quad\Hb\to\Hb,
\quad
t_{s_\a}\mapsto -t_{s_\a}^{-1},
\quad
x_\l\mapsto x_\l^{-1},
\quad
q\mapsto q^{-1},
$$
and a $\CC_{q,t}$-algebra anti-automorphism
$$OP:\quad\Hb\to\Hb,
\quad
t_{s_\a}\mapsto t_{s_\a},
\quad
x_\l\mapsto x_\l.$$

Consider the torus $T_d$ whose character group is $X_0\oplus\ZZ\delta$.
Set $H=T_d\times\CC^\times$.
For each $h=(s,\tau,\zeta)\in H$, with $s\in T_0$,
we put $$h^\dag=(s^{-1},\tau^{-1},\zeta),
\quad
h^\ddag=(s,\tau,\zeta^{-1}).$$
In the whole paper 
we'll assume that $\tau,\zeta$ are not root of unity.

For any $\Hb$-module $M$ which is locally finite over $\CC X_0$ we set
$\widehat M=\prod_{h\in H}M_h$ where $M_h$ is the $h$-weight space
$$M_h=\bigcap_{\l\in X^*(H)}\bigcup_{r\ge 0}\{m\in M;
(x_\l-\l(h^\dag))^rm=0\}.$$

Let $\mod(\Hb)_\adm\subset\mod(\Hb)$ be the full subcategory consisting 
of the modules which are locally finite over $\CC X_0$ 
with finite dimensional weight spaces.
Objects in this category are called admissible modules.
If $M$ is admissible we equip $\widehat M$ with the product topology, i.e., 
the coarser topology such that the projections to the weight
subspaces are continuous (each weight subspace being given the
discrete topology).
Since $M=\bigoplus_hM_h$,
the module $M$ is a dense subset of $\widehat M$.

Let $O(\Hb)\subset\mod(\Hb)$ be the full subcategory consisting 
of the finitely generated modules 
which are locally finite over $\CC X_0$ 
and such that $q,t$ act by multiplication by some
complex number. 

Let $O_h(\Hb)\subset O(\Hb)$ be the full subcategory consisting of
the modules $M$ such that $q=\tau^{-1}$, $t=\zeta$, and
$M_{h'}=\{0\}$ if $h'$ is not in the orbit of $h$,
relatively to the $W$-action on $H$ such that
$${}^wh=(s\lambdav(\tau),\tau,\zeta)
\ \roman{if}\ w=\xi_{\lambdav},
\,\roman{and}\ ({}^w s,\tau,\zeta)
\ \roman{if}\ w\in W_0.
\leqno(2.1.1)$$

For a future use let us quote the following basic facts.

\proclaim{2.1.2. Lemma}
(a)
Finite dimensional modules belong to $O(\Hb)$.

(b)
The category $O(\Hb)$ is a Serre subcategory of $\mod(\Hb)$,
and any object in $O(\Hb)$ has a finite length.

(c)
We have $O(\Hb)=\bigoplus_{h}O_h(\Hb)$, where $h$ varies 
in a set of representatives of the $W$-orbits in $H$.
In particular each category $O_h(\Hb)$ is closed under taking subquotients
and extensions,
and any simple module in $O(\Hb)$ belongs to $O_h(\Hb)$ for some $h$.
\endproclaim

\noindent{\sl Proof :}
Part (a) is obvious, and (c) also because any object $M$ in $O(\Hb)$
is the direct sum of its weight subspaces and $\bigoplus_{h\in O}M_h$
is a submodule for each $W$-orbit $O\subset H$.

Now, let us concentrate on part (b).
The first claim is proved as in \cite{VV, prop.~3.1}.
For the second one, recall that $M$ has a finite length iff
it is Artinian and Noetherian, and that if $N\subset M$ then
$M$ is Artinian and Noetherian iff both $N$ and $M/N$ have that
property.
Thus, since objects of $O(\Hb)$ are finitely generated and locally finite
over $\CC[x_\lambda]$, we are reduced to prove that modules induced from
finite dimensional $\Hb_X$-modules (see below for notation) are Artinian and
Noetherian. This is proved as in loc.~cit.~ prop.~3.4.

\qed

\vskip3mm

Let $\Hb_X$ be the $\CC_t$-subalgebra generated by
$t_w,x_{\l}$ with $w\in W_0$, $\l\in X_0$.
It is isomorphic to the Iwahori-Hecke algebra of
the group $\check G_\CC$.
Let $\sgnb,\trivb$ be the Steinberg and trivial representations of $\Hb_X$,
i.e., the projective one-dimensional $\CC_t$-modules
such that $x_{\l}$ acts by the scalars $t^{-2\l\cdot\rhov}$,
$t^{2\l\cdot\rhov}$, $t_{s_i}$ by the scalars $-t^{-1}$, $t$,
and $t_\pi$ by the scalar 1.

There is a $\CC_{q,t}$-algebra embedding
$\CC_{q,t}Y_0\to\Hb$ taking
$\lambdav\in Y_0$ to $x_{\lambdav}=t_{\xi_{\lambdav_1}}t_{\xi_{\lambdav_2}}^{-1}$,
with $\lambdav=\lambdav_1-\lambdav_2$ and
$\lambdav_1,\lambdav_2\in Y_0$ dominant.
The multiplication in $\Hb$ yields an isomorphism
$\CC_{q,t}Y_0\otimes_{\CC_t}\Hb_X\to\Hb.$
See \cite{C4, thm.~3.2.1} for details.
%Let $\trivb$ be the trivial representation of $\Hb_X$,
%i.e., the projective one dimensional $\CC_t$-module
%such that $x_\l$, $t_{s_i}$ act by the scalars
%$t^{2\l\cdot\rhov}$, $t$
%respectively.
The induced module $\Pb=\Pb_Y=\Hb\otimes_{\Hb_X}\trivb$
is identified with $\CC_{q,t}Y_0$ as a $\CC_{q,t}$-module.
It is called the polynomial representation.
In the same way, setting $\Hb_Y$ equal to the $\CC_t$-subalgebra
generated by $t_w$, $w\in W$, we get the
induced module $\Pb_X=\Hb\otimes_{\Hb_Y}\trivb$.
Recall that $c\in\QQ^\times$.

A $\Hb$-module is ($X$-)spherical if it is a quotient of $\Pb$,
and $Y$-spherical if it is a quotient of $\Pb_{X}$.
We write $\Hb_c$, $\Pb_c$, $\Pb_{X,c}$
for the specialization of $\Hb$, $\Pb$, $\Pb_X$ both at the ideals
$(t-q^{-c/2})$ and $(q-\tau^{-1},t-\tau^{c/2})$,
hopping it will not create any confusion. 
Put $\Pb_{X}^\dag={}^{IM}\Pb_{X}$,
$\Pb^\dag={}^{IM}\Pb$,
$\Pb_{X,c}^\dag={}^{IM}\Pb_{X,-c}$, and 
$\Pb_{c}^\dag={}^{IM}\Pb_{-c}$.
%and $\Hb_{X,c}$ for the specialization of $\Hb_X$
%at the maximal ideal $(t-\tau^{c/2})\subset\CC_t$.
If $\zeta=\tau^{c/2}$ we may write $O_{\tilde s}(\Hb_c)$ for $O_h(\Hb)$,
with $\tilde s=(s,\tau)$.

\vskip3mm

\noindent{\bf 2.1.3. Examples.}
(a)
The element $x_\l$ acts as $t^{2\rhov\cdot\l}$ on $\trivb$.
Thus it acts also as $t^{2\rhov\cdot\l}$ on the obvious
generator of $\Pb$.
Hence $x_\l,q,t$ act on the obvious generator of the
$\Hb_{c}$-module $\Pb_{c}$ as 
$\tau^{c\rhov\cdot\l},\tau^{-1},\tau^{c/2}$ respectively.
It is generated by its $(h_c)^{-1}$-weight subspace,
where $$h_c=(\rhov(\tau^{c}),\tau^{-1},\tau^{-c/2}).$$

(b)
The $\Hb_c$-module $\Pb_{c}^\dag$
is generated by its $(h_c^\dag)^{-1}$-weight subspace,
where $$h^\dag_{c}=(\rhov(\tau^{c}),\tau,\tau^{c/2}).$$
The notation for $h_c^\dag$ may be confusing because we have
$h^\dag_{c}=(h_{-c})^\dag\neq(h_{c})^\dag$.

\vskip3mm

Since the choice of the non-zero complex number $\tau$ 
is indifferent, under substituting $\tau$ by $\tau^{-1}$ everywhere
we may identify
$h_c$, $h^\dag_{c}$ with 
$(h_c)^{-1}$, $(h^\dag_{c})^{-1}$ respectively. 
In particular this yields an identification of the elements
$h^\dag_{c}$, ${}^{w_0}h_{c}$,
and of the categories $O_{h^\dag_c}(\Hb)$, $O_{h_c}(\Hb)$. 
To simplify notations we'll also write $c$ for $h_c^\dag$ 
if it creates no confusion.
For instance, we write $O_c(\Hb)$ for both
$O_{h_c^\dag}(\Hb)$ and $O_{h_c}(\Hb)$.

\proclaim{2.1.4. Lemma}
The $\Hb_c$-modules $\Pb_{X,c}$, $\Pb_{X,c}^\dag$ are not admissible.
The $\Hb_c$-modules $\Pb_c$, $\Pb_c^\dag$ belong to $O_c(\Hb)$.
\endproclaim

Let $J^h,I^h\subset\CC_{q,t}X_0$ be the ideals 
$$J^h=\Bigl(f-f(h^\dag);f\in\CC_{q,t}X_0\Bigr),
\quad
I^h=\Bigl(J^h\cap\CC_{q,t}X_0^{Z_W(h)}\Bigr).$$
For each $E\subset{}^Wh$ we set $I^E=\bigcap_{h'\in E}I^{h'}$.
Notice that ${}^w(I^E)=I^{{}^wE}$.

\proclaim{2.1.5. Lemma}
(a)
If $E\subset{}^Wh$ is finite, there is a finite subset
$E'\subset{}^Wh$ such that $E\subset E'$ and
$I^{E'}t_{s_\a}\subset t_{s_\a}I^{E}+I^{E}$,
$t_{s_\a}I^{E'}\subset I^{E}t_{s_\a}+I^{E}$ 
in $\Hb$ for each $\a\in\Pi$.

(b)
We have $O(\Hb)\subset\mod(\Hb)_\adm.$ 

(c)
The $\Hb$-action 
on an admissible module $M$ 
extends uniquely to
a continuous action on $\widehat M$. 
\endproclaim

\noindent{\sl Proof :}
Then we prove (a). The proof is modelled on \cite{VV, lem.~3.2}.
We claim that the set $E'=E\cup\bigcup_{\a\in\Pi}{}^{s_\a}E$ is a solution.
To prove this it is enough to check that if $F=E\cup{}^{s_\a}E$ then 
$I^{F}t_{s_\a}\subset t_{s_\a}I^{F}+I^{F}$ in $\Hb$.
For each $p\in\CC X_0$ we have 
$$pt_{s_\a}-t_{s_\a}{}^{s_\a}p=(t-t^{-1})\vartheta_\a(p),$$
in $\Hb$, where $\vartheta_\a(p)=(p-{}^{s_\a}p)/(1-x_{-\a})$.
Hence it is enough to prove that
$\vartheta_\a(I^{F})\subset I^{F}$, i.e., that 
$\vartheta_\a(I^{h'}\cap I^{{}^{s_\a}h'})\subset I^{h'}\cap I^{{}^{s_\a}h'}$
for each $h'\in E$.
If $h'={}^{s_\a}h'$ this is obvious, because
$$\vartheta_\a(pp')=p\,\vartheta_\a(p')+\vartheta_\a(p)\,{}^{s_\a}p',
\quad
\vartheta_\a(\{f-f(h');f\in(\CC X_0)^{Z_W(h')}\})=0.$$
Assume now that $h'\neq{}^{s_\a}h'$.
We must check that $\vartheta_\a(p)\in I^{h'}$
for each $p\in I^{h'}\cap I^{{}^{s_\a}h'}.$
We have $(1-x_{-\a})\vartheta_\a(p)=p-{}^{s_\a}p\in I^{h'}.$
So, since $\CC X_0/I^{h'}$ is a local ring 
and $(1-x_{-\a})$ is not in the maximal
ideal, we get that $\vartheta_\a(p)\in I^{h'}.$

Now we prove (b).
%The first claim is \cite{VV1, prop.~3.8(ii)}.
It is enough to check that the weights spaces of
the module $M=\Hb/\Hb(J^h)^r$ are finite dimensional for all $r\ge 0$.
The PBW theorem for $\Hb$ yields an isomorphism
$M=\bigoplus_{v\in W}t_v\CC_{q,t}X_0/(J^h)^r$.
For each $w\in W$ we set 
$M_{\le w}=\bigoplus_{v\le w}t_v\CC_{q,t}X_0/(J^h)^r$.
Since
$$x_\l t_w\in t_wx_{{}^w\l}\oplus\bigoplus_{v<w}t_v\CC_{q,t}X_0$$ in $\Hb$
for each $\l$, the subspace $M_{\le w}$ is preserved by left multiplication
by elements of $\CC_{q,t}X_0$.
Since $\dim(M_{\le w})<\infty$ the quotient
$M_w=M_{\le w}/\bigoplus_{v<w}M_{\le v}$ is finite dimensional
and supported on $\{{}^wh^\dag\}$ as a $\CC_{q,t}X_0$-module.
In the Grothendieck group of $\CC_{q,t}X_0$-modules we have
$M=\bigoplus_{w\in W}M_w$.
Hence the dimension of $M_{h'}$
is equal to $\sum_{{}^wh=h'}\dim M_w$.
Therefore $\dim M_{h'}<\infty$, because $Z_W(h)$ is finite.

Finally we prove (c).
We must prove that the action of any element $x\in\Hb$ extends
to a continuous operator $\widehat M\to\widehat M$.
The unicity is obvious because $M$ is a dense subset of $\widehat M$.
Given a finite subset $E\subset{}^Wh$ 
we must check that there is
a finite subset $E'\subset{}^Wh$ such that 
$xM_{E'}\subset M_E$, where
$M_E=\prod_{h\notin E}M_h$.
We have
$(I^h)^rx=0$ if $x\in M_h$ and $r\gg 0$.
Further $I^hM_{h'}=M_{h'}$ if $h\neq h'$. 
Therefore, since $M$ has finite dimensional
weight subspaces and $E$ is finite
there is an integer $r_E$ such that 
$M_E=(I^E)^{r}M$ for each $r\ge r_E$.
Let $E'$ be as in part (a) above.
We have $xI^{E'}\subset I^Ex+I^E$
and $I^{E'}\subset I^E$.
Thus
$$xM_{E'}=x(I^{E'})^{r}M\subset (I^E)^{r}M=M_E$$
if $r\ge r_E,r_{E'}$.
We are done.

\qed

\vskip3mm

In the rest of the section let $G_\CC$ be a Chevalley group.
Thus $W$ is the affine Weyl group of the root system $\Delta_0$.
Let $\la\tau,\zeta^2\ra\subset\CC^\times$ be the subgroup generated
by $\tau$, $\zeta^2$.
For each $h$ the set
$$\Delta_{0,(h)}=\{a\in\Delta_0;a(s)\in\la\tau,\zeta^2\ra\}$$
is a root system.
Let $\Delta_{\re,(h)}=\Delta_{0,(h)}\times\ZZ$ be the set of affine real roots,
$\Delta_{(h)}=\Delta_{\re,(h)}\cup(\{0\}\times\ZZ)$ the set of affine roots,
$\check\Delta_{\re,(h)}$ the set of affine real coroots,
and $W_{0,(h)}$, $W_{(h)}$ the Weyl group and the affine Weyl group associated 
to the root system $\Delta_{0,(h)}$.
Thus $W_{(h)}$ is a subgroup of $W_{0,(h)}\ltimes Y_0$ which is
equal to the subgroup of $W=W_0\ltimes Y_0$
generated by the affine reflections $s_\a$ with $\a\in\Delta_{\re,(h)}$.
See section 1.

Let $\Hb_{(h)}$ be the subalgebra of the DAHA associated to
$(X_0,\Delta_{0,(h)},Y_0,\check\Delta_{0,(h)})$ 
which is generated by the elements
$x_\l,$ $t_w$ with $\l\in X_0$, $w\in W_{(h)}$.
Let $O_h(\Hb_{(h)})$, $O(\Hb_{(h)})$, $\mod(\Hb_{(h)})$ be the corresponding
categories of modules.

\vskip3mm

\proclaim{2.1.6. Proposition}
Assume that $Z_W(h)$ is a reflection group.

(a)
There is an equivalence of categories
$O_h(\Hb)\to O_h(\Hb_{(h)})$.

(b)
It preserves the finite dimensional modules iff $W/W_{(h)}$ is finite.
It preserves the dimension of finite dimensional modules iff $W=W_{(h)}$.
\endproclaim

Before the proof we introduce one more notation.
For each subset $F\subset {}^Wh$ and each $\CC_{q,t}X_0$-module $M$ we set
$${}^FM=\pro_{E}\pro_nM/(I^E)^nM,$$
where $E$ is any finite subset of $F$.
To simplify we'll write
${}^hM$ for ${}^FM$ if $F={}^Wh$.
We equip ${}^FM$ with the limit topology, i.e.,
the coarsest topology which makes the obvious maps 
${}^FM\to M/(I^E)^nM$ continuous.
Since $\CC_{q,t}X_0=\CC[H]$ we may write
${}^F\CC[H]$ for ${}^F\CC_{q,t}X_0$.
It is a topological algebra 
equal to the product of the topological algebras 
$${}^{\{h'\}}\CC_{q,t}X_0=\pro_n\CC_{q,t}X_0/(I^{h'})^n$$
with $h'\in F$.
The topology is the product topology, where each factor
is equipped the $I^{h'}$-adic topology.

A topological $\CC$-algebra $\Ab$ is an algebra with a $\CC$-linear topology, 
i.e.,
a topology for which there is a basis of neighborhoods of 0
consistings of vector subspaces, 
such that the multiplication $\Ab\times\Ab\to\Ab$ is jointly continuous.
We do not impose the topology to be $\Ab$-linear.
A topological module over $\Ab$
is a module $M$ with a $\CC$-linear topology
%such that for each open submodule $N\subset M$
%the ideal $\{a\in\Ab; aM\subset N\}$ is open.
such that the action map $\Ab\times M\to M$ is jointly continuous.
So, if $\Ab$ is discrete this means that the action of any $a\in\Ab$
is a continuous operator $M\to M$, while if $M$ is discrete this means that
the annihilator of any element of $M$ is an open subset of $\Ab$. 

\proclaim{2.1.7. Lemma}
(a)
For each $F$ we have $\Hb\subset{}^F\Hb$ 
and there is an unique topological algebra structure on ${}^F\Hb$
extending the multiplication on $\Hb$.
The assignement $M\mapsto{}^FM$ yields a functor
from the category of $\Hb$-modules to the category of 
continuous ${}^F\Hb$-modules. 

(b)
We have ${}^hM=\widehat M$ if $M\in O_h(\Hb)$.
\endproclaim

\noindent{\sl Proof :}
We have $\Hb\subset{}^F\Hb$ by the PBW theorem for $\Hb$.
To prove the second claim we must check that for each $E,n$ and each $x\in\Hb$
there are $E',n'$ such that
$x(I^{E'})^{n'}\Hb\subset (I^{E})^{n}\Hb$.
Let $E'$ be as in 2.1.5(a) and $n'=n$. Then we have
$$t_{s_i}(I^{E'})^{n'}\Hb\subset I^Et_{s_i}(I^{E'})^{n-1}\Hb+
(I^{E})^{n}\Hb\subset\cdots\subset (I^{E})^{n}\Hb,\quad\forall i.$$
We are done.
The rest of (a) is left to the reader.
Compare 2.1.5(c).

Part (b) is obvious, because $I^hM_{h'}=M_{h'}$ if $h\neq h'$ and
$(I^h)^nM_h=\{0\}$ if $n\gg 0$, because $\dim(M_h)$ is finite by 2.1.5(a).

\qed

\vskip3mm

\noindent{\sl Proof of 2.1.6 :}
First, we prove (a).
The proof is modelled after \cite{L1},
and was sketched in \cite{VV, sec.~3.5}.
For the comfort of the reader, let us recall a few facts.

We have ${}^Wh=\bigsqcup_{h'\in\Pc_h}(h')$, where
$(h')={}^{W_{(h')}}h'$, for some subset $\Pc_h\subset{}^Wh$,
see \cite{VV, p.~1309}.
The $W$-action on ${}^Wh$ factors to a $W$-action on $\Pc_h$.
Using the fact that $Z_W(h)$ is generated by reflections,
one gets $\Pc_h\simeq W/W_{(h)}$ as $W$-sets. 
See loc.~cit.
In particular, we have $|\Pc_h|=1$ iff $W=W_{(h)}$.

Since ${}^h\CC_{q,t} X_0$ is isomorphic to the product
$\prod_{h'\in\Pc_h}{}^{(h')}\CC_{q,t} X_0,$ 
the identity in ${}^{(h')}\CC_{q,t} X_0$ 
may be viewed as an indempotent $1_{(h')}$ in ${}^h\CC_{q,t} X_0$.
For all $(h')\in\Pc_h$,
one can check as in \cite{L1, sec.~8.8} that if 
$(h')\neq{}^{s_i}(h')$
then we have
the following equality in $1_{{}^{s_i}(h')}{}^h\Hb 1_{(h')}$ 
$$1_{{}^{s_i}(h')}\psi_{s_i}=\psi_{s_i}1_{(h')}.$$
For all $(h')$ we fix an element $w_{(h')}\in W$ 
with a reduced decomposition 
$w_{(h')}=s_{i_1}s_{i_2}...s_{i_p}$ 
such that $$(h)\neq {}^{s_{i_p}}(h)\neq...\neq
{}^{s_{i_2}...s_{i_p}}(h)\neq{}^{w_{(h')}}(h)=(h').$$
We have
$1_{(h')}\psi_{w_{(h')}}\in{}^h\Hb\, 1_{(h)}$.
Thus each element in
$1_{(h')}{}^h\Hb\, 1_{(h'')}$
is of the form
$\psi_{w_{(h')}}y\psi_{w_{(h'')}}^{-1}$
for some $y\in 1_{(h)}{}^h\Hb\,1_{(h)}$.

Let
$A:{}^h\Hb\to\roman{Mat}_{\Pc_h}(1_{(h)}{}^h\Hb\,1_{(h)})$
be the unique topological algebra homomorphism such that
$\psi_{w_{(h')}}y\psi_{w_{(h'')}}^{-1}\mapsto E_{(h'),(h'')}(y)$,
where $E_{(h'),(h'')}(y)$ is the matrix with $y$ in the entry
$(h'),(h'')$ and 0 elsewhere.
The image of $A$ is dense.

There is an unique topological algebra isomorphism
$B:{}^h\Hb_{(h)}\to 1_{(h)}{}^h\Hb\, 1_{(h)}$
such that $x_\l\mapsto x_\l$, $t_w\mapsto t_w^{(h)}$, where
$\l\in X_0$, $w\in W_{(h)}$, and $t_w^{(h)}$ is as in \cite{L1, sec.~8.10}.

Composing $A$, $B$, and the Morita equivalence we get the required 
equivalence of categories $O_h(\Hb)\to O_h(\Hb_{(h)})$.

Now, claim (b) is obvious from the proof of (a).

\qed

\vskip3mm

For each $h\in H$ we may write $\Delta_h$ for the root system 
$$\Delta_{\tilde s}=\{\a\in\Delta;\alpha(\tilde s)=1\}$$
(it is independant on $\zeta$).
Notice that $\Delta_h\subset\Delta_\re$, because $\tau$ is not a root of unity.
%Since $\tau$ is not a root of unity, the connected component
%$G^{h,\circ}$ of the fixed point subset $G^h$ is the set
%of $\CC$-points of a connected reductive group.
%See 2.4.1(a) below.
Let $W_h$ be the Weyl group of $\Delta_h$.
It is a subgroup of $Z_W(h)$.

Recall that $X$ is the lattice $X_0\oplus\ZZ\delta\oplus\ZZ\omega_0$.
Consider the tori $\hat T_0$, $\tilde T$ whose character groups are
$X_0\oplus\ZZ\omega_0$, $X$ respectively.
Set $\hat H=\tilde T\times\CC^\times$.
The action 2.1.1 lifts to a $W$-action on $\hat H$ such that,
for each $h=(s,\tau,\zeta)$, we have
$$\aligned
&{}^wh=
(s\lambdav(\tau)\check\delta(\tau\lambda(s))^{-\lambdav\cdot\lambdav/2},
\tau,\zeta)
\ \roman{if}\ w=\xi_{\lambdav},\hfill\\
&{}^wh=({}^w s,\tau,\zeta)
\ \roman{if}\ w\in W_0.
\endaligned
$$
Let $\hat\Hb$ be the DAHA considered in \cite{V1}.
It is the algebra generated by $\Hb_Y$ and $x_\l$, $\l\in X$,
modulo the same relations as above. 
Hence the element $x_{\o_0}$ commutes with $\Hb_X$ and we have 
$$x_{\o_0} t_{s_0}x_{\o_0}^{-1}=x_{\a_0}t_{s_0}^{-1}.$$
So the algebra $\hat\Hb$ is a semi-direct product
$\CC[x_{\o_0}^{\pm 1}]\ltimes\Hb$.
The categories $O(\hat\Hb)$, $O_h(\hat\Hb)$, with $h\in\hat H$, 
are defined in the obvious way.
We will use the same notation for the polynomial representations of
the algebras $\Hb$ and $\hat\Hb$, hopping it will not create any confusion.
The equivalence 2.1.6 has an obvious analogue for $\hat\Hb$.
Notice that the assumption in 2.1.6 becomes void for $\hat\Hb$,
because $Z_W(h)$ is a reflection group by 2.1.8(b) below.
Indeed, if $h\in\hat H$ we define $\Delta_h, W_h$ as above
(with $\hat H$ instead of $H$).
Then $W_h=Z_W(h)$. See below.

\vskip3mm

\noindent{\bf 2.1.8. Remarks.}
(a)
A finite subgroup of $W$ is conjugated into a parabolic subgroup.
This is well-known.
We recall the proof for the comfort of the reader.
Fix a triple $(\lambdav,\ell,z)$ in the Tits cone, i.e., such that $\ell>0$.
If $W'\subset W$ is a finite subgroup the sum
$\muv=\sum_{w\in W'}{}^w(\lambdav,\ell,z)$
is again in the Tits cone. Thus there is an element $w\in W$
such that ${}^w\muv$ belongs to the closed fundamental alcove 
$\bar A=\{\check\nu\in\check V_{0,\RR};
a_i\cdot\check\nu\ge 0,\theta\cdot\check\nu\le 1\}$.
Thus $W'$ is conjugated into the parabolic subgroup
$Z_W({}^w\muv)$. See \cite{K2, prop.~3.12}.

(b)
If $h\in\hat H$ then we have $W_h=Z_W(h)$ (a reflection group).
Indeed, the element $\xi_{\lambdav}w$ belongs to $Z_W(h)$ iff 
$\lambdav(\tau)={}^wss^{-1}$.
Thus the map $\xi_{\lambdav}w\mapsto w$ yields a group isomorphism
$$Z_W(h)\to\{w\in W_0; {}^wss^{-1}\in Y_0(\tau)\},$$
where $Y_0(\tau)=\{\lambdav(\tau);\lambdav\in Y_0\}$.
So the group $Z_W(h)$ is finite.
Hence it is conjugate into the proper
parabolic subgroup $W_J$ generated by $\{s_i;i\in J\}$
for some subset $J\subsetneq I$ by part (a).
So we may assume that
$$W_h\subset Z_W(h)\subset W_J.$$
Thus we are reduced to check that $(W_J)_h=Z_{W_J}(h)$
(with the same notation for $(W_J)_h$ as for $W_h$ above).
This is a standard result because $\tilde T$ is the maximal torus
of a connected reductive group with simply connected derived subgroup 
whose Weyl group is isomorphic to $W_J$ (a Levi subgroup in $\tilde G$).

Alternatively, the axioms of a split BN pair imply that the group
$\pi_0(\tilde G^h)$ is isomorphic to $Z_W(h)/W_h$.
See \cite{S6, rem.~8.3(c)} for the finite type case.
Then, apply 2.4.1(a) below.
%Since $W_J$ is a finite Coxeter group, by \cite{S6, 4.4} the 
%group $Z_{W_J}(h)$ is the product of the group generated by the reflections
%in $W_J$ which fix $h$ and a finite subgroup in 
%$$Y_0\oplus\ZZ\check\omega_0\oplus\ZZ\check\delta/\sum_{i\in J}\ZZ\alphav_i.$$
%We are done, because this quotient is torsion free since $G_\CC$ is simply
%connected.

\subhead 2.2. Reminder on degenerate DAHA's \endsubhead

The degenerate DAHA associated to the Chevalley group $G_\CC$ is
the unital associative $\CC_\kappa$-algebra $\Hb'$ generated by
the group algebra $\CC W$ and symbols $\xi_\l$, $\l\in X_0\times\ZZ$,
with defining relations given by
$$
\xi_{(0,1)}=1,\quad \xi_\l+\xi_\mu=\xi_{\l+\mu},\quad
\xi_\l\xi_\mu=\xi_\mu \xi_\l,\quad
\xi_\l s_i-s_i\xi_{{}^{s_i}\l}=-\kappa \l\cdot\alphav_i.
$$
For each $i\in I$ we set
$$\psi'_{s_i}=(s_i+1)\xi_{\a_i}(\xi_{\a_i}-\kappa)^{-1}-1
=\xi_{\a_i}(\xi_{\a_i}+\kappa)^{-1}(s_i-1)+1.$$
It is the (normalized) intertwiner operator, 
and it belongs to a ring of quotients of $\Hb'$.
The $\psi'_{s_i}$'s satisfy the same relations as the 
$s_i$'s do. 
In particular the element $\psi'_w$ is well-defined for each $w\in W$.
Recall that $\psi'_w\,\xi_\l=\xi_{{}^w\l}\,\psi'_w.$
There is a $\CC$-algebra automorphism
$$IM:\quad\Hb'\to\Hb',
\quad
s_i\mapsto -s_i,
\quad
\xi_\l\mapsto\xi_\l,
\quad
\kappa\mapsto-\kappa,$$
and a $\CC_\kappa$-algebra anti-automorphism
$$OP:\quad\Hb'\to\Hb',
\quad
s_i\mapsto s_i,
\quad
\xi_\l\mapsto\xi_\l.
$$

For each $M\in\mod(\Hb'_c)$ and 
$\lambdav=(\muv,m)$ with $\muv\in\check X_0$ and $m\neq 0$, we set
$$M_{\lambdav}=
\bigcap_{\l\in X_0\times\ZZ}\bigcup_{r\ge 0}
\{x\in M;(\xi_\l-\l\cdot\lambdav/m)^rx=0\}.$$

Let $O(\Hb')\subset\mod(\Hb')$ be the full subcategory consisting 
of the finitely generated modules 
which are locally finite over $\CC[\xi_\l]$ and such that
$\kappa$ acts by multiplication by some complex number.
Let $O_\lambdav(\Hb')\subset O(\Hb')$ be the full subcategory consisting of
the modules $M$ such that
$M_{\lambdav'}=\{0\}$ if $\lambdav'$ is not in the
$W$-orbit of $\lambdav$.

\proclaim{2.2.1. Lemma}
(a)
Finite dimensional modules belong to $O(\Hb')$.

(b)
We have $M=\bigoplus_\lambdav M_\lambdav$ for each $M\in O(\Hb')$.

(c)
We have $O(\Hb')=\bigoplus_\lambdav O_\lambdav(\Hb').$ 
\endproclaim

The $\CC$-subalgebra $\CC[\xi_\l]$ generated by the $\xi_\l$'s is isomorphic
to $\CC[\check V_0]$.
The $\CC$-subalgebra $\Hb'_X$ generated by $\CC W_0$, $\CC[\xi_\l]$
is isomorphic to the degenerate affine Hecke algebra.
Put $\rhov_\kappa=(\kappa\rhov,1).$
Let $\trivb$ denote both the trivial $W$-module and
the $\Hb'_X$-module such that
$\xi_\l$, $s_i$ act by the scalars $\l\cdot\rhov_\kappa$, $1$ respectively.
The induced $\Hb'$-module $\Pb'=\Pb'_Y=\Hb'\otimes_{\Hb'_X}\trivb$
is identified with $\CC_\kappa Y_0$ as a $\CC_\kappa$-module.
We set also $\Pb'_X=\Hb'\otimes_{\CC W}\trivb$.
We write $\Hb'_c$, $\Pb'_c$, $\Pb'_{X,c}$
for the specialization of $\Hb'$, $\Pb'$, $\Pb'_X$
at the maximal ideal $(\kappa-c)$.
A $\Hb'_c$-module is ($X$-)spherical if is a quotient of $\Pb'_c,$
and $Y$-spherical if is a quotient of $\Pb'_{X,c}$.
Write $O(\Hb'_c)$, $O_{\lambdav}(\Hb'_c)$
for the full subcategory of $O(\Hb')$,
$O_\lambdav(\Hb')$ consisting of
$\Hb'_c$-modules.

We'll write $c$ for $\rhov_c$ and $c\rhov$ if it creates no confusions.

\proclaim{2.2.2. Lemma}
The modules $\Pb'_c$, 
${\Pb'_{c}}^\1dag={}^{IM}\Pb'_{-c}$
belong to the category $O_c(\Hb'_c)$.
\endproclaim

For each subset $F\subset{}^W\lambdav$ and each $\CC_k[\xi_\l]$-module $M$
we define ${}^FM$ as in 2.1.7(a). 
In particular, we get a topological algebra
${}^F\Hb'$ and, for each $\Hb'$-module $M$, a continuous ${}^F\Hb'$-module
${}^FM$.

\proclaim{2.2.3. Lemma}
Assume that $\zeta=\tau^{c/2}$
and $\Delta_0=\Delta_{0,(h)}$.

(a)
There is an element $\muv\in\check X_0$ such that
$\muv(\tau_m)=s\, \mod\, Z(G_0)$.

(b)
Put $\lambdav=(\muv,m)$ where $\muv$ is as above.
There is an isomorphism of topological algebras
${}^{\lambdav}\CC[\check V_0\times\CC]\to{}^{h}\CC[H]/(t-q^{-c/2})$.
\endproclaim

\noindent{\sl Proof :}
Part (a) is obvious, because $\la\tau,\zeta^2\ra=\la\tau_m\ra$. 
Let us concentrate on (b).
Set $\nuv=(\muv,m,k/2)$.
We have $Z_W(h)=Z_W(\nuv)$.
Here the $W$-action on $\nuv$ is as in 1.0.2
(it is trivial on the third component).
Compare 2.1.1.
So there is an unique isomorphism of $W$-sets 
$\iota:{}^W\nuv\to{}^Wh$
such that $\nuv\mapsto h$.
Let $\een:\CC\to\CC^\times$ be the exponential map,
normalized so that $\een(1)=\tau_m$.
It yields a $W$-equivariant group homomorphism
$$\een:\ \check V_0\times\CC^2=(Y_0\times\ZZ^2)\otimes\CC\to H,\ 
(\lambdav',m',k')\otimes z\mapsto(\lambdav'(\een z),\een (m'z),\een(k'z)).$$
For each $\nuv'\in{}^W\nuv$ set $h_0=\iota(\nuv')\een(-\nuv')$.
The element $h_0\in H$ is indendent of the choice of $\nuv'$ because
$\iota$, $\een$ are $W$-equivariant.
A direct computation shows that $h_0=(s_0,1,1)$ with
$a(s_0)=1$ for each $a\in\Delta_0$, i.e., with $s_0\in Z(G_0)$.
So the holomorphic map 
$$\Psi:\check V_0\times\CC^2\to H,\quad x\mapsto h_0\een(x)$$
is $W$-equivariant and takes $\nuv$ to $h$.
We claim that the comorphism of $\Psi$ yields
an isomorphism of topological algebras
${}^{h}\CC[H]/(t-q^{-c/2})\to{}^{\lambdav}\CC[\check V_0\times\CC]$.

To prove this, observe first that for each $\nuv'\in{}^W\nuv$,
setting $h'=\Psi(\nuv')$, the comporphism of $\Psi$ yields  
a topological algebra isomorphism of the completions of the algebras
$\CC[\check V_0\times\CC^2]$, $\CC[H]$ at the points $\muv'$, $h'$ respectively.
Thus it yields also an isomorphism
$$
{}^{h}\CC[H]/(t-q^{-c/2})=
\prod_{h'\in {}^Wh}
{}^{\{h'\}}\CC[H]/(t-q^{-c/2})\to 
\prod_{\lambdav'\in {}^W\lambdav}
{}^{\{\lambdav'\}}\CC[\check V_0\times\CC]
={}^{\lambdav}\CC[\check V_0\times\CC]
.$$
We are done.

\qed

\proclaim{2.2.4. Proposition}
Assume that $\zeta=\tau^{c/2}$ and $\Delta_0=\Delta_{0,(h)}$.
Fix $\lambdav$ as in 2.2.3(a). 

(a)
There is an equivalence of categories
$O_{\tilde s}(\Hb_c)\to O_{\lambdav}(\Hb'_c)$
which preserves the dimension of finite dimensional modules.

(b)
It yields an equivalence
from the full subcategories of spherical, $Y$-spherical modules to the
full subcategories of spherical, $Y$-spherical ones respectively.
\endproclaim

\noindent{\sl Proof :}
First we prove (a).
The proof is modelled after \cite{L1},
and was sketched in \cite{VV, sec.~3.5}.
Let $\muv$ be as above.
There is an isomorphism of topological algebras
$$C:{}^{\lambdav}\Hb'_c\to{}^{h}\Hb_c.$$ 
It is uniquely determined by the following properties 

\vskip1mm

\itemitem{(i)}
it takes $\psi'_{s_i}$ to $\psi_{s_i}$ for each $i\in I$,

\vskip1mm

\itemitem{(ii)}
it restricts to the topological algebra isomorphism in 2.2.3
$${}^{\lambdav}\CC[\check V_0\times\CC]\to{}^{h}\CC_{q,t} X_0/(t-q^{-c/2}).$$

\vskip1mm

\noindent
Observe that 
%the element $x_{\a_i}^{1/2}$ is well defined in ${}^h\CC_{q,t}$,
%because $x_{\a_i}\in\a_i(h)+J^h$ and $\a_i(h)\neq 0$.
%However
$\psi_{s_i}\notin{}^h\Hb$ and $\psi'_{s_i}\notin{}^\lambdav\Hb'$.
Claim (i) should be understood in the following way :
we have
$$C(s_i+1)=(t_{s_i}+t^{-1})g_{+,i},
\quad
C(s_i-1)=g_{-,i}(t_{s_i}-t),
$$
where 
$$g_{+,i}={x_{\a_i}-1\over t^{-1}x_{\a_i}-t}\,
C\biggl({\xi_{\a_i}-\kappa\over \xi_{\a_i}}\biggr),
\quad
g_{-,i}={x_{\a_i}-1\over tx_{\a_i}-t^{-1}}\,
C\biggl({\xi_{\a_i}+\kappa\over \xi_{\a_i}}\biggr).$$ 
Indeed, one check easily as in \cite{L1, lem.~9.5} that $g_{\pm,i}$
is a well-defined unit in ${}^h\CC[H]$ 
(although both factors of $g_{\pm,i}$ are not in general).
Write $\mod({}^{h}\Hb)\to\mod({}^{\lambdav}\Hb')$, $M\mapsto{}^CM$
for the twist by $C$.
It restricts to an equivalence
$O_{\tilde s}(\Hb_c)\to O_{\lambdav}(\Hb'_c)$ which
preserves the dimension of the finite dimensional modules.

Part (b) is obvious. 

\qed

\vskip3mm

For a future use let 
$\mod(\Hb'_c)_{uni}\subset\mod(\Hb'_c)$,
$O(\Hb')_{uni}\subset O(\Hb')$
be the full subcategories consisting of the modules which are locally
unipotent as $Y_0$-modules.

\vskip3mm

\noindent{\bf 2.2.5. Examples.}
(a)
If $c=1/h$
the one dimensional representations such that
$$\aligned
\Lb_c^\dag:&
\quad
t\mapsto\tau^{c/2},
\quad
q\mapsto\tau^{-1},
\quad
x_\l\mapsto\tau^{-\rhov_c\cdot\l},
\quad
t_i\mapsto -t^{-1},
\hfill\cr
{\Lb'}_c^\dag:&
\quad
\kappa\mapsto c,
\quad
\xi_\l\mapsto \rhov_c\cdot\l,
\quad
s_i\mapsto -1
\hfill
\endaligned$$
correspond to each other under 2.2.4.
We set $\Lb_{-c}={}^{IM}\Lb_c^\dag$ and $\Lb'_{-c}={}^{IM}{\Lb'}_c^\dag$.
The $\Hb_{-c}$-modules $\Lb_{-c},$ $\Lb'_{-c}$ are spherical and $Y$-spherical.

(b)
If $(G_\CC,c)=(C_2,1/2)$
there are one dimensional $\Hb_c$-modules
such that
$$\aligned
\Sb_c^\dag:&
\quad
t\mapsto\tau^{c/2},
\quad
q\mapsto\tau^{-1},
\quad
x_\l\mapsto\tau^{-\rhov_c\cdot\l},
\quad
t_{s_i}\mapsto -t^{-1}\ \roman{if}\ i\neq 0,
\quad
t_{s_0}\mapsto t,
\hfill\cr
\bar\Sb_c^\dag:&
\quad
t\mapsto\tau^{c/2},
\quad
q\mapsto\tau^{-1},
\quad
x_\l\mapsto\tau^{-{}^{s_2s_0}\rhov_c\cdot\l},
\quad
t_{s_i}\mapsto -t^{-1}\ \roman{if}\ i\neq 2,
\quad
t_{s_2}\mapsto t.
\hfill
\endaligned$$
The corresponding $\Hb'_c$-modules are
$$\aligned
{\Sb'_c}^\1dag:&
\quad
\kappa\mapsto c,
\quad
\xi_\l\mapsto\rhov_c\cdot\l,
\quad
s_i\mapsto -1\ \roman{if}\ i\neq 0,
\quad
s_0\mapsto 1,
\hfill\cr
{\bar{\Sb'}_c}^\dag :&
\quad
\kappa\mapsto c,
\quad
\xi_\l\mapsto {}^{s_2s_0}\rhov_c\cdot\l,
\quad
s_i\mapsto -1\ \roman{if}\ i\neq 2,
\quad
s_2\mapsto 1.
\hfill\endaligned$$
Set $\Sb_{-c}={}^{IM}\Sb_c^\dag$, $\bar\Sb_{-c}={}^{IM}\bar\Sb_c^\dag$, etc.
The $\Hb_{-c}$-modules $\Sb_{-c},$ $\Sb'_{-c}$ are spherical but not
$Y$-spherical.
The $\Hb_{-c}$-modules $\bar\Sb_{-c},$ $\bar\Sb'_{-c}$ are neither
spherical nor $Y$-spherical. 

\vskip3mm

\noindent{\bf 2.2.6. Remarks.}
(a)
Fix $N$-th root $q_N$ of $q$.
The group $\Omega$ acts on $X_0\oplus(1/N)\ZZ\delta$ as in 1.0.2.
Let $\HHc$ be the $\CC_{q_{_N},t}$-algebra generated by $\Hb$, $\Omega$
with the relations 
$$
\pi_rt_{s_i}\pi_r^{-1}=t_{s_j},
\ 
\pi_rx_\lambda\pi_r^{-1}=x_\mu
\quad \roman{if}\ {}^{\pi_r}\!\lambda=\mu,
{}^{\pi_r}\a_i=\alpha_j.
$$
See \cite{C4, def.~2.12.1, 2.12.3}.
We call $\HHc$ the extended double affine Hecke algebra
(written EDAHA for short).
The degenerate EDAHA, denoted by $\HHc'$, is defined
in the same way.
For each $\muv\in Y_0$, $r\in O$ we set
$x_{\ov_r+\muv}=\pi_rt_{w_r}x_{\muv}\in\HHc.$
Recall that $\check X_0=Y_0+\sum_r(\ov_r+Y_0).$
Then $\HHc$ is generated by the elements
$x_\lambda,$ $x_{\lambdav}$, $t_w$ with
$\lambdav\in\check X_0,$
$\lambda\in X_0,$
and $w\in W_0$.

Write $\check\HHc$ for the EDAHA associated to $\check G_\CC$.
The Cherednik-Fourier transform is the unique 
$\CC$-algebra isomorphism 
such that
$$CF:
\quad
\HHc\to\check\HHc,
\quad
x_{o_i}\mapsto x_{\ov_i},
\quad
x_{\ov_i}\mapsto x_{o_i},
\quad
t_{s_i}\mapsto t_{s_i}^{-1},
\quad
q_{_N}\mapsto q_{_N}^{-1},
\quad 
t\mapsto t^{-1},
$$
with $i\in I_0$.
The dual here is due to the fact that $\HHc$ is built over the Hecke algebra
of $W_0$ and the lattices $X_0$, $\check X_0$, while in loc.~cit.~
the EDAHA is built over the Hecke algebra of $W_0$ and 2 copies of $X_0$.

(b)
The equivalence 2.2.4 has also an obvious analogue for $\hat\Hb$.

(c)
By A.1.1(b), A.4.1, and the analogue statement for $\Hb'$ it is easy to
see that the diagrams below commute
(here we set $\tilde s^\dag=(s^{-1},\tau^{-1})$) 
$$\matrix
\Irr(O_{\tilde s}(\Hb_c))&\to&\Irr(O_{\lambdav}(\Hb'_c))
\cr
{\ss IM}\downarrow&&\downarrow{\ss IM}
\cr
\Irr(O_{\tilde s^\dag}(\Hb_{-c}))&\to&\Irr(O_{\lambdav}(\Hb'_{-c})),
\endmatrix
\quad
\matrix
\Irr(O_h(\Hb))&\to&\Irr(O_h(\Hb_{(h)}))
\cr
{\ss IM}
\downarrow&&\downarrow
{\ss IM}
\cr
\Irr(O_{h^\dag}(\Hb))&\to&\Irr(O_{h^\dag}(\Hb_{(h^\dag)})).
\endmatrix
$$
Recall that
the equivalence between graded and affine Hecke algebras commutes with IM up to 
conjugation by an invertible holomorphic function on the torus.
For DAHA's one must consider, instead of an holomorphic function,
an infinite product of holomorphic functions, 
yielding delicate questions of convergence.

(d)
For each $\lambdav=(\muv,m)$ with $\muv\in\check X_0$,
setting $s=\muv(\tau_m)$ we get $\Delta_0=\Delta_{0,(h)}$.
Thus 2.2.4(a) yields an equivalence 
$O_{\tilde s}(\Hb_c)\to O_{\lambdav}(\Hb'_c)$
which preserves the dimension of finite dimensional modules.

(e)
In the paper we'll only use degenerate DAHA associated to Chevalley groups.
A degenerate DAHA can be associated to any 
connected reductive group $G_\CC$,
or to any root datum $(X_0,\Delta_0,Y_0,\check\Delta_0)$,
by mimicking the definition of $\Hb$ in section 2.1.

\subhead 2.3. Reminder on rational DAHA's \endsubhead

The rational
DAHA associated to $G_\CC$ is the unital associative
$\CC_\kappa$-algebra $\Hb''$ generated by $\check V_0$, $V_0$, $W_0$
with defining relations
$$\aligned
&w\lambdav w^{-1}={}^w\lambdav,\quad
w\l w^{-1}={}^w\l,\quad
\lambdav \muv=\muv \lambdav,\quad
\l\mu=\mu\l,\quad
\cr
&\qquad\lambdav\l-\l\lambdav=
\l\cdot\lambdav+\kappa\sum_{a\in\Delta_0^+}(a\cdot\lambdav)(\l\cdot \av)s_a,
\endaligned$$
for all $\l,\mu\in V_0,$ $\lambdav,\muv\in\check V_0,$ $w\in W_0.$
There are $\CC$-algebra isomorphisms
$$\aligned
&IM:\quad\Hb''\to\Hb'',
\quad
s_i\mapsto -s_i,
\quad
\l\mapsto \l,
\quad
\lambdav\mapsto\lambdav,
\quad
\kappa\mapsto-\kappa,
\cr
&CF:\quad\Hb''\to\Hb'',
\quad
w\mapsto w,
\quad
a_i\mapsto-2\av_i/(\av_i,\av_i),
\quad
\av_i\mapsto 2a_i/(a_i,a_i),
\quad
\kappa\mapsto\kappa.
\endaligned$$

The subalgebras generated by
$V_0,\check V_0$ are isomorphic to the polynomial algebras
$\CC_\kappa[\check V_0]$, $\CC_\kappa[V_0]$ respectively.
Recall the subcategory $O(\Hb'')\subset\mod(\Hb'')$ from
\cite{GGOR}, 
which consists of the finitely generated modules which are locally finite over
$\CC_\kappa[\check V_0]$ and such that $\kappa$ acts by multiplication by
a complex number.

The multiplication yields an isomorphism
$\CC_\kappa[V_0]\otimes\CC W_0\otimes\CC[\check V_0]\to\Hb''.$
See \cite{EG} for details.
Let $\trivb$ be the one dimensional $W_0\ltimes\CC[\check V_0]$-module such that
$V_0$, $W_0$ act trivially.
The induced $\Hb''$-module
$\Pb''=\Pb''_Y=\Hb''\otimes_{W_0\ltimes\CC[\check V_0]}\trivb$
is equal to $\CC_\kappa[V_0]$ as a $\CC_\kappa$-module.
We set also
$\Pb''_X=\Hb''\otimes_{W_0\ltimes\CC[V_0]}\trivb$.
We write $\Hb''_c$, $\Pb''_c$, $\Pb''_{X,c}$
for the specialization at the maximal ideal $(\kappa-c)$.
A $\Hb''_c$-module is ($X$)-spherical if it is a quotient of $\Pb''_c,$
and $Y$-spherical if it is a quotient of $\Pb''_{X,c}$.

Let $\mod(\Hb''_c)_{nil}\subset\mod(\Hb''_c)$,
$O(\Hb'')_{nil}\subset O(\Hb'')$
be the full subcategories consisting of the modules which are locally
nilpotent over $\CC_\kappa[V_0].$
Finite dimensional $\Hb''_c$-modules belong to
$\mod(\Hb''_c)_{nil}$ by 2.3.1(d) below.

%The obvious $\ZZ_{\ge 0}$-grading on $\CC[X_\CC]$
%induces a grading on $\Pb''_b$.
%Let $\Pb''_b(n)$ be the degree $n$ subspace.
%We write $1\in\Pb''_b(0)$ for the obvious generator of $\Pb''_b$.
%We equip $\Hb''_b$, $\Pb''_b$ with the $\CC[X_\CC]_0$-adic topology.
%Let $\hat\Hb''_b$, $\hat\Pb''_b$ be the corresponding completions.
%We have $\hat\Pb''_b=\prod_{n\ge 0}\Pb''_b(n)$.
%Equip $V$ with the $W_0$-invariant non degenerate symmetric pairing such that
%$B(x,y)=\sum_{a\in\Delta_0}(a^\vee\cdot x)(a^\vee\cdot y)$.
%The pairing $``\cdot"$
%yields an isomorphism of $\CC W_0$-modules $*:V\to V^\vee$ such that
%$a^*={1\over 2}B(a,a)a^\vee$.
%The inverse map is denoted by $*$ as well.
%We extend $*$ to a $\CC_\kappa$-linear anti-involution $\Hb''\to\Hb''$
%such that $w\mapsto w^*=w^{-1}$ for each $w\in W_0$.
%There is an unique $W_0$-invariant symmetric pairing on
%$\Pb''_b$ such that $(1:1)=1$ and
%$(x\bullet p:q)=(p:x^*\bullet q)$
%for each $p,q\in\Pb''_b$, $x\in\Hb''$.
%We use the symbol $\bullet$ to denote the $\Hb''$-action,
%hopping it will not create any confusion,
%and we write 1 for the obvious generator of the polynomial representation.
%See \cite{DO, sec.~2.4}, \cite{D2, sec.~3} for details.

\proclaim{2.3.1. Proposition}
(a)
The $\Hb''_c$-module $\Pb''_c$ has a unique simple quotient $\Lb''_c$.
%The radical of $(\ :\ )$ is the unique maximal
%proper $\Hb''_c$-submodule of $\Pb''_c$.

(b)
There is an equivalence of categories
$\mod(\Hb''_c)_{nil}\to\mod(\Hb'_c)_{uni}$
which preserves the dimension of finite dimensional modules.

(c)
It restricts to an equivalence
$O(\Hb''_c)_{nil}\to O(\Hb'_c)_{uni}$,
and to an equivalence
from the full subcategories of spherical, $Y$-spherical modules to the
full subcategories of spherical, $Y$-spherical ones respectively.
%Further it commutes with the involution $IM$.

(d)
Each finite dimensional $\Hb''_c$-module $M$ is isomorphic to ${}^{CF}M$.
In particular it is spherical iff it is $Y$-spherical. 
Further the subalgebras $\CC[V_0]$, $\CC[\check V_0]$ 
act nilpotently on any finite dimensional module. 

%(e)
%If $M,M'$ are simple modules in $O(\Hb'')$ with the same weights,
%then they are isomorphic.
\endproclaim

\noindent{\sl Proof :}
Claim (a) is proved in \cite{DO, cor.~2.28, prop.~2.34}, \cite{D2, prop.~3.5}.

Part (b) is due to Cherednik.
For any $\lambdav\in\check V_0$ consider the formal series
$f(\lambdav)=(1-\exp(-\lambdav))^{-1}-\lambdav^{-1}$ in $\CC[[V_0]]$.
The action of $\lambdav$ on $M\in\mod_{nil}(\Hb''_c)$
is nilpotent.
We assign to $M$ the $\Hb'_c$-module with the same underlying vector space
and the action such that
$$w\mapsto w,\quad\xi_\l\mapsto\l-\kappa\l\cdot\rhov+
\kappa\sum_{a\in\Delta_0^+}(\av\cdot\l)f(\av)(1-s_a),\quad
\xi_{\lambdav}\mapsto \exp(\lambdav),$$
with $\l\in X_0$, $\lambdav\in Y_0$, and $w\in W_0$.
See \cite{BEG, prop.~7.1}, \cite{C4, sec.~2.12.4} for details.
This functor yields obviously an equivalence
from the full subcategories of spherical, $\check V_0$-spherical modules to the
full subcategories of spherical, $Y$-spherical ones respectively.

Part (c) is obvious, while Part (d) is proved in \cite{EG, sec.~1.5}, 
\cite{BEG1, sec.~3}.
Indeed, the elements
$2\sum_ia_io_i/(a_i,a_i)$,
$2\sum_i\av_i\ov_i/(\av_i,\av_i)$,
and $\sum_i(o_i\av_i+\av_io_i)/2$ form a $\sen\len_2$-triple in $\Hb''_c$.
Under the adjoint action $\Hb''_c$ is a locally finite $\sen\len_2(\CC)$-module.
Thus this action exponentiate in a compatible way, both on $\Hb''_c$ and
the finite dimensional module.
The action of the element $(\smallmatrix 0&1\cr-1&0\endsmallmatrix)$ in 
$\SL_2(\CC)$
yields a linear automorphism $CF$ of $M$ such that $CF(xy)=CF(x)CF(y)$
for all $x\in\Hb''_c$, $y\in M$.
Thus ${}^{CF}M\simeq M$ as a $\Hb_c''$-module.
\qed

\vskip3mm

\noindent{\bf 2.3.2. Example.}
If $c=-1/h$ the $\Hb''_c$ module 
$\Lb''_c$ is particularly simple, since
it is one dimensional such that 
$\kappa\mapsto c,$
$\l,\lambdav\mapsto 0,$
$s_i\mapsto 1.$

\vskip3mm

If $\dim(\Lb''_c)<\infty,$ let $\Lb'_c$ be the $\Hb'_c$-module
associated to $\Lb''_c$ via 2.3.1(b),(d).
The $\Hb'_c$-module $\Lb'_c$ is spherical by 2.3.1(c).
Thus it belongs to $O_c(\Hb')$ by 2.2.2, because this category is stable under
taking quotients.
By 2.2.6(d) there is an equivalence 
$O_c(\Hb'_c)\to O_c(\Hb_c)$
which preserves the dimension of finite dimensional modules.
Let $\Lb_c$ be the image of $\Lb'_c$ by this equivalence.
The modules $\Lb_c$, $\Lb'_c$ are both
spherical, finite dimensional, and simple.

We'll say that a finite dimensional $\Hb_c$-module $M$ comes from $\Hb''_c$ if
$M\in O_h(\Hb)$ for some $h$ such that $\Delta_0=\Delta_{0,(h)}$
and if the corresponding $\Hb'_c$-module belongs to
the subcategory $\mod(\Hb'_c)_{uni}$.
We'll also say that a finite dimensional $\hat\Hb_c$-module 
comes from $\Hb''_c$ if its restriction to $\Hb_c$ comes from $\Hb''_c$.

Set $\Lb_{c}^\dag={}^{IM}\Lb_{-c}$,
${\Lb'_{c}}^\1dag={}^{IM}\Lb'_{-c}$, 
and ${\Lb''_{c}}^\dag={}^{IM}\Lb''_{-c}$.
For a future use, let us quote the following basic facts.

\proclaim{2.3.3. Proposition}
(a)
One simple finite dimensional quotient of $\Pb_c$ at most
comes from $\Hb''_c$.
Further, a finite dimensional quotient of $\Pb_c$
which comes from $\Hb''_c$ is isomorphic to $\Lb_c$.

(b)
The modules $\Lb_c$, $\Lb'_c$ are
spherical and $Y$-spherical.

(c)
The equivalence 2.2.4 takes ${\Lb'_{c}}^{\1dag}$ to $\Lb_{c}^\dag$,
and 2.3.1(b) takes ${\Lb''_{c}}^\dag$ to ${\Lb'_{c}}^\1dag$.
\endproclaim

\noindent{\sl Proof :}
Part (a) is obvious, because 
if a simple finite dimensional quotient of $\Pb_c$ comes from the
$\Hb''_c$-module $M''$, 
then $M''$ is spherical and finite dimensional, thus it is
isomorphic to $\Lb''_c$ by 2.3.1(a). 
Part (b) follows from 2.2.3(b), 2.3.1(c),(d).

Let us concentrate on Part (c).
The first claim is obvious by A.4.1. 
The automorphism $IM$ of $\Hb'$ takes $\xi_\lambdav$ to itself.
Hence, since $\Lb'_{-c}\in O_{-c}(\Hb')_{uni}$, we have
${\Lb'_c}^\1dag\in O_{c}(\Hb')_{uni}$. 
Thus ${\Lb'_c}^\1dag$ comes from a simple $\Hb''_c$-module $M$. 
Since $\Lb'_{-c}$ is $Y$-spherical by Part (b),
the $\Hb'_c$-module ${\Lb'_{c}}^\1dag$ is a quotient of
${\Pb'_{X,c}}^\7dag=\Hb'_c\otimes_{\CC W}\sgnb$,
where $\sgnb$ is the signature of $W$.  
Thus $M$ is a simple quotient of
the $\Hb''_c$-module
${\Pb''_{X,c}}^\5dag=\Hb''_c\otimes_{W_0\ltimes\CC[V_0]}\sgnb$
by definition of 2.3.1(b).
So ${}^{IM}M$ is a simple quotient of
the $\Hb''_{-c}$-module
$\Pb''_{X,-c}$, i.e., it is spherical by 2.3.1(d),
yielding $M={}^{IM}\Lb''_{-c}={\Lb''_c}^\dag$.
We are done.

\qed

\vskip3mm

\proclaim{2.3.4. Question}
Is it true that the equivalence 2.3.1(b) commutes with $IM$
(at least when restricted to irreducible modules)?
\endproclaim

\vskip3mm

\noindent{\bf 2.3.5. Remarks.}
(a)
In type $A_n$ the finite dimensional $\Hb''_c$-modules are classified in
\cite{BEG2}.
We have $\dim(\Lb''_c)<\infty$ iff $m=h$, $k<0$.

(b)
The category of finite dimensional $\Hb'_c$-modules is not equivalent
to the category of finite dimensional $\Hb''_c$-modules.
Compare \cite{BEG2, prop.~7.1}.
If $(G_\CC,c)=(C_2,1/2)$
then ${\Sb'_c}^\1dag$, ${\bar{\Sb'}_c}^\dag$ do not come from $\Hb''_c$, because
$\xi_{\av_2}\mapsto-1$ on both.

(c)
If $(G_\CC,c)=(C_2,1/2)$
then ${\Sb'_c}^\1dag$, ${\bar{\Sb'}_c}^\dag$ do not extend to $\HHc'$-modules.
However, since $O=\{2\}$ and ${}^{\pi_2}\rhov_c={}^{s_2s_0}\rhov_c$
the sum ${\Sb'_c}^\1dag\oplus{\bar{\Sb'}_c}^\dag$ 
is a two dimensional simple $\HHc'$-module
such that $\pi_2\mapsto(\smallmatrix 0&1\cr 1&0\endsmallmatrix)$.
It does not come from $\Hb''$ either.

(d)
If $G_\CC=G_2$ then $\Hb'=\HHc'$.
If $c=1/3$ the assignement
$${\Sb'_c}^\1dag:\quad
\xi_\l\mapsto
\Bigl(\smallmatrix\rhov_c\cdot\l&0\cr0&{}^{s_0}\rhov_c\cdot\l
\endsmallmatrix\Bigr),\quad
s_1\mapsto -\Id,\quad
s_2\mapsto\Bigl(\smallmatrix -1&0\cr0&1\endsmallmatrix\Bigr),\quad
s_0\mapsto{1\over 2}\Bigl(\smallmatrix 1&3\cr 1&-1\endsmallmatrix\Bigr),\quad
\kappa\mapsto c$$
yields a two dimensional simple $\Hb'_c$-module which does not come
from $\Hb''_c$, because
$\xi_{\av_2}\mapsto{1\over 2}
\Bigl(\smallmatrix -1&-3\cr 1&-1\endsmallmatrix\Bigr)$
which is not unipotent.
If $c=1/2$ the assignement
$${\Sb'_c}^\1dag:\quad
\xi_\l\mapsto
\Biggl(\matrix\rhov_c\cdot\l&0&0\cr0&{}^{s_0}\rhov_c\cdot\l&0\cr
0&0&{}^{s_2s_0}\rhov_c\cdot\lambda
\endmatrix\Biggr),\quad
s_1\mapsto\Biggl(\matrix -1&0&0\cr0&-1&0\cr0&0&1\endmatrix\Biggr),$$
$$s_2\mapsto\Biggl(\matrix -1&0&0\cr0&1/2&3\cr0&1/4&-1/2\endmatrix\Biggr),
\quad
s_0\mapsto\Biggl(\matrix 1/3&8&0\cr 1/9&-1/3&0\cr 0&0&-1\endmatrix\Biggr),\quad
\kappa\mapsto c$$
yields a three dimensional simple $\Hb'_c$-module which does not come
from $\Hb''_c$.
In both cases the module $\Sb'_{-c}={}^{IM}{\Sb'_c}^\1dag$ is spherical.

(e)
The modules $\Pb'_c$, $\Pb_c$ do not have a simple top in general.
They may even admit several finite dimensional simple quotients.
For instance,  if $(G_\CC,c)=(C_2,-1/2)$
then $\Sb'_{c}$ is the quotient of
$\Pb'_{c}$ by the ideal generated by
$\xi_{\av_1}-1$, $\xi_{\av_2}+1$,
while $\Lb'_{c}$ is the quotient
by the ideal generated by
$(\xi_{\av_1}-1)(\xi_{\av_2}-1)+(\xi_{\av_2}-1)^2$,
$2(\xi_{\av_1}-1)(\xi_{\av_2}-1)+(\xi_{\av_1}-1)^2$.

(f)
In this paper we consider $\Hb''_c$-modules with $c\in\QQ$ only.
There are no finite dimensional simple $\Hb''_c$-modules
if $c$ is not rational. See \cite{BEG1, cor.~2.11} for details.

(g)
In type $A_n$ all simple finite dimensional $\Hb$-modules come from 
$\Hb''$ by \cite{FG}. This does not follow from \cite{BEG2}.

%The central character of $\Sb'_c$, $\bar\Sb'_c$
%as $W$-modules is the $W_0$-orbit
%$\{a(-1);a\in\Delta_{0,s}\}$ in $X\otimes\CC^\times$.
%$\xi_{\av_1}\mapsto 1,$
%Its central character,
%as a $W$-module, is the $W_0$-orbit
%$\{a(\varpi_3);a\in\Delta_{0,\ell}\}$ in $X\otimes\CC^\times$.

\subhead 2.4. DAHA's and affine Springer fibers \endsubhead

Let $\hat G$ be the maximal (=Kac-Moody)
central extension of $G$ by $\CC^\times$.
Consider the groups
$\tilde G=\hat G\rtimes\CC^\times_\delta$,
$G_d=G\rtimes\CC_\delta^\times$.
Here $\CC^\times_\delta$ denotes the group $\CC^\times$ acting by
`rotating the loops'.
Write $\CC^\times_{\omega_0}$ for the image of $\CC^\times$ in
the center $Z(\hat G)$.
View $T_0$, $T_d$ as closed tori in $G$, $G_d$ in the obvious way.
Then $\hat T_0=T_0\times\CC^\times_{\omega_0}$, 
$\tilde T=\hat T_0\times\CC^\times_\delta$
are the inverse image of $T_0$, $T_d$ in $\hat G$, $\tilde G$ respectively.
We will also view $H$, $\hat H$ as closed tori in
$G_d\times\CC^\times$, $\tilde G\times\CC^\times$.
Let $\tilde\gen$, $\tilde\ten$, $\hat\hen$, $\hen$ be the Lie algebras
of $\tilde G$, $\tilde T$, $\hat H$, $H$.

For each $w\in W$ let $\ben_w\subset\tilde \gen$ be the Iwahori Lie subalgebra
containing $\tilde\ten$ whose set of real roots is equal to $w(\Delta^+_\re)$.
Let $\Bc$ be the set of all Iwahori Lie subalgebras in $\tilde\gen$
(i.e., the set of Lie subalgebras $\tilde G$-conjugate into $\ben_1$).
Set $\uen_w$ equal to the pronilpotent radical of $\ben_w$,
and $B_w=N_{\tilde G}(\ben_w)$ (an Iwahori subgroup).

The parahoric Lie subalgebras of $\tilde\gen$ are those Lie subalgebras
which contain an Iwahori Lie subalgebra.
The conjugacy classes of parahoric Lie subalgebras correspond bijectively
to the subsets of $I$.
A parahoric Lie subalgebra is said to be of type $i$
if it is conjugated to $\ben_{1,i}$.
Set $\ben_{w,i}=\ben_w+\ben_{ws_i}$ and $B_{w,i}=N_{\tilde G}(\ben_{w,i})$.
So $\ben_{w,i}$ is the unique parahoric Lie subalgebra of type $i$ containing
$\ben_w$.
More generally, for each $\ben\in\Bc$ let $\ben_i\subset\tilde\gen$
be the unique parahoric Lie subalgebra of type $i$ containing $\ben$,
and $\uen_i$ be its pronilpotent radical.
So $\uen_{w,i}$ denotes the pronilpotent radical of $\ben_{w,i}$.
Let $\Bc_i$ be the set of parahoric Lie subalgebras of $\gen$ of type $i$.
A parahoric subgroup of $\tilde G$ 
is the normalizer of a parahoric Lie subalgebra of $\tilde \gen$.

Let $\Nc\subset\gen$ be the set of TN elements,
and $$\dot\Nc=\{(x,\ben)\in\gen\times\Bc;x\in\uen\},$$
where $\uen$ is the pronilpotent radical of $\ben$.
The first projection $\pi:\dot\Nc\to\gen$ 
takes $\dot\Nc$ onto $\Nc$ by \cite{KL2, sec.~2, lem.~1}. 
There is an unique $\tilde G\times\CC^\times$-action on
$\Nc$, $\Bc$, and $\dot\Nc$ such that
$\CC_\delta^\times$ acts by `rotating the loops' and
$\CC^\times$ by scalar multiplication on $\gen$ and trivially on $\Bc$.
For each $h=(s,\tau,\zeta)\in\tilde G\times\CC^\times$, 
we set $\tilde s=(s,\tau)$
and $(\ad\tilde s)=(\ad s)\circ F_\tau$.
For each $\ben\in\Bc$, $x\in\gen$ we have
$$
h\cdot x=\zeta^{-2}(\ad\tilde s)(x),
\quad
h\cdot\ben=(\ad \tilde s)(\ben).
$$

Now let $h\in\hat H$, i.e., we assume that $s\in\hat T_0$. 
Set $\tilde G^{\tilde s}=Z_{\tilde G}(\tilde s)$. 
Since $\zeta$ does not act on $\tilde G$ we may also write
$\tilde G^h$ for $\tilde G^{\tilde s}$. 
Consider also the fixed points subsets 
$\gen^h\subset\gen$ and $G^h\subset G$.
Let $G^{h,\circ}\subset G^h$ be the connected component containing the unit.
The set of affine roots of $\gen$ which appear in $\gen^h$ is 
$$\Den_h=\{\alpha\in\Delta;\alpha(\tilde s)=\zeta^2\}.$$ 
Notice that $\Den_h\subset\Delta_\re$ if $\zeta^2\notin\la\tau\ra$,
and that $\Delta_h\subset\Delta_\re$ because $\tau$ is not a root of unity.

\proclaim{2.4.1. Lemma}
(a)
The groups $\tilde G^h$, $G^{h,\circ}$ 
are the sets of $\CC$-points of connected reductive groups,
and $\gen^h$ is a finite dimensional space.

(b)
The root systems of the pairs 
$(\tilde G^h,\tilde T)$,
$(G^{h,\circ},T_0)$ are both isomorphic to $\Delta_h$. 
The $G^{h,\circ}$-orbits and the $\tilde G^h$-orbits in the nilpotent cone
of $\gen^h$ are the same.
\endproclaim

\noindent{\sl Proof :}
Part (a) is proved in \cite{V1, lem.~2.13}.
The connected group $G^{h,\circ}$ is generated by $T_0$ and the root subgroups
associated to the roots in $\Delta_h$ by the axioms of a split BN pair.
See \cite{SS, sec.~II.4.1}.
The $\CC_\delta^\times\times\CC^\times_{\omega_0}$-torsor 
$\tilde G\to G$ factors to 
a $\CC_\delta^\times\times\CC^\times_{\omega_0}$-torsor 
$\tilde G^h\to G^{h,\circ}$.
To prove (b) it is enough to check that any nilpotent
$G^{h,\circ}$-orbit $O\subset\gen^h$ is preserved by the 
$(\ad\CC_\delta^\times)$-action.
This is obvious.
Indeed, the actions of $(1,\tau^{-1},1)$ and $(s,1,\zeta)$ on $\gen^h$ 
are the same.
Further $O$ is preserved by the action of $(s,1,1)$,
because $s\in G_0^{h,\circ}$ since $G_0$ is simply connected, 
and $O$ is preserved by the action of $(1,1,\zeta)$ because $O$ is a cone.

\qed

\vskip3mm

Let $n_h$ be the cardinal of the set
$\Delta^+_h=\Delta_h\cap\Delta^+$.
Let $B^h\subset \tilde G^h$ be the Borel subgroup containing $\tilde T$
associated to $\Delta_h^+$.
Recall that $W_h=Z_W(h)$ by 2.1.8(b).
The set $W^h=\{w\in W;B_w\cap \tilde G^h=B^h\}$
is a set of representatives of the right $W_h$-cosets in $W$.

The element $h$ acts on $\Nc,\Bc$, $\dot\Nc$.
The fixed points set
$\Nc^h$ is the set of $\CC$-points of a $\CC$-scheme of finite type,
$\dot\Nc^h$ is the set of $\CC$-points of a smooth
$\CC$-scheme locally of finite type,
and $\pi$ yields a
$\tilde G^h$-equivariant proper morphism
$\pi^h:\dot\Nc^h\to\Nc^h$.
Since the $h$-action on $\Bc$ depends only on $\tilde s$ we may write
$\Bc^h=\Bc^{\tilde s}$.
The set $\Bc$ has the structure of an ind-projective ind-$\CC$-scheme,
such that $\Bc^h$ is a closed $\CC$-subscheme locally of finite type.
So are also the affine Springer fibers
$\Bc_x^h=(\pi^h)^{-1}(x)$ for any $x$,
or $\Bc_x=\pi^{-1}(x)$ if $x$ is RS.
It is known that $\Bc_x$ is connected, and that it is a 
(projective) variety iff $x$ is ERS.
See \cite{KL2, sec.~3, cor.~2, sec.~4, lem.~2}
for details.
Let us quote the following basic facts for a future use.

\proclaim{2.4.2. Lemma}
(a)
For each $x\in\Nc^h$ we have
$\dim\,H_*(\Bc_x^h,\CC)<\infty$ iff $x\in\Nc^h_\ERS$.

(b)
If $\tau^k\neq\zeta^{2m}$ whenever $m>0>k$, then
$\Nc^h=\gen^h$ and it may have an infinite number of $\tilde G^h$-orbits.

(c)
If $\tau^k\neq\zeta^{2m}$ whenever $m,k>0$, then
$\Nc^h$ consists of nilpotent elements
and has a finite number of $\tilde G^h$-orbits.

(d)
If $\tau^k\neq\zeta^{2m}$ whenever $m,k\neq 0$, then
$\Nc^h=\gen^h$ consists of nilpotent elements
and has a finite number of $\tilde G^h$-orbits.
\endproclaim

\noindent{\sl Proof :}
Part (a) follows from \cite{KL2}.
Indeed, if $x\in\Nc^h\setminus\Nc^h_\reg$
then $\Bc_x\subset\Bc$ is a closed
infinite dimensional sub-ind-scheme by loc.~cit.,
sec.~2, lem.~6. Hence the fixed points set $\Bc^h_x$, which is locally of finite
type, has necessarily an infinite number of irreducible components.
Thus $H_*(\Bc^h_x,\CC)$ is infinite dimensional.
Further, if $x\in\Nc^h_\reg$ then the dimension of the homology is finite iff
$x\in\Nc^h_\ERS$ by loc.~cit., sec.~3.

Parts (b) to (d) are proved in \cite{V1, lem.~2.14}.
We give another proof of (b), (c) for the comfort of the reader.
Observe that (d) is a direct consequence of (b), (c).
First recall the following :

\vskip1mm

\itemitem{(i)}
we have $x\in\Nc$ iff $\varphi_i\eta(x)\in \eps A,$

\vskip1mm

\itemitem{(ii)}
an element $x\in\gen$ is nilpotent iff 
$\varphi_i\eta(x)=0$ for each $i$, 

\vskip1mm

\itemitem{(iii)}
if $x\in\gen^h$ we have
$F_\tau\varphi_i\eta(x)=\zeta^{2d_i}\varphi_i\eta(x)$.

\vskip1mm

Therefore,
if $\tau^k\neq\zeta^{2m}$ for each $m,k\neq 0$ 
and $x\in\gen^h$ then 
$\varphi_i\eta(x)=0$ by (iii), hence
$x$ is nilpotent by (ii).
Further, if $\tau^k\neq\zeta^{2m}$ for each $m,k>0$ 
and $x\in\Nc^h$ then
$\varphi_i\eta(x)\in\eps A$ by (i), hence
$\varphi_i\eta(x)=0$ by (iii), so 
is $x$ nilpotent by (ii).
Finally, if $\tau^k=\zeta^{2m}$ for some $m,k>0$ and
$x\in\gen^h$ then
$\varphi_i\eta(x)\in \eps A$ by (iii), hence
$x\in\Nc^h$ by (i).

Thus we are reduced to check the second claim of (c).
For the comfort of the reader we reproduce here the proof from \cite{V1}.
Let $L\subset X_0$ be the subgroup generated by the set 
$\Delta_{0,(h)},$
and let $e$ be the cardinal of the maximal torsion subgroup of
$X_0/L$.

Assume first that $\tau^k\neq\zeta^{2m}$ for each $m,k\neq 0$.
There is an unique group homomorphism $\g:L\to\ZZ$ such that
$a(s)\tau^{\gamma(a)}\in\la\zeta^2\ra$ for each $a\in\Delta_{0,(h)}$.
Fix an element $\lambdav\in Y_0$ such that
$\gamma(a)=a\cdot\lambdav/e$ for each $a\in\Delta_{0,(h)}$.
We have
$$\aligned
&\Delta_h=\{(a,\ell)\in\Delta_{\re,(h)};
a(s_1)=a(s_2)=1,a\cdot\lambdav=e\ell\},
\\
&\Den_h=\{(a,\ell)\in\Delta_{\re,(h)};
a(s_1)=1,a(s_2)=\zeta^2,a\cdot\lambdav=e\ell\}.
\endaligned
$$
where $s_1=\lambdav(\varpi_e),$ $s_2=s\lambdav(\tau_e)$.
Notice that $\Den_h\subset\Delta_\re$, because $\zeta^2\notin\la\tau\ra$.
Therefore, since $\Delta_h$ is the root system of $G^{h,\circ}$ 
the conjugation by $\lambdav(\eps_e)^{-1}$ 
yields a group isomorphism
$G^{h,\circ}\to G_0^{S,\circ}$ with $S=\{s_1,s_2\}$.
For a similar reason, it yields a linear isomorphism 
$$\gen^h\to\gen_0^S=\{x\in\gen_0;(\ad s_1)(x)=\zeta^{-2}(\ad s_2)(x)=x\}$$
which is compatible with the isomorphism 
$G^{h,\circ}\simeq G_0^{S,\circ}$ 
under the adjoint action.
Now, by \cite{R6, thm.~B} the $G_0^{S,\circ}$-orbits on
$\gen_0^S$ are the connected components of the traces in $\gen_0^S$ of the
nilpotent $G_0$-orbits in $\gen_0$. 
By 2.4.1(b) the $\tilde G^h$-orbits in $\Nc^h$ 
are the same as the $G^{h,\circ}$-orbits in $\gen^h$. 
So there are only finitely many of them.

Now, assume that $\tau^k=\zeta^{2m}$ with $m>0>k$ relatively prime.
There is an unique group homomorphism $\g:L\to\ZZ$ such that
$a(s)=\tau^{-\gamma(a)/m}$ for each $a\in\Delta_{0,(h)}$.
Fix an element $\lambdav\in Y_0$ such that
$\gamma(a)=a\cdot\lambdav/e$ for each $a\in\Delta_{0,(h)}$.
%Assume that $m\neq 1$.
We have
$$\aligned
&\Delta_h=\{(a,\ell)\in\Delta_{\re,(h)};
a(s_1)=a(s_2)=1,a\cdot\lambdav=em\ell\},
\\
&\Den_h=\{(a,\ell)\in\Delta_{(h)};
a(s_1)=\varpi^c,a(s_2)=1,a\cdot\lambdav=em\ell-ek\},
\endaligned$$
where $s_1=\lambdav(\varpi_{em})$ and $s_2=s\lambdav(\tau_{em})$.
%Notice that $\Den_h\subset\Delta_\re$, because $\zeta^2\notin\la\tau\ra$.
Therefore
the conjugation by $\lambdav(\eps_{em})^{-1}$ 
yields a group isomorphism
$G^{h,\circ}\to G_0^{S,\circ}$ with $S=\{s_1,s_2\}$
and a compatible linear isomorphism 
$$\gen^h\to\gen_0^S\otimes\eps^c=
\{x\in\gen_0;(\ad s_1)(x)=\varpi^cx,(\ad s_2)(x)=x\}\otimes\eps^c.$$
Further, the latter restricts to an isomorphism
$\Nc^h\to\Nc_0^S\otimes\eps^c.$
Since the $\tilde G^h$-orbits in $\Nc^h$ 
are the same as the $G^{h,\circ}$-orbits in $\Nc^h$ by 2.4.1(b), 
there are only finitely many of them by loc.~cit. 
%If $m=1$ then $\Den_h$ contains imaginary roots, but none of them
%contribute to $\Nc^h$ because $k<0$.
%To do that, observe that
%$G^{h,\circ},\gen^h$ are contained into
%$G_\CC(\CC[\eps,\eps^{-1}])$,
%$\gen_0\otimes\CC[\eps,\eps^{-1}]$
%respectively.
%Compare section 3.1.
%Thus the assignement $\eps\mapsto 1$ yields a closed embedding 
%$ev:G^{h,\circ}\to G_0$ and a linear embedding $ev:\gen^h\to\gen_0$.
%Set $H_0=ev(G^{h,\circ})$ and $\Xc_0=ev(\Nc^h)$.
%Notice that $H_0$ is a closed subgroup of $G_0$ and that
%$\Xc_0$ is a closed subset of the nilpotent cone $\Nc_0\subset\gen_0$.
%Further, the map $ev$ takes the $\tilde G^h$-orbits in $\Nc^h$ 
%to the $H_0$-orbits in $\Xc_0$.
%Now, we can assume that $\Delta_0=\Delta_{0,(h)}$.
%Fix an element $\lambdav\in\check X_0$ as in 2.2.3(a).
%Then $H_0=Z_{G_0}(\lambdav_m)$
%and $\Xc_0=\{x\in\Nc_0;(\ad \lambdav_m)(x)=\varpi_c x\}.$
%Hence 
%there is a finite number of $H_0$-orbits in $\Xc_0$
%by \cite{R6, thm.~A}.

\qed

\vskip3mm

The connected components of $\Bc^h$ are precisely the subsets
$\Bc_w^h=(\ad \tilde G^h)(\ben_w)$ with $w\in W^h$.
See \cite{R5, prop.~2.3} for the finite type case.
So we have 
$$\Bc_x^h=\bigsqcup_w\Bc_{x,w}^h,
\quad
\dot\Nc^h=\bigsqcup_w\dot\Nc_{w}^h$$
with $\Bc_{x,w}^h=\Bc_x\cap\Bc_w^h$ and
$\dot\Nc_w^h=\{(x,\ben)\in\dot\Nc^h;\ben\in\Bc_w^h\}$.
For a future use, recall that $\tilde G^h$ is connected and that
there is an isomorphism of $\tilde G^h$-varieties
$$\tilde G^h/B^h\to\Bc^h_{w},\quad
gB^h\mapsto(\ad g)(\ben_w).$$
Set
$$H_*(\Bc^h_x,\CC)=\bigoplus_{w\in W^h}H_*(\Bc^h_{x,w},\CC),
\quad
\widehat H_*(\Bc^h_x,\CC)=\prod_{w\in W^h}H_*(\Bc^h_{x,w},\CC),$$
and 
$$K(\Bc^h_x)=\bigoplus_{w\in W^h}K(\Bc^h_{x,w}),
\quad
\widehat K(\Bc^h_x)=\prod_{w\in W^h}K(\Bc^h_{x,w}).$$
The following is proved in \cite{V1, sec.~4.9, 6.2}.

\proclaim{2.4.3. Lemma}
There are representations of $\hat\Hb$ on the vector spaces
$H_*(\Bc^h_x,\CC)$,
$K(\Bc^h_x),$
$H_*(\dot\Nc^h,\CC)$,
$K(\dot\Nc^h).$
\endproclaim

For future use let us recall briefly the proof of 2.4.3.
Set
$$\ddot\Nc=\{(x,\ben,\ben')\in\Nc\times\Bc^2;x\in\ben\cap\ben'\},
\quad
\ddot\Nc^h=\ddot\Nc\cap(\dot\Nc^h)^2, 
\quad
\ddot\Nc^h_{v,w}=\ddot\Nc\cap(\dot\Nc^h_v\times\dot\Nc^h_w).$$
The scheme $\dot\Nc^h$ is smooth locally of finite type,
and $\ddot\Nc^h$ is a closed subscheme of $(\dot\Nc^h)^2$.
Thus there is a topological algebra structure on the vector space
$$\widehat K(\ddot\Nc^h)=\prod_v\bigoplus_wK(\ddot\Nc^h_{v,w}),$$
such that
$$\Ec\star\Ec'=Rp_{13,*}(p_{12}^*\Ec\otimes^Lp_{23}^*\Ec').$$
Here $p_{ij}$ is the obvious projection
$(\dot\Nc^h)^3\to(\dot\Nc^h)^2$,
and the Tor product is relative to the smooth scheme
$(\dot\Nc^h)^3$ 
and its closed subschemes $\ddot\Nc^h\times\dot\Nc^h$,
$\dot\Nc^h\times\ddot\Nc^h$ with intersection contained into 
$p_{13}^{-1}(\ddot\Nc^h)$. 
The topology is the product topology.
A similar construction yields a topological algebra structure on the 
vector space
$$\widehat H_*(\ddot\Nc^h,\CC)=\prod_v\bigoplus_wH_*(\ddot\Nc^h_{v,w},\CC)$$
and right continuous representations of 
$\widehat K(\ddot\Nc^h)$, $\widehat H_*(\ddot\Nc^h,\CC)$
on the discrete vector spaces
$K(\Bc^h_x),$ $H_*(\Bc^h_x,\CC)$,
$K(\dot\Nc^h),$ $H_*(\dot\Nc^h,\CC)$ respectively.

Now, there is a chain of algebra homomorphisms
$$\hat\Hb\,{\buildrel\Phi\over\to}\,\widehat K(\ddot\Nc^h)
\,{\buildrel ch\over\to}\,
\widehat H_*(\ddot\Nc^h,\CC),
\leqno(2.4.4)$$
where $ch$ is the Chern character.
So we get right representations of $\hat\Hb$ 
on $K(\Bc^h_x),$ $H_*(\Bc^h_x,\CC)$.
From now on we will view any right $\hat\Hb$-module as a left module,
up to twisting the action by the anti-automorphism $OP$.
Notice that there is no algebra homomorphism
$\hat\Hb\to\bigoplus_{v,w}K(\ddot\Nc^h_{v,w})$,
because the rhs has no unit.

\proclaim{2.4.5. Lemma}
(a)
The modules $K(\Bc^h_x),$ $H_*(\Bc^h_x,\CC)$ are admissible with 
${}^{w^{-1}\!}h^\dag$-weight subspaces 
equal to $H_*(\Bc^h_{x,w},\CC)$, $K(\Bc^h_{x,w})$.
Both $\hat\Hb$-modules
extend uniquely to topological representations on
$\widehat K(\Bc^h_x),$ $\widehat H_*(\Bc^h_x,\CC)$.
The same holds also with the vector spaces
$\widehat K(\dot\Nc^h),$ $\widehat H_*(\dot\Nc^h,\CC)$.

(b)
The $\hat\Hb$-modules $K(\dot\Nc^h_x),$ $H_*(\dot\Nc^h_x,\CC)$
are isomorphic (via the Chern character).
If $x$ is nilpotent 
the $\hat\Hb$-modules $K(\Bc^h_x),$ $H_*(\Bc^h_x,\CC)$
are also isomorphic.
\endproclaim

\noindent{\sl Proof :}
The first claim of part (a) follows from 
the definition of $\Phi$ below.
So $H_*(\Bc^h_x,\CC)$ is admissible.
The other cases are identical.
The second claim follows from by 2.1.5(c).
The first claim of part (b) is obvious, while the second
one follows from \cite{V1, lem.~6.8}.

\qed

\vskip3mm

We'll need more details on the definition of the map $\Phi$.
Each character $\l\in X=X^*(\tilde T)$ extends uniquely
to a character of the group $B_1$ which is trivial on the prounipotent radical.
Let $\Oc_\Bc(\l)$ be the $\tilde G$-equivariant line bundle on $\Bc$
whose sections consist of the functions $f$ on $\tilde G$
such that $f(gb)=\l(b)f(g)$ for each $g\in\tilde G$, $b\in B_1$.
Let $\Oc_{\dot\Nc}(\l)$ be its pull-back to $\dot\Nc$ by the obvious projection.
For each subvariety $\Yc\subset\dot\Nc$ we write
$\Oc_{\Yc}(\l)$ for $\Oc_{\dot\Nc}(\l)|_{\Yc}$.
If $\Yc\subset\dot\Nc^h$ we write $\Oc^h_{\Yc}(\l)$ 
for the virtual bundle equal to
$\prod_{w\in W^h}{}^w\l(\tilde s)\Oc_{\Yc\cap\dot\Nc^h_w}(\l)$
(it depends on the choice of $h$).
We'll use the same notation for a sheaf, or a virtual vector bundle,
and its class in the Grothendieck group.

Since $\dot\Nc^h$ is smooth, the Tor product of coherent sheaves
equip $\widehat K(\dot\Nc^h)$ with a commutative algebra structure. 
Given $z\in\CC^\times$
let $\Lambda(z\Oc^h_{\dot\Nc^h}(\l))$
be the invertible element 
whose component in $K(\dot\Nc^h_w)$ is equal either to
$1-z\Oc^h_{\dot\Nc^h_w}(\l)$ if ${}^w\l(\tilde s)z\neq 1$ or to 1 else.

Consider the closed subsets
$$\aligned
&\ddot\Nc_1=\{(x,\ben,\ben')\in\Nc\times\Bc^2;x\in\ben,\ben=\ben'\},
\\
&\ddot\Nc_{s_i}=\{(x,\ben,\ben')\in\Nc\times\Bc^2;x\in\uen_i,\ben_i=\ben'_i\}.
\endaligned$$
There is an obvious isomorphism
$\ddot\Nc_1\simeq\dot\Nc$.
Set $\ddot\Nc_1^h=\ddot\Nc_{1}\cap\ddot\Nc^h$ and  
$\ddot\Nc_{s_i}^h=\ddot\Nc_{s_i}\cap\ddot\Nc^h$. 

The map $\Phi$ is defined as follows.
First, we have $\Phi(t)=\zeta$, $\Phi(q)=\tau^{-1}$.
Then $\Phi(x_\l)$ equals 
the image of the virtual sheaf $\Oc^h_{\ddot\Nc^h_1}(-\l)$
by the direct image 
$\widehat K(\ddot\Nc^h_1)\to \widehat K(\ddot\Nc^h).$
Finally, given $\l',\l''\in X$ such that
$\l'\cdot\alphav_i=\l''\cdot\alphav_i=1$ and
$\l'+\l''=\a_i$,
the element $\Phi(tt_{s_i}+1)$ is the image of 
$$-\Oc^h_{\dot\Nc^h}(\l')\Lambda(\Oc^h_{\dot\Nc^h}(\a_i))^{-1}\boxtimes
\Oc^h_{\dot\Nc^h}(\l'')\Lambda(\zeta^2\Oc^h_{\dot\Nc^h}(-\a_i))
\leqno(2.4.6)$$
by the restriction
$\widehat K(\dot\Nc^h\times\dot\Nc^h)\to \widehat K(\ddot\Nc^h_{s_i})$
and the direct image $\widehat K(\ddot\Nc^h_{s_i})\to \widehat K(\ddot\Nc^h).$
The same argument as in \cite{L2} shows that the expression 2.4.6 does not 
depend on the choice of $\l'$, $\l''$.

The action of $\hat\Hb$ on $\widehat K(\dot\Nc^h)$ in 2.4.4 can be computed
explicitely.
Set $p_i:(\dot\Nc^h)^2\to\dot\Nc^h$ 
equal to the $i$-th projection.
By definition, the action of $x$ takes $\Ec$ to 
$$x\Ec=Rp_{2,*}(\Phi(x)\otimes^L(\Ec\boxtimes\Oc_{\dot\Nc^h})).$$
The Tor product is relative to the smooth scheme
$(\dot\Nc^h)^2$
and its closed subsets $\ddot\Nc^h$, $(\dot\Nc^h)^2$ 
with intersection $\ddot\Nc^h$.
Recall that there is an unique 
continuous $W$-action such that 
${}^w\Oc_{\dot\Nc^h}(\l)=\Oc_{\dot\Nc^h}({}^w\l)$ for each $\l$.

\proclaim{2.4.7. Proposition}
For each vector bundle $\Ec$ over $\dot\Nc^h$ and each 
$i\in I$, $\l\in X$ we have 
(tensor product is omitted)
$$
t_{s_i}\Ec=
{\zeta^{-1}({}^{s_i}\Ec-\Ec\Oc^h_{\dot\Nc^h}(-\alpha_i))
+\zeta(\Ec-{}^{s_i}\Ec)\Oc^h_{\dot\Nc^h}(-\alpha_i)
\over\Oc^h_{\dot\Nc^h}(-\alpha_i)-1},
\quad
x_{\l}\Ec=\Oc^h_{\dot\Nc^h}(-\l)\Ec.$$
\endproclaim

\noindent{\sl Proof :}
This formula can be proved either by the same technics as in \cite{V1},
using a version of Thomason concentration theorem, or by an explicit
computation using 2.4.6.
We'll use the second method. We must prove that
$$Rp_{2,*}(\Phi(tt_{s_i}+1)\otimes^L(\Ec\boxtimes\Oc_{\dot\Nc^h}))=
{(1-\zeta^2\Oc^h_{\dot\Nc^h}(-\alpha_i))
(\Ec-{}^{s_i}\Ec)
\over 1-\Oc^h_{\dot\Nc^h}(-\alpha_i)}.
\leqno(2.4.8)
$$

Set
$\Bc_{w,(i)}=\{\ben\in\Bc;\ben_i=\ben_{w,i}\}$
and $\dot\Nc_{w,(i)}=\dot\Nc\cap(\ben_{w,i}\times\Bc_{w,(i)}).$
We have
$$\dot\Nc\simeq\tilde G\times_{B_{w,i}}\dot\Nc_{w,(i)},
\quad
\ddot\Nc_{s_i}\simeq\tilde G\times_{B_{w,i}}(\uen_{w,i}\times\Bc_{w,(i)}^2).$$
Notice that
$\dot\Nc_{w,(i)},$
$\Bc_{w,(i)}$,
$\ben_{w,i}$,
and $\uen_{w,i}$ 
depend only on the class of $w$ in $W/\{1,s_i\}$.

Now, let us consider $h$-fixed points subsets.
For each $u\in W_h\setminus W/\{1,s_i\}$ and each $w\in u$,
let $\Bc_{u,i}^h$ be the $\tilde G^h$-orbit of $\ben_{w,i}$.
It depends only on the class $u$.
We have 
$\Bc_i^h=\coprod_u\Bc_{u,i}^h.$
Compare \cite{R5, prop.~2.3}.
Set 
$$\aligned
\Bc_{u,(i)}^h&=\{\ben\in\Bc^h;\ben_i\in\Bc_{u,i}^h\},
\hfill\\
\dot\Nc^h_{u,(i)}&=\{(x,\ben)\in\dot\Nc^h;\ben\in\Bc_{u,(i)}^h\},
\hfill\\
\ddot\Nc^h_{s_i,u,(i)}&=
\{(x,\ben,\ben')\in\ddot\Nc_{s_i}^h;\ben,\ben'\in\Bc_{u,(i)}^h\}.
\hfill
\endmatrix$$
We have
$\Bc_{u,(i)}^h
%=\coprod_{w\in u}\Bc_w^h
\simeq
\tilde G^h\times_{B_{w,i}^h}\Bc_{w,(i)}^h,$
with $B_{w,i}^h=B_{w,i}\cap\tilde G^h.$
Therefore
$$\dot\Nc_{u,(i)}^h\simeq \tilde G^h\times_{B_{w,i}^h}\dot\Nc_{w,(i)}^h,
\quad
\ddot\Nc_{s_i,u,(i)}^h\simeq
\tilde G^h\times_{B_{w,i}^h}(\uen_{w,i}^h\times\Bc_{w,(i)}^h\times\Bc_{w,(i)}^h).$$

We must compute
the direct image by $p_2$ of the restriction to
$\ddot\Nc_{s_i,u,(i)}^h$ 
of the virtual bundle
$$-\Oc^h_{\dot\Nc^h}(\l')\Lambda(\Oc^h_{\dot\Nc^h}(\a_i))^{-1}\Ec\boxtimes
\Oc^h_{\dot\Nc^h_w}(\l'')\Lambda(\zeta^2\Oc^h_{\dot\Nc^h_w}(-\a_i)).$$
The projection 
$p_{2,w}:\uen_{w,i}^h\times\Bc_{w,(i)}^h\times\Bc_{w,(i)}^h\to
\ben_{w,i}^h\times\Bc_{w,(i)}^h$,
$(x,\ben,\ben')\mapsto(x,\ben')$
maps into
$\dot\Nc_{w,(i)}^h$.
Since the forgetful map 
$K^{\tilde G^h}(\dot\Nc^h_w)\to K(\dot\Nc^h_w)$ is surjective,
by base change it is enough to compute the direct image by 
$p_{2,w}$
of the virtual bundle
$$-\Oc_{\uen^h_{w,i}}\boxtimes
\Oc^h_{\Bc_{w,(i)}^h}(\l')\Lambda(\Oc^h_{\Bc_{w,(i)}^h}(\a_i))^{-1}\Fc\boxtimes
\Oc^h_{\Bc_{w,(i)}^h}(\l'')\Lambda(\zeta^2\Oc^h_{\Bc_{w,(i)}^h}(-\a_i))
\leqno(2.4.9)$$
%over $\uen_{w,i}^h\times\Bc_{w,(i)}^h\times\Bc_{w,(i)}^h$
for any vector bundle $\Fc$ over $\Bc_{w,(i)}^h$.
Notice that $\Bc_{w,(i)}\simeq B_{w,i}/B_w$ as a $B_{w,i}$-variety, and that
$\Oc_{\Bc_{w,(i)}}(\l)$ is the line bundle induced from the character 
${}^w\l$ of the torus $\tilde T\subset B_w$.
So $h$ acts on $\Oc_{\Bc^h_{w,(i)}}(\l)$ 
by fiberwise multiplication by the scalar ${}^w\l(\tilde s)$.

Let $p$ be the projection $\Bc_{w,(i)}^h\to\{point\}$.
If ${}^w\alpha_i(\tilde s)=1$ then
$\Bc_{w,(i)}^h=\Bc_{w,(i)}\simeq\PP^1$ 
and $\Lambda(\Oc^h_{\Bc^h_{w,(i)}}(\a_i))=1$, hence
$$p^*Rp_*(
\Lambda(\Oc^h_{\Bc^h_{w,(i)}}(\a_i))^{-1}\Fc)=
{\Fc-{}^{s_i}\Fc\Oc^h_{\Bc^h_{w,(i)}}(\alpha_i)
\over 1-\Oc^h_{\Bc^h_{w,(i)}}(\alpha_i)}.$$
Else, we have
$\Bc_{w,(i)}^h=\{\ben_w,\ben_{ws_i}\}$ and 
$\Lambda(\Oc^h_{\Bc^h_{w,(i)}}(\a_i))=
1-\Oc^h_{\Bc^h_{w,(i)}}(\alpha_i),$ hence
$$p^*Rp_*(
\Lambda(\Oc^h_{\Bc^h_{w,(i)}}(\a_i))^{-1}\Fc)=
{\Fc\over 1-\Oc^h_{\Bc^h_{w,(i)}}(\alpha_i)}+
{{}^{s_i}\Fc\over 1-\Oc^h_{\Bc^h_{w,(i)}}(-\alpha_i)}=
{\Fc-{}^{s_i}\Fc\Oc^h_{\Bc^h_{w,(i)}}(\alpha_i)
\over 1-\Oc^h_{\Bc^h_{w,(i)}}(\alpha_i)}.$$
Further, if ${}^w\alpha_i(\tilde s)\neq\zeta^2$ then 
$\uen_{w,i}^h=\uen_w^h$, hence
$\dot\Nc_{w,(i)}^h=\uen_{w,i}^h\times\Bc^h_{w,(i)}$,
while if
${}^w\alpha_i(\tilde s)=\zeta^2$ then $\uen_{w,i}^h\subsetneq\uen_w^h$ and
$\Lambda(\zeta^2\Oc^h_{\Bc^h_{w,(i)}}(-\a_i))=1.$
In both cases we have
$$\Oc_{\uen^h_{w,i}}\boxtimes\Lambda(\zeta^2\Oc^h_{\Bc^h_{w,(i)}}(-\a_i))=
1-\zeta^2\Oc^h_{\dot\Nc^h_{w,(i)}}(-\a_i)$$
in $K(\dot\Nc_{w,(i)}^h)$.
Therefore, the direct image by $p_{2,w}$ of 2.4.9 is
$$
-{(\Oc^h_{\dot\Nc^h_{w,(i)}}(\lambda'')-
\zeta^2\Oc^h_{\dot\Nc^h_{w,(i)}}(\lambda''-\a_i))
(\Fc\Oc^h_{\dot\Nc^h_{w,(i)}}(\lambda')
-{}^{s_i}\Fc\Oc^h_{\dot\Nc^h_{w,(i)}}(\alpha_i+{}^{s_i}\lambda'))
\over 1-\Oc^h_{\dot\Nc^h_{w,(i)}}(\alpha_i)}=
$$
$$
=-\Oc^h_{\dot\Nc^h_{w,(i)}}(\a_i){(1-\zeta^2\Oc^h_{\dot\Nc^h_{w,(i)}}(-\a_i))
(\Fc-{}^{s_i}\Fc)
\over 1-\Oc^h_{\dot\Nc^h_{w,(i)}}(\alpha_i)} 
$$
(we use the same notation for the pull-back of $\Fc$ to 
$\dot\Nc_{w,(i)}^h$).
Formula 2.4.8 is proved.

\qed

\vskip3mm

\noindent{\bf 2.4.10. Remarks.}
(a)
The map $\Phi$ is the composition of $IM$ and the homomorphism
$\hat\Hb\to \widehat K(\ddot\Nc^h)$ in \cite{V1}.
In loc.~cit.~ the element
$\Phi(t_{s_i}+t^{-1})$ is defined via a 
version of Thomason concentration map, see loc.~cit., sec.~A.3.3.
Formula 2.4.6 is obtained by computing explicitely this element.
Details are left to the reader.

(b)
Assume that $m=h$, $k>0$.
Lusztig's element $e_{cox}=\sum_{\a}e_\a$,
with $\a\in\{(a,\ell)\in\Delta^+_\re;\rhov\cdot a+h\ell=k\}$,
is HERS of type $\cox$ and degree $c$.
The vector space $H_*(\Bc^c_{e_{cox}},\CC)$ is an irreducible 
$\hat\Hb_c$-module by \cite{V1}.
It restricts to an irreducible $\Hb_c$-module.
If $k=1$ this $\Hb_c$-module is equal to $\Lb_c^\dag$.
It has dimension $k^n$ by \cite{F, prop.~1}.
The dimension of $H_*(\Bc^c_{x},\CC)$
for other types is considered in the last section.

(c)
The group $\tilde G^h$ may not be a Levi factor of any 
parahoric subgroup of $\tilde G$.
It is an affine analogue of the pseudo-Levi subgroups of a reductive group.
See \cite{MS, sec.~9} for instance.

\subhead 2.5. Isogenies and simple modules of DAHA's \endsubhead

The geometric description of the DAHA in \cite{V1} 
involves the algebra $\hat\Hb$.
In this paper we are interested by $\Hb$ and its rational degeneration
$\Hb''$. The purpose of this section is to recall the classification
of elements of $\Irr(O(\hat\Hb))$, proved in \cite{V1}, and to explain
how it yields a classification of $\Irr(O(\Hb))$.
The way to pass from $\hat\Hb$ to $\Hb$ is rather classical, and is
closely related to Clifford theory. In the affine Hecke algebra case a similar
construction was done in \cite{RR}, \cite{R7}.

First, recall the classification of $\Irr(O_h(\hat\Hb))$, for $h\in\hat H$.
Let 
$$\pi^h_w:\dot\Nc_w^h\to\Nc^h$$ 
be the restriction of the map $\pi^h$ to $\dot\Nc_w^h$.
Put $d_{h,w}=\dim(\dot\Nc^h_w)$.
Since $\dot\Nc_w^h$ is a smooth variety,
the complex $\CC_{\dot\Nc^h_w}[d_{h,w}]$ is perverse.
Put $$\Lc_{h,w}=\pi^h_{w,*}(\CC_{\dot\Nc^h_w}[d_{h,w}]).$$
By Gabber's theorem there are finite dimensional vector spaces $L_{h,w,\chi,i}$
such that
$$\Lc_{h,w}\simeq\bigoplus_{\chi,i}L_{h,w,\chi,i}\otimes IC(\chi)[i],
$$
with $\chi$ a $\tilde G^h$-equivariant irreducible local system supported on
a locally closed subset of $\Nc^h$, and $IC(\chi)$ the simple perverse sheaf
over $\Nc^h$ associated to $\chi$.
%$\sup(L_{h,\chi})=\{w;L_{h,w,\chi}\neq 0\},$
Let $\Xc_{h,w}$ be the set of simple perverse sheaves
$IC(\chi)$ which occur in $\Lc_{h,w}$.
We may identify $\Xc_{h,w}$ with a set of 
irreducible $\tilde G^h$-equivariant local systems
$\{\chi\}$ such that $\Xc_{h,w}=\{IC(\chi)\}$, 
hopping it will not create any confusion.
Put $\Xc_h=\bigcup_w\Xc_{h,w}$ and $\Lc_h=\bigoplus_w\Lc_{h,w}$.
Put also $L_{h,w,\chi}=\bigoplus_i L_{h,w,\chi,i}$
and $L_{h,\chi}=\bigoplus_w L_{h,w,\chi}$.

For each $x\in\Nc^h$ set 
$$A(h,x)=\pi_0(Z_{\tilde G^h}(x)),
\quad
A(h,x)^+=\pi_0(Z_{G^{h}_d}(x)).$$ 
They are finite groups containing the center of $G_0$.
Since the groups $A(h,x)$, $A(h,x)^+$ depend only on the image of
$h$ in $\tilde T$, $T_d$ respectively, we may abreviate
$A(\tilde s,x)$, $A(\tilde s,x)^+$.
For a future use we set $$A(h,x)^\circ=A(h,x)/Z(G_0).$$

If $\tau^k\neq\zeta^{2m}$ for each $m,k>0$ the set $\Nc^h$ 
consists of a finite number of nilpotent $\tilde G^h$-orbits by 2.4.2(c).
For each orbit $O$,
taking the stalk at any element $x\in O$ identifies
the $\tilde G^h$-equivariant
irreducible local system $\chi$ on $O$
with the irreducible $A(h,x)$-module $\chi_x$.
The pair $(h,\chi)$ is called the Langlands parameter of $L_{h,\chi}$.
The group $A(h,x)$ acts obviously on the space $H_*(\Bc^h_x,\CC).$
The actions of $\hat\Hb$, $A(h,x)$ centralize each other.

We'll say that the Langlands parameters $(h,\chi)$,
$(h',\chi')$ are $\tilde G$-conjugate if
$\tau=\tau'$, $\zeta=\zeta'$, and there is an element
$g\in\tilde G$ such that
$s'=(\ad g)(s)$,
$\chi=(\ad g)^{-1}(\chi')$ 
(where $(\ad g)^{-1}$ is the pull back of constructible sheaves
by the map $\ad g$).
Notice that, since $s,s'\in T_0$ the elements $h,h'$ are in the same
$W$-orbit relatively to the action 2.1.1, by the axioms
of split BN-pairs.

For each $x\in\Nc^h$
let $IC(\chi)_x$ be the stalk at $x$ of the cohomology sheaves
of the complex $IC(\chi)$.
It is a graded vector space.
The following is proved in \cite{V1, thm.~7.6}.

\proclaim{2.5.1. Proposition}
(a)
The vector space $L_{h,\chi}$
has the structure of a simple module in $O_h(\hat\Hb)$.
We have $\Irr(O_h(\hat\Hb))=\{L_{h,\chi};\chi\in\Xc_h\}$,
and two simple modules are isomorphic iff 
their Langlands parameters are $\tilde G$-conjugate.

(b)
If $\tau^k\neq\zeta^{2m}$ for each $m,k>0,$ 
then $\Xc_h$ is identified with the set of $\tilde G^h$-conjugacy classes
of representations
in $\Irr(A(h,x))$, where $x$ varies in $\Nc^h$,
which are Jordan-H\"older factors
of the module $H_*(\Bc^h_x,\CC)$.

(c)
The $\hat\Hb$-module $H_*(\Bc^h_x,\CC)$ has a finite length, 
and its class in the
Grothendieck group is $\sum_{\chi}L_{h,\chi}\otimes IC(\chi)_x$.
\endproclaim

Further, 
if $\tau^k\neq\zeta^{2m}$ for each $m,k>0$,
given any $w$ there are surjective maps
$\Hom_{A(h,x)}(\chi,H_*(\Bc^h_{x,w},\CC))\to L_{h,w,\chi}.$
It yields a surjective $\Hb$-modules homomorphism
$$\Hom_{A(h,x)}(\chi,H_*(\Bc^h_{x},\CC))\to L_{h,\chi}.
\leqno(2.5.2)$$

Now, 
the set $\Irr(O(\Hb))$
may be recovered from
$\Irr(O(\hat\Hb))$
as in
\cite{RR, appendix}, \cite{R7} for the affine Hecke algebras.
Given $s\in T_0$, $z\in\CC^\times_{\omega_0}$
we write 
$$s_z=(s,z)\in\hat T_0,
\quad
\tilde s_z=(\tilde s,z)\in\tilde T,
\quad
h_z=(\tilde s_z,\zeta)\in\hat H.$$
We have the following
(see also \cite{V1, sec.~8.3}).

\proclaim{2.5.3. Theorem}
Assume that $\tau^k\neq\zeta^{2m}$ for each $m,k>0.$ 
There is a bijection between $\Irr(O_{h}(\Hb))$
and the set of $G^h$-conjugacy classes of representations 
$\chi^+\in\Irr(A(h,x)^+)$, with $x\in\Nc^h$, such that
$\chi^+$ appears in $H_*(\Bc_x^h,\CC)$. 
\endproclaim
 
\noindent{\sl Proof :}
%$$x_{2\o_0} t_{s_0}x_{2\o_0}^{-1}=x_{\theta}^{-1}t_{s_0}x_{\theta},
%\quad
%x_{2\o_0} t_ix_{2\o_0}^{-1}=t_i,\ \forall i\neq 0.$$
Write $\tilde G^{\tilde s}$ for $\tilde G^{\tilde s_z}$
(it is independent of the choice of $z$).
The group
$$Z_{\tilde G}(\tilde s)^+=
\{g\in\tilde G;(\ad g)(\tilde s_z)\in\tilde s_z\CC^\times_{\omega_0}\}$$
is well-defined (it is independent of the choice of $z$).
The axioms of a split BN pair imply that
$Z_{\tilde G}(\tilde s)^+$ is generated by $\tilde T$,
the root subgroups associated to the roots in $\Delta_{\tilde s}$,
and the $n_w$'s with $w\in W$ such that
${}^w\tilde s_z\in\tilde s_{zz'}$ for some $z,z'\in\CC^\times_{\omega_0}$.
Thus it is a linear group which contains
$\tilde G^{\tilde s}$
as a closed normal subgroup.
Further, since the latter is generated by $\tilde T$
and the root subgroups associated to the roots in 
$\Delta_{\tilde s_z}$ by 2.4.1(a),
the quotient
$Z_{\tilde G}(\tilde s)^+/\tilde G^{\tilde s}$
is a finite group.

For each $x\in\Nc^h$ we set
$Z_{\tilde G}(\tilde s,x)^+=
Z_{\tilde G}(\tilde s)^+\cap Z_{\tilde G}(x).$
Observe that 
$$A(\tilde s,x)^+\simeq\pi_0(Z_{\tilde G}(\tilde s,x)^+),$$
because
$Z_{\tilde G}(\tilde s,x)^+$
is the inverse image of $Z_{G_d^{\tilde s}}(x)$
by the obvious projection $\tilde G\to G_d$.
Thus there is an exact sequence of groups
$$1\to A(\tilde s,x)\to A(\tilde s,x)^+\to B(\tilde s,x)\to 1,$$
where $B(\tilde s,x)$ is the image of $Z_{\tilde G}(\tilde s,x)^+$ 
in the finite group 
$Z_{\tilde G}(\tilde s)^+/\tilde G^{\tilde s}.$
In particular the group $A(\tilde s,x )^+$ acts by conjugation on
$A(\tilde s,x)$, yielding an action $\chi\mapsto{}^g\chi$ of
$B(\tilde s,x)$ on the set $\Irr(A(\tilde s,x))$.
For each $\chi$ let $E(\tilde s,\chi)$ be the endomorphism ring of the
induced $A(\tilde s,x)^+$-module
$$I(\chi)=\Ind_{A(\tilde s,x)}^{A(\tilde s,x)^+}(\chi).$$
Observe that if $g\in B(\tilde s, x)$ then $I({}^g\chi)\simeq I(\chi)$ as
$A(\tilde s,x)^+$-modules.
Thus
$$E(\tilde s,{}^g\chi)\simeq E(\tilde s,\chi)$$
as algebras, yielding a bijection
$\Irr(E(\tilde s,\chi))\simeq\Irr(E(\tilde s,{}^g\chi)),$
$\psi\mapsto{}^g\psi$.
Equip the $\hat\Hb$-module $L_{h_z,\chi}$ with the restriction of
the action to $\Hb$.

\proclaim{2.5.4. Lemma}
(a)
There is an algebra isomorphism
$E(\tilde s,\chi)=\End_{\Hb}(L_{h_z,\chi}).$

(b)
The restriction of
$L_{h_z,\chi}$ to $\Hb$
is a sum 
$\bigoplus_\psi L_{h_z,\chi}^\psi\otimes\psi$
with $(\psi,L_{h_z,\chi}^\psi)\in
\Irr(E(\tilde s,\chi))\times\Irr(O_h(\Hb)).$ 

(c)
Each module in $\Irr(O_h(\Hb))$ arises in this way.
Further 
$L_{h_z,\chi_1}^{\psi_1}\simeq L_{h_z,\chi_2}^{\psi_2}$
iff 
there is an element $g\in B(\tilde s,x)$ such that
$\chi_2={}^g\chi_1$,
$\psi_2={}^g\psi_1$.
\endproclaim

To finish the proof we must check that the assignement 
$$(z,x,\chi,\psi)\mapsto (x,\chi^\psi),
\leqno(2.5.5)$$
where $\chi^\psi$ is as in 2.5.10 below,
factors to a bijection from the set of quadruples
$(z,x,\chi,\psi)\in\CC^\times\times\Nc^h\times\Irr(A(\tilde s,x))\times
\Irr(E(\tilde s,\chi))$ 
such that $\chi$ appears in $H_*(\Bc_x^h,\CC)$ 
modulo the equivalence 
$$(z_1,x_1,\chi_1,\psi_1)\sim
(z_2,x_2,\chi_2,\psi_2)
\iff
\exists g\in G^{\tilde s}_d
\ \roman{s.t.}\ 
\left\{\aligned
&x_2=(\ad g)(x_1),
\hfill\\
&\chi_2={}^g\chi_1, 
\hfill\\
&\psi_2={}^g\psi_1\hfill\endaligned\right.$$ 
to the set of $G^{\tilde s}$-conjugacy classes of pairs
$(x,\chi^+)\in\Nc^h\times\Irr(A(\tilde s,x)^+)$ 
such that $\chi^+$ appears in $H_*(\Bc_x^h,\CC)$. 
Here, when writing $\chi_2={}^g\chi_1$,
$\psi_2={}^g\psi_1$
we are identifying
$A(\tilde s,x_1)$ with its image in $A(\tilde s,x_1)^+$ and we 
consider the bijection 
$A(\tilde s,x_1)^+\to A(\tilde s,x_2)^+$ induced by $(\ad g)$.
This notation is coherent with the one before 2.5.4.
Indeed, if $x_1=x_2$ then $g\in Z_{G^{\tilde s}_d}(x_1)$.
So, identifying $g$ with its image by the chain of maps
$$
Z_{G^{\tilde s}_d}(x_1)
\to
A(\tilde s,x_1)^+=\pi_0(Z_{G^{\tilde s}_d}(x_1))
\to
B(\tilde s,x_1),$$
the elements ${}^g\chi_1$, ${}^g\psi_1$ are the same as above.

To prove that 2.5.5 is one-to-one, 
observe that both 2.5.5 and the assignement 
$(x,\chi^+)\mapsto(z,x,\chi,\psi)$,
where $\chi$ is any composition factor in the restriction of $\chi^+$ to
$A(\tilde s,x)$ and $\psi=\Hom_{A(\tilde s,x)}(\chi,\chi^+)$,
respect the equivalence relations. 
Therefore, it is enough to prove that
$\chi^+$ appears in $H_*(\Bc_x^h,\CC)$ 
iff $\chi^+=\chi^\psi$ for some pair
$(\chi,\psi)$ where $\chi$ appears in $H_*(\Bc_x^h,\CC)$.
This follows from the proof of 2.5.4(a).

\qed

\vskip3mm

Before to prove 2.5.4 we recall some standard material.
See \cite{RR}.
The map
$$f:
\
Z_{\tilde G}(\tilde s)^+/\tilde G^{\tilde s}\to
\CC^\times_{\omega_0}\simeq\CC^\times,
\
g\mapsto (\ad g)(\tilde s_z)\tilde s_{z}^{-1}
$$
is a closed embedding of groups.
Its image, $C(\tilde s)$,
is a finite subgroup of $\CC^\times$
which does not depend on $z$.
The group $\CC^\times$ acts by algebra automorphisms on $\hat\Hb$ 
in such a way that the action on $\Hb$ is trivial and 
$z'\cdot x_{\o_0}=z' x_{\o_0}$ for each $z'\in\CC^\times$.
The algebra $\Hb$ is the fixed points subset $\hat\Hb^{\CC^\times}$.
Twisting a $\hat\Hb$-module by the automorphism $z'$
yields a map 
$$\Irr(O_{h_z}(\hat\Hb))\to\Irr(O_{h_{zz'}}(\hat\Hb)).$$
So for any module
$M\in\Irr(O_{h_z}(\hat\Hb))$
the inertia group
$$C(M)=\{z'\in\CC^\times;{}^{z'}\!M\simeq M\}$$
is contained into the finite group $C(\tilde s).$
Choosing $\hat\Hb$-module isomorphisms
$$\phi_{z'}:M\to{}^{z'}\!M,\quad\forall z'\in C(M),$$ 
we define a cocycle
$\eta:\ C(M)\times C(M)\to\CC^\times$
such that
$\phi_{z_1}\circ\phi_{z_2}=\eta(z_1,z_2)\phi_{z_1z_2}$.
Here ${}^{z'}\!M$
is canonically identified with $M$ as a vector space.
Let $C(M)_\eta$ be the $\CC$-algebra spanned by linearly independent
elements $e^{z'}$, $z'\in C(M)$, with the multiplication rule
$$e^{z_1}e^{z_2}=\eta(z_1,z_2)e^{z_1z_2}.$$
Its isomorphism class 
is independent of the choice of the isomorphisms
$\phi_{z'}$.

\proclaim{2.5.6. Lemma}
The restriction of each $M\in\Irr(O_{h_z}(\hat\Hb))$ to $\Hb$ is a sum 
$\bigoplus_\psi M^\psi\otimes\psi,$
where $\psi\in\Irr(C(M)_\eta)$ and 
$M^\psi\in\{0\}\cup\Irr(O_h(\Hb))$.
%The modules in $\Irr(O_h(\Hb))$ are of the form $M^\psi$, and
%$M_1^{\psi_1}\simeq M_2^{\psi_2}$ iff there is $z'\in\CC^\times$ such that
%$M_1={}^{z'}\!M_2$ and $\psi_2=\psi_1$.
\endproclaim

\noindent{\sl Proof :}
The algebra $\Hb$ acts on $M$ by restriction, 
and $C(M)_\eta$ via the $\Hb$-modules isomorphisms 
$\phi_{z}:M\to {}^zM\simeq M$.
Since the group $C(M)$ is finite, the algebra $C(M)_\eta$
is semi-simple.
Thus 
$M=\bigoplus_\psi M^\psi\otimes\psi$
with $\psi\in\Irr(C(M)_\eta)$
and $M^\psi\in\mod(\Hb).$
We must prove that 
$M^\psi\in\Irr(O_{h}(\Hb))$.
For all $h'\in{}^ Wh$ the set
$\{w\in W;{}^wh_z\in h'\CC^\times_{\omega_0}\}$
is finite.
So the weight spaces of the restriction of $M$ to $\Hb$ are finite dimensional.
Hence the $\Hb$-module $M^\psi$ is admissible.
Thus it is enough to check that it is irreducible.

The proof is modelled after \cite{RR}
(we can't apply loc.~cit.~ because $\CC^\times$ is not a finite group).
%Taking $\CC^\times$-invariants in a rational $\CC^\times$-module
%is an exact functor.
Let $\psi^*\in\Irr(C(M)_{\eta^{-1}})$ be the module dual to $\psi$.
The tensor product $M\otimes\psi^*$
is a $\hat\Hb\rtimes C(M)$-module such that
$$xz(m\otimes t)=x\phi_z(m)\otimes zt,
\quad
\forall (x,z)\in\hat\Hb\times C(M).$$
%Taking $C(M)$-invariants is a quotient functor (in the sense of Gabriel)
%$$\mod(\hat\Hb\rtimes C(M))\to\mod(\hat\Hb^{C(M)}).$$
Let $\mod(\hat\Hb\rtimes\CC^\times)_\rat$ be the category of 
$\CC^\times$-equivariant $\hat\Hb$-modules
which are rational as
$\CC^\times$-modules.
The algebra $\hat\Hb$ is $\ZZ$-graded,
because the $\CC^\times$-action is rational.
Let $\rho$ be the linear character of $\CC^\times$.
The induced module 
(from the finite group $C(M)$ to the algebraic group $\CC^\times$)
$$\Ind(M\otimes\psi^*)=
\bigoplus_{i\in\ZZ}(\rho^{i}\otimes M\otimes\psi^*)^{C(M)},$$
equipped with the $\hat\Hb\rtimes\CC^\times$-action such that
$$xz(\rho^i\otimes m\otimes t)=z^{-i}\rho^{i-j}\otimes xm\otimes t,$$
for each $x\in\hat\Hb$
which is homogeneous of degree $j$ and for each
$z\in\CC^\times$,
belongs to $\mod(\hat\Hb\rtimes\CC^\times)_\rat.$
By Frobenius reciprocity we have
$$
\Ind(M\otimes\psi^*)^{\CC^\times}
=(M\otimes\psi^*)^{C(M)}=M^\psi.
$$
Therefore, we are reduced to prove the following

\itemitem{(i)}
the $\hat\Hb\rtimes\CC^\times$-module 
$\Ind(M\otimes\psi^*)$
is irreducible,

\itemitem{(ii)}
taking $\CC^\times$-invariants is a quotient functor
(in the sense of Gabriel)
$$\mod(\hat\Hb\rtimes\CC^\times)_\rat\to\mod(\Hb).$$
%So it takes a simple modules to a simple module or 0.

Let us prove (i).
Let $n_M$ be the order of the group $C(M)$.
Then $\Ind(M\otimes\psi^*)$ is a free $\CC[\rho^{\pm n_M}]$-module.
For each $z'\in\CC^\times$ the evaluation map $\rho^{n_M}\mapsto 1/z'$ factors
to a surjective $\hat\Hb\rtimes C(M)$-module homomorphism
$$ev_{z'}:\Ind(M\otimes\psi^*)\to{}^{z'}\!M\otimes\psi^*.$$
The rhs is irreducible, see \cite{M3, p.~203} for instance.
Thus if
$N\subset\Ind(M\otimes\psi^*)$ is a $\hat\Hb\rtimes\CC^\times$-submodule then 
$ev_{z'}(N)$ equals either $\{0\}$ or ${}^{z'}\!M\otimes\psi^*$.
If $ev_{z'}(N)=\{0\}$ for some $z'$ then
$ev_{z'}(N)=\{0\}$ for each $z'$, because $N$ is stable under the
$\CC^\times$-action, yielding
$$N\subset\bigcap_{z'}(\rho^{n_M}-z')\Ind(M\otimes\psi^*)=\{0\}.$$ 
Thus $ev_{z'}(N)={}^{z'}\!M\otimes\psi^*$ for each $z'$,
yielding $N=\Ind(M\otimes\psi^*)$.

Let us prove (ii).
Let $B$ denote $\hat\Hb$ with its natural
$(\Hb,\hat\Hb\rtimes\CC^\times)$-bimodule structure.
%and $C$ denote $\hat\Hb$ with its natural
%$(\Hb,\hat\Hb\rtimes\CC^\times)$-bimodule structure.
The functor
$$F_*:\mod(\Hb)\to\mod(\hat\Hb\rtimes\CC^\times)_\rat,\quad
N\mapsto\Hom_\Hb(B,N)$$
is right adjoint to the functor
$$F^*:\mod(\hat\Hb\rtimes\CC^\times)_\rat\to\mod(\Hb),\quad
N\mapsto N^{\CC^\times}.$$
Further we have 
$$F^*F_*(N)=\Hom_{\Hb}(B,N)^{\CC^\times}=\Hom_{\Hb}(\Hb,N)=N.$$
So $\Ker(F^*)$ is a localizing subcategory and $F^*$ yields an equivalence
$$\mod(\hat\Hb\rtimes\CC^\times)_\rat/\Ker(F^*)\simeq\mod(\Hb).$$

\qed

\vskip3mm

\noindent{\sl Proof of 2.5.4 :}
It is enough to prove part (a), because (b), (c) follow from it.
It is easy to prove that
$$H_*(\Bc^{\tilde s}_{x},\CC)\simeq {}^{f(g)}H_*(\Bc^{\tilde s}_x,\CC),
\quad\forall g\in Z_{\tilde G}(\tilde s,x)^+,$$
as $\hat\Hb$-modules.
Compare \cite{R7, prop.~2.6.1}.
Thus the group $A(\tilde s,x)^+$ acts on the space $H_*(\Bc^{\tilde s}_x,\CC)$.
Notice that this 
action does not centralizes the $\hat\Hb$-action, but it centralizes
the $\Hb$-action.
Further, 2.5.2 and section A.3
imply that $L_{h,\chi}$ is the top of
the $\hat\Hb$-module
$$\Hom_{A(\tilde s,x)^+}(I(\chi),H_*(\Bc^{\tilde s}_{x},\CC))\simeq
\Hom_{A(\tilde s_z,x)}(\chi,H_*(\Bc^{\tilde s}_{x},\CC)).$$
Therefore there is an algebra homomorphism
$$E(\tilde s,\chi)\to\End_{\Hb}(L_{h_z,\chi}).
\leqno(2.5.7)$$
The group $B(\tilde s,x)$ acts on 
$\Irr(A(\tilde s,x))$ in the obvious way.
Let $B(\tilde s,\chi)$ be the isotropy subgroup of $\chi$.
We have $|B(\tilde s,\chi)|=\dim E(\tilde s,\chi)$
by Mackey's theorem.
It is easy to check that
$${}^{f(g)}\!L_{h_z,\chi}\simeq L_{h_z,{}^g\chi},
\quad
\forall g\in B(\tilde s,x).$$
Compare \cite{R7, prop.~ 2.7.2}.
Thus we have
$$C({L_{h_z,\chi}})=f(B(\tilde s,\chi)),
\quad
C({L_{h_z,\chi}})_\eta\simeq E(\tilde s,\chi).$$
Therefore 2.5.6 yields 
$$\dim \End_\Hb(L_{h_z,\chi})\le
|\Irr(C(L_{h_z,\chi})_\eta)|=
\dim E(\tilde s,\chi).$$
Hence it is enough to prove that 2.5.7 is injective.

Recall the decomposition
$\Bc_x^{\tilde s}=\coprod_w\Bc^{\tilde s}_{x,w}$ in section 2.4.
%Each subvariety $\Bc^{\tilde s}_{x,w}$ is stable under 
%the action of the group $Z_{\tilde G^{\tilde s}}(x)$. 
Each subspace $H_*(\Bc^{\tilde s}_{x,w},\CC)$ of 
$H_*(\Bc^{\tilde s}_{x},\CC)$ is stable under the action of the group
$A(\tilde s,x)$. 
Fix $w$ such that $L_{h_z,w,\chi}\neq\{0\}$.
Then $\chi$ appears in $H_*(\Bc^{\tilde s}_{x,w},\CC)$.
Now the action of $Z_{\tilde G}(\tilde s,x)^+$ on $\Bc_x^{\tilde s}$
permutes the components $\Bc^{\tilde s}_{x,w}$. 
We claim that the stabilizer of any of them is $Z_{\tilde G^{\tilde s}}(x)$. 
Thus, the $A(\tilde s,x)^+$-submodule of
$H_*(\Bc_x^{\tilde s},\CC)$ generated by $H_*(\Bc_{x,w}^{\tilde s},\CC)$ 
is equal to 
$$\Ind_{A(\tilde s,x)}^{A(\tilde s,x)^+}(H_*(\Bc_{x,w}^{\tilde s},\CC)).$$
It follows that
$H_*(\Bc_x^{\tilde s},\CC)$ contains $I(\chi)$.
Thus 2.5.7 is injective.

It remains to verify the claim.
Assume that the element $g\in Z_{\tilde G}(\tilde s)^+$ 
preserves the component $\Bc_w^{\tilde s}$.
Modifying $g$ by an element of $\tilde G^{\tilde s}$ we may assume
that it normalizes $B^{\tilde s}_w$ and $\tilde T$.
Thus $B_w$ and $(\ad g)(B_w)$ are both fixed point of $B^{\tilde s}_w$
acting in $\Bc^{\tilde s}_w$ via the restriction of the 
$\tilde G^{\tilde s}$-action.
Hence $B_w=(\ad g)(B_w)$.
So $g\in B_w\cap N_{\tilde G}(T)=T$.
Therefore $g\in\tilde G^{\tilde s}$.
We are done.

\qed

\proclaim{2.5.8. Corollary}
(a)
In type $A$ the restriction yields a bijection 
$\Irr(O_{h_z}(\hat\Hb))\to\Irr(O_h(\Hb))$.

(b)
If $\zeta=\tau^{c/2}$ and $\Delta_0=\Delta_{0,(h)}$ 
the restriction yields a bijection 
$\Irr(O_{h_z}(\hat\Hb))\to\Irr(O_h(\Hb))$ 
for any $z$.

(c)
Any simple spherical $\Hb_c$-module 
is the restriction of a module in
$\Irr(O_{h^\dag_{c,z}}(\hat\Hb))$ for any $z$.
\endproclaim

\noindent{\sl Proof :}
Part (a) is obvious, because no local system enter in the classification of
$\Irr(O(\hat\Hb))$ in type $A$. See \cite{V1}.
Let us concentrate on part (b).
It is enough to check that
$A(\tilde s_z,x)=A(\tilde s,x)^+$ (see, f.i., 2.5.9 below).
We'll prove that $Z_{\tilde G}(\tilde s)^+=\tilde G^{\tilde s}$.
To do this, according to the description of the quotient
$Z_{\tilde G}(\tilde s)^+/\tilde G^{\tilde s}$
in the proof of 2.5.3, it is enough to check that 
$Z_W(\tilde s)$ is generated by the reflections with respect to the roots
in $\Delta_{\tilde s}$.
Fix $\muv\in\check X_0$ as in 2.2.3(a).
We have $Z_W(\tilde s)=Z_W(\check\mu,m)$, because $\tau$ is not a root of unity.
Formula 1.0.2 yields
$$\aligned
Z_W(\check\mu,m)
&=\{\xi_{\lambdav}w; \lambdav\in Y_0, w\in W_0,\muv-{}^w\muv=m\lambdav\}
\hfill\\
&\simeq\{w\in W_0;\muv-{}^w\muv\in m Y_0\}
\hfill\\
&\simeq Z_{W_0}(\muv_m).
\endaligned
$$
Further, since $G_0$ is simply connected, the group
$Z_{W_0}(\muv_m)$ is generated by the reflections with respect to the roots in
$$\{a\in\Delta_0;a(\muv_m)=1\}.$$
On the other hand, we have
$$\aligned
\Delta_{\tilde s}&=\{(a,\ell)\in\Delta;a\cdot\muv+m\ell=0\}
\hfill\cr
&\simeq\{a\in\Delta_0;a(\muv_m)=1\}.
\endaligned$$

Part (c) is a direct consequence of (b), because
$\Delta_0=\Delta_{0,(c)}$ and a spherical $\Hb_c$ belongs to $O_c(\Hb)$
by 2.1.4.

\qed

\vskip3mm

\noindent{\bf 2.5.9. Remark.}
The bijection in 2.5.3 is given as follows.
For each $z\in\CC^\times_{\omega_0}$,
$x\in\Nc^h$,
$\chi\in\Irr(A(\tilde s,x))$
the induced $A(\tilde s,x)^+$-module has the form
$$I(\chi)=
\bigoplus_\psi\chi^\psi\otimes\psi,
\leqno(2.5.10)$$
with
$\psi\in\Irr(E(\tilde s,\chi))$,
$\chi^\psi\in\Irr(A(\tilde s,x)^+).$
Now, given 
$\chi^+\in\Irr(A(\tilde s,x)^+)$,
we fix $z$,
$\chi$,
$\psi$ 
such that $\chi^+=\chi^\psi$.
Then the module
$L_{h,\chi^+}\in\Irr(O_h(\Hb))$ with Langlands parameters
$(h,\chi^+)$ is equal to
$L^\psi_{h_z,\chi}$.

\subhead 2.6. Fourier-Sato transform\endsubhead

In this section we assume that $\zeta=\tau^{c/2}$ to simplify.
Fix $h\in\hat H$.
Recall that $\Irr(O_h(\hat\Hb))=\{L_{h,\chi};\chi\in\Xc_h\}$ by 2.5.1(a).
If $k<0$ then $\Xc_h$ is identified with
the set of $\tilde G^h$-conjugacy classes of Jordan-H\"older factors of the
$A(h,x)$-modules $H_*(\Bc^h_x,\CC)$, with $x\in\Nc^h$, by 2.5.1(b).
We'll need the sets $\Xc_h$ with $k>0$.
Since there may be an infinite number of
$\tilde G^h$-orbits in $\Nc^h$, the set $\Xc_h$ is
more difficult to describe than for negative $k$'s.
This is done via the Fourier-Sato transform.

Write $\CC_{h,0},\CC_{h,\reg}$ for the constant sheaves 
$\CC_{\{0\}},\CC_{\Nc^h_\reg}$
on $\{0\},\Nc^h_\reg$ respectively.
Write $L_{h,0}$, $L_{h,\reg}$ 
for the $\hat\Hb$-modules 
$L_{h,\CC_{h,0}}$, $L_{h,\CC_{h,\reg}}$
whenever they are defined.

Fix a nondegenerate $(\ad \tilde G)$-invariant 
symmetric pairing $(\ :\ )$ on $\gen$.
We have $\tilde G^h=\tilde G^{h^{\ddag}}$.
So the bilinear form $(\ :\ )$ restricts to a nondegenerate
$(\ad \tilde G^h)$-invariant bilinear form $\gen^h\times\gen^{h^{\ddag}}\to\CC$.
Let $D^b_{\RR^+}(\gen^h)$ be the derived category of
bounded complexes of constructible sheaves over $\gen^h$ with conic cohomology
(i.e., the cohomology sheaves are constant on the $\RR_+$-orbits).
The Fourier-Sato transform yields an equivalence
$D^b_{\RR^+}(\gen^h)\to D^b_{\RR^+}(\gen^{h^\ddag}).$
Since $\tilde G^h\times\CC^\times$-equivariant sheaves are conic,
it yields also an equivalence
$FS:D^b_{\tilde G^h\times\CC^\times}(\gen^h)\to
D^b_{\tilde G^{h^{\ddag}}\times\CC^\times}(\gen^{h^{\ddag}}).$
See \cite{B3, sec. 6}, \cite{KS1, sect.~3.7} for details on the functor $FS$.
The functor $FS$ is difficult to compute in general.
The following particular case will be important for us. 

Set $h=h_c^\dag$.
Since $(h_c^\dag)^\ddag=h_{-c}^{-1}$,
we have $\tilde G^{(h_c^\dag)^\ddag}=\tilde G^{h_{-c}}$ and
$\gen^{(h_c^\dag)^\ddag}=\gen^{h_{-c}}$.
Thus the functor $FS$ yields an equivalence 
$$D^b_{\tilde G^c\times\CC^\times}(\gen^c)
\to D^b_{\tilde G^{h_{-c}}\times\CC^\times}(\gen^{h_{-c}}).$$

\proclaim{2.6.1. Lemma}
Let $k>0$.

(a)
There is a bijection
$\Xc_{c,w}\to\Xc_{h_{-c},w}$, $\chi\mapsto\chi^\dag$
such that $FS(IC(\chi))$, $IC(\chi^\dag)$ are isomorphic.

(b)
For each $\chi\in\Xc_c$ we have
$L_{c,\chi}\simeq{}^{IM}L_{h_{-c},\chi^\dag}$
as $\hat\Hb_{c}$-modules.
\endproclaim

\noindent{\sl Proof :}
Claim (a) is proved in \cite{V1, sec.~7.4}.
See also section A.1.
Part (b) is conjectured in loc.~cit.
The modules ${}^{IM}L_{c,\chi}$, $L_{h_{-c},\chi^\dag}$
belong to $O(\hat\Hb_{-c})$.
They have the same weights by A.1.1(b).
Thus they are isomorphic by A.4.1.

\qed

\vskip3mm

\proclaim{2.6.2. Lemma}
Let $k>0$.

(a)
We have $\CC_{h_{-c},0}\in\Xc_{h_{-c}}$.
If $\CC_{h_{-c},0}\in\Xc_{h_{-c},w}$ then 
$\CC_{\Nc^c}[\dim\Nc^c]$ belongs to $\Xc_{c,w}$.
If $m\in\RN$ then
$\CC_{c,\reg}\in\Xc_c$
and
$(\CC_{c,\reg})^\dag=\CC_{h_{-c},0}$.

(b)
If $m\in\RN$ then
$L_{c,\reg}$ is a Jordan-H\"older
composition factor of
$H_*(\Bc^c_x,\CC)$ for each $x\in\Nc^c_\reg.$
In particular
$\dim(L_{c,\reg})<\infty$ if $m\in\EN$. 

(c)
We have 
$\dim(L_{h_{-c},0})<\infty$ iff $m\in\EN$. 
\endproclaim

\noindent{\sl Proof :}
The first claim of (a) follows from 2.5.1(b). 
Since $\Nc^c$ is a vector space we have
$$FS(\CC_{\Nc^c}[\dim(\Nc^c)])=\CC_{h_{-c},0}=IC(\CC_{h_{-c},0}),$$ 
yielding the 2-nd claim.
If $m\in\RN$ then
$\Nc^c_\reg\neq\emptyset$ by 1.3.3.
Thus $IC(\CC_{c,\reg})=\CC_{\Nc^c}[\dim\Nc^c]$,
and
$\CC_{c,\reg}\in\Xc_c$,
$(\CC_{c,\reg})^\dag=\CC_{h_{-c},0}$ by 2.6.1(a).

Now, let us prove part (b).
The first claim is a general fact which follows from the sheaf theoretical
approach to convolutions algebras, see \cite{CG, sec.~8.5}.
The second one is obvious because
$\dim H_*(\Bc^c_x,\CC)<\infty$ if $m\in\EN$. 

Finally we prove (c).
It is enough to prove that if $\dim(L_{h_{-c},0})<\infty$ then $m\in\EN$. 
By 2.6.1(b) and part (a) we have
$\dim(L_{c,\reg})=\dim(L_{h_{-c},0})$.
Thus we must prove that 
$L_{c,\reg}$ is infinite dimensional if $m\notin EN$.
It is enough to check that 
if $m\notin\EN$ then
$\CC_{\Nc^{c}}[\dim\Nc^{c}]\in\Xc_{c,w}$
for an infinite number of $w$'s.
To do that, first notice the following basic fact. 

\proclaim{2.6.3. Lemma}
Let $\pi:\Xc\to\Yc$ be a surjective projective morphism of smooth connected
complex algebraic varieties.
The constant sheaf $\CC_\Yc$ is a direct summand of the complex
$R\pi_*(\CC_\Xc)$.
\endproclaim

The variety $\Nc^{c}$ is smooth and irreducible by 2.4.2(b).
Therefore it is enough to check that if $m\notin\EN$ the map 
$\pi^{c}_w$ is surjective 
for an infinite number of $w$'s.
In the rest of the proof we'll
use the terminology of section 3, to which we refer for more details.
There, we introduce a partition of $W^c$ into clans such that

\itemitem{(i)}
there is a finite number of clans,

\itemitem{(ii)}
if $w,w'$ are in the same clan the morphisms
$\pi_w^c$, $\pi^c_{w'}$ are the same by 3.4.1(a).

Set $E=\{w\in W^{c};\pi^{c}_w(\dot\Nc^c_w)=\Nc^c\}$.
The set $S=\bigcup_{w\notin E}\pi^{c}_w(\dot\Nc^c_w)$ 
is a proper closed subset of $\Nc^c$ such that
$$x\in\Nc^c\setminus S\Rightarrow
\{w\in W^{c};x\in\pi^{c}_w(\dot\Nc_{w}^{c})\}\subset E.$$ 
Further $H_*(\Bc^{c}_x,\CC)$ is finite dimensional iff 
$x\in\Nc^{c}_{\ERS}$ by 2.4.2(a),
and $\Nc^c_\ERS\neq\emptyset$ iff $m\in\EN$ by 1.3.3.  
Thus, if $m\notin\EN$ the set
$\{w\in W^{c};x\in\pi^{c}_w(\dot\Nc_{w}^{c})\}$ is infinite
for each $x\in\Nc^{c}$. 
Hence, if $m\notin\EN$ the set $E$ is infinite. 
We are done.

\qed

\vskip3mm

\noindent{\sl Proof of 2.6.3 :}
By complex we mean a bounded complex of sheaves of
$\CC$-vector spaces with constructible cohomology.
Fix a non-empty open subset $\Yc'\subset\Yc$ such that
$\pi$ restricts to a locally trivial fibration 
$\pi':\Xc'=\pi^{-1}(\Yc')\to\Yc'$.
Let $i$ be the inclusion of $\Yc'$ in $\Yc$.

First observe that if $F$ is any semi-simple complex over $\Yc$
with locally constant cohomology sheaves over $\Yc'$ 
such that $i^*(F)$ contains $\CC_{\Yc'}$ 
as a direct summand,
then the complex $IC(\CC_{\Yc'})[-\dim(\Yc)]$ is a direct summand of $F$.

Now, set $F=R\pi_*(\CC_\Xc)$.
Since $\Xc$ is smooth the sheaf $\CC_\Xc$ is a semi-simple complex,
hence $F$ is semi-simple by Gabber's theorem.
Since $(\pi')^*,R\pi'_*$ are adjoint functors, 
there is a canonical map
$\CC_{\Yc'}\to R\pi'_*(\CC_{\Xc'})=i^*(F)$.
Since $\pi'$ is a locally trivial fibration, the complex
$i^*(F)$ is a direct sum of shifts of irreducible local systems
over $\Yc'$, i.e., $i^*(F)=\bigoplus_k\Lc_k[d_k]$ where
$\Lc_k$ is an irreducible local system and $d_k$ 
is a $\le 0$ integer.
By construction, the canonical map factors to an injective sheaf homomorphism 
$\CC_{\Yc'}\to\bigoplus_{d_k=0}\Lc_k$.
Thus $\Lc_k=\CC_{\Yc'}$ for at least one integer $k$ with $d_k=0$, i.e.,
$i^*(F)$ contains $\CC_{\Yc'}$ as a direct summand. 

Finally we have
$IC(\CC_{\Yc'})[-\dim(\Yc)]=\CC_\Yc$ because $\Yc$ is smooth. 
Therefore $\CC_\Yc$ is a direct summand of the complex $R\pi_*(\CC_\Xc)$.
We are done.

\qed

\vskip3mm

\noindent{\bf 2.6.4. Example.}
The case $m=h$, $k>0$ was already considered in \cite{V1, sec.~9.3}.
The vector space $\gen^c$ is
linearly spanned by the root vectors
$e_{\a_1},...e_{\a_n},e_{-\theta}\otimes\eps^k$.
Let $x_1,...x_{n+1}$ be the corresponding coordinates.
We have
$\gen^c_\reg=\{x_1...x_{n+1}\neq 0\}$,
$\sen^c=\{x_1=x_2=...x_n=1\},$
$G^{c,\circ}\times\CC_\delta^\times=(\CC^\times)^{n+1}$, and $A_m=Z(G_0)$.
Here $A_m$ is the group introduced in sections 3.1, 3.2 below.
The weights of the
$(\CC^\times)^{n+1}$-action on $\gen^c$ are
$(1,0,...0,-\theta\cdot\ov_1),$
$(0,1,0,...0,-\theta\cdot\ov_2),$
...
$(0,...0,1,-\theta\cdot\ov_n),$
and $(0,...0,k).$
The semi-simple complex $\Lc_c$
is an infinite sum of the simple extensions of
the trivial local systems on each orbit except $\{0\}$
(up to some shift).

\vskip3mm

%\noindent{\bf 2.6.5. Remark.}
%We have
%$W^c=W^{h_{-c}}$,
%$\tilde G^c=\tilde G^{h_{-c}}$,
%and
%$\Bc^c=\Bc^{h_{-c}}$.
%So the action of $\hat\Hb_{-c}$ on $H_*(\Bc^{h_{-c}},\CC)$ in
%2.4.3 (for $h=h_{-c}$, $x=0$) may be viewed as 
%an action of $\hat\Hb_{-c}$ on $H_*(\Bc^{c},\CC)$.
%It is different from the $\hat\Hb_{c}$-action on $H_*(\Bc^{c},\CC)$
%obtained by taking $h=h_c^\dag$, $x=0$.
%
%\vskip3mm

\subhead 2.7. Geometrization of the polynomial representation\endsubhead

First, let us recall the Borel map for affine flag manifolds.
Let $\Xc_\CC$ be the affine flag manifold of type $G_\CC$ studied in \cite{K2}.
It is an infinite dimensional (not quasi-compact)
$\tilde T$-equivariant scheme over $\CC$
paved by finite codimensional affine cells $\Xc_{\CC,w}$ with $w\in W$.
We have $\Bc\subset\Xc_0$.
Further the fixed-points subscheme $\Xc_\CC^h$ is 
locally of finite type, and its set of $\CC$-points is equal to $\Bc^h$. 
See \cite{V1, lem.~2.13} for details.
We have the Borel map
(see also \cite{KS2})
$$\quad R(\tilde T)\otimes_{\CC[q^{\pm 1}]}\CC X\to K_{\tilde T}(\Xc_\CC),
\quad
V\otimes x_\l\mapsto V\otimes\Oc_{\Xc_\CC}(-\l).$$ 
%Set $R(\tilde T)_{\ell oc}=R(\tilde T)\Sigma^{-1}$, where
%$\Sigma$ is the multiplicative set generated by
%$\{q^n-1\}_{n>0}$.
%The Borel map factors to an injection
%$$
%R(\tilde T)_{\ell oc}\otimes_{\CC[q^{\pm 1}]}\CC_q X\to
%R(\tilde T)_{\ell oc}\otimes_{R(\tilde T)}K_{\tilde T}(\Xc_\CC)\simeq
%\prod_{w\in W}
%R(\tilde T)_{\ell oc}\otimes\Oc_{\bar\Xc_{\CC,w}}
%$$
%whose image contains $\{\Oc_{\bar\Xc_{\CC,w}}\}$
%by \cite{KS2, thm.~4.4}.

Set $h=(s,\tau,\zeta)\in\hat H$ with $\zeta=\tau^{c/2}$.
Composing the Borel map with
the evaluation at $\tilde s$ and 
the restriction to the fixed points subset 
we get a linear map 
$\CC X\to\widehat K(\Bc^{\tilde s}).$
Composing it with the Chern character
and the pull-back by the vector bundle $\dot\Nc^h\to\Bc^{\tilde s}$
we get a map
$$\CC X\to\widehat H_*(\dot\Nc^h,\CC).
\leqno(2.7.1)$$
For each $w\in W^h$ we consider the ideal $I^{h,w}=I^{{}^{w^{-1}\!\!}h}$.
Write $I^{h,F}=\bigcap_{w\in F}I^{h,w}$ for each $F\subset W^h$,
and
$$\widehat\Pb_{X,h}^\dag=\pro_E(\Pb_{X}^\dag/I^{h,E}),
\quad
\widehat\Pb_{X,h}={}^{IM}\widehat\Pb_{X,h^\dag}^\dag
=\pro_E(\Pb_{X}/I^{h,E}),$$
where $E\subset W^h$ is any finite subset,
with the limit topology.

\proclaim{2.7.2. Lemma}
Assume that $\zeta=\tau^{c/2}$.

(a)
The $\hat\Hb$-actions on $\Pb_{X}$, $\Pb^\dag_{X}$,
$H_*(\Bc^{h},\CC)$, $H_*(\dot\Nc^h,\CC)$ 
extend uniquely to continuous 
actions on $\widehat\Pb_{X,h}$, 
$\widehat\Pb_{X,h}^\dag$, 
$\widehat H_*(\Bc^{h},\CC)$,
$\widehat H_*(\dot\Nc^h,\CC)$.

(b)
There is an isomorphism of topological $\hat\Hb$-modules
$\widehat\Pb_{X,h}^\dag\simeq\widehat H_*(\dot\Nc^h,\CC)$.

(c)
If $k<0$
there is an isomorphism of topological $\hat\Hb$-modules
$\widehat\Pb_{X,h}\simeq\widehat H_*(\Bc^{h},\CC)$.

\endproclaim

\noindent{\sl Proof :}
Part (a) is obvious.
The case of $H_*(\Bc^{h},\CC)$, $H_*(\dot\Nc^h,\CC)$
follows from 2.1.5(c), because both modules are admissible. 
The case of $\Pb_{X}$, $\Pb^\dag_{X}$
is proved as in 2.1.7(a).

Let us concentrate on Part (b).
For each $w$ set
$\Pb_{X,h,w}^\dag$ equal to $\Pb_{X}^\dag/I^{h,w}\Pb_{X}^\dag$.
We have
$\widehat\Pb_{X,h}^\dag=\prod_w\Pb_{X,h,w}^\dag$,
with the product topology.
The natural map $\CC X\to\Pb_{X}^\dag$
factors to an isomorphism
$$\CC X/I^{h,w}\to\Pb_{X,h,w}^\dag.$$
Further, the map 2.7.1 yields an isomorphism
$$\CC X/I^{h,w}\to H_*(\dot\Nc^h_w,\CC),$$
by 2.7.4(b) below.
Therefore, there is an isomorphism of topological vector spaces
$$\widehat\Pb_{X,h}^\dag\to\widehat H_*(\dot\Nc^h,\CC).$$
See also \cite{V1, lem.~4.8(ii)}.
It is an $\hat\Hb$-module isomorphism by 2.4.7.

Now we prove (c).
By part (b) there is an isomorphism of topological $\hat\Hb_{}$-modules
$${}^{IM}\widehat\Pb_{X,h}=\widehat\Pb_{X,h^\dag}^\dag
\simeq\widehat H_*(\dot\Nc^{h^\dag},\CC).$$
Further
$H_*(\dot\Nc^{h^\dag},\CC)\simeq{}^{IM}H_*(\Bc^{h},\CC)$ by A.2.10.
Thus there is an isomorphism of topological $\hat\Hb$-modules
$$\widehat\Pb_{X,h}\simeq\widehat H_*(\Bc^{h},\CC).$$

\qed

\vskip3mm

\proclaim{2.7.3. Lemma}
Assume that $\zeta=\tau^{c/2}$ and $\Delta_{0,(h)}=\Delta_0$.

(a)
Any $Y$-spherical module in $\Irr(O_h(\hat\Hb))$ or any 
simple finite dimensional quotient of
$\widehat\Pb_{X,h}$ comes from a simple $\Hb''_c$-module.

(b)
If $k<0$
there is a finite dimensional $Y$-spherical module in $\Irr(O_h(\hat\Hb))$
iff $\dim(L_{h,0})<\infty$.
Further this $\hat\Hb_c$-module is isomorphic to $L_{h,0}$, 
and its restriction to $\Hb_c$ is isomorphic to $\Lb_c$.
\endproclaim

\noindent{\sl Proof :}
First we prove (a).
Since the $\Hb'_c$-module
$\Pb'_{X,c}$ is free of rank one over $\CC[\xi_\l]$,
it is $\ZZ_{\ge 0}$-graded by the total degree in the $\xi_\l$'s.
The $\CC W$-action is locally finite and
for each $i\in I$, $\l_1, \l_2,...\,\l_r\in X_0$
we have
$$s_i(\xi_{\l_1}\xi_{\l_2}...\xi_{\l_r})=
\xi_{{}^{s_i}\l_1}\xi_{{}^{s_i}\l_2}...\xi_{{}^{s_i}\l_r}$$
modulo terms of lower degree.
Thus the action of $\xi_{\check\theta}$ is locally finite and
$$\xi_{\check\theta}(\xi_{\l_1}\xi_{\l_2}...\xi_{\l_r})=
\xi_{\l_1}\xi_{\l_2}...\xi_{\l_r}$$
modulo terms of lower degree by formula 1.0.2.
Therefore $\xi_{\lambdav}$ is locally unipotent on $\Pb'_{X,c}$
for each $\lambdav\in Y_0$, because $Y_0$ is spanned by the $W_0$-orbit of
$\check\theta$.
See \cite{C3, prop.~17.23} for instance.
So, any quotient of $\Pb'_{X,c}$
comes from $\Hb''_{c}$ by 2.3.1(b).

Now, let $M\in \Irr(O_h(\hat\Hb))$ be a $Y$-spherical module. 
The restriction of $M$ to $\Hb$ is simple and $Y$-spherical by 2.5.8(b).
Further $M$ comes from a simple quotient of $\Pb'_{X,c}$ by 2.2.4(a),(b). 
So $M$ comes from a simple $\Hb''_{c}$-module by the discussion above.

Let $M$ be a finite dimensional $\hat\Hb$-module.
Notice that $M$ belongs to $O(\hat\Hb)$ and is a finite sum of weights spaces.
%Thus there is a finite subset $E\subset W^h$ such that $M_{{}^wh}=\{0\}$
%for each $w\notin E$.
%Therefore we have $I^{h,E}M=\{0\},$
Further $M=\widehat M$ with the discrete topology.

Finally, let $M$ be a finite dimensional simple quotient of $\widehat\Pb_{X,h}$.
Set $\Pb_{X,h,w}={}^{IM}\Pb_{X,h^\dag,w}^\dag$,
where the rhs is as in the proof of 2.7.2.
We have $\widehat\Pb_{X,h}=\prod_w\Pb_{X,h,w}$,
with the product topology.
Since $M$ is a discrete quotient of $\widehat\Pb_{X,h}$, 
there is a finite subset $E\subset W^h$ such that 
$\bigoplus_{w\in E}\Pb_{X,h,w}$ maps onto $M$.
Since the ideal $I^{h,w}$ acts trivially on
$\Pb_{X,h,w}$, we have also $I^{h,E}M=\{0\}.$ 
Therefore $M\in O_h(\hat\Hb)$.
Let $M'$ be the image of $\Pb_X$ by 
the surjective $\hat\Hb_c$-module homomorphism $\widehat\Pb_{X,h}\to M$.
We have
$\widehat M'=\widehat M$.
So $M'=M$, because $M$ is discrete.
Hence $M\in \Irr(O_h(\hat\Hb))$ and it is $Y$-spherical.

Now, let us concentrate on part (b).
If $\dim(L_{h,0})<\infty$
then $L_{h,0}=\widehat L_{h,0},$  with the discrete topology,
and it is a quotient of $\widehat H_*(\Bc^{h},\CC)$ by 2.5.2.
Hence it is a simple finite dimensional quotient of
$\widehat\Pb_{X,h}$ 
by part 2.7.2(c). 
Hence it comes from a simple
finite dimensional spherical $\Hb''_c$-module by part (a)
and 2.2.4(b), 2.3.1(c),(d).
Hence its restriction to $\Hb_c$ is isomorphic to $\Lb_c$. 
So we must prove the following facts :

\vskip1mm

\itemitem{(i)}
any quotient $M$ of $\Pb_X$ 
which belongs to $\Irr(O_h(\hat\Hb))$
is also a quotient of the $\hat\Hb_c$-module $H_*(\Bc^{h},\CC)$,

\vskip1mm

\itemitem{(ii)}
the $\hat\Hb_c$-module $H_*(\Bc^{h},\CC)$ has a simple top equal to $L_{h,0}$.

\vskip1mm

Let $M\in\Irr(O_h(\hat\Hb))$.
Then $M$ is isomorphic to $L_{h,\chi}$ for some $\chi\in\Xc_h$
by 2.5.1(a).
So it is a subquotient of the module $H_*(\Bc^{h}_x,\CC)$ for some
$x\in\Nc^h$ by loc.~cit. 
Recall that
$$H_*(\Bc^{h}_{w,x},\CC)=
H^*(\Bc^{h}_{w},\Bc^{h}_w\setminus\Bc^{h}_{w,x},\CC)$$
by Poincar\'e duality.
Further, the $\CC X$ action on
$H_*(\Bc^{h}_{w,x},\CC)$
factors through the Borel map
$\CC X\to H^*(\Bc^{h}_{w},\CC)$
and the cup product in cohomology.
Therefore, since the ideal $I^{h,w}$ maps to 0 in
$H^*(\Bc^{h}_{w},\CC)$
we have $I^{h,w}M_{{}^{w^{-1}\!\!}h^\dag}=\{0\}$.
In particular we have
$$\widehat M=\pro_E(M/I^{h,E})$$
%\leqno()$$
by 2.1.7(b).
Hence, if $M$ is $Y$-spherical then the
surjective map $\Pb_X\to M$
factors to a surjective map 
$\widehat\Pb_{X,h}\to\widehat M$.

%If $M\in\Irr(O_h(\hat\Hb))$ is finite dimensional,
%then there is a finite subset $E\subset W^h$ such that $M_{{}^wh}=\{0\}$
%for each $w\notin E$, i.e., we have 
%$I^{h,E}M=\{0\}.$

Now, assume that $M$ is a quotient of $\Pb_X$ in $\Irr(O_h(\hat\Hb))$.
There is a surjective map $\widehat\Pb_{X,h}\to \widehat M$.
Hence there is a surjective $\hat\Hb$-homomorphism
$\widehat H_*(\Bc^{h},\CC)\to\widehat M,$ by 2.7.2(c).
Let $M'\subset\widehat M$ be the image of $H_*(\Bc^{h},\CC).$
Since $H_*(\Bc^{h},\CC)$ is admissible, $M'$ is contained into the sum
of the weight subspaces of $\widehat M$.
Thus $M'\subset M$.
Further $M'\neq\{0\}$, because $H_*(\Bc^{h},\CC)$ is dense in 
$\widehat H_*(\Bc^{h},\CC)$.
Hence $M'=M$, because $M$ is simple.
This yields (i).

Part (ii) is a particular case of a more general result saying that
if $k<0$ then standard modules have a simple top.
A proof is given in section A.3.

\qed

\vskip3mm

\noindent{\bf 2.7.4. Remarks.}
(a)
For future applications of 2.7.3, observe that $\Delta_{0,(c)}=\Delta_0$.

(b)
Fix a connected reductive group $G'$.
Let $B'$, $T'$, $h'$, $W'$ be a Borel subgroup, a maximal torus,
a central element, and the Weyl group respectively.
Then the Borel map gives an isomorphism
$$\CC[T']/(f-f(h');f\in\CC[T']^{W'})\to K(G'/B').$$
If $G'$ has a simply connected derived subgroup this is well-known.
Else, fix a connected reductive group $G''$ such that
$G'=G''/Z$ for some finite subgroup $Z\subset Z(G'')$
and $G''$ has a simply connected derived subgroup.
Fix also a maximal torus $T''\subset G''$ such that $T'=T''/Z$,
and an element $h''\in T''$ which maps to $h'$.
Consider the following ideals
$$J'=(f-f(h'))\subset\CC[T']^{W'},
\quad
J''=(f-f(h''))\subset\CC[T'']^{W'}.$$
The inclusion $\CC[T']\subset\CC[T'']$ factors to a map
$$\CC[T']/J'\CC[T']\to\CC[T'']/J''\CC[T'']\simeq K(G'/B').$$
Since the quotient morphism $T'\to T'/W'$ is finite and generically 
a Galois cover with group $W'$, the dimension of the lhs is 
larger or equal to
$$\dim_{\CC(T')^{W'}}(\CC(T'))=|W'|=\dim(K(G'/B')).$$
To prove that the left arrow is injective
we must check that 
$$(J''\CC[T''])^Z=J'\CC[T'].$$
Observe that $J''\CC[T'']$ is the ideal of the
scheme-theoretic fiber $\xi$ of the quotient
$T''\to T''/W''$ over the closed point $\{h''\}$ of $T''/W'$.
Further, since the ideal of the closed subset $h''Z\subset T''/W'$
is equal to $J'\CC[T'']^{W'},$ 
the ideal of the closed subscheme $Z\xi\subset T''$ 
is equal to $J'\CC[T'']$.
Taking $Z$-invariants we get the equality above.
%Since $z_1\xi\cap z_2\xi=\emptyset$ for each $z_1\neq z_2\in Z$, we have

\vskip3mm

\subhead 2.8. Classification of spherical simple finite dimensional modules
coming from $\Hb''$
\endsubhead

The main result of the paper is the following.

\proclaim{2.8.1. Theorem}
(a)
If $m\in\EN$, $k<0$ then $L_{h_c,0}$
is spherical finite dimensional and simple.
Its restriction to $\Hb_c$ is isomorphic to $\Lb_c$.

(b)
If $M$ is a spherical finite dimensional simple $\Hb_c$-module 
which comes from $\Hb''_c$ then 
$m\in\EN$, $k<0$ and $M$ is isomorphic to the restriction of
$L_{h_c,0}$ to $\Hb_c$.
\endproclaim

\noindent{\sl Proof :}
If $m\in\EN$, $k<0$, then
we must check that $L_{h_c,0}$ finite dimensional.
We have $L_{h_c,0}\simeq{}^{IM}L_{-c,\reg}$ by 2.6.1(b), 2.6.2(a),
and $\dim(L_{-c,\reg})<\infty$ by 2.6.2(b).
So the restriction of $L_{h_c,0}$ to $\Hb_c$ is isomorphic to
$\Lb_{c}$ by 2.7.3(b), proving claim (a).

Let us concentrate on claim (b).
If $M$ is a simple finite dimensional spherical
$\Hb_c$-module which comes from a $\Hb''_c$-module $M''$,
then $M''$ is also simple and spherical.
Thus $M''=\Lb''_c$.
Hence $M=\Lb_c$ and it is 
a simple finite dimensional quotient of $\Pb_{X,c}$
by 2.3.3(b), which belongs to $O_c(\Hb)$ by 2.1.4.
So $M$ is the restriction to $\Hb_c$ of 
a finite dimensional $Y$-spherical module 
in $\Irr(O_{h^\dag_{c,z}}(\hat\Hb))$ by 2.5.8(b), for some $z$.
Now we study the cases $k<0$ and $k>0$ separatly.

If $k<0$ then $\Lb_c$ is the restriction of $L_{h_c,0}$ to $\Hb_c$ by 2.7.3(b).
So $L_{h_c,0}$ is finite dimensional.
Hence $m\in\EN$ by 2.6.2(c).

If $k>0$ then a simple computation using \cite{GGOR, thm.~2.19} implies that
$\Pb''_c$ is simple. This was also proved before in \cite{DDO} 
by other technics.
Hence no finite dimensional
spherical $\Hb_c$-module comes from $\Hb''_c$.
Another argument is indicated in 2.8.2(g) below.
We are done.

\qed

\vskip3mm

\noindent{\bf 2.8.2. Remarks.}
(a)
Notice that 2.8.1(b) means precisely that 
the top of the $\Hb''_c$-module $\Pb_c''$ is finite dimensional iff
$m\in EN$, $k<0$. 

(b)
There are finite dimensional $\Hb''_c$-modules with $m\notin EN$.
They are not spherical.
For instance, if $(G_\CC,c)=(B_n,1/(2n-4))$ with $n\ge 3$
there is an unique one-dimensional $\Hb''_c$-module such that
$$s_1,s_2,...s_{n-1}\mapsto -1,\quad s_n\mapsto 1,\quad \l,\lambdav\mapsto 0.$$
Further $2n-4\notin\EN$ if $n\ge 5$.
There is also a non-spherical finite dimensional module
if $(G_\CC,c)=(D_4,1/2)$ by \cite{BEG2, sec.~ 6.6}, and another one if 
$(G_\CC,c)=(E_8,1/15)$ by \cite{BE}.
We'll come back to this elsewhere.

(c)
The theorem implies that there is a $Y$-spherical module in
$\Irr(O(\Hb_c))$ 
iff 
$m\in\EN$, $k<0$, and, then, this module is the restriction of $L_{h_c,0}$
to $\Hb$.
Indeed, observe first that a $Y$-spherical module in $\Irr(O(\Hb_c))$ 
is both locally finite and finitely generated over $\CC X_0$.
Hence it is necessarily finite dimensional.
Thus it is enough to prove that if
$M\in\Irr(O_{h_{c,z}}(\hat\Hb))$ 
is $Y$-spherical finite dimensional then 
$m\in\EN$, $k<0$, and $M\simeq L_{h_c,0}$.
The module
$M$ comes from a $\Hb''_c$-module $M''$
by 2.7.3(a).
Thus we have $M''=\Lb''_c$ and the restriction of $M$
to $\Hb_c$ is isomorphic to $\Lb_c$.
The rest of the proof is as above.

(d)
A $Y$-spherical finite dimensional simple $\Hb_c$-module 
may not belong to $O_c(\Hb)$.
Counter examples can be constructed by 
twisting a spherical finite dimensional simple $\Hb_c$-module
which does not come from $\Hb''_c$ by the Cherednik-Fourier transform.

(e)
By part (c) any finite dimensional module in $\Irr(O(\Hb_c))$ which
is $Y$-spherical is also spherical. The reverse is false by 2.2.5(b).

(f)
If $m\in\EN$, $k>0$
then $\Lb_{-c}$ is isomorphic to the restriction of
$L_{h_{-c},0}$ to $\Hb_{-c}$ by 2.8.1.
Thus $\Lb_{c}^\dag$ is isomorphic to the restriction of
$L_{c,\reg}$ to $\Hb_c$ by 2.6.1, 2.6.2.
Thus $\Lb_{c}^\dag$ is a Jordan-H\"older composition factor
of the $\Hb_{c}$-module $H_*(\Bc^c_x,\CC)^{A(c,x)}$.
Computations in low rank suggest that, as $\Hb_{c}$-modules,
$$m\in\EN, k>0, x\in\Nc^c_\ERS\Rightarrow
H_*(\Bc^c_x,\CC)^{A(c,x)}\simeq\Lb_{c}^\dag.$$ 

(g)
If $k>0$ then any quotient of
the $\Hb_c$-module $\Pb_{c}$ is infinite dimensional.
Equivalently, 
by 2.5.8(b), 
so is any simple quotient of the $\hat\Hb_{-c}$-module 
$\Pb_{-c}^\dag$.
Indeed, since $e_R\in\gen_0\cap\Nc^{-c}$
we can equip the space $H_*(\Bc_{0,e_R}^{-c},\CC)$ 
with the representation of $\Hb_{X,-c}$ in 
\cite{KL1}, \cite{L2}.
It is isomorphic to the Steinberg representation $\sgnb$.
The intersection of the Kostant slice $\sen$ with
the subset $\Nc^{-c}$ is reduced to $\{e_R\}$,
because $\Nc^{-c}$ consists of nilpotent elements by 2.4.2(c).
Hence the induction theorem
\cite{V, sec.~6.6} (see also A.2.4)
yields an isomorphism of $\hat\Hb_{-c}$-modules
$$\Pb_{-c}^\dag=
\hat\Hb_{-c}\otimes_{\Hb_{X,-c}}\sgnb\to H_*(\Bc^{-c}_{e_R},\CC).$$
According to 2.5.1(c) the simple quotients of 
$H_*(\Bc^{-c}_{e_R},\CC)$ 
are labelled by a subset of $\Xc_{-c}$ consisting of
irreducible representations of the group $A(h_{-c}^\dag,e_R)$.
Fix a simple quotient $L_{-c,\chi}$
with $\chi\in\Xc_{-c}$.
The $\hat\Hb_{-c}$-module $H_*(\Bc_{e_R}^{-c},\CC)$
is generated by the weight subspace $H_*(\Bc_{e_R,1}^{-c},\CC)$,
because it is equal to $\Pb_{-c}^\dag$
and the latter is generated by its
$h_{-c}^\dag$-weight subspace.
So $L_{-c,1,\chi}\neq\{0\},$ 
because $L_{-c,\chi}\neq\{0\}.$ 
Thus $\chi\in\Xc_{-c,1}$.
Now, we have
$\Xc_{-c,1}=\Xc_{-c,w}$
for all $w$ such that the alcove $\Ac_w$ is in the dominant Weyl chamber
by 3.4.1(a),(c) below.
Hence $\chi\in\Xc_{-c,w}$ for all $w$ as above.
Since there is an infinite number of such $w$'s, we have
$\dim(L_{-c,\chi})=\infty$.

\qed

\head 3. Finite dimensional representations of DAHA's and affine Springer
fibers\endhead

After our paper was writen we had some discussion with P. Etingof
who found how to classify
all finite dimensional irreducible $\Hb'$-modules
form our theorem 2.8.1.
More precisely, the set of all finite dimensional irreducible $\Hb'$-modules
is the union of the sets of all finite dimensional irreducible
modules of rational DAHA's associated to the maximal subgroups of $W_0$.
Further this decomposition is compatible with the 
subsets consisting of the modules which are spherical. 
See \cite{E} for more details.

In the rest of the paper we concentrate on a different problem :
classify all (spherical) Jordan-H\"older factors of the homology
of the elliptic homogeneous affine Springer fibers.
As mentioned in the introduction, this yields interesting combinatorics
which appear already in the representation theory of $p$-adic groups.

Sections 3.1, 3.2 contain technical facts.
The results are in section 3.3, while 3.4, 3.5 contain 
complementary facts on the dimension of finite dimensional modules.

\subhead 3.1. The evaluation map\endsubhead

Equip $\gen_0$ with the automorphism $\ad\rhov_m$.
In this section we consider the graded Lie algebra $(\gen_0,\ad\rhov_m)$.
We'll assume that $m>1$ and $k$ is arbitrary.
The root system of $G^{c,\circ}$ is the subset $\Delta_c\subset\Delta_\re.$
The set of affine roots of $\gen^c$ is the subset $\Den_c\subset\Delta_\re.$
See section 2.4.

\proclaim{3.1.1. Lemma}
We have
$\Delta_c^+=\{\a\in\Delta^+_\re;\ell=-c\rhov\cdot a, \rhov\cdot a\in m\ZZ\}$
and
$\Den_c=\{\a\in\Delta_\re;\ell=c(1-\rhov\cdot a),\rhov\cdot a\in1+m\ZZ\}.$
%and $\Den_{c,w}=\{\a\in\Den_c;A_w\subset V_\a^-\}=\Den_c\cap w(\Delta^-)$.
\endproclaim

Since both sets are finite, there are inclusions 
$G^{c,\circ}\subset G_\CC(\CC[\eps,\eps^{-1}])$
and
$\gen^c\subset \gen_0\otimes\CC[\eps,\eps^{-1}].$
Set 
$$G_0^m=Z_{G_0}(\rhov_m), 
\quad
\gen_0^m=\{x\in\gen_0;(\ad\rhov_m)(x)=\varpi_mx\}.$$
The assignement $\eps\mapsto 1$ yields a group isomorphism
$ev:G^{c,\circ}\to G_0^m$ and a linear isomorphism $ev:\gen^c\to\gen_0^m$.
For each subset $\Sc$ of $\gen_0$, $G_0$ we abbreviate 
$\Sc^m=\gen_0^m\cap \Sc$, $G_0^m\cap \Sc$ respectively.
We have $ev(\sen^c)=\sen_0^m$.

Now, fix $m\in\RN$ and $x\in\sen^c_\reg$.
We set $x_0=ev(x)$
and $A_0=Z_{G_0}(x_0)$.
Notice that $x_0\in\sen^m_{0,\reg}$.
Since $G_0$ is simply connected $A_0$
is a maximal torus by \cite{S6, thm.~8.1}.
Let $\aen_0$ be its Lie algebra.
It is the set of $\CC$-points of a Lie algebra $\CC$-scheme $\aen_\CC$.

The element $\rhov_m$ belongs to $N_{G_0}(A_0)$.
We'll identify 
$N_{G_0}(A_0)/A_0$ with $W_0$,
hopping it will not create any confusion.
Let $w_m$ be the image of $\rhov_m$ in $N_{G_0}(A_0)/A_0$.
The corresponding element of $W_0$ is uniquely determined up to conjugacy.
We'll still write $w_m$ for it.
We have ${}^{w_m}x_0=\varpi_m x_0$ and $x_0\in\aen_{0,\reg}$.
Thus $w_m$ belongs to $W_0[m]$.
Further, the conjugacy class $w_{m,*}$ 
is the type of $x$ by 1.3.2.

Let $W_m$ be the centralizer of $w_m$ in $W_0$.
The isomorphism class of the group $W_m$ depends only on the integer $m$.
It is a complex reflection group in the space $\aen_0^m$ by 1.1.1. 

Set $A(m,x_0)=Z_{G_0^m}(x_0)$ and $A(m,x_0)^\circ=A(m,x_0)/Z(G_0)$.
If there is no confusion we'll writte
$A_m$, $A_m^\circ$ for $A(m,x_0),$ $A(m,x_0)^\circ$.
Notice that $A_m$ is equal to the fixed points subgroup $A_0^{w_m}$.
In particular, its isomorphism class does not depend on the choice of $x$.
Thus $W_m$ acts on $A_m$, $A_m^\circ$ in the obvious way.
Write $N_m$ for $N_{G_0^m}(A_0)$.

\proclaim{3.1.2. Proposition} 
(a) 
If $m\in\RN$
then $\aen_0^m$ is a Cartan subspace of $(\gen_0,\ad\rhov_m)$.

(b) 
If $m\in\EN$ the Weyl group of $(\gen_0,\ad\rhov_m)$
is isomorphic to $W_m$ and $N_m/A_m$.
\endproclaim

\noindent{\sl Proof :} According to \cite{V2, sec.~3.1} a Cartan subspace
of $(\gen_0,\ad\rhov_m)$ is a maximal Abelian subspace of $\gen_0^m$ which
consists of semi-simple elements of $\gen_0$. 
So we must prove that $\aen_0^m$ is a maximal Abelian subspace of $\gen_0^m$.
Let $y_0\in\zen_{\gen_0^m}(\aen_0^m)$. 
Since $x_0\in\aen_0^m$ we have $y_0\in\zen_{\gen_0^m}(x_0)$.
Hence we have $y_0\in\aen_0^m$, because $\zen_{\gen_0}(x_0)=\aen_0$.
We have proved claim (a).

According to \cite{V2, sec.~3.4} 
the Weyl group of $(\gen_0,\ad\rhov_m)$
is $$W'_m=N_{G_0^m}(\aen_0^m)/Z_{G_0^m}(\aen_0^m).$$
Identify the group $W_0$ with $N_{G_0}(\aen_0)/A_0$.
Since the group
$Z_{W_0}(\aen_0^m)$ is generated by complex reflections by \cite{S6, 1.20}
and since $\aen_{0,\reg}^m\neq\emptyset$, we have
$Z_{W_0}(\aen_0^m)=\{1\}.$ 
For each element $w\in N_{W_0}(\aen_0^m)$ we have 
$ww_mw^{-1}w_m^{-1}\in Z_{W_0}(\aen_0^m)$.
Hence $ww_m=w_mw$ in $W_0$.
Therefore we have
$N_{W_0}(\aen_0^m)\subset Z_{W_0}(w_m)=W_m.$ 
The reverse inclusion is obvious.
%See \cite{M1, rem.~7.7} for more details.
So, to prove (b) we must check that
$$W'_m\simeq N_{W_0}(\aen_0^m).$$
Since $\aen_{0,\reg}^m\neq\emptyset$ we have 
$N_{G_0^m}(\aen_0^m)=N_{G_0^m}(\aen_0)$.
We have also $Z_{G_0}(\aen_0^m)=A_0$. 
Thus 
$$W'_m=(G_0^m\cap N_{G_0}(\aen_0))/(G_0^m\cap A_0)\subset W_m.$$
Notice that the equality above implies that
$W'_m=N_m/A_m.$
The group-scheme homomorphism 
$1-w_m:A_0\to A_0$ yields an exact sequence
$$1\to A_m\to A_0\to A_0.$$ 
Since $m\in\EN$ the element $w_m$ is elliptic,
hence the group $A_m$ is finite.
Thus
$1-w_m$ is a proper map $A_0\to A_0$.
Since $A_0$ is an irreducible variety, 
we have $A_0=(1-w_m)(A_0)$. 
Any element in $W_m$ is the $A_0$-coset
of an element $g\in N_{G_0}(\aen_0)$ such that
$(\ad\rhov_m)(g)\in gA_0$.
For each such $g$ there is an element $h\in A_0$ such that
$$h^{-1}(\ad\rhov_m)(h)=g^{-1}(\ad\rhov_m)(g).$$
Therefore $gh^{-1}\in (G_0^m\cap N_{G_0}(\aen_0))$.
Hence $W'_m=W_m$.

\qed

\vskip3mm

The group $W_m$ acts on the set $G_0^m/A_m$ by right multiplication
by 3.1.2(b).
Let $$(G_0^m/A_m)\times_{W_m}\aen_{0}^m$$ be the quotient 
relative to the $W_m$-action such that
$w(g,y_0)=(gw^{-1},{}^wy_0)$.

\proclaim{3.1.3. Corollary} 
Assume that $m\in\EN$.

(a)
The restriction
$\CC[\gen_\CC^m]^{G_0^m}\to\CC[\aen_\CC^m]^{W_m}$ is an algebra isomorphism.

(b) 
The map $q:(G_0^m/A_m)\times\aen_{0,\reg}^m\to\gen_{0,\reg}^m$,
$(gA_m/A_m,x)\mapsto(\ad g)(x)$ factors to an isomorphism of varieties
$(G_0^m/A_m)\times_{W_m}\aen_{0,\reg}^m\to\gen_{0,\reg}^m$.

(c)
The primitive $W_m$-invariants yield an isomorphism $\aen_0^m/W_m\to\sen_0^m$.

(d)
Any $(\ad G_0^m)$-orbit in $\gen_{0,\reg}^m$ meets
$\sen_{0,\reg}^m$ exactly once.
\endproclaim

\noindent{\sl Proof :} 
Claim (a) follows from \cite{V2, thm.~7}.
The proof of 3.1.2(a) 
implies that $\zen_{\gen_0^m}(y_0)$ is a Cartan subspace of
$(\gen_0,\ad\rhov_m)$ for each $y_0\in\gen_{0,\reg}^m$.
In particular $y_0$ is $G_0^m$-conjugate into $\aen_{0,\reg}^m$ by
\cite{V2, thm.~1}.
Thus the map $q$ is surjective.
On the other hand, given $g,g'\in G_0^m$ such that
$q(gA_m,y_0)=q(g'A_m,y_0')$, we get 
$y_0=(\ad g^{-1}g')(y_0')$.
Hence there is an element $g''\in N_m$ such that
$y_0'=(\ad g'')(y_0)$ and
$g^{-1}g'g''\in Z_{G_0^m}(y_0)$, by 3.1.2(b) and \cite{V2, thm.~2}.
Recall that $A_m=Z_{G_0^m}(y_0)$, because $y_0\in\aen^m_{0,\reg}$.
Thus 
$g'A_mg''=g'g''A_m=gA_m$.
Let $w$ be the image of $g''$ in $W_m$.
We have $y_0'={}^wy_0.$ 
Thus $q$ factors to an injective map
$$G_0^m/A_m\times_{W_m}\aen_{0,\reg}^m\to\gen_{0,\reg}^m.$$
Claim (b) is proved.

Now, let us concentrate on (c). 
We abbreviate $\varphi_m=(\varphi_i)_{i\in I_m}$.
Observe that $\varphi_i\eta(\gen_0^m)=\{0\}$ for all $i\notin I_m$,
because $\varphi_i\eta$ is a homogeneous invariant
polynomial of degree $d_i$.
Therefore, by 1.1.1 
the obvious inclusions $\aen_0^m\subset\aen_0$, $\CC^{I_m}\subset\CC^I$
yield a commutative square
$$\matrix
\aen_0^m/W_m&{\buildrel\varphi_m\eta\over\lra}&\CC^{I_m}\cr
\downarrow&&\downarrow\cr
\aen_0/W_0&{\buildrel\varphi\eta\over\lra}&\CC^{I_0}\cr
\endmatrix$$
with invertible horizontal maps.
Since $\varphi\eta$ yields also an isomorphism $\sen_0\to\CC^{I_0}$,
for each $y_0\in\sen_0$ we have
$$\aligned
y_0\in\sen^m_0
&\iff
(\ad\rhov_m)(y_0)=\varpi_my_0
\hfill\cr
&\iff
\varphi_i\eta(y_0)=\varpi^{d_i}_m\varphi_i\eta(y_0),\ \forall i
\hfill\cr
&\iff
\varphi_i\eta(y_0)=0,\ \forall i\notin I_m.
\hfill
\endaligned$$
The claim follows.

Finally we prove (d).
By part (b) the intersection of any $(\ad G_0^m)$-orbit in $\gen_{0,\reg}^m$
with $\aen_{0,\reg}^m$ is a $W_m$-orbit.
So the claim follows from part (c).

\qed

\vskip3mm

For a future use, let $\nu:\gen_{0,\reg}^m\to\sen_{0,\reg}^m$
be the unique $(\ad G_0^m)$-invariant map which restricts to the identity 
of $\sen_{0,\reg}^m$.
Via the evaluation map $ev$ it yields a morphism
$\nu:\gen_{\reg}^c\to\sen_{\reg}^c$.

\subhead 3.2. Homology of affine Springer fibers
\endsubhead

The purpose of this section is to gather a few basic facts on
the homology of affine Springer fibers.
For each $h=(s,\tau,\tau^{c/2})$ and $x\in\Nc^h_\ERS$ 
the $\hat\Hb_c$-module $H_*(\Bc^h_x,\CC)$ is finite dimensional.
We'll prove in particular that it is
independent of the choice of $x$.
See 2.4.2(a) and 3.2.1(a) below.
From now on we'll assume that $h=h_c^\dag$.

For each $w\in W_0$ set $X_w=X_0/(1-w)X_0$ and $A_w=T_0^w$.
Recall that $A_w=\Spec(\CC X_w)$.
See \cite{SS, sec.~II.1.7} for instance.
If $m\in\EN$, $w\in W_0[m]$
the Abelian groups $A_w$, $X_w$, $A_m$ are finite and isomorphic.
See section 3.1.
Further $X_w=\prod_j(\ZZ/e_j)$, where the $e_j$'s are the elementary
divisors of $1-w$ acting on $X_0$.

As in the previous section we'll assume that $m>1$.

\proclaim{3.2.1. Proposition}
Let $m\in\RN$, $k>0$.

(a)
The set $\Nc^c_\reg=\gen^c_\reg$ 
is affine, connected, non empty,
and $\pi^c_w$ is a locally trivial fibration over it.
The representation of $\hat\Hb_c$ on $H_*(\Bc^c_x,\CC)$
is independent of the choice of $x\in\Nc^c_\reg$.

(b)
If $x\in\Nc^c_\reg$ then $\Bc^c_x$ is a smooth projective
locally finite scheme.

(c)
If $x\in\Nc^c_\ERS$ 
the group of connected components of $Z_{\hat G^c}(x)$
is a normal subgroup of
$A(c,x)$ with complement $Z_{\CC^\times_\delta}(\nu(x))$,
and it is isomorphic to $A_m$. 
\endproclaim

Recall that the $\tilde G^c$-variety
$\Bc^c_{w}$ is isomorphic to
$\tilde G^c/B^c$.
The map $\pi^c_w:\dot\Nc_w^c\to\Nc^c$ is isomorphic to the map
$$\tilde G^c\times_{B^c}\uen_w^c\to\Nc^c,
\quad
(g,x)\mapsto(\ad g)(x).$$
We'll abbreviate $\dot\Nc_{w,\reg}^c=\dot\Nc_w^c\cap(\pi^c_w)^{-1}(\Nc_\reg^c)$ 
and $\pi^c_{w,\reg}=\pi^c_w|_{\dot\Nc_{w,\reg}^c}$.
Consider the vector bundle
$\Ec_{c,w}=\tilde G^c\times_{B^c}(\gen^c/\uen_w^c)$ over $\Bc^c_w$.
For any element $x\in\gen^c$ the map 
$$\tilde G^c\to\gen^c/\uen_w^c,\quad g\mapsto (\ad g)^{-1}(x)+\uen_w^c$$
factors to a global section $\een_{c,w}$ of $\Ec_{c,w}$ (depending on $x$).

\proclaim{3.2.2. Lemma}
(a)
If $x\in\Nc^c_\reg$ then $\een_{c,w}$ is transverse to the zero section
of $\Ec_{c,w}$.

(b)
The map $\pi^c_{w,\reg}$ is a proper submersion.

(c)
We have $\Bc^c_{x,w}=\een_{c,w}^{-1}(0)$.
\endproclaim

\noindent{\sl Proof :}
Part (c) is obvious.
The proof of (a), (b) is modelled after \cite{GKM1, sec.~2.5}.
Fix $x\in\Nc^c_\reg$.
It is enough to prove that the map
$\pi^c_w$ is a submersion.
We are reduced to check the following :
for each $g\in \tilde G^c$ such that $(\ad g)^{-1}(x)\in\uen_w^c$ we have
$$[\gen,x]^c+(\ad g)(\uen^c_w)=\gen^c.
\leqno(3.2.3)$$
To prove 3.2.3 we must check that any linear form
$\varphi$ on $\gen^c$ which annihilates both
$[\gen,x]^c$ and $(\ad g)(\uen^c_w)$ is zero.
The group $\tilde G^c$ acts on the dual space $(\gen^c)^*$ by the coadjoint
representation.
Since $\gen^c$ is finite-dimensional its dual space is the direct sum
of the weight subspaces with respect to the torus $T_0\times\CC_\delta^\times$.

We use the terminology of section 3, to which we refer for more details.
Fix $\check\mu\in\Ac_w$ (an alcove in $\check V_{0,\RR}$).
Set $\lambdav=(\muv,1)$.
The weights of
$\gen^c$, $(\gen^c)^*$ are real affine roots, because $m>1$. 
Since the linear form $\varphi$ annihilates $(\ad g)(\uen^c_w)$ and 
$\uen^c_w$ is the sum of the weight subspaces of $\gen^c$ associated 
to the affine roots $\alpha$ such that $\alpha\cdot\lambdav>0$,
the element $(\ad g)^{-1}(\varphi)$ belongs to the sum
of weight subspaces of $(\gen^c)^*$ associated 
the affine roots $\alpha$ such that $\alpha\cdot\lambdav\ge 0$. 
We have $\alpha\cdot\lambdav\neq 0$ if $\a$ is real. 
Thus $(\ad g)^{-1}(\varphi)$ is a sum of linear forms
$\varphi_i$ of weight $\beta_i$ with 
$\beta_i\cdot\lambdav=i>0$.
Hence $\varphi$ is an unstable vector
in the $\tilde G^c$-module $(\gen^c)^*$ by the Hilbert-Mumford criterium, i.e.,
we have $0\in\overline{(\ad \tilde G^c)(\varphi)}$.

Now, view $\varphi$ as an element $v\in\gen$ via the Killing form.
Since $\varphi$ vanishes on $[\gen,x]^c$ we have $v\in\zen_\gen(x)$.
Hence $v$ is semisimple, because $x\in\gen_\reg$. 
It is unstable by the discussion above.
We are done since Chevalley's restriction theorem implies that there are no
semisimple unstable elements.

\qed

\vskip3mm

\noindent{\sl Proof of 3.2.1 :}
By 3.2.2(b) the map $\pi^c_{w,\reg}$
is a locally trivial fibration. This proves claim (a).
The smoothness in (b) follows from the smoothness of $\dot\Nc^c$ and 3.2.2(b).
Now we concentrate on part (c).
There are group isomorphisms
$$A(c,x)\simeq A(c,\nu(x)),\quad
Z_{\hat G^c}(x)\simeq Z_{\hat G^c}(\nu(x)),$$
because $x$, $\nu(x)$ are 
$\hat G^c$-conjugate by definition of the map $\nu$.
Hence we may assume that $x\in\sen^c_\reg$.

For each element $hg\in Z_{\tilde G^c}(x)$, 
with $g\in\hat G^c$ and $h\in\CC^\times_\delta$,
we have $g,h\in Z_{\tilde G^c}(x)$ because the action of $h$ 
preserves the subset $\sen_\reg^c$
and two distinct elements of $\sen^c$ are not $\hat G^c$-conjugate.
Thus $Z_{\tilde G^c}(x)$ is the semi-direct product
$Z_{\hat G^c}(x)\rtimes Z_{\CC^\times_\delta}(x)$.
The group $Z_{\CC^\times_\delta}(x)$ is finite.
We must prove that $A_m$ is isomorphic to the
group of connected components of $Z_{\hat G^c}(x)$.
The obvious projection $\hat G^c\to G^{c,\circ}$ is a 
$\CC^\times_{\omega_0}$-torsor.
So is also the induced map
$Z_{\hat G^c}(x)\to Z_{G^{c,\circ}}(x)$.
The evaluation map yields a group isomorphism
$Z_{G^{c,\circ}}(x)\to Z_{G_0^m}(x_0).$
Further $Z_{G_0^m}(x_0)=A_m.$

\qed

\vskip3mm

\noindent{\bf 3.2.4. Remarks.}
(a)
The finite group $A(c,x)$ depends on the choice of the 
$\tilde G^c$-conjugacy class of the element $x\in\Nc^c_\ERS$.
For instance,
if $(G_\CC,c)=(C_2,1/2)$ then
$$\sen^c_\reg=\{e_R+af_R\otimes\eps^k+bf_\theta\otimes\eps^{2k};
16a^2+b\neq 0\}.$$
Further, if $x_i=e_R+f_R\otimes\eps^k$ then 
$A(c,x_i)=(\ZZ/2)^2\rtimes(\ZZ/ik)$ for $i=1,2$.  

(b)
The smoothness statement in 3.2.1(b) is mentioned
just for the sake of completness. We'll not use it.

\subhead 3.3. Spherical simple finite dimensional modules\endsubhead

Using sections 3.1, 3.2 we now describe the spherical 
Jordan-H\"older composition factors of the $\Hb$-module 
${}^{IM}H_*(\Bc^c_{x},\CC)$, with $x\in\Nc^h_\ERS$, $m\in\EN$, and $k>0$.
It involves some nice combinatoric which appears also in the 
tamely ramified Langlands correspondence (see the introduction).

Fix the following notation.
Let $\Xc_{c,\reg}\subset\Xc_c$ be the subset
of all local systems $\chi$
such that $\chi|_{\Nc^c_\ERS}\neq\{0\}$.
Set
$$\Lc_{c,1,\reg}=\Lc_{c,1}|_{\Nc^c_\ERS},
\quad
\Xc_{c,1,\reg}=\Xc_{c,1}\cap\Xc_{c,\reg}.$$
Thus $\Xc_{c,1,\reg}$
is the set of non-isomorphic simple direct summands of the 
semi-simple local system $\Lc_{c,1,\reg}$.

\proclaim{3.3.1. Theorem}
Let $m\in\EN$, $k>0$, and $x\in\Nc^c_\ERS$. 
%
%(a)
The non-isomorphic 
spherical Jordan-H\"older composition factors of
${}^{IM}H_*(\Bc^c_{x},\CC)$
are in one-to-one correspondence with
elements of the set $\Xc_{c,1,\reg}$.
The multiplicity of
${}^{IM}L_{c,\chi}$ in
${}^{IM}H_*(\Bc^c_{x},\CC)$
is equal to the rank of the local system $\chi$.
%
%(b)
%Any spherical finite dimensional simple $\Hb_{-c}$-module is
%a Jordan-H\"older composition factors of
%${}^{IM}H_*(\Bc^c_{x},\CC)$. 
\endproclaim

\noindent{\sl Proof :}
Let $m\in\EN$, $k>0$, and $x\in\Nc^c_\ERS$. 
By 2.5.1(c) we have
$$H_*(\Bc^c_{x},\CC)\dot=
\bigoplus_{\chi\in\Xc_{c,\reg}} L_{c,\chi}\otimes\chi_x
\leqno(3.3.2)$$ 
in the Grothendieck group of $\hat\Hb_c$.
Here, the symbol $\dot=$ denotes an equality of graded objects 
which do not preserve the grading.
Thus non-isomorphic Jordan-H\"older composition factors of
the $\Hb_{-c}$-module
${}^{IM}H_*(\Bc^c_{x},\CC)$ are in one-to-one correspondence with
$\Xc_{c,\reg}$ by 2.5.8(b).

This correspondence restricts to a bijection
between the set of non-isomorphic Jordan-H\"older composition factors of
the $\Hb_{-c}$-module
${}^{IM}H_*(\Bc^c_{x},\CC)$ whose
$h_{-c}$-weight subspace is non zero and
$\Xc_{c,1,\reg}$.

A spherical Jordan-H\"older composition factor of
${}^{IM}H_*(\Bc^c_{x},\CC)$ has a non zero $h_{-c}$-weight subspace,
because it is a quotient of the polynomial representation.
Conversely, 
let $M$ be a Jordan-H\"older composition factor of
the $\Hb_{-c}$-module
${}^{IM}H_*(\Bc^c_{x},\CC)$ whose
$h_{-c}$-weight subspace is non zero.
We have $H_*(\Bc^c_{x,w},\CC)=\{0\}$
for each $1\neq w\in W_0$ by 3.3.5(a) and 3.2.2(a),(c).
Therefore, for each element $v\in H_*(\Bc^c_{x,1},\CC)$ we have
$\psi_{s_i}(v)=0$ for all $i\in I_0$.
Hence, the formula for $\psi_{s_i}$ in section 2.1 yields
$$t_{s_i}v=-\tau^{-c/2}v,\quad\forall i\in I_0.$$
Thus $t_{s_i}$ acts by scalar multiplication by $\tau^{c/2}$ on 
${}^{IM}H_*(\Bc^c_{x,1},\CC)$, and a fortiori on the weight subspace 
$M_{h_{-c}}$.
So $M$ is a spherical module, because it is simple.

Finally, 
the multiplicity of
the Jordan-H\"older composition factor ${}^{IM}L_{c,\chi}$ in 
the $\Hb_{-c}$-module ${}^{IM}H_*(\Bc^c_{x},\CC)$ is
the rank of $\chi$ by 3.3.2.
We are done.

\qed

\vskip3mm

By 3.3.1, to compute the spherical composition factors of
$H_*(\Bc_x^c,\CC)$ with $x\in\Nc^c_\ERS$
we must first compute 
the set $\Xc_{c,1,\reg}$.
This can be done explicitely.
First, let us introduce more notations.
Since $1\notin\EN$, 
in the rest of the section we'll assume that $m\neq 1$
to simplify.
Thus, for each $h=(s,\tau,\tau^{c/2})$
the set $\Den_h$ consists of real affine roots. 

Given $x$ as above, we set $x_0=ev(x)$.
To $x_0$ we associate the groups
$A_0$, $A_m$, $A_m^\circ$, $N_m$, and $W_m$ as in section 3.1.
The following technical conjecture is important for the rest of the paper.

\proclaim{3.3.3. Conjecture}
If $m\in\EN$, $k>0$, and $x\in\Nc_\ERS^c$ then
$|\Bc^c_{x,1}|=|A_m^\circ|.$
\endproclaim

For each $w$ 
let $\Den_{h,w}\subset\Den_h$ be such that
$$\uen_w^h=\bigoplus_{\a\in\Den_h\setminus\Den_{h,w}}\CC e_\a.$$
It is proved in 3.3.5(b) below that the set 
$\Bc^c_{x,1}$ is finite of cardinal $\ge|A^\circ_m|$. 
Set 
$$p=\prod_{\a\in\Den_{c,1}}a,
\quad
eu=\prod_{\a\in \Delta_c^+}a.$$
To prove the equality, by 3.2.2(a),(c) it is enough to check that the image of
the polynomial $p$ by the Borel map 
$\CC[\check V_0]\to H_*(\Bc_1^c,\CC)$ is
$|A_m^\circ|$ times the fundamental class.
Using the localization theorem in equivariant cohomology
we are reduced to prove the following formula 
$$(-1)^{n_c}|A_m^\circ|\, eu
=\sum_{w\in W_c}(-1)^ {\ell(w)}\,{}^w\!p.\leqno(3.3.4)$$ 
See \cite{BGV, thm.~7.13} for instance.
We have checked it by direct computations in a lot of cases.
See section 4 below for more details.

\proclaim{3.3.5. Lemma}
Let $k>0$.

(a)
We have $|\Den_{c,w_0}|\ge n_{c}$.
If $m\in\EN$ then
$|\Den_{c,w}|\ge n_{c}$ for each $w\in W_0$,
with an equality iff $w=1$.

(b)
If $m\in\EN$, $x\in\Nc_\ERS^c$, and 3.3.3 holds then
$A(c,x)$ acts transitively on $\Bc^c_{x,1}$.
In particular the map $g\mapsto(\ad g)(\ben_0^-)$ yields a bijection 
$A_m^\circ\simeq A(m,x_0)^\circ\simeq\Bc^c_{x,1}.$
\endproclaim

\noindent{\sl Proof :}
First, we prove part (a).
Notice that $\uen_w^c=\ben_w^c$, because $\zeta\neq 1$.
Thus, 
$$\Den_{c,w}=\Den_c\cap w(\Delta^-_\re)=
\{\a\in\Den_c;\ell<0\}\cup(\Pi_0\cap w(\Delta_0^-)),\quad
\forall w\in W_0.$$
The $\slen_2(\CC)$-representation theory, applied to a triple
containing $e_R$, implies that
$[\uen_0,e_R]$ is a complement of $\bigoplus_i\CC e_{a_i}$ in $\uen_0.$
Thus we have 
$$n_c-\dim\,Z_{U_0^m}(e_R)=\dim(\uen_0^m)-n.$$
Using 3.1.1, it yields
$$\aligned
|\Den_{c,w}|&=\dim(\uen_0^m)-n+|\Pi_0\cap w(\Delta_0^-)|
\\
&=n_c-\dim\,Z_{U_0^m}(e_R)+|\Pi_0\cap w(\Delta_0^-)|.
\endaligned$$
Further $n\ge\dim\,Z_{U_0^m}(e_R)$, because $e_R$ is a regular element.
Thus $|\Den_{c,w_0}|\ge n_{c}$.

Now, assume that $m\in\EN$.
Recall that
$\zen_{\uen_0}(e_R)$
is spanned by vectors $e_i$, $i\in I_0$, such that
$[\rhov,e_i]=(d_i-1)e_i$, see 1.3.1(c).
Thus the weights of $(\ad\rhov_m)$ in
$\zen_{\uen_0}(e_R)$
belong to the set $\{\varpi^{(d_i-1)/m}\}$.
Now, a case by case checking shows that an exponent,
i.e., an integer of the form $d_i-1$, can not be a multiple
of an elliptic number.
Therefore $\dim\,Z_{U_0^m}(e_R)=0$.

Now, we prove part (b).
Let $m\in\EN$ and $x\in\Nc_\ERS^c$.
We may assume that $x\in\sen^c_\reg$ by section 3.2.
So, in particular, we have $\ben_1\in\Bc^c_{x,1}$.
The set $\Bc^c_{x,1}$ is finite by 3.2.2(a),(c), because $|\Den_{c,1}|=n_c$
by part (a).
Identify $A_m$ with the group of connected components
of $Z_{\hat G^ c}(x)$ as in 3.2.1(c).
It acts on $\Bc^c_{x,1}$ by left multiplication,
because $\Bc^c_{x,1}$ is finite.
The complement of $A_m$ in $A(c,x)$ fixes $\ben_1$.
So the $(\ad A_m)$-orbit of $\ben_1$ and the
$(\ad A(c,x))$-orbit of $\ben_1$ coincide.
%To prove that they are equal to $\Bc^c_{x,1}$
%we'll check that both sets have
%$|A_m^\circ|$ elements.

Now, let $\Bc_0$ be the flag variety of $\gen_0$
and
$\Bc^m_{0,1}\subset\Bc_0$ be the $G_0^m$-orbit of $\ben_0^-$.
The bijection 
$G_0/B_0^-\to\Bc_0,$ $gB_0^-/B_0^-\mapsto(\ad g)(\ben_0^-)$
identifies $G_0^mB_0^-/B_0^-$ with $\Bc^m_{0,1}$.
%Since $G_0^mB_0^-/B_0^-\simeq G_0^m/B_0^{-,m}$,
%the set $\Bc^m_{0,1}$ may also be identified with the flag variety of $G_0^m$.
The evaluation map $ev:G^{c,\circ}\to G_0^m$
factors to an isomorphism of varieties
$\Bc^c_1\to\Bc_{0,1}^m.$

Finally, the $(\ad U_0^-)$-orbit of $x_0$ is contained into
$e_R+\ben^ -_0$ and meets $\sen_0$ exactly once (at $x_0$) by 1.3.1(d).
So $Z_{B_0^-}(x_0)=Z_{T_0}(e_R)=Z(G_0)$.
Since $A_0$ centralizes $x_0$, this implies that $B_0^-\cap A_0=Z(G_0)$.
Thus the $(\ad A_0)$-orbit of $\ben_0^-$ in $\Bc_0$ is equal to $A_0/Z(G_0)$.
Hence the $(\ad A_m)$-orbit of $\ben_0^-$ in $\Bc_{0,1}^m$ is 
equal to $A_m^\circ$.
Therefore, applying the evaluation map, we get that the
$(\ad A_m)$-orbit of $\ben_1$ in $\Bc^c_{x,1}$ is isomorphic to $A_m^\circ$.
So 3.3.3 implies that $\Bc^c_{x,1}\simeq A_m^\circ$.

\qed

\proclaim{3.3.6. Theorem}
Let $m\in\EN$, $k>0$.
If 3.3.3 holds there is a one-to-one correspondence between
elements of the set $\Xc_{c,1,\reg}$ and 
$W_m$-orbits in $\Irr(A_m^\circ)$.
The multiplicity of
${}^{IM}L_{c,\chi}$ in
${}^{IM}H_*(\Bc^c_{x},\CC)$
is the cardinal of the corresponding orbit. 
\endproclaim

\noindent{\sl Proof :}
Set $V=ev(\ben_1^c)$.
We have
$V=\ben_0^{-,m}\oplus\bigoplus_i\CC e_{a_i}$.
Let $\dot\Nc_{0,1}^m\subset\Bc^m_{0,1}\times\gen_0^m$ be the set of couples 
$((\ad g)(\ben^{-}_0),y)$ such that $g\in G_0^m$,
$(\ad g^{-1})(y)\in V$,
and let $\pi_{0,1}^m:\dot\Nc_{0,1}^m\to\gen_0^m$ be the second projection.
The fibers of $\pi_{0,1}^m$ are generalized Hessenberg varieties, 
see \cite{DPS}, \cite{GKM1}.

To simplify we'll abbreviate
$\bar\gen=\gen_{0,\reg}^m$,
$\bar\sen=\sen_{0,\reg}^m$,
$\bar\Nc=(\pi_{0,1}^m)^{-1}(\bar\gen)$,
$\bar G=G_0^m/Z(G_0),$
and
$\bar\pi=\pi_{0,1}^m|_{\bar\Nc}:\bar\Nc\to\bar\gen.$
Since the evaluation map $ev$ identifies the morphisms
$\pi^c_1$ and $\pi_{0,1}^m$, by 3.3.1 we must classify the
non-isomorphic simple direct summands of the 
$\bar G$-equivariant local
system $\Lc=\bar\pi_*(\CC_{\bar\Nc})$
(the center $Z(G_0)$ acts trivially on $\CC_{\bar\Nc}$, hence on $\Lc$).

For any subscheme $\Xc\subset\bar\gen$ we consider the group-$\Xc$-schemes
$$\bar G_\Xc=\bar G\times\Xc,
\quad
\bar A_\Xc=\{(g,y)\in \bar G_\Xc;(\ad g)(y)=y\}.$$
By 3.1.3(d) the map
$$\beta:\bar G_{\bar\sen}\to\bar\gen,
\quad
(g,y)\mapsto(\ad g)(y)$$
factors to an isomorphism
$\bar G_{\bar\sen}/\bar A_{\bar\sen}\to\bar\gen.$
Further, the fiber of $\bar A_{\bar\sen}$ over any element $y\in\bar\sen$
is equal to $A(m,y)^\circ$.
So 3.3.5(b) yields a $\bar\sen$-scheme isomorphism
$$\bar A_{\bar\sen}\to\bar\pi^{-1}(\bar\sen),
\quad
(g,y)\mapsto((\ad g)(\ben_0^-),y).$$
Notice that
$\{\ben_0^-\}\times\bar\sen\subset\bar\Nc$,
because $\bar\sen\subset V$.
Thus we have a scheme isomorphism
$$\bar G_{\bar\sen}\to\bar\Nc,
\quad
(g,y)\mapsto((\ad g)(\ben_0^-),(\ad g)(y))$$
such that $\bar\pi$ is the quotient
$\bar G_{\bar\sen}\to\bar G_{\bar\sen}/\bar A_{\bar\sen}$ 
(it is an \'etale cover of smooth $\bar\sen$-schemes).
Now, we must compute the $\bar G_{\bar\sen}$-equivariant local
system $\Lc=\pi_*(\CC_{\bar G_{\bar\sen}})$.
%Further, the fiber of $\pi$ over the $\CC$-point $y_0\in\sen$
%is equal to $Z_{G^m_0}(y_0)/Z(G_0)$, which is a finite group
%isomorphic to $A^\circ_m$ (non canonically).

Recall that we have fixed an element $x_0\in\bar\sen$.
Set $\bar\aen=\zen_{\bar\gen}(x_0).$ 
By 3.1.3(a),(c) the primitive $W_m$-invariants yield a map
$$f:\bar\gen\to\bar\aen/W_m\simeq\bar\sen,$$
which restricts to a Galois cover 
$f|_{\bar\aen}:\bar\aen\to\bar\sen$ with group $W_m$.
The same argument as in \cite{N, prop.~3.2} yields a
group-$\bar\gen$-scheme isomorphism 
$f^{-1}(\bar A_{\bar\sen})\to\bar A_{\bar\gen}.$
It restricts to a group-$\bar\aen$-scheme isomorphism 
$$(f|_{\bar\aen})^{-1}(\bar A_{\bar\sen})\to\bar A_{\bar\aen}.$$
A little attention shows that 
$\bar A_{\bar\sen}$ is indeed the quotient of
$\bar A_{\bar\aen}=A^\circ_m\times\aen$
by the obvious diagonal $W_m$-action.
Further, under taking base change with respect to $f|_{\bar\aen}$, 
the map $\bar\pi$ is again taken to the quotient 
$\bar G_{\bar\aen}\to \bar G_{\bar\aen}/\bar A_{\bar\aen}$.
So $(f|_{\bar\aen})^{-1}(\Lc)$ 
is the constant sheaf $(\CC A^\circ_m)_{\bar G_{\bar\aen}/\bar A_{\bar\aen}}$.

Since the group $\bar G$ is connected, 
a $\bar G_{\bar\aen}$-equivariant local system
over the homogeneous space
$\bar G_{\bar\aen}/\bar A_{\bar\aen}$ is the same as a representation 
of the finite group $A^\circ_m$.
Thus, by faithfully flat descent, a
$\bar G_{\bar\sen}$-equivariant local system
over $\bar G_{\bar\sen}/\bar A_{\bar\sen}$ is the same as a 
representation of the semi-direct product 
$W_m\ltimes A_m^\circ$.
Further, $\Lc$ is identified to the 
obvious representation of $W_m\ltimes A^\circ_m$
on $\CC A^\circ_m$.
Thus the theorem is a consequence of
the following basic result in Clifford theory.

\proclaim{3.3.7. Lemma}
Let $N=W\ltimes A$, where $W,A$ are finite groups and $A$ is Abelian.

(a)
There is a one-one correspondence between irreducible representations of
$N$ and isomorphism classes of pairs $(O,\pi_O)$,
where $O\subset \Irr(A)$ is a $W$-orbit and
$\pi_O$ is an irreducible representation of the isotropy subgroup 
$W_O$ in $W$ of a character $\chi_O\in O$. 
The representation associated to the pair $(O,\pi_O)$ 
is the induced module 
$L_{(O,\pi_O)}=\Ind^N_{W_O\ltimes A}(\pi_O\otimes\chi_O)$.

(b)
There is a one-one correspondence between the non-isomorphic Jordan-H\"older 
composition factors of the $N$-module $\CC A$ and the
$W$-orbits in $O\subset \Irr(A)$.
The representation associated to $O$ 
is $L_{(O,triv)}$.
We have $\dim(L_{(O,triv)})=|O|$.

\endproclaim

\qed

\vskip3mm

\noindent{\bf 3.3.8. Example.}
Let $m\in\EN$, $k>0$, and $x\in\Nc^c_\ERS$.
By 3.3.1, 3.3.6 the non-isomorphic 
spherical Jordan-H\"older composition factors of
the $\Hb_{-c}$-module ${}^{IM}H_*(\Bc^c_{x},\CC)$
are in one-to-one correspondence with
the $W_m$-orbits in $A_m^\circ\simeq\Irr(A_m^\circ)$.
Among those $\Hb_{-c}$-modules, 
there is the unique spherical finite dimensional simple module
$\Lb_{-c}$ which comes from $\Hb''_{-c}$.
It is easy to see that $\Lb_{-c}$ corresponds to the orbit
of the unit in $A_m^\circ$.
In particular $\Lb_{-c}$ has multiplicity one in the homology.

\vskip3mm

\noindent{\bf 3.3.9. Remarks.}
(a)
We collect here two remarks concerning 3.3.4. 
First, since $A_m^\circ\simeq\check Y_0/(1-w_m)X_0$,
the integer $|A_m^\circ|$ is
the square of the defect of $w_m$ in the terminology in \cite{KP, sec.~10}.
Second,
to prove that the embedding 
$A_m^\circ\subset\Bc^c_{x,1}$
is onto we must check that each element in $\Bc_{x,1}^c$ is
$Z_{\tilde G^c}(x)$-conjugate to $\ben_1$ for each $x\in\sen^c_\ERS$.
Further, by \cite{B1, prop.~(a)} the set $\Pc_x$ of 
parahoric Lie subalgebras $\qen\subset\gen$ of type $I_0$ containing $x$ 
such that the image of $x$ in the quotient $\qen/\eps\qen$ 
is a regular nilpotent element is an orbit of the group
$Z_{\tilde G}(x)$.
Therefore, one should check
that for any $\ben\in\Bc^c_{x,1}$
the unique parahoric Lie subalgebra of type $I_0$ containing $\ben$ 
belongs to $\Pc_x$.

(b)
A list of elements $e_c\in\sen^c_\ERS$, for each
$m\in\EN$, $k>0$, is given in section 4 below.
Notice that $e_z\in\sen^c_\reg$ iff $e_R+\sum_iz_if_i\in\sen_{0,\reg}$,
because $e_z=\eps^{c/2}\bullet(e_R+\sum_iz_if_i)$,
and that $e_R+f_n\in\sen_{0,\reg}$,
because $f_n=f_\theta$ up to a scalar.
In the computations we'll order the simple roots as in 
\cite{B2, planches I-IX}.

(c)
If $2\in\EN$ then $w_2=-1$.
Hence $W_2=W_0$,
$A_2\simeq X_0/2X_0$, and
$A_2^\circ\simeq\check Y_0/2X_0$.
The group $W_2$ acts on $A_2^\circ$ in the obvious way.

(d)
Since the $h_{-c}$-weight subspace of a simple spherical $\Hb_{-c}$-module
is one, the formula 
$\Lc_{c,1,\reg}=
\bigoplus_{\chi\in\Xc_{c,1,\reg}}({}^{IM}L_{c,\chi})_{h_{-c}}
\otimes\chi$
implies that the multiplicity in 
$\Lc_{c,1,\reg}$ of each irreducible local system is also one.

%Further $\dim\,A_0B_0/B_0=\dim\,X_{B_0,\Hen,x_0}$.
%Thus $\overline{A_0B_0/B_0}$ is a connected component of $X_{B_0,\Hen,x_0}.$
%Any connected component of $X_{B_0,\Hen,x_0}$ contains an element of
%$(\Bc_0)^{A_0}$.
%This fixed points set is isomorphic to $W_0$.
%Further, it is known that
%$(\Bc_0)^{A_0}\subset\overline{A_0B_0/B_0}$.
%See \cite{DPS, sec.~IV}.
%Thus $\overline{A_0B_0/B_0}=X_{B_0,\Hen,x_0}$.
%The map $\g$ factors to an isomorphism
%$\Bc[c]_{e_c,1}\to X_{B_0,\Hen,x_0}\cap\Bc_0[m]_1.$
%Since $\rhov_m\in N_{G_0}(A_0)$ it yields an element in $W_0$.
%The $A_m$-action on
%$H_*(\Bc[c]_{e_c,1},\CC)$ is induced by the
%left multiplication on $\Bc_0[m]_1$.

\subhead 3.4. Dimension\endsubhead

The dimension of the simple finite dimensional $\Hb_c$-modules and 
the dimension of
the homology of affine Springer fibers are both difficult 
to compute in general. 
In this section we collect some partial result.

For each real affine root $\a$ let
$V_\a, V^+_\a, V^-_\a$
be the set of elements $\lambdav\in\check V_{0,\RR}$
such that
$a\cdot\lambdav+\ell=0$,
$a\cdot\lambdav+\ell>0$,
and
$a\cdot\lambdav+\ell<0$ respectively.
Let $\Aen$ be the set of alcoves in $\check V_{0,\RR}$,
and $\Ac$ be the unique alcove contained in the Weyl chamber
containing 0 in its closure.
We may identify an element $w\in W$ with the alcove $\Ac_w=w(\Ac)$,
hopping it will not create any confusion.
Then, the set $W^h$ is identified with
$\Aen^h=\{\Ac;\Ac\subset\bigcap_{\a\in\Delta^+_h}V_\a^+\}$.

The open subset $V_{h,\reg}=\check V_{0,\RR}\setminus\bigcup_{\a\in\Den_h}V_\a$ 
consists
of a finite number of connected components.
We'll say that two elements $w,w'\in W^h$ are in the same
clan, or $h$-clan to avoid confusions, iff
$\Ac_w$, $\Ac_{w'}$ are contained in the same connected component.
Let $\Cen_A$ be the set of $w$'s such that
$\Ac_w$ belongs to the antidominant Weyl chamber,
and $\Cen_D$ be the set of $w$'s such that
$\Ac_w$ belongs to the dominant Weyl chamber.

%Now, let $h=h_c^\dag$.
%For any element $x\in\Nc^c$ the section $\een_{c,w}$ 
%vanishes at $(\ad g)(\ben_w)$ iff $x\in(\ad g)(\uen_w)$.
%Recall that $\Bc^h_{x,w}=\een_{h,w}^{-1}(0)$
%for each $x\in\Nc^h$.

\proclaim{3.4.1. Proposition}
(a)
If $w,w'$ are in the same clan then
$\pi^h_w$, $\pi^h_{w'}$
are isomorphic. 
Hence $\Xc_{h,w}=\Xc_{h,w'}$
and the $\tilde G^c$-varieties $\Bc^h_{x,w},\Bc^h_{x,w'}$ are equal 
for each $x\in\Nc^h$.

(b)
Let $h=(\lambdav(\tau^c),\tau,\tau^{c/2})$
and $h'=(\lambdav(\tau^c),\tau^k,\tau^{c/2})$ with $\lambdav\in Y_0$.
Then $\Nc^h\simeq\Nc^{h'}$, $\Xc_h\simeq\Xc_{h'}$, 
$\dim(L_{h,\chi})=k^n\dim(L_{h',\chi})$,
and $\dim H_*(\Bc^h_x,\CC)=k^n\dim H_*(\Bc^{h'}_x,\CC)$ 
for each $x,\chi$.

(c)
If $k<0$ 
then $\Cen_D$
is a $h_{c}^\dag$-clan and a $h_{-c}$-clan.
If $k>0$ 
then $\Cen_A$
is a $h_{c}^\dag$-clan and a $h_{-c}$-clan.
\endproclaim

\noindent{\sl Proof :}
First, we prove (a).
We have $\Ac_w\subset V_\a^+$ iff $e_\a\in\uen_w$,
%$\a\in\Den_h$ iff $e_\a\in\gen^h$,
and $\a\in\Den_h\setminus\Den_{h,w}$ iff $e_\a\in\uen^h_w$.
Therefore $\Den_{h,w}=\{\a\in\Den_h;\Ac_w\in V_\a^-\}.$ 
So $\Den_{h,w}=\Den_{h,w'}$ iff 
$\uen_w^h=\uen_{w'}^h$ iff $w,w'$ are in the same clan.
The claim follows, because 
$\dot\Nc_w^h=\tilde G^h\times_{B^h}\uen_w^h$
for each $w\in W^h$.

Now we prove (b).
The real affine root $\a=(a,\ell)$ belongs to
$\Den_h$ iff $\ell=c(1-\lambdav\cdot a)$.
Thus, since $\lambdav\cdot a\in\ZZ$ and $(k,m)=1$, we have
$\ell\in k\ZZ$ for each $\a\in\Den_h$.
So
$\Den_{h}$ consist of the real affine roots
$(a,k\ell)$ such that $(a,\ell)\in\Den_{h'}$.
Hence the assignement 
$\eps^\ell\mapsto \eps^{k\ell}$
yields isomorphisms
$\gen^{h'}\simeq\gen^{h}$,
$\Nc^{h'}\simeq\Nc^{h}$, 
and
$G^{h',\circ}\simeq G^{h,\circ}$.

Further, since scalar multiplication by $k$ identifies the
set $V_{h',\reg}$ with $V_{h,\reg}$,
it yields a bijection from the set of $h$-clans
to the set of $h'$-clans taking any finite
clan $\Cen$ to a clan $\Cen'$ such that 
$|\Cen|=k^n|\Cen'|$.
We claim that the morphisms
$\pi^h_w$, $\pi^{h'}_{w'}$
are the same whenever $\Ac_w\in k\Ac_{w'}$.
This will imply part (b).
To prove this it is enough to check that 
$\Den_{h,w}=\{(a,k\ell);\a\in\Den_{h',w'}\}$. 
This is obvious. Indeed, given an element $\lambdav\in\Ac_{w'}$
we have
$$\Den_{h',w'}=\{\a\in\Den_{h'};a\cdot\lambdav+\ell<0\}.$$
Since all the alcoves in $k\Ac_{w'}$ belong to the same $h$-clan
and the set $\Den_{h,w}$ 
is independent of the choice of $\Ac_w$ in its $h$-clan,
we have also
$$\aligned
\Den_{h,w}
&=\{\a\in\Den_{h};ka\cdot\lambdav+\ell<0\}
\hfill\cr
&=\{(a,k\ell);\a\in\Den_{h'},\,a\cdot\lambdav+\ell<0\}.
\hfill
\endaligned$$

(c) 
We'll assume that $k<0$ and $h=h_{c}^\dag$, the other cases beeing similar.
If $a\in\Delta_0^+$ and $\ell<0$ then $\ell<c(1-\rhov\cdot a)$,
hence $\a\notin\Den_c$.
The affine roots $\a$ such that $V_\a$ intersects the dominant Weyl 
chamber (i.e., the polyhedral cone generated by the fundamental coweights) 
are all of the above form.
Thus the dominant Weyl chamber is contained in a clan.
On the other hand we have $\alpha_i\in\Den_c$ for each $i\in I_0$ by 3.1.1.
Thus $\Cen_D$ is a clan.

\qed

\proclaim{3.4.2. Corollary}
If $m\in\EN$, $k<0$ then $\dim\Lb''_{k/m}=(-k)^n\dim\Lb''_{1/m}$.
\endproclaim

\vskip3mm

\noindent{\bf 3.4.3. Remarks.}
(a)
If $k>0$, $m\in\EN$, $x\in\Nc^c_{ERS}$, 
the $A(c,x)$-action on $H_*(\Bc^c_{x,1},\CC)$
is known by 3.3.5(b).
Another subspace 
$H_*(\Bc_{x,w}^c,\CC)\subset H_*(\Bc_x^c,\CC)$
is easy to compute.
There is an unique alcove $\Ac_{w_c}$ in $\Aen^c$
containing $c\rhov$ in its closure.
%Thus $Z_W(\lambdav_c)$ is a parabolic subgroup of $W$.
If $\a\in\Den_c$ then $c\rhov\in V_\a^+$,
because $\a\cdot\rhov_c=c>0$.
Hence $\Ac_{w_c}\subset\bigcap_{\a\in\Den_c}V_\a^+$, 
i.e., $\Den_{c,w_c}=\emptyset$.
Therefore we have $\Bc^c_{x,w_c}=\Bc^c_{w_c}$ for each $x\in\Nc^c$ by 3.2.2(c).

In particular the simple perverse sheaf
$\CC_{\Nc^c}[\dim\Nc^c]$ belongs to $\Xc_{c,w_c}$.
Hence the weight subspace 
$L_{h_{-c},w_c,0}$ of $L_{h_{-c},0}$
is always $\neq 0$.

(b)
An algebraic counterpart to part (a) above is the following.
For any $\lambdav\in\check X\cap\bar\Ac$
let $W_\lambdav$ be the group
generated by $\{s_\a;\a\cdot(\lambdav,1)=0\}$. 
It is a parabolic subgroup of the Coxeter group $W$.
Therefore, Kato's criterium applied to the
affine Hecke algebra $\Hb_\lambdav$ associated to 
couple $(X,W_\lambdav)$
implies that the $\lambdav$-weight space of any $\hat\Hb$-module 
%There is an obvious inclusion $\Hb'_{\lambdav}\subset\Hb'.$
%is preserved by $\Hb'_{\lambdav}.$
%Since the group $Z_W(\lambdav)$ fixes $\lambdav$, 
%applied to the $\Hb'_{\lambdav}$-module $M_{\lambdav}$ implies that
is a multiple of the regular representation as a module
over the Hecke algebra of $W_\lambdav$.
In particular, since the element $\lambdav_c={}^{w_c^{-1}}(c\rhov)$ 
belongs to $\bar\Ac$ and since $W_{\lambdav_c}=w_c^{-1}W_cw_c$,
the dimension of the $\lambdav_c$-weight subspace
is a multiple of $|W_c|$.

(c)
Let $w\in W^h$.
Then $ws_{\a_i}\in W^h$ iff ${}^w\a_i(h)\neq 1,$
and $ws_{\a_i}$, $w$ are in the same $h$-clan iff
${}^w\a_i(h)\neq\zeta^{\pm 2}.$
Compare \cite{R1, lem.~12.6}.
Thus the clans are an affine analogue of the connected components
of the weight diagrams in \cite{R1}.
Indeed, 3.4.1(a) may be viewed as an affine analogue of
\cite{loc.~cit., cor.~3.6}.

(d)
The same argument as in 3.4.1(b) proves that the Jordan-H\"older multiplicities
of the standard modules in $O_h(\Hb)$, $O_{h'}(\Hb)$ are the same,
with $h$, $h'$ as in loc.\ cit. Probably this implies that the corresponding
categories are equivalent. 
This can also probably be proved using the KZ functor
as in \cite{VV}.

\subhead 3.5. The Coxeter case and the sub-Coxeter case\endsubhead

In this section let $k>0$.
If $m=h$ then the dimension of
$H_*(\Bc^c_{e_c},\CC)$ is $k^n$
by \cite{F, prop.~1}.
In this case the $\Hb_c$-module
$\Lb_c^\dag$ has a particular nice algebraic description, 
see \cite{C4, thm.~3.7.2, thm.~3.10.6}.
This is called the Coxeter case.
In this case the equality 3.3.4 is obvious.
Indeed, it is enough to assume that $k=1$, and in this case we have
$\Delta_c=\emptyset$, 
$\Den_c=\Pi$, $\Den_{c,1}=\emptyset$, and $A^\circ_m=\{1\}$.
Since $hn$ is equal to the number of roots,
using the dimension formula of affine Springer fibers
in \cite{KL2}, \cite{B1} it is easy to see that
$\Bc_{e_c}$ is zero dimensional iff $c=1/h$.
Thus the Coxeter case is the set of pairs
$(G_\CC,m)$ such that $\Bc_{e_c}$ is finite for some $k$.

A pair $(G_\CC,m)$ is of sub-Coxeter type if it is not
of Coxeter type and
if $\Bc_{e_c}$ is a curve for some $k$.
In this case the homology of $\Bc_{e_c}$ is given in 
\cite{KL2, prop.~7.7}.
From loc.~cit.~ there is no odd homology,
and $h_c^\dag$ is a generator of
the one parameter group $z\mapsto(\rhov(z^{2k}),z^{2m}, z^k)$.
Therefore we have 
$$\dim\,H_*(\Bc^c_{e_c},\CC)=\dim\,H_*(\Bc_{e_c},\CC)$$
by \cite{B4}.
Using the dimension formula of affine Springer fibers
we get that $\Bc_{e_c}$ is a curve iff
$k=1$ and $m$ is in the following list. Therefore,
in the sub-Coxeter we have 
$$\vbox{
\offinterlineskip

\def\tablerule{\noalign{\hrule}}

\def\tableskip{&&&&&&&\cr}

\halign{
\tabskip= .3em#\vrule 
&\it#&\vrule#
&\vtop{\hsize=2pc\it #}&#\vrule\tabskip=-.1em
&\vtop{\hsize=7pc\it #}&#\vrule\tabskip=-.1em
&\vtop{\hsize=4pc\it #}&#\vrule\tabskip=-.1em\cr
\tablerule
\tableskip
&$G_\CC$&
&$m$&
&$\dim\,H_*(\Bc^c_{e_c},\CC)$&
&$A_m^\circ$&
\cr
\tablerule
\tableskip
&$C_2$&
&$2$&
&$6k^2$&
&$\ZZ/2$&
\cr
\tablerule
\tableskip
&$D_4$&
&$4$&
&$6k^4$&
&$\{1\}$&
\cr
\tablerule
\tableskip
&$E_6$&
&$9$&
&$8k^6$&
&$\{1\}$&
\cr
\tablerule
\tableskip
&$E_7$&
&$14$&
&$9k^7$&
&$\{1\}$&
\cr
\tablerule
\tableskip
&$E_8$&
&$24$&
&$10k^8$&
&$\{1\}$&
\cr
\tablerule
\tableskip
&$F_4$&
&$8$&
&$9k^4$&
&$\ZZ/2$&
\cr
\tablerule
\tableskip
&$G_2$&
&$3$&
&$8k^2$&
&$\ZZ/3$&
\cr
\tablerule
}}
$$
In this case we expect that the $\Hb_c$-module
$\Lb_c^\dag$ has a particular nice algebraic description. 
In particular the following should be true.
%In this case, the $\Hb_c$-module
%$H_*(\Bc^c_{e_c},\CC)$ is not irreducible in general.

\proclaim{3.5.1. Conjecture}
In the sub-Coxeter case we have $\dim(\Lb_c^\dag)=(n+2)k^n$.
\endproclaim

To prove 3.5.1 it is enough to assume that $k=1$ by 3.4.1.
Since 3.3.4 is true in the sub-Coxeter case (see section 4),
the $\Hb$-module ${}^{IM}H_*(\Bc^c_{e_c},\CC)$
has a unique spherical Jordan H\"older factor, except for 
$(G_\CC,m)=(C_2,2)$, $(G_2,3)$, $(F_4,8)$
by 3.3.1, 3.3.6.
Away of these cases 3.5.1 reduces to prove that
the $\Hb$-module ${}^{IM}H_*(\Bc^c_{e_c},\CC)$
is indeed irreducible.
If $(G_\CC,m)=(C_2,2)$, $(G_2,3)$
then 3.5.1 is true.
See sections 4.3, 4.9.

\head 4. Computations\endhead

In the rest of the paper we list case by case computations.
We'll write $\Pi_c$ for the basis of $\Delta^+_c$,
and $\underline{i}$ for the set $\{1,2,...i\}$.
Let $\Cen_c$ be the set of clans consisting of the $w$'s
such that $\Bc^c_{x,w}\neq\emptyset$ for some $x\in\Nc^c_\reg$
(then $\Bc^c_{x,w}\neq\emptyset$ for all $x\in\Nc^c_\reg$).
%Observe that if $m\in\RN$, $k>0$, $x\in\Nc^c_\reg$ then $\Bc^c_{x,w}$
%is either smooth equidimensional or empty by 3.2.2(b).
In the computations below we have checked equality 3.3.4 if
$(G_\CC,m)$ is either of sub-Coxeter type or equal to
$(G_2,2),$ $(F_4,6),$ 
$(E_6,6),$ $(E_8,15)$, $(E_8,20).$
In the Coxeter case it is obvious.

\subhead 4.1. Type $A_n$\endsubhead

We have
%$h=n+1$, $\{d_i\}=\{2,3,...n+1\}$, and
$\EN=\{n+1\}$, and
$A_{n+1}=Z(G_0)\simeq\ZZ/(n+1)$.
Set $e_c=e_R+f_\theta\otimes\eps^{(n+1)c}$.
Then $\dim\,H_*(\Bc^c_{e_c},\CC)=k^n$.
This is the Coxeter case.

\subhead 4.2. Type $B_n$\endsubhead

Let $n\ge 3$.
The sets EN, $I_m$, the groups $A_m$, $Z(G_0)$,
and the element $e_c$ are as in section 4.3 below.

\subhead 4.3. Type $C_n$\endsubhead

Let $n\ge 2$.
We have
%$h=2n$, $\{d_i\}=\{2,4,...2n\}$,
$\EN=\{m=2r; r|n\}$,
$I_m=\{i\in I_0;r|i\}$,
$A_m\simeq(\ZZ/2)^{n/r},$
and $Z(G_0)=A_{2n}.$
Set $e_c=e_R+f_\theta\otimes\eps^{2nc}$.

If $(n,m)=(2,2)$ then 
%\Den_c=\{(\theta,-1),\a_1,\a_2,(-a_1,1),(-a_2,1),(-\theta,2)\},
$\Den_{c,1}=\{(\theta,-k)\},$
$\Pi_c=\{(-o_2,k)\}.$
Thus 3.3.4 is $2o_2=\theta-{}^{s_{o_2}}\theta$.
The set $\Cen_c$ consists of the clans $\{s_0\},$
$\{1\},$ $\{s_2s_0\}$.
Further
$\Den_{c,s_0}=\emptyset$, $\Bc^c_{e_c,s_0}=\PP^1$,
and $|\Den_{c,w}|=1$, $|\Bc^c_{e_c,w}|=2$ if $w=1,s_2s_0$.
Thus we have
$$\Bc^c_{e_c}=\PP^1\sqcup\{4\ \roman{points}\},
\quad
\dim H_*(\Bc^c_{e_c},\CC)=6.$$
This is the sub-Coxeter case.

%(a)
%If $(n,m)=(3,2)$ then 
%$$\aligned
%&\Den_c=
%\{(\theta,-2),(2a_2+a_3,-1),(a_1+a_2+a_3,-1),
%\a_1,\a_2,\a_3,(-a_1,1),(-a_2,1),
%\hfill\\
%&\qquad\qquad(-a_3,1),
%(-a_1-a_2-a_3,2),
%(-2a_2-a_3,2),(-\theta,3)\},
%\hfill\\
%&\Delta^+_c=
%\{(-a_1-a_2,1),(-a_2-a_3,1),(-a_1-2a_2-a_3,2)\}.
%\hfill
%\endaligned$$
%Set $\b=(-o_2,1)$.
%$\Den[c]_1=\{(\theta,-1)\}$,
%and $\deg(\Ec_{c,w})=2$
%We have $s_\b=s_0s_1s_0$.
%Fix $n_{s_\b}\in N_{G[c]^\circ}(T_0)$, $g\in G[c]^\circ$, such that
%$n_{s_\b}^2=1$, $(\ad g)e_{(-\theta,2)}=e_{\a_2}+e_{(-\theta,2)}.$
%We have
%$\Gamma_2=\{\pm 1,\pm n_{s_\b}\},$
%$\een_1^{-1}(0)=(\ad\Gamma_2)(\ben_1)$,
%and
%$\een_{s_0s_2}^{-1}(0)=(\ad\Gamma_2 g)(\ben_{s_0s_2})$.
%Thus
%$\dim H_*(\Bc[c]_{e_c},\CC)^{\Gamma_2}=4$.
%A computation with rational Dunkl operators shows that
%$\Lb''_b=\CC[a_1^\vee,a_2^\vee]/(a_1^\vee a_2^\vee+(a_2^\vee)^2,
%(a_1^\vee)^2-2(a_2^\vee)^2).$

\vskip3mm

\noindent{\bf 4.3.1. Remark.}
The group $W_0$ is the dihedral group $I_2(4)$.
Set $(n,c)=(2,1/2)$.
Let $\Sb_c^\dag$, $\bar\Sb_c^\dag$ be as in 2.2.5(b).
We have 
%$A_2\simeq(\ZZ/2)^2$, 
$A_2^\circ\simeq\ZZ/2$,
and $W_2\simeq W_0$ by 3.3.9(c).
The group $W_2$ has two one-element orbits in $A_2^\circ$. 
So there are two non-isomorphic spherical Jordan-H\"older
composition factors in the $\Hb_{-c}$-module ${}^{IM}H_*(\Bc^c_{e_c},\CC)$, 
which are $\Lb_{-c}$ and $\Sb_{-c}$
(with multiplicity one).
Let $triv,sgn$ be the trivial character and the signature of $A_2^\circ$.
As a $\Hb_c\otimes A(c,e_c)^\circ$-module, we have
$$H_*(\Bc^c_{e_c},\CC)=
(\Lb_c^\dag\otimes triv)\oplus
(\Sb_{c}^\dag\otimes sgn)\oplus(\bar\Sb_{c}^\dag\otimes sgn).$$
The $\Hb_{-c}$-module ${}^{IM}\bar\Sb_{c}$ is not spherical.
The character is of ${\Lb_c''}^\dag$ is
$$ch({\Lb''_{c}}^\dag)=t^{-1}\cdot sgn+std+t\cdot sgn$$ 
(see \cite{C5} for details).
As $\CC[x_\l]$-modules, we have
$$\Lb_c^\dag=\CC_{h_{c}^\dag}\oplus V_{{}^{s_0}h_{c}^\dag}
\oplus\CC_{{}^{s_2s_0}h_{c}^\dag},\quad
\Sb_c^\dag=\CC_{h_{c}^\dag},\quad
\bar\Sb_c^\dag=\CC_{{}^{s_2s_0}h_{c}^\dag},
$$ 
where
$\CC_{h}$ is the unique one dimensional
$\CC[x_\l]$-module such that $x_\l\mapsto\l(h^\dag)$ and 
$V_{{}^{s_0}h_{c}^\dag}$
is a two-dimensional module equal to 
$\CC_{{}^{s_0}h_{c}^\dag}^{\oplus 2}$ up to a filtration.
%$\CC_{\lambdav}$ is the one dimensional
%$\CC[\xi_\l]$-module such that $x_\l\mapsto\l\cdot\lambdav$ and 
%$V_{{}^{s_0}\rhov_{c}}$
%is a two dimensional module equal to 
%$\CC_{{}^{s_0}\rhov_{c}}^{\oplus 2}$ up to a filtration.

%$\Bc[c]_{E,\Ien_w}=\{(\ad g_1)\Ien_{s_0s_2}, (\ad g_2)\Ien_{s_0s_2}\}$
%for some $g_1,g_2\in Z_{G[c]^\circ}(n_{s_{(-o_2,2c)}})$
%So, the representation of $\Hb''_{1/2}\times Z_{\tilde G[1/2]}(E)$ on
%$H_*(\Bc[c]_E,\CC)$ is irreducible.
%However, observe that
%$\tilde G[1/2]$ acts on $\gen[1/2]$ with an infinite number of orbits,
%because the orbit of $E$ has dimension $5$ and $\dim(\gen[1/2])=6$.

\subhead 4.4. Type $D_n$\endsubhead

%We have $h=2n-2$ and $\{d_i\}=\{2,4,...2n-4,2n,2n-2\}$.
Let $n\ge 4$.

If $n$ is odd then
$\EN=\{m=2r;r|n-1\ \roman{and}\ 2r\nmid n-1\}$,
$I_m=\{i=1,2,...n-2;r|i\}\cup\{n\}$,
$A_m\simeq(\ZZ/2)^{(n-1-r)/r}\times\ZZ/4$,
and $Z(G_0)=A_{2n-2}$.
We set $e_c=e_R+f_\theta\otimes\eps^{(2n-2)c}$.

If $n$ is even then
$\EN=\{m=2r; r|n-1 \ \roman{or}\ 2r|n\}$,
$I_m=\{i=1,2,...n-2;r|i\}
\cup\{n-1\ \roman{if}\ 2r|n\}\cup\{n\ \roman{if}\ r|n-1\}$,
$A_m\simeq(\ZZ/2)^{(n-1+r)/r}$ if $r|n-1$,
$A_m=(\ZZ/2)^{n/r}$ if $2r|n$,
and
$Z(G_0)=A_{2n-2}$.
We set $e_c=e_R+f_R\otimes\eps^{2c}$.

If $(n,m)=(4,4)$ then 
%\Den_c=\{(\theta,-1),\a_1,\a_2,\a_3,\a_4,(c_1,1),(c_2,1),(c_3,1)\},
$\Den_{c,1}=\{(\theta,-k)\},$
$\Pi_c=\{(b,k)\}$
where $b=-a_1-a_2-a_3-a_4$.
% $c_1=-a_1-a_2-a_3$, $c_2=-a_2-a_3-a_4$, $c_3=-a_1-a_2-a_4$.
Thus 3.3.4 is $-b=\theta-{}^{s_b}\theta$.
This is the sub-Coxeter case.

\subhead 4.5. Type $E_6$\endsubhead

We have
%$h=12$, $\{d_i\}=\{2,5,6,8,9,12\}$,
%$\{\theta_i\}=\{1,2,2,3,2,1\}$,
$\EN=\{3, 6, 9, 12\}$.
The elliptic regular conjugacy classes in $W_0$ are
$A_2^3$, $E_6(a_2)$, $E_6(a_1)$, $E_6$
respectively.
We have $I_3=\{3,5,6\}$,
$I_6=\{3,6\}$,
$I_9=\{5\}$,
$I_{12}=\{6\}$.
We have $A_3\simeq(\ZZ/3)^3$,
and $A_m=Z(G_0)\simeq\ZZ/3$ else.
We set $e_c=e_R+f_5\otimes\eps^{9c}$ if $m=9$ and
$e_c=e_R+f_\theta\otimes\eps^{12c}$ else.
The vector $f_5$ is the sum of a root vector
of weight $-a_1-a_2-2a_3-2a_4-a_5-a_6$
and a root vector of weight 
$-a_1-a_2-a_3-2a_4-2a_5-a_6$.

%The vector space $\gen[c]$ is a polar $G^\sim[c]$-module of rank
%$\sharp I_m=\dim(\sen[c])$.

(a) If $m=3$ the semisimple type of $\tilde G^c$ is $A_2^3$.
Further
$$\Den_{c,1}=\{(c_1,-3k), (c_2,-2k),(c_3,-2k),(c_4,-2k),(c_{i+4},-k)\},$$
$\Pi_c=\{(b_i,k), (b_6,k)\}$ with $i\in\underline{5}$,
%\a_1, \a_2, \a_3, \a_4, \a_5, \a_6, (c_4+b_1+b_6,1), (c_8+b_1,1),
%(c_8+b_1+b_6,1), \hfill\\ &\hskip1cm (c_6+b_2+b_5,1), (c_4+b_4+b_6,1),
%(c_4+b_1+b_2+b_6,2), \hfill\\ &\hskip1cm (c_6+b_2+b_3+b_5,2),
%(c_4+b_1+b_4+b_6,2), (c_4+b_3+b_4+b_6,2), \hfill\\ &\hskip1cm
%(c_4+b_1+b_3+b_4+b_6,3), (c_4+b_1+b_2+b_4+b_6,3), (-\theta,4)
$b_1=-a_1-a_3-a_4$,
%=-o_1+o_2-o_4+o_5$, 
$b_2=-a_2-a_4-a_5$,
%=-o_2+o_3-o_5+o_6$, 
$b_3=-a_2-a_3-a_4$,
%=o_1-o_2-o_3+o_5$, 
$b_4=-a_4-a_5-a_6$,
%=o_2+o_3-o_4-o_6$, 
$b_5=-a_3-a_4-a_5$,
%=o_1+o_2-o_3-o_5+o_6$, 
$b_6=-a_1-a_2-a_3-a_4-a_5-a_6$,
%=-o_1-o_2+o_4-o_6$, 
and
$c_1=a_1+a_2+2a_3+3a_4+2a_5+a_6$,
$c_2=a_2+a_3+2a_4+2a_5+a_6$,
$c_3=c_1+b_4$,
$c_4=c_1+b_5$,
$c_5=c_3+b_5$,
$c_6=c_2+b_4$,
$c_7=c_4+b_1$,
$c_8=c_3+b_3$,
$c_9=c_1+b_1+b_2$.
%$a_1=c_8+b_6$,
%$a_2=c_6+b_5$,
%$a_3=c_9+b_4$,
%$a_4=c_4+b_6$,
%$a_5=c_8+b_1$,
%$a_6=c_7+b_2$.
Thus 3.3.4 is
$$-9(b_1+b_2)(b_3+b_4)(b_5+b_6)\prod_{i=1}^6b_i=
\sum_{w\in W_c}(-1)^ {\ell(w)}\prod_{i=1}^9{}^wc_i.$$

(b) If $m=6$ the semisimple type of $\tilde G^c$ is $A_1^3$.
Further
$\Den_{c,1}=\{(c_i,-k)\},$
%\a_1,\a_2,...\a_6,(a_1+b_1,1),(a_1+b_2,1),(a_2+b_2,1),(a_6+b_2,1),(-\theta,2)
$\Pi_c=\{(b_i,k)\}$ with $i\in\underline{3},$
$b_1=-a_1-a_2-a_3-2a_4-a_5$,
%=-o_1+o_3-o_4+o_6$, 
$b_2=-a_1-a_2-a_3-a_4-a_5-a_6$,
%=-o_1-o_2+o_4-o_6$, 
$b_3=-a_2-a_3-2a_4-a_5-a_6$,
%=o_1-o_4+o_5-o_6$, 
and
$c_1=a_1+a_2+2a_3+2a_4+a_5$,
$c_2=a_1+a_2+a_3+2a_4+a_5+a_6$,
$c_3=a_2+a_3+2a_4+2a_5+a_6$.
Thus 3.3.4 is
$$-b_1b_2b_3=(1-s_{b_1})(1-s_{b_2})(1-s_{b_3})(c_1c_2c_3).$$

(c) If $m=9$ the semisimple type of $\tilde G^c$ is $A_1$.
Further
$\Den_{c,1}=\{(c,-k)\},$
%\a_1,\a_2,\a_3,\a_4,\a_5,\a_6, (a_3+b,1), (a_5+b,1)
$\Pi_c=\{(b,k)\}$
where $b=-a_1-a_2-2a_3-2a_4-2a_5-a_6$,
%=-o_3+o_4-o_5$, 
$c=a_1+a_2+2a_3+3a_4+2a_5+a_6$.
Thus 3.3.4 is $-b=c-{}^{s_b}c.$

%The elements $f_i$ can be computed explicitely.
%We have
%$$\aligned
%&
%f_1=16 f_{a_1}+22 f_{a_2}+30 f_{a_3}+42 f_{a_4}+30 f_{a_5}+16 f_{a_6},\quad
%f_2=f_{b_1}+f_{b_2}+f_{b_3}+f_{b_4},
%\cr
%&
%\qquad
%f_3=f_{b_5}+f_{b_6}+f_{b_7}+f_{b_8},\quad
%f_4=f_{b_9}+f_{b_{10}}+f_{b_{11}},\quad
%f_5=f_{b_{12}}+f_{b_{13}},
%\endaligned$$
%where the $f_b$ is a root vector of weight $-b$,
%and $b_1,b_2,...b_{13}$ are equal to
%$a_1+a_2+a_3+a_4,$
%$a_1+a_3+a_4+a_5,$
%$a_2+a_4+a_5+a_6,$
%$a_3+a_4+a_5+a_6,$
%$a_1+a_2+a_3+a_4+a_5,$
%$a_1+a_3+a_4+a_5+a_6,$
%$a_2+a_3+2a_4+a_5,$
%%$a_2+a_3+a_4+a_5+a_6,$
%$a_1+a_2+2a_3+2a_4+a_5,$
%$a_1+a_2+a_3+2a_4+a_5+a_6,$
%$a_2+a_3+2a_4+2a_5+a_6,$
%$a_1+a_2+2a_3+2a_4+a_5+a_6,$
%and
%$a_1+a_2+a_3+2a_4+2a_5+a_6$
%respectively.

\subhead 4.6. Type $E_7$\endsubhead

We have
%$h=18$, $\{d_i\}=\{2,6,8,10,12,14,18\}$,
%$\{\theta_i\}=\{2,2,3,4,3,2,1\}$,
$\EN=\{2,6,14,18\}$.
The elliptic regular conjugacy classes in $W_0$ are
$A_1^7$, $E_7(a_4)$, $E_7(a_1)$, $E_7$
respectively.
We have $I_2=I_0$,
$I_6=\{2,5,7\}$,
$I_{14}=\{6\}$,
$I_{18}=\{7\}$.
We have $A_2\simeq(\ZZ/2)^7$,
and $A_m=Z(G_0)\simeq\ZZ/2$ else.
We set $e_c=e_R+f_6\otimes\eps^{14c}$ if $m=14$ and
$e_c=e_R+f_\theta\otimes\eps^{18c}$ else.

(a) If $m=2$ the semisimple type of $\tilde G^c$ is $A_7$.
Further
$$\aligned
&\Den_{c,1}=\{
(c_{28},-13k),
(c_{27},-11k),
(c_{26},-9k),
(c_{25},-9k),(c_{21+i},-7k)\}_{i\in\underline{3}}\cup
\hfill\\
&\hskip12mm
\cup\{(c_{17+i},-5k)\}_{i\in\underline{4}}
\cup\{(c_{12+i},-3k)\}_{i\in\underline{5}}
\cup\{(c_{6+i},-2k),(c_i,-k)\}_{i\in\underline{6}},
\hfill\\
&\Pi_c=\{(b_1,2k),(b_{1+i},k)\}_{i\in\underline{6}},
\hfill
\endaligned$$
where
$b_1=-a_2-a_3-a_4-a_5$,
$b_2=-a_6-a_7$,
$b_3=-a_4-a_5$,
$b_4=-a_1-a_3$,
$b_5=-a_2-a_4$,
$b_6=-a_5-a_6$,
$b_7=-a_3-a_4$,
and
$c_1=a_1+a_3+a_4$,
$c_2=a_2+a_3+a_4$,
$c_3=a_2+a_4+a_5$,
$c_4=a_3+a_4+a_5$,
$c_5=a_4+a_5+a_5$,
$c_6=a_5+a_6+a_7$,
$c_7=a_1+a_2+a_3+a_4+a_5$,
$c_8=a_1+a_3+a_4+a_5+a_6$,
$c_9=a_2+a_3+a_4+a_5+a_6$,
$c_{10}=a_2+a_3+2a_4+a_5$,
$c_{11}=a_2+a_4+a_5+a_6+a_7$,
$c_{12}=a_3+a_4+a_5+a_6+a_7$,
$c_{13}=a_1+a_2+a_3+a_4+a_5+a_6+a_7$,
$c_{14}=a_2+a_3+2a_4+a_5+a_6+a_7$,
$c_{15}=a_1+a_2+2a_3+2a_4+a_5$,
$c_{16}=a_1+a_2+a_3+2a_4+a_5+a_6$,
$c_{17}=a_2+a_3+2a_4+2a_5+a_6$,
$c_{18}=a_1+a_2+2a_3+2a_4+2a_5+a_6$,
$c_{19}=a_1+a_2+2a_3+2a_4+a_5+a_6+a_7$,
$c_{20}=a_1+a_2+a_3+2a_4+2a_5+a_6+a_7$,
$c_{21}=a_2+a_3+2a_4+2a_5+2a_6+a_7$,
$c_{22}=a_1+2a_2+2a_3+3a_4+2a_5+a_6$,
$c_{23}=a_1+a_2+2a_3+2a_4+2a_5+2a_6+a_7$,
$c_{24}=a_1+a_2+2a_3+3a_4+2a_5+a_6+a_7$,
$c_{25}=a_1+a_2+2a_3+3a_4+3a_5+2a_6+a_7$,
$c_{26}=a_1+2a_2+2a_3+3a_4+2a_5+2a_6+a_7$,
$c_{27}=a_1+2a_2+2a_3+4a_4+3a_5+2a_6+a_7$,
$c_{28}=2a_1+2a_2+3a_3+4a_4+3a_5+2a_6+a_7$.
Thus 3.3.4 is
$$-2^6\prod_{i=1}^{28}b_i=
\sum_{w\in W_c}(-1)^{\ell(w)}\prod_{i=1}^{28}{}^wc_i.$$
where $b_1,...b_{28}$ are the positive roots of the system of type $A_7$ with
basis $b_1,...b_7$.

(b) If $m=6$ the semisimple type of $\tilde G^c$ 
is $A_2^2\times A_1$.
Further we have
$\Den_{c,1}=\{(c_i,-k),(c_6,-2k),(c_7,-2k)\},$
$\Pi_c=\{(b_i,k)\}$ with $i\in\underline{5}$,
$b_1=-a_1-a_2-a_3-2a_4-a_5$,
$b_2=-a_2-a_3-a_4-a_5-a_6-a_7$,
$b_3=-a_2-a_3-2a_4-a_5-a_6$,
$b_4=-a_1-a_3-a_4-a_5-a_6-a_7$,
$b_5=-a_1-a_2-a_3-a_4-a_5-a_6$,
and
$c_1=a_1+a_2+a_3+a_4+a_5+a_6+a_7$,
$c_2=a_2+a_3+2a_4+a_5+a_6+a_7$,
$c_3=a_1+a_2+2a_3+2a_4+a_5$,
$c_4=a_1+a_2+a_3+2a_4+a_5+a_6$,
$c_5=a_2+a_3+2a_4+2a_5+a_6$,
$c_6=a_1+a_2+2a_3+3a_4+3a_5+2a_6+a_7$,
$c_7=a_1+2a_2+2a_3+3a_4+2a_5+2a_6+a_7$.
Thus 3.3.4 is
$$-(b_1+b_2)(b_3+b_4)\prod_{i=1}^5b_i=
\sum_{w\in W_c}(-1)^{\ell(w)}\prod_{i=1}^7{}^wc_i.$$

(c) If $m=14$ the semisimple type of $\tilde G^c$ is $A_1$.
Further
$\Den_{c,1}=\{(c,-k)\},$
$\Pi_c=\{(b,k)\}$
where $b=-a_1-2a_2-2a_3-3a_4-3a_5-2a_6-a_7$,
$c=a_1+2a_2+2a_3+4a_4+3a_5+2a_6+a_7$.
Thus 3.3.4 is $-b=c-{}^{s_b}c.$

\subhead 4.7. Type $E_8$\endsubhead

We have
%$h=30$, $\{d_i\}=\{2,8,12,14,18,20,24,30\}$,
%$\{\theta_i\}=\{2,3,4,6,5,4,3,2\}$,
$\EN=\{2, 3, 4, 5, 6, 8, 10, 12, 15, 20, 24, 30\}$.
The elliptic regular conjugacy classes in $W_0$ are
$A_1^8$, $A_2^4$, $D_4(a_1)^2$, $A_4^2$, $E_8(a_8)$, $D_8(a_3)$,
$E_8(a_6)$, $E_8(a_3)$, $E_8(a_5)$, $E_8(a_2)$, $E_8(a_1)$, $E_8$
respectively.
We have $I_2=I_0$,
$I_3=I_6=\{3,5,7,8\}$,
$I_4=\{2,3,6,7\}$,
$I_8=\{2,7\}$,
$I_5=I_{10}=\{6,8\}$,
$I_{12}=\{3,7\}$,
$I_{20}=\{6\}$,
$I_{24}=\{7\}$,
$I_{15}=I_{30}=\{8\}$.
We have
$A_2\simeq(\ZZ/2)^8$,
$A_3\simeq(\ZZ/3)^4$,
$A_4\simeq(\ZZ/2)^4$,
$A_5\simeq(\ZZ/5)^2$,
$A_8\simeq(\ZZ/2)^2$,
and
$A_m=Z(G_0)=\{1\}$ else.
We set $e_c=e_R+f_6\otimes\eps^{20c}$ if $m=20$,
$e_c=e_R+f_7\otimes\eps^{24c}$ if $m=4,8,12,24$,
and
$e_c=e_R+f_\theta\otimes\eps^{30c}$ else.

(a) If $m=2$ the semisimple type of $\tilde G^c$ is $D_8$.

(b) If $m=3$ the semisimple type of $\tilde G^c$ is $A_7\times A_1^8$.

(c) If $m=4$ the semisimple type of $\tilde G^c$ 
is $D_4\times A_2^2\times A_1^8$.

(d) If $m=5$ the semisimple type of $\tilde G^c$ is $A_4\times A_3\times A_1^4$.

(e) If $m=6$ the semisimple type of $\tilde G^c$ is $A_4\times A_3$.

(f) If $m=8$ the semisimple type of $\tilde G^c$ is $A_3\times A_2\times A_1^2$.

(g) If $m=10$ the semisimple type of $\tilde G^c$ is $A_2^2\times A_1^2$.

(h) If $m=12$ the semisimple type of $\tilde G^c$ is $A_2\times A_1^3$.

(i) If $m=15$ the semisimple type of $\tilde G^c$ is $A_1^4$.
Further
$\Den_{c,1}=\{(c_i,-k)\},$
$\Pi_c=\{(b_i,k)\}$
with $i\in\underline 4$, 
$b_1=-a_1-2a_2-2a_3-3a_4-3a_5-2a_6-a_7-a_8$,
$b_2=-a_1-2a_2-2a_3-3a_4-2a_5-2a_6-2a_7-a_8$,
$b_3=-a_1-a_2-2a_3-3a_4-3a_5-2a_6-2a_7-a_8$,
$b_4=-a_1-2a_2-2a_3-4a_4-3a_5-2a_6-a_7$,
$c_1=a_1+2a_2+2a_3+4a_4+3a_5+2a_6+a_7+a_8$,
$c_2=a_1+2a_2+2a_3+3a_4+3a_5+2a_6+2a_7+a_8$,
$c_3=a_1+a_2+2a_3+3a_4+3a_5+3a_6+2a_7+a_8$,
$c_4=a_1+2a_2+3a_3+4a_4+3a_5+2a_6+a_7$.
Thus 3.3.4 is 
$\prod_ib_i=\bigl(\prod_i(1-s_{b_i})\bigr)\bigl(\prod_ic_i\bigr).$

(j) If $m=20$ the semisimple type of $\tilde G^c$ is $A_1^2$.
Further $\Den_{c,1}=\{(c_i,-k)\},$
$\Pi_c=\{(b_1,k), (b_2,k)\}$
with $i\in\underline 2$,
$b_1=-2a_1-2a_2-3a_3-4a_4-3a_5-3a_6-2a_7-a_8$,
$b_2=-a_1-2a_2-3a_3-4a_4-4a_5-3a_6-2a_7-a_8$,
$c_1=a_1+2a_2+3a_3+5a_4+4a_5+3a_6+2a_7+a_8$,
$c_2=2a_1+2a_2+3a_3+4a_4+4a_5+3a_6+2a_7+a_8$.
Thus 3.3.4 is $b_1b_2=(1-s_{b_1})(1-s_{b_2})(c_1c_2).$

(k) If $m=24$ the semisimple type of $\tilde G^c$ is $A_1$.
Further
$\Den_{c,1}=\{(c,-k)\},$
$\Pi_c=\{(b,k)\}$
where $b=-2a_1-3a_2-4a_3-5a_4-4a_5-3a_6-2a_7-a_8$,
$c=2a_1+3a_2+4a_3+6a_4+4a_5+3a_6+2a_7+a_8=o_4-o_5$.
Thus 3.3.4 is $-b=c-{}^{s_b}c.$

\subhead 4.8. Type $F_4$\endsubhead

We have
%$h=12$, $\{d_i\}=\{2,6,8,12\}$,
%$\{\theta_i\}=\{2,3,4,2\}$,
$\EN=\{2, 3, 4, 6, 8, 12\}$.
The elliptic regular conjugacy classes in $W_0$ are
$A_1^4$, $A_2\times\tilde A_2$, $D_4(a_1)$, $F_4(a_1)$, $B_4$, $F_4$
respectively.
We have $I_2=I_0$,
$I_3=I_6=\{2,4\}$,
$I_4=\{3,4\}$,
$I_8=\{3\}$,
$I_{12}=\{4\}$.
We have $A_2\simeq(\ZZ/2)^4$,
$A_3\simeq(\ZZ/3)^2$,
$A_4\simeq(\ZZ/2)^2$,
$A_8\simeq\ZZ/2$,
and $A_m=Z(G_0)=\{1\}$ else.
We set $e_c=e_R+f_3\otimes\eps^{8c}$ if $m=8$ and
$e_c=e_R+f_\theta\otimes\eps^{12c}$ else.

(a) If $m=2$ the semisimple type of $\tilde G^c$ is $A_3\times A_1^2$.

(b) If $m=3$ the semisimple type of $\tilde G^c$ is $A_2\times A_1^2$.

(c) If $m=4$ the semisimple type of $\tilde G^c$ is $A_2\times A_1$.

(d) If $m=6$ the semisimple type of $\tilde G^c$ is $A_1^2$.
Further
$\Den_{c,1}=\{(c_i,-k)\},$
$\Pi_c=\{(b_i,k)\}$ with $i\in\underline 2$,
$b_1=-a_1-2a_2-2a_3-a_4$,
$b_2=-a_1-a_2-2a_3-2a_4$,
$c_1=a_1+2a_2+3a_3+a_4$,
$c_2=a_1+2a_2+2a_3+2a_4$.
Thus 3.3.4 is $b_1b_2=(1-s_{b_1})(1-s_{b_2})(c_1c_2).$

(e) If $m=8$ the semisimple type of $\tilde G^c$ is $A_1$.
Further
$\Den_{c,1}=\{(c,-1)\},$
$\Pi_c=\{(b,1)\}$
where $b=-a_1-2a_2-3a_3-2a_4=-o_4$,
$c=a_1+2a_2+4a_3+2a_4$.
Thus 3.3.4 is $-2b=c-{}^{s_b}c.$

\subhead 4.9. Type $G_2$\endsubhead

We have
%$h=6$, $\{d_i\}=\{2,6\}$, $\{\theta_i\}=\{3,2\}$,
%$f_1=6f_{a_1}+10f_{a_2}$, $f_2=f_\theta$,
$\EN=\{2, 3, 6\}$.
The elliptic regular conjugacy classes in $W_0$ are
$A_1\times\tilde A_1$, $A_2$, $G_2$
respectively.
We have $I_2=\{1,2\}$, $I_3=I_6=\{2\}$.
We have $A_2\simeq(\ZZ/2)^2$,
$A_3\simeq\ZZ/3$,
and
$A_6=Z(G_0)=\{1\}$.
We set $e_c=e_R+f_\theta\otimes\eps^{6c}$.

\proclaim{4.9.1. Proposition}
We have $\dim\,H_*(\Bc^c_{e_c},\CC)=k^2, 8k^2, 18k^2$
if $m=6,3,2$ respectively, and
$\dim\,H_*(\Bc^c_{e_c},\CC)^{A(c,e_c)}=36c^2$.
\endproclaim

\noindent{\sl Proof :}
(a)
If $m=3$ the semisimple type of $\tilde G^c$ is $A_1$.
Further
$\Den_{c,1}=\{(3a_1+a_2,-k)\},$
%\a_1,\a_2,(-a_1-a_2,1),(-\theta,2)\},
$\Pi_c=\{(-o_1,k)\}.$
Thus 3.3.4 is $3o_1=(1-s_{o_1})(3a_1+a_2).$
Now, set $k=1$ to simplify.
The set $\Cen_c$ consists of the clans $\{1,s_0\},$ $\{s_0s_2\}$.
Further we have
$\Den_{c,s_0s_2}=\emptyset$, $\Bc^c_{e_c,s_0s_2}=\PP^1$,
and $|\Den_{c,w}|=1$, $|\Bc^c_{e_c,w}|=3$ if $w=1,s_0$.
Thus 
$$\Bc^c_{e_c}=\PP^1\sqcup\{6\ \roman{points}\},
\quad
\dim H_*(\Bc^c_{e_c},\CC)=8,
\quad
\dim H_*(\Bc^c_{e_c},\CC)^{A(c,e_c)}=4.$$
This is the sub-Coxeter case.

(b)
If $m=2$ the semisimple type of $\tilde G^c$ is $A_1^2$.
Further
$\Den_{c,1}=\{(c_i,-ik)\},$
%,\a_1,\a_2,(-a_1,1),(-a_2,1),(-2a_1-a_2,2),(-\theta,3)\},
$\Pi_c=\{(b_i,ik)\}$
with $i\in\underline 2$,
$b_1=-a_1-a_2$, $b_2=-3a_1-a_2$,
$c_2=\theta$, $c_1=2a_1+a_2.$
%We have $\check b_1=\ov_1-3\ov_2$, $\check b_2=-\ov_1+\ov_2$.
Thus 3.3.4 is $4b_1b_2=(1-s_{b_1})(1-s_{b_2})(c_1c_2).$
Now, set $k=1$ to simplify.
The set $\Cen_c$ consists of the clans
$$\{1,s_0,s_0s_2\},\ \{s_0s_2s_1\},\ \{s_0s_2s_1s_2\}.$$
We have
$\Den_{c,s_0s_2s_1s_2}=\emptyset$ and $\Bc^ c_{e_c,s_0s_2s_1s_2}=(\PP^1)^2$,
$|\Den_{c,s_0s_2s_1}|=1$ and $\Bc^c_{e_c,s_0s_2s_1}=\PP^1$,
$|\Den_{c,w}|=2$ and $|\Bc^c_{e_c,w}|=4$ if $w=1,s_0,s_0s_2$
because the second Chern number of the vector bundle
$\Ec_{c,w}$ over $\Bc^c_w$ is $4$.
We have also
$$\Bc^c_{e_c}=(\PP^1)^2\sqcup\PP^1\sqcup\{12\ \roman{points}\},
\quad
\dim H_*(\Bc^c_{e_c},\CC)=18,
\quad
\dim H_*(\Bc^c_{e_c},\CC)^{A(c,e_c)}=9.
$$

\qed
%$\Den_{c,s_0s_2s_1}=\{(\theta,-2)\}$,
%$\Den_{c,s_2s_0s_2s_1}=\{(\theta,-2),\a_2\}$,
%$\Den_{c,w}=\{(\theta,-2),(o_1,-1)\}$
%(use Demazure's operators)
%$w=s_2s_0s_2s_1$
%(observe that $\theta a_2$ is zero in
%$H^*(\Bc^c_w,\CC)
%\simeq\CC[\ten_0]/((a_2+\theta)^2,(a_2-\theta)^2)$).

\vskip3mm

\noindent{\bf 4.9.2. Remarks.}
(a)
Since $W_0$ is the dihedral group $I_2(6)$,
we have $\dim({\Lb''_{c}}^\dag)=36c^2$
if $m=2,3,6$, $k>0$ 
by \cite{C5}.

(b)
Assume that $c=1/3$ (the sub-Coxeter case).
Let $\Sb_c^\dag$ be the $\Hb_c$-module associated to the
$\Hb'_c$-module ${\Sb'_c}^\1dag$ in 2.3.5(d),
and $\Sb_{-c}={}^{IM}\Sb_c^\dag$.
We have
$A_3=A_3^\circ\simeq\ZZ/3$
and
$W_3/\la w_3\ra\simeq\ZZ/2$.
The group $W_3$ has one one-element orbit and one two-elements
orbits in $A_3^\circ$.
So there are two non-isomorphic spherical Jordan-H\"older composition factors
in the $\Hb_{-c}$-module ${}^{IM}H_*(\Bc^c_{e_c},\CC)$, 
which are $\Lb_{-c}$ (with multiplicity one)
and $\Sb_{-c}$ (with multiplicity two).
We have $A(c,e_c)\simeq I_2(3)$.
Let $triv,$ $std$ denote the trivial and the standard 
representation of $I_2(3)$.
As a $\Hb_c\otimes A(c,e_c)$-module, we have
$$H_*(\Bc^c_{e_c},\CC)=
(\Lb_c^\dag\otimes triv)\oplus(\Sb_{c}^\dag\otimes std).$$
The group $I_2(6)$ has four irreducible one-dimensional representations
and two two-dimensional ones.
Let $triv$, $sgn$, $std$ denote the trivial representation, 
the signature, and the standard representation.
The character of ${\Lb''_{c}}^\dag$ is 
$$ch({\Lb''_{c}}^\dag)=t^{-1}\cdot sgn+std+t\cdot sgn$$
(see loc.~cit.~ for details).
As $\CC[x_\l]$-modules we have
$$\Lb_c^\dag=\CC_{h_c^\dag}\oplus\CC_{{}^{s_0}h_c^\dag}\oplus
V_{{}^{s_2s_0}h_c^\dag},\quad
\Sb_{c}^\dag=\CC_{h_c^\dag}\oplus\CC_{{}^{s_0}h_c^\dag}.$$

(c)
Assume that $c=1/2$.
Let $\Sb_c^\dag$ be the $\Hb_c$-module associated to the
$\Hb'_c$-module $\Sb'_c$ in 2.3.5(d),
and $\Sb_{-c}={}^{IM}\Sb_c^\dag$.
We have
$A_2^\circ\simeq(\ZZ/2)^2$
and
$W_2=W_0$.
The group $W_2$ has one one-element orbit and one three-elements
orbits in $A_2^\circ$.
So there are two non-isomorphic spherical Jordan-H\"older composition factors
in the $\Hb_{-c}$-module ${}^{IM}H_*(\Bc^c_{e_c},\CC)$, 
which are $\Lb_{-c}$ (with multiplicity one)
and $\Sb_{-c}$ (with multiplicity three).
We have $A(c,e_c)\simeq \ZZ/3\ltimes(\ZZ/2)^2$.
As a $\Hb_c\otimes A(c,e_c)$-module, we have
$$H_*(\Bc^c_{e_c},\CC)=
(\Lb_c^\dag\otimes triv)\oplus(\Sb_c^\dag\otimes rep),$$
where $rep$ is a three dimensional irreducible representation of
$A(c,e_c)$.
%The character is (see loc.~cit.~ for details)
As $\CC[x_\l]$-modules we have
$$\aligned
&\Lb_c^\dag=
\CC_{h_c^\dag}\oplus\CC_{{}^{s_0}h_c^\dag}\oplus\CC_{{}^{s_2s_0}h_c^\dag}\oplus
V_{{}^{s_1s_2s_0}h_c^\dag}\oplus V_{{}^{s_2s_1s_2s_0}h_c^\dag},
\\
&\Sb_{c}^\dag=
\CC_{h_c^\dag}\oplus\CC_{{}^{s_0}h_c^\dag}\oplus\CC_{{}^{s_2s_0}h_c^\dag}.
\endaligned$$

(d)
Let $O_{e_c}$
be the $\tilde G^c$-orbit of $e_c.$
If $m=3,6$ then
$O_{e_c}$ is open,
while
if $m=2$ then
$\dim(O_{e_c})=\dim(\gen^c)-1$.
Recall that if $O_{e_c}$ is open in $\gen^c$,
then there is a finite number of $\tilde G^c$-orbits in $\gen^c$
and a standard argument implies that the $A(c,e_c)$-isotypic components
in $H_*(\Bc^c_{e_c},\CC)$
are irreducible as $\Hb_c\otimes A(c,e_c)$-submodules.

\head A. Appendix\endhead

\subhead A.1. Fourier-Sato transform\endsubhead

Fix $h=(s,\tau,\zeta)\in\hat H$, 
with $\tau$ not a root of unity and $\zeta=\tau^{c/2}$.
Set $\tilde s=(s,\tau)$. 
Recall the Fourier-Sato transform 
$FS:D^b_{\tilde G^h\times\CC^\times}(\gen^h)\to
D^b_{\tilde G^{h^{\ddag}}\times\CC^\times}(\gen^{h^{\ddag}}).$

\proclaim{A.1.1. Proposition}
(a)
There is a bijection
$\Xc_{h,w}\to\Xc_{h^{\ddag},w}$, $\chi\mapsto\chi^\ddag$
such that the complexes $FS(IC(\chi))$, $IC(\chi^\ddag)$ are isomorphic.

(b)
If $\chi\in\Xc_h$ then the $\hat\Hb_{-c}$-modules
${}^{IM}L_{h,\chi}$, $L_{h^\ddag,\chi^\ddag}$
have the same weights.
%If $\chi\in\Xc_h$ then ${}^{IM}L_{h,\chi}\simeq L_{h^\ddag,\chi^\ddag}$
%as $\Hb_{-c}$-modules.
\endproclaim

\noindent{\sl Proof :}
The proof of part (a) is sketched in \cite{V1, sec.~7.4}.
Let us give more details.
It is modelled after \cite{L3, cor.~10.5}, \cite{M1}.
We have $W_h=W_{h^{\ddag}}$.
We must prove that if $IC(\chi)$ is a direct summand in
the semisimple complex
$\Lc_{h,w}$ then $FS(IC(\chi))$ is a direct summand in
$\Lc_{h^{\ddag},w}$.
Thus, it is enough to check that
for each $w\in W$ we have 
$$FS(\Lc_{h,w})\simeq\Lc_{h^{\ddag},w}.\leqno(A.1.2)$$

For any finite dimensional representation $V$ of
$\tilde G^h\times\CC^\times$ let
$$F_h: D_{\tilde G^h\times\CC^\times}^b(V)\to 
D_{B^h\times\CC^\times}^b(V)$$
be the forgetful functor.
Consider the diagram
$$V{\buildrel \a^h\over\lla}\,
\tilde G^h\times V
\,{\buildrel \b^h\over\lra}\,
\tilde G^h\times_{B^h} V
\,{\buildrel\gamma^h\over\lra}\,
V,$$
where $\a^h$ is the second projection,
$\b^h$ is the obvious quotient map (a $B^h$-torsor),
and $\g^h$ is the action map $(g,x)\mapsto (\ad g)(x)$.
We have $\dim(\tilde G^h/B^h)=n_h$.
Thus $F_h$ admits a left adjoint equal to
$\Gamma_h[n_h]$, where $\Gamma_h$ is Bernstein's induction functor
$$\Gamma_h :
D_{B^h\times\CC^\times}^b(V)\to D_{\tilde G^h\times\CC^\times}^b(V),
\quad
\Fc\mapsto \g^h_*(\b^{h,*})^{-1}\a^{h,*}(\Fc).$$
See \cite{M1, sec.~1.7}, \cite{MV}.
The functor $\Gamma_h$ commutes with all standard functors, including
$FS$, because the same is true for $F_h$.

Since $\tau$ is not a root of unity and $c\neq 0$ we have $\zeta\neq 1$.
Hence the fixed points set $\hen^h$ is $\{0\}$.
For each $w\in W$ we consider the diagram
$$\{0\}{\buildrel a^h\over\lla}\,\ben_w^h
\,{\buildrel b^h\over\lra}\,\gen^h,$$
where $b^h$ is the obvious inclusion.
Observe that $D^b_{B_w^h\times\CC^\times}(\{0\})=D^b_{\hat H}(\{0\})$.
The restriction functor is
$$\Res_{h,w}:D_{\tilde G^h\times\CC^\times}^b(\gen^h)\to D_{\hat H}^b(\{0\}),
\quad
\Fc\mapsto a^h_*b^{h,!}F_h(\Fc).$$
It admits a left adjoint, which is the induction functor
$$\Ind_{h,w}:D_{\hat H}^b(\{0\})\to D_{\tilde G^h\times\CC^\times}^b(\gen^h),
\quad
\Fc\mapsto
\Gamma_hb^h_!a^{h,*}(\Fc[n_h])=
\Gamma_hb^h_*a^{h,!}(\Fc[n_h]).$$
Now, consider the Cartesian square
$$\matrix
\ben^h_w&{\buildrel a^h\over\lra}&\{0\}\cr
\hskip-5mm{\ss b^h}\downarrow&&\hskip-5mm{\ss d^h}\downarrow\cr
\gen^h&{\buildrel c^h\over\lra}&\gen^h/\uen^h_w.
\endmatrix$$
The base change formula yields
$$\Res_{h^{\ddag},w}\circ FS=
d^{h^{\ddag},!}\circ c^{h^{\ddag}}_*\circ F_h\circ FS.$$
The functor $FS$ commutes with $F_h$.
Since the dual of $\gen^{h^{\ddag}}/\uen^{h^{\ddag}}_w$ 
is identified with $\ben_w^h$,
by \cite{KS1, prop.~3.7.14} we have 
$$FS\circ b^{h,!}=c^{h^{\ddag}}_*\circ FS$$ 
as functors
$D^b_{B_w^h\times\CC^\times}(\gen^{h})\to 
D^b_{B_w^h\times\CC^\times}(\gen^{h^{\ddag}}/\uen^{h^{\ddag}}_w)$.
For the same reason we have
$$a^h_*=d^{h^{\ddag},!}\circ FS$$
as functors $D^b_{B_w^h\times\CC^\times}(\ben^h_w)\to D^b_{\hat H}(\{0\})$.
Therefore
$$\Res_{h^{\ddag},w}\circ FS=a_*^h\circ b^{h,!}\circ F_h=\Res_{h,w}.
$$
See also \cite{M1, lem.~4.2}.
Since $\Ind_{h,w}$ is left adjoint to $\Res_{h,w}$ we have also
$$FS\circ\Ind_{h,w}=\Ind_{h^{\ddag},w}.\leqno(A.1.3)$$
Finally we compute the complex $\Ind_{h,w}(\CC_{\{0\}})$.
Observe that $\dot\Nc_w^h=\tilde G^h\times_{B_w^h}\ben^h_w$,
because $\hen^h=\{0\}$.
The diagram
$$\matrix
\gen^h&{\buildrel \a^h\over\lla}&
\tilde G^h\times\gen^h&{\buildrel \b^h\over\lra}&
\tilde G^h\times_{B_w^h}\gen^h&{\buildrel\gamma^h\over\lra}&\gen^h\cr
\hskip-5mm{\ss b^h}\uparrow&&
\hskip-5mm{\ss f^h}\uparrow&&
\hskip-5mm{\ss g^h}\uparrow&&\parallel\cr
\ben^h_w&\lla&\tilde G^h\times\ben^h_w&\lra&
\dot\Nc_w^h&{\buildrel\pi^h_w\over\lra}&\gen^h
\endmatrix$$
is formed of Cartesian squares.
Thus, the base change formula yields
$$\aligned
\Ind_{h,w}(\CC_{\{0\}})
&=\g^h_*(\b^{h,*})^{-1}\a^{h,*}b^h_!a^{h,*}(\CC_{\{0\}})
[n_h]
\cr
&=\g^h_*(\b^{h,*})^{-1}\a^{h,*}b^h_!(\CC_{\ben^h_w})
[n_h]
\cr
&=\g^h_*(\b^{h,*})^{-1}f^h_!(\CC_{\tilde G^h\times\ben^h_w})
[n_h]
\cr
&=\g^h_*g^h_!(\CC_{\dot\Nc^h_w})
[n_h].
\endaligned$$
Therefore we have
$$\Ind_{h,w}(\CC_{\{0\}})=\Lc_{h,w}[n_h-d_{h,w}].
\leqno(A.1.4)$$
Applying the isomorphism of functors A.1.3 to the object $\CC_{\{0\}}$,
formula A.1.4 yields a canonical isomorphism
$$\nu_{h,w}\ :\ \Lc_{h^\ddag,w}\to FS(\Lc_{h,w}).$$

Now we prove (b).
Consider the algebra
$$\Extb_h=\prod_v\bigoplus_w\Ext_{D^b(\Nc^h)}^*(\Lc_{h,v},\Lc_{h,w}),$$
with the Yoneda product,
where $v,w$ belong to a set of representatives of the right $W_h$-cosets in $W$.
The functoriality of $\pi^h_{w,*}$ 
yields an algebra homomorphism
$$H^*(\dot\Nc^h_w,\CC)=
\Ext_{D^b(\dot\Nc^h_w)}^*(\CC_{\dot\Nc^h_w},\CC_{\dot\Nc^h_w})\to
\Ext_{D^b(\Nc^h)}^*(\Lc_{h,w},\Lc_{h,w})$$
for each $w$.
Recall that
$$\CC_\kappa[\xi_\l]=\CC[\hen]\subset\CC[\hat\hen],$$
where $\l$ varies in $X_0\times\ZZ$.
Set $\dot y_{w,\l}\in H^*(\dot\Nc^h_w,\CC)$ equal to
$\exp(-c_1(\Oc_{\dot\Nc^h_w}(\l)))$,
and $\ddot y_{w,\l}$ to the image of
the Poincar\'e dual of $\dot y_{w,\l}$
by the chain of maps
$$H_*(\dot\Nc^h_w,\CC)\to
H_*(\ddot\Nc_1^h\cap\ddot\Nc_{w,w}^h,\CC)\to H_*(\ddot\Nc^h_{w,w},\CC).$$
By \cite{CG, lem.~8.6.1, 8.9.5}
there is a $\ZZ$-graded algebra isomorphism
$$\prod_v\bigoplus_w H_{d_{h,v}+d_{h,w}-*}(\ddot\Nc^h_{v,w},\CC)\to\Extb_h
$$
taking 
$(\ddot y_{w,\l})$
to
$(\pi^h_{w,*}(\dot y_{w,\l})).$
Composing it with the map $ch\circ\Phi$ we get an algebra homomorphism
$$\Psi_h:\hat\Hb_c\to\Extb_h.$$
By definition of $\Phi$ the element $ch\circ\Phi(x_\l)$ 
is equal to the family 
$({}^w\l(\tilde s)^{-1}\ddot y_{w,\l})$.
Therefore, we have
$$\Psi_h(x_\l)=
\Bigl({}^w\l(\tilde s)^{-1}\pi^h_{w,*}(\dot y_{w,\l})\Bigr).$$
The functoriality of $\Ind_{h,w}$
yields an algebra homomorphism
$$\CC[\hat\hen]=\Ext^*_{D^b_{\hat H}(\{0\})}(\CC_{\{0\}},\CC_{\{0\}})\to
\Ext_{D^b_{\tilde G^h\times\CC^\times}(\Nc^h)}^*(\Lc_{h,w},\Lc_{h,w}).$$
Composing it with the forgetful map we get a morphism,
again denoted by $\Ind_{h,w}$,
$$\CC[[\hat\hen]]\to\Ext_{D^b(\Nc^h)}^*(\Lc_{h,w},\Lc_{h,w}),
\quad
\exp(-\xi_\l)\mapsto
\pi^h_{w,*}(\dot y_{w,\l}).$$
So, taking the product over all $w$'s we get the following formula
$$\Psi_h(x_\l)=
\Bigl({}^w\l(\tilde s)^{-1}\Ind_{h,w}(\exp(-\xi_\l))\Bigr).$$

Composing the isomorphisms $\nu_{h,w}$ and the functor $FS$
we get an algebra isomomorphism
$$\aligned
\Extb_{h^\ddag}
&\simeq
\prod_v\bigoplus_w\Ext_{D^b(\Nc^{h^\ddag})}^*(\Lc_{h^\ddag,v},\Lc_{h^\ddag,w})
\cr
&\simeq
\prod_v\bigoplus_w\Ext_{D^b(\Nc^{h^\ddag})}^*(FS(\Lc_{h,v}),FS(\Lc_{h,w}))
\cr
&\simeq
\prod_v\bigoplus_w\Ext_{D^b(\Nc^h)}^*(\Lc_{h,v},\Lc_{h,w})
\cr
&\simeq
\Extb_{h}.
\endaligned
\leqno(A.1.7)$$
A little attention shows that we must prove that the composed maps
$$\hat\Hb_{-c}{\buildrel IM\over\lra}\hat\Hb_c
{\buildrel \Psi_h\over\lra}\Extb_h,
\quad
\hat\Hb_{-c}{\buildrel \Psi_{h^\ddag}\over\lra}\Extb_{h^\ddag}
{\buildrel A.1.7\over\lra}\Extb_h$$
coincide on $\CC_{q,t} X_0$.
Notice that $\Psi_h\,IM(t)=\tau^{c/2}$ and
$\Psi_{h^\ddag}(t)=\tau^{-c/2}$.
So, up to substituting $\tau$ with $\tau^{-1}$ (see the remark after 2.13),
we must check that $FS$ takes
$$\Psi_{h^\ddag}(x_\l)=
\Bigl({}^w\l(\tilde s)^{-1}\Ind_{h^\ddag,w}(\exp(-\xi_\l))\Bigr)$$
to
$$\Psi_h(x_{\l})=\Bigl({}^w\l(\tilde s)^{-1}\Ind_{h,w}(\exp(-\xi_\l))\Bigr).$$
This follows from $A.1.3$, because $FS^2=\Id$.

\subhead A.2. Induction\endsubhead

Fix $h=(s,\tau,\zeta)\in\hat H$ with $\zeta=\tau^{c/2}$.
Recall that $\tau$ is not a root of unity,
that $\tilde s=(s,\tau)$, and that $(\ad\tilde s)=(\ad s)\circ F_\tau$.
Fix a group homomorphism $v:\CC^\times\to\QQ$ such that $v(\tau)<0.$
View $\tilde G$ as the set of $\CC$-points of an
ind-affine group-ind-scheme in the usual way. 
For any 
group-ind-scheme homomorphism $\phi:\SL_2(\CC)\to\tilde G$, 
let $\phi':\slen_2(\CC)\to\tilde\gen$ denote the
differential at 1.
Fix a nilpotent element $e\in\Nc^h$. 

\proclaim{A.2.1. Lemma}
There is a homomorphism $\phi$ as above such that

(a)
$\phi'(\smallmatrix 0&1\cr 0&0\endsmallmatrix)=e,$

(b)
$\phi_\zeta=\phi(\smallmatrix\zeta&0\cr 0&\zeta^{-1}\endsmallmatrix)
\in\tilde T$,

(c)
$(\ad \tilde s)(x)=(\ad\phi_{\zeta})(x)$ for each $x\in\phi'(\slen_2(\CC))$.
\endproclaim

Set $\tilde s_\phi=\tilde s\phi_{\zeta}^{-1}$,
$\tilde\gen_{t,i}=\{x\in\tilde\gen;(\ad\tilde s_\phi)(x)=tx,
\,(\ad\phi_\zeta)(x)=\zeta^ix\}$,
and
$$\qen=\bigoplus_{v(t)\le 0}\tilde\gen_{t,i},
\quad
\nen=\bigoplus_{v(t)<0}\tilde\gen_{t,i},
\quad
\len=\bigoplus_{v(t)=0}\tilde\gen_{t,i}.
\leqno(A.2.2)$$

\proclaim{A.2.3. Lemma}
(a)
The space $\qen$ is a parahoric Lie subalgebra of $\tilde\gen$, $\nen$ is its
pronilpotent radical, and $\len$ is isomorphic to the quotient Lie algebra
$\qen/\nen$.

(b)
There are unique connected closed group-subschemes
$Q,N,L\subset\tilde G$ 
whose Lie algebras are equal to $\qen$, $\nen$, $\len$ respectively, such that
$Q$ is parahoric, $N$ is the prounipotent radical of $Q$,
and $L\simeq Q/N$.

(c)
We have $\tilde s_\phi=(s_\phi,\tau)$ for some $s_\phi\in T_0$.
Further $\qen\neq\tilde\gen$, the group $L$ is connected reductive 
with simply connected derived subgroup, and the fixed points subgroups
$L^h$, $Q^h$ are connected. 
\endproclaim

\vskip3mm

\noindent{\sl Proof :}
Lemmas A.2.1 and  A.2.3(a),(b) are proved in \cite{V1, sec.~6}.
For instance, the group $Q$ is equal to the normalizer of $\qen$ in $\tilde G$.
Let us concentrate on A.2.3(c).
The element $\delta\in X$ may be viewed as a character of $\tilde G$ such that
$\delta(\tilde s)=\tau$.

To prove the first claim we must check that
$\delta(\tilde s_\phi)=\tau$.
This is obvious, because $\phi_\zeta\in\phi(\SL_2(\CC))$,
the group $\phi(\SL_2(\CC))$
is contained in the derived subgroup $\hat G$ of $\tilde G$, and
the character $\delta$ is trivial on $\hat G$.

To prove that $\qen\neq\tilde\gen$, notice that, else, 
we have $\len=\tilde\gen$.
Thus, if $0\neq x\in\tilde\gen$ is such that $(\ad\tilde s_\phi)(x)=tx$ 
then $v(t)=0$.
Now, since $\qen=\tilde\gen$,
we can choose $x$ so that it is a non-zero vector of weight $\delta$.
Then we get
$v(\delta(\tilde s_\phi))=0$. This is absurd, because $v(\tau)<0$.
Since $\qen\neq\tilde\gen$ the group $L$ is reductive.
It is obviously connected.
It is not difficult to prove that its derived subgroup 
is simply connected. Compare \cite{SS, cor.~II.5.4} for the
finite type case.
The other claims are standard by Steinberg's theorem.

\qed

\vskip3mm

For a future use observe that

\vskip1mm

\itemitem{(i)}
the elements $e$ and $f=\phi'(\smallmatrix 0&0\cr 1&0\endsmallmatrix)$
belong to $\len$,

\vskip1mm

\itemitem{(ii)}
the eigenvalues of $(\ad \tilde s)$ acting on
$\zen(e)$, $\zen(f)$, $\qen$, $\nen$ are in 
$\{t\zeta^{i};i\ge 0\}$,
$\{t\zeta^{i};i\le 0\}$,
$\{t\zeta^{i};v(t)\le 0\}$,
$\{t\zeta^{i};v(t)<0\}$ respectively.

\vskip1mm

\noindent
For part (ii) notice that $(\ad \tilde s)$ 
acts on $\tilde\gen_{t,i}$ by multiplication by the scalar $t\zeta^i$,
and that by $\slen_2(\CC)$-representation theory the eigenvalues of
$(\ad\,\phi_\zeta)$ acting on $\zen(e)$ are contained into 
$\{\zeta^{i};i\ge 0\}$.

Let $\Bc_L$ be the flag variety of $\len$
and $\Bc_{L,e}=\{\ben\in\Bc_L;e\in\ben\}$.
Let $I_L$  be the set of elements $i\in I$ such that
$\qen$ contains a parahoric Lie subalgebra of type $i$.
Let $\Hb_L\subset\hat\Hb$ be the subalgebra generated by the elements
$x_\l$, $t_{s_i}$, with $\l\in X$, $i\in I_L$.
It is isomorphic to the affine Hecke algebra of the reductive group $L$.
Thus it acts on the space $H_*(\Bc_{L,e}^h,\CC)$ as in \cite{L2} 
(which is a modified version of the action in \cite{KL1}, \cite{CG}).

Let $W_L\subset W$ be the Weyl group of $L$.
For each right coset $v\in W_L\!\setminus\! W$ let ${}^L\Bc_v\subset\Bc$ 
be the $Q$-orbit containing the set $\{\ben_w;w\in v\}$.
Write ${}^L\Bc$ for ${}^L\Bc_{W_L\!\setminus W_L}$.
Set 
$$\Mc_e=\{(x,\ben)\in\dot\Nc;x\in e+\zen_\nen(f)\}.$$
The same construction as in 2.4.3, using the homomorphism 2.4.4, yields
a $\hat\Hb$-module structure on the vector space $H_*(\Mc_e^h,\CC)$.

\proclaim{A.2.4. Proposition}
There is an isomorphism of $\hat\Hb$-modules
$$\hat\Hb\otimes_{\Hb_L}H_*(\Bc^h_{L,e},\CC)\simeq H_*(\Mc_e^h,\CC).$$
\endproclaim

\noindent{\sl Proof :}
The proof is modelled after \cite{L2, thm.~7.11}.
See also \cite{V1, thm.~5.8} for another version. 
We'll give a direct proof, without using the concentration map as in
\cite{V1}.

Put ${}^L\dot\Nc_v=\dot\Nc\cap(\qen\times{}^L\Bc_v)$,
$\dot\Nc_L=\{(x,\ben)\in\len\times\Bc_L;x\in\ben\}$,
and ${}^L\Mc_{e,v}=\Mc_e\cap{}^L\dot\Nc_v$.
Identify $\len$ with the quotient $\qen/\nen$.
There is the vector bundle (of infinite rank)
$$p_v:{}^L\dot\Nc_v\to\dot\Nc_L,\quad
(x,\ben)\mapsto(x+\nen,(\ben\cap\qen+\nen)/\nen).$$
We'll prove the following  

\vskip1mm

\itemitem{(iii)}
the map $p_v$ restricts to a vector bundle
$r_v:{}^L\Mc_{e,v}^h\to\Bc_{L,e}^h$,

\vskip1mm

\itemitem{(iv)}
we have $H_*(\Mc_e^h,\CC)=\bigoplus_vH_*({}^L\Mc_{e,v}^h,\CC),$ 

\vskip1mm

\itemitem{(v)}
the $\hat\Hb$-action on $H_*(\Mc_e^h,\CC)$
yields a surjective $\CC$-linear map 
$$\hat\Hb\otimes H_*(\Bc^h_{L,e},\CC)\to H_*(\Mc_e^h,\CC),
\quad
x\otimes y\to x\cdot r_{W_L}^*(y),$$

\vskip1mm

\itemitem{(vi)}
the map in (v) factors to an isomorphism
$\hat\Hb\otimes_{\Hb_L}H_*(\Bc^h_{L,e},\CC)\to H_*(\Mc_e^h,\CC)$.

\vskip1mm

Let us concentrate on (iii).
First, we prove that the variety
$${}^L\widetilde\Mc_{e,v}^h=\{(x,\ben)\in{}^L\dot\Nc_v^h;x\in e+\zen_\qen(f)\}$$
is smooth.
Set ${}^L\dot\Vc^h_{v}=
\{(x,\ben)\in\qen^h\times{}^L\Bc^h_v;x\in\ben\}$,
$\dot\Vc^h_{L}=\{(x,\ben)\in\len^h\times\Bc_L^h;x\in\ben\}$.
The map
$$p_v:{}^L\dot\Vc_v^h\to\dot\Vc^h_L,\quad
(x,\ben)\mapsto(x+\nen,(\ben\cap\qen+\nen)/\nen)$$
is a vector bundle.
Let $\cen$ be the canonical Cartan Lie subalgebra.
We have $\cen=\ben/[\ben,\ben]$ for each $\ben\in\Bc$. 
Compare \cite{CG, lem.~3.1.26}.
The map
$$\delta:{}^L\dot\Vc^h_{v}\to\cen^h,
\quad
(x,\ben)\mapsto (x+[\ben,\ben])/[\ben,\ben]$$ 
is smooth, because it is the composition of the smooth map 
$p_v:{}^L\dot\Vc^h_{v}\to\dot\Vc_L^h$
and of the smooth map
$\dot\Vc_L^h\to\cen^h$, $(x,\ben)\mapsto (x+[\ben,\ben])/[\ben,\ben]$.
Set
$${}^L\dot\Vc^h_{e,v}=\{(x,\ben)\in{}^L\dot\Vc_v^h;x\in e+\zen_\qen(f)\}.$$
Thus
$${}^L\widetilde\Mc_{e,v}^h=\{(x,\ben)\in {}^L\dot\Vc^h_{e,v};
\delta(x,\ben)=0\}.$$
So, to prove that the variety
${}^L\widetilde\Mc_{e,v}^h$ is smooth, 
it is enough to check that the restriction of $\delta$ to 
${}^L\dot\Vc^h_{e,v}$
is again smooth.
The $\tilde G^h$-action on $\dot\Nc^h$ restricts to a
$Q^h$-action on ${}^L\dot\Vc^h_{v}$.
The map $\delta$ is constant along the $Q^h$-orbits.
Hence it is sufficient to prove that the map
$$Q^h\times{}^L\dot\Vc^h_{e,v}\to{}^L\dot\Vc_v^h,
\quad
(g,x,\ben)\mapsto(\ad g)(x,\ben)\leqno(A.2.5)$$
is smooth.
The representation theory of $\slen_2(\CC)$ implies that
$\qen=[\qen,e]\oplus\zen_\qen(f).$
Hence
$$\qen^h=T_e((\ad Q^h)e)\oplus\zen_\qen(f)^h.$$
Thus the set $e+\zen_\qen(f)^h$ is a transversal slice 
to the $Q^h$-orbit of $e$
in $\qen^h$ locally near $e$ for the analytic topology.
Further,
%$$(e+\zen_\qen(f)^h)\cap((\ad Q^h)e)=\{e\},$$
due to the
contracting $\CC^\times$-action
$$\CC^\times\times(e+\zen_\qen(f)^h)\to e+\zen_\qen(f)^h,
\quad
(z,x)\mapsto z^{-2}(\ad\phi_z)(x),$$
the action map
$Q^h\times(e+\zen_\qen(f)^h)\to\qen^h$
is smooth.
The smoothness of $A.2.5$ follows.

Now, we have the following commutative diagram
$$\matrix
{}^L\Mc_{e,v}^h&\subset& {}^L\widetilde\Mc_{e,v}^h
&\subset&{}^L\dot\Vc^h_v
\cr
\hskip3mm\downarrow{\ss r_v}
&&\hskip3mm\downarrow{\ss p_v}
&&\hskip3mm\downarrow{\ss p_v}
\cr
\Bc_{L,e}^h&\subset&\dot\Nc_L^h&\subset&\dot\Vc^h_L.
\endmatrix$$
Since the left square is Cartesian,
to prove that the left vertical map is a vector bundle it is enough 
to prove that so is the middle one.
Now, recall the following basic fact, which follows from \cite{BH, thm.~9.1}.

\proclaim{A.2.6. Lemma}
If $r:\Yc\to\Xc$ is a vector bundle over a smooth variety with a fiber 
preserving linear $\CC^\times$-action with positive weights, and 
$\Yc'\subset\Yc$ is a $\CC^\times$-stable smooth closed subvariety, then
$r|_{\Yc'}$ is a subbundle of $\Yc$ restricted to $r(\Yc')$.
\endproclaim

\noindent
Fix a group homomorphism $$\gamma:\CC^\times\to Z(L^h)^\circ$$ such that
$(\ad\gamma)$ has positive weights on $\nen^h$.
The $L^h$-action on $\dot\Vc_L^h$ restricts to the trivial $\gamma$-action.
Apply A.2.6 to $r=p_v$, $\Yc={}^L\dot\Vc^h_v$,
$\Yc'={}^L\widetilde\Mc_{e,v}^h$,
$\Xc=\dot\Vc_L^h$,
and to the $\CC^\times$-action induced by $\gamma$.

Then we prove (iv).
Let $\le$ be the partial order on $W_L\!\setminus\!W$ such that
$v'\le v$ whenever ${}^L\Bc_{v'}$ is contained in the closure of
${}^L\Bc_v$.
Set 
$${}^L\Mc_{e,<v}^h=\coprod_{v'<v}{}^L\Mc_{e,v'}^h,
\quad
{}^L\Mc_{e,\le v}^h=\coprod_{v'\le v}{}^L\Mc_{e,v'}^h.$$
Since both are closed subsets of $\Mc_e^h$, for each $i$
we have an exact sequence
$$H_i({}^L\Mc_{e,<v}^h,\CC)\to H_i({}^L\Mc_{e,\le v}^h,\CC)\to
H_i({}^L\Mc_{e}^h,\CC).
\leqno(A.2.7)$$
Let $\gamma$ be as above.
Let ${}^L\Mc_{e,v}^{h,\g},$ $\Mc_{e,v}^{h,\g}$
be the fixed points subsets.
The vector bundle $r_v$ restricts to an isomorphism 
$${}^L\Mc_{e,v}^{h,\g}\to\Bc_{L,e}^h.$$
Thus, the pieces of the partition 
$\Mc_e^{h,\g}=\coprod_v{}^L\Mc_{e,v}^{h,\g}$ 
are both open and closed, yielding exact sequences
$$0\to H_i({}^L\Mc_{e,<v}^{h,\g},\CC)\to 
H_i({}^L\Mc_{e,\le v}^{h,\g},\CC)\to
H_i({}^L\Mc_{e}^{h,\g},\CC)\to 0.$$
So, a standard argument in equivariant homology
using A.2.7 implies that  
$$H_*({}^L\Mc_{e,\le v}^h,\CC)=\bigoplus_{v'\le v}H_*({}^L\Mc_{e,v'}^h,\CC).$$
We are done, because
$H_*(\Mc_{e}^h,\CC)$ is the direct limit of the subsets
$H_*({}^L\Mc_{e,\le v}^h,\CC)$.

Now, let us prove (v).
Formula 2.4.6 yields
$$t_{s_i}H_*({}^L\Mc^h_{e,v},\CC)\subset 
H_*({}^L\Mc^h_{e,\le vs_i},\CC)$$
for each $i\in I$ such that $v\le vs_i$, because the following holds :
$$(x,\ben)\in{}^L\Mc^h_{e,v},\ (x,\ben,\ben')\in\ddot\Nc_{s_i}^h\Rightarrow
(x,\ben')\in{}^L\Mc^h_{e,\le vs_i}.$$
Compare \cite{L2, sec~4.3, lem.~7.2}. 
So the map in (v) is well-defined, and we must check surjectivity.
Let
$$\sigma_i:H_*(\Bc^h_{L,e},\CC)\to H_*({}^L\Mc^h_{e,v},\CC)\to 
H_*({}^L\Mc^h_{e,vs_i},\CC)$$
be the composition of $r_v^*$,
the action of $-tt_{s_i}-1$, and the restriction to ${}^L\Mc^h_{e,vs_i}$.
By induction it is enough to check that $\sigma_i$ is surjective 
if $v<vs_i$.

Consider the open subset
$\ddot\Uc_{s_i}^h\subset\ddot\Nc_{s_i}^h$
consisting of all $(x,\ben,\ben')$ 
such that $(\ben,\ben')$ is $\tilde G$-conjugate to $(\ben_1,\ben_{s_i})$,
and the open subset $\Uc\subset(\dot\Nc^h)^2$
consisting of all $(x,\ben,x',\ben')$ such that
$(x,\ben)\in{}^L\Mc_{e,\ge v}^h$ and
the relative position of $(\ben,\ben')$ contains the simple reflection 
$s_i$ in any reduced expression. 
Here ${}^L\Mc_{e,\ge v}$ is defined in the obvious way. 
We have
$$
\aligned
\Uc\cap\ddot\Nc_{s_i}^h&=\Uc\cap\ddot\Uc^h_{s_i},
\\
\Uc\cap({}^L\Mc_{e,\le v}^h\times\dot\Nc^h)&=
\Uc\cap({}^L\Mc_{e,v}^h\times\dot\Nc^h),
\\
\Uc\cap\ddot\Nc_{s_i}^h\cap({}^L\Mc_{e,\le v}^h\times\dot\Nc^h)&=
\ddot\Uc_{s_i}^h\cap({}^L\Mc_{e,v}^h\times\dot\Nc^h).
\endaligned
$$ 
The axioms of a Tits system
imply that the map $p_2$ in the commutative diagram
$$\matrix
\ddot\Uc_{s_i}^h\cap({}^L\Mc_{e,v}^h\times\dot\Nc^h) 
&{\buildrel p_1\over\lra}&
{}^L\Mc_{e,v}^h 
\cr
\hskip3mm\downarrow{\ss p_2}&&\hskip3mm\downarrow{\ss r_v}
\cr
{}^L\Mc_{e,vs_i}^h 
&{\buildrel r_{vs_i}\over\lra}&
\Bc_{L,e}^h
\endmatrix$$
is invertible.
Further, there are invertible elements 
$a,b\in \widehat H_*(\dot\Nc^h,\CC)$ such that
the image of 2.4.6 in $\widehat H_*(\ddot\Uc_{s_i},\CC)$ is equal to
the restriction of $a\otimes b$.
Finally, the varieties 
$\Uc\cap\ddot\Uc_{s_i}^h$ and
$\Uc\cap({}^L\Mc^h_{e,v}\times\dot\Nc^h)$ 
are smooth and intersect transversaly along
$\ddot\Uc_{s_i}^h\cap({}^L\Mc^h_{e,v}\times\dot\Nc^h)$.
Therefore, given any invertible element $a_L\in H_*(\Bc^h_L,\CC)$ such that
$a^{-1}r_v^*(x)=r_v^*(a^{-1}_Lx)$ for each $x\in H_*(\Bc_{L,e}^h,\CC)$, we get
$$\sigma_i(a^{-1}_Lx)=p_{2,*}(p_1^*r_v^*(x)p_2^*(b))=
p_{2,*}p_2^*(r_{vs_i}^*(x)b)=r_{vs_i}^*(x)b.$$ 
We are done, because $r_{vs_i}$ is a vector bundle and $b$ is invertible. 

Finally we prove (vi).
The map in (v) factors to a surjective map 
$$\hat\Hb\otimes_{\Hb_L}H_*(\Bc^h_{L,e},\CC)\to H_*(\Mc_e^h,\CC).
\leqno(A.2.8)$$
We must check injectivity.
The right $\CC X$-submodule $\Hb_{\le v}\subset\hat\Hb$ spanned by the $t_w$'s
with $w\in v'$ and $v'<v$ is a free right $\Hb_L$-submodule
of rank equal to the minimal length of an element in $v$.
We claim that $A.2.8$ restricts to a surjective map
$$\Hb_{\le v}\otimes_{\Hb_L}H_*(\Bc^h_{L,e},\CC)\to H_*(\Mc^h_{e,\le v},\CC).
\leqno(A.2.9)$$
If $v=W_L$ this is obvious. 
So this holds for any $v$ by induction on the length, using the surjectivity of
$\sigma_i$.
The map $A.2.9$ is invertible, 
because both sides are finite dimensional vector spaces of the same rank
by (iii), (iv).

\qed

\vskip3mm

\proclaim{A.2.10. Corollary}
Assume that $k<0$. 
We have an isomorphism of $\hat\Hb$-modules
$$H_*(\dot\Nc^{h^\dag},\CC)={}^{IM}H_*(\Bc^{h},\CC).$$ 
\endproclaim

\noindent{\sl Proof :}
Fix $k<0$ and $e=0$. 
Then $\tilde s_\phi=\tilde s$.
Recall that $\gen^{h^\dag}$ is spanned by the $e_\a$'s 
such that $a(s)\tau^\ell=\tau^{-c}$, 
and $\gen^{h}$ by the $e_\a$'s 
such that $a(s)\tau^\ell=\tau^{c}$. 
Hence, since $v(\tau)<0$ and the Lie algebra $\nen$ is spanned by the $e_\a$'s 
such that $v(a(s)\tau^\ell)<0$, we have
$\nen^{h^\dag}=\gen^{h^\dag}$ and $\nen^{h}=\{0\}$.
Therefore 
$\Mc^{h^\dag}_e=\dot\Nc^{h^\dag}$ and 
$\Mc^{h}_e=\Bc^{h}$.
So A.2.4 yields the following isomorphisms 
$$H_*(\dot\Nc^{h^\dag},\CC)\simeq
\hat\Hb_c\otimes_{\Hb_L}H_*(\dot\Nc^{h^\dag}_{L},\CC),
\quad
H_*(\Bc^{h},\CC)\simeq\hat\Hb_{-c}\otimes_{\Hb_{L}}
H_*(\Bc^{h}_{L},\CC).$$
Further, since $\Hb_{L}$ is an affine Hecke algebra,
it is well known that
$$H_*(\dot\Nc^{h^\dag}_L,\CC)\simeq {}^{IM}H_*(\Bc^{h}_L,\CC)$$
as $\Hb_L$-modules, see \cite{L4, sec.~ 7.20}.

\qed

\vskip3mm

\noindent{\bf A.2.11. Remark.}
Let $e=0$ and $h=h_{-c}^{-1}$.
Hence $\tilde s_\phi=(\rhov(\tau^c),\tau)$.
Therefore $L=\tilde G^c$ and
$\qen$ is spanned by the $e_\a$'s 
such that $ca\cdot\rhov+\ell\ge 0$.
See $A.2.2$.
For each $\lambdav\in\check V_{0,\RR}$ the subspace of
$\tilde\gen$ spanned by the $e_\alpha$'s such that 
$a\cdot\lambdav+\ell\ge 0$ is the sum
$\sum_{\lambdav\in\bar\Ac_w}\ben_w.$
Let $w_c$ be as in 3.4.3(a).
We have 
$$\qen=\sum_{c\rhov\in\bar\Ac_w}\ben_w=
\sum_{w\in W_cw_c}\ben_w.$$
Recall that $W_cw_c=w_cW_{\lambdav_c}$ with
$\lambdav_c={}^{w_c^{-1}}(c\rhov)$,
and that $W_{\lambdav_c}$ is a parabolic subgroup of $W$. 
Hence $Q$ is the unique parahoric group of $\tilde G$ of type $W_{\lambdav_c}$
containing $B_{w_c}.$ 
In particular we have $\Hb_L=\Hb_{\lambdav_c}$ and $\Bc^c_{L,e}=\Bc^c_{w_c}$.

\subhead A.3. Standard modules\endsubhead

We'll use the same notations as in section A.2.
Further, for all $\ben\in\Bc$ we write $B$ for $N_{\tilde G}(\ben)$, 
and $\uen$ for the pronilpotent radical of $\ben$.
Set
$$\aligned
X_{L,+}&=\{\l\in X_+;
\l\cdot\alphav_i=0,\,\forall i\in I_L\}
=\bigoplus_{i\notin I_L}\ZZ_{\ge 0}\omega_i\oplus\ZZ\delta,
\\
X_{L,++}&=\{\l\in X_{L,+};
\l\cdot\alphav_i>0,\,\forall i\notin I_L\}
=\bigoplus_{i\notin I_L}\ZZ_{>0}\omega_i\oplus\ZZ\delta.
\endaligned$$

For each $\lambda\in X_+$ let
$E(\l)$ be the integrable simple $\tilde\gen$-module with highest weight $\l$.
Notice that the $\tilde\gen$-action lifts to a $\tilde G$-action.
The space $E(\l)^\nen$ is one-dimensional iff $\l\in X_{L,+}$,
and the normalizer of $E(\l)^\nen$ in $\tilde\gen$ is $\qen$
iff $\l\in X_{L,++}$. 

If $\ben\in\Bc$ then $B$ acts on the line $E(\l)^\uen$.
If $\ben\in\Bc^{\tilde s}$ then $\tilde s\in B$, hence the numbers 
$$r_\l(\tilde s,\ben)=\roman{tr}_{E(\l)^\uen}(\tilde s)\in\CC^\times,
\quad
v_\l(\tilde s,\ben)=v(r_\l(\tilde s,\ben))\in\QQ$$
are well-defined. 
We'll abbreviate $v(\tilde s,\ben)$, $r(\tilde s,\ben)$
when no confusion is possible.
Notice that $r(\tilde s,\ben)$ is also the trace of the action of $\tilde s$
on the fiber of $\Oc_\Bc(\l)$ at $\ben$.

\proclaim{A.3.1. Lemma}
Let $\l\in X_+$ and $k<0$.

(a)
The number $v_\phi=v(\tilde s_\phi,\ben_\phi)$ is independent on the choice of
$\ben_\phi\in{}^L\Bc^{\tilde s_\phi}$.

(b)
If $\ben\in\Bc^{\tilde s}_e$ then
$v(\tilde s,\ben)\ge v_\phi$.

(c)
If $\l\in X_{L,+}$ and 
$\ben\in{}^L\Bc^{\tilde s}_e$
then
$v(\tilde s,\ben)=v_\phi$.

(d)
If $\l\in X_{L,++}$, $\ben\in\Bc^{\tilde s}_e$,
and $v(\tilde s,\ben)=v_\phi$ then $\ben\in{}^L\Bc^{\tilde s}_e$.
\endproclaim

\noindent{\sl Proof :}
Setting 
$E(\l)_t=\Ker(\tilde s_\phi-t\,\Id_{E(\l)})$, we get
$$E(\l)=\bigoplus_{t\in\CC^\times}E(\l)_t, 
\quad
\nen(E(\l)_t)\subset\bigoplus_{v(t')<v(t)}E(\l)_{t'},
\quad
\qen(E(\l)_t)\subset\bigoplus_{v(t')\le v(t)}E(\l)_{t'}.$$
Thus the set $\{v(t);E(\l)_t\neq 0\}$ is bounded by below because 
$E(\l)$ is highest weight.
Let $v'_\phi$ be the minimum, and 
$E(\l)_\phi=\bigoplus_{v(t')=v_\phi}E(\l)_{t'}$.
Observe that
$$
E(\l)^\uen\subset E(\l)_\phi\subset E(\l)^\nen,\quad
\forall \ben\subset\qen.$$
For each $y\in E(\l)$ we set $y=\sum_ty_t$ with $y_t\in E(\l)_t$.

Now we prove (a).
For each $\ben_\phi\in{}^L\Bc^{\tilde s_\phi}$ we have 
$\ben_\phi\subset\qen$, hence
$E(\l)^{\uen_\phi}\subset E(\l)_\phi$.
Thus $v_\phi=v'_\phi$. 

Next we prove (b).
For each $\ben\in\Bc^{\tilde s}_e$, $y\in E(\l)^\uen$ we have 
$$\sum_t\tilde s(y_t)=\tilde s(y)=r(\tilde s,\ben)y
=\sum_tr(\tilde s,\ben)y_t.$$
Since $\tilde s$ commutes with $\tilde s_\phi$, this yields
$\tilde s(y_t)=r(\tilde s,\ben)y_t.$
Hence 
$$r(\tilde s,\ben)y_t=\tilde s_\phi\phi_\zeta(y_t)
=t\phi_\zeta(y_t).$$
%Since $e$, $\phi_\zeta$, $\tilde s$ commute to $\tilde s_\phi$ they act 
%on $E(\l)_t$.
Since $e\in\uen$ we have $e(y)=0.$ 
Since $e$ commutes with $\tilde s_\phi$, this yields $e(y_t)=0.$ 
Thus, if $y_t\neq 0$ the representation theory of $\sen\len_2(\CC)$ implies that
$r(\tilde s,\ben)=t\zeta^i$ for some $i\ge 0$.
We have $v(\zeta)>0$, because $k,v(\tau)<0$.
So $v(\tilde s,\ben)\ge v(t)\ge v_\phi$.
We have proved (b).

Next we prove (c).
Assume that $\l\in X_{L,+}$.
Fix $\ben\in{}^L\Bc^{\tilde s}_e$ and 
$\ben_\phi\in{}^L\Bc^{\tilde s_\phi}$.
We have $E(\l)^\uen\subset E(\l)^\nen$, because $\ben\subset\qen$.
So $E(\l)^\uen=E(\l)^\nen$, because $E(\l)^\nen, E(\l)^\uen$ are both
one-dimensional.
Hence $r(\tilde s,\ben)=\l(\tilde s)$.
For the same reason we have also
$r(\tilde s_\phi,\ben_\phi)=\l(\tilde s_\phi).$ 
Since $\phi(\SL_2(\CC))$ is contained into the derived subgroup of $Q$ 
and $\l$ is the restriction of the character 
of $Q$ acting on the line $E(\l)^\nen$,
we have
$\l(\phi_\zeta)=1$.
So $\l(\tilde s)=\l(\tilde s_\phi).$ 
Therefore $v(\tilde s,\ben)=v(\tilde s_\phi,\ben_\phi)=v_\phi.$ 

Finally we prove (d).
Assume that $\l\in X_{L,++}$
and that $\ben\in\Bc^{\tilde s}_e$ is such that $v(\tilde s,\ben)=v_\phi$.
The proof of (b) implies that for each $y\in E(\l)^\uen$ and each $t$ such that
$y_t\neq 0$ we have 
$v(\tilde s,\ben)\ge v(t)\ge v_\phi$.
Since $v(\tilde s,\ben)=v_\phi$ this yields
$v(t)=v_\phi$, hence $y\in E(\l)_\phi$.
Therefore, we have $E(\l)^\uen\subset E(\l)_\phi$.
Since $E(\l)_\phi\subset E(\l)^\nen$
and $\dim E(\l)^\uen=\dim E(\l)^\nen=1$,
we have $E(\l)^\uen=E(\l)^\nen$.
Therefore $\ben\subset\qen$,
because $\l\in X_{L,++}$. 

\qed

\vskip3mm

\proclaim{A.3.2. Corollary}
If $k<0$ the subset ${}^L\Bc_{e}^{\tilde s}\subset\Bc^{\tilde s}_e$ is 
$Z_{\tilde G^h}(e)$-stable, open, and closed.
\endproclaim

\noindent{\sl Proof :}
The Lie algebra of $Z_{\tilde G^h}(e)$ is contained into
$\bigoplus_{i\ge 0}\tilde\gen_{\zeta^{-i},i}$ by A.2(ii).
Since $k<0$ we have $v(\zeta)>0$.
Thus the Lie algebra of $Z_{\tilde G^h}(e)$ is also contained into $\qen$.
Thus ${}^L\Bc_{e}^h$ is $Z_{\tilde G^h}(e)$-stable.
Choose an element $\l\in X_{L,++}$.
Since the function
$\ben\mapsto v_\l(\tilde s,\ben)$
is constant on the connected components of $\Bc^{\tilde s}_e$,
each component is either contained in
${}^L\Bc^{\tilde s}_e$ or disjoint from it by A.3.1. 

\qed

\vskip3mm

For each representation $\chi$ of the finite group $A(h,e)$
which is a Jordan-H\"older factor of $H_*(\Bc^{\tilde s}_e,\CC)$ we set
$$N_{h,\chi}=\Hom_{A(h,e)}(\chi,H_*(\Bc^{\tilde s}_e,\CC)).$$
The algebra $\hat\Hb$ acts on the space $N_{h,\chi}.$

In the rest of this section we'll assume that $k<0$.
Then, the group $Z_{\tilde G^h}(e)$ acts on ${}^L\Bc^{\tilde s}_e$ by A.3.2. 
Hence the finite group
$A(h,e)$ acts on $H_*({}^L\Bc^{\tilde s}_{e},\CC)$ and we may set
$${}^LN_{h,\chi}=\Hom_{A(h,e)}(\chi,H_*({}^L\Bc^{\tilde s}_{e},\CC)).$$
Recall that each irreducible 
$\tilde G^h$-equivariant local system $\chi\in\Xc_h$
is supported on a nilpotent orbit $O\subset\Nc^h$,
and, taking the stalk at an element $e\in O$,
yields an irreducible 
representation of the finite group $A(h,e)$
which is a Jordan-H\"older factor of $H_*(\Bc^{\tilde s}_e,\CC)$.
We'll call $N_{h,\chi}$ a standard module (if $k<0$).

\proclaim{A.3.3. Proposition}
If $k<0$ then ${}^LN_{h,\chi}$ is a simple $\Hb_L$-submodule of
$N_{h,\chi}$. Further there is an isomorphism of $\hat\Hb$-modules
$$N_{h,\chi}\simeq\hat\Hb\otimes_{\Hb_L}{}^LN_{h,\chi}.$$
\endproclaim

\noindent{\sl Proof :}
Consider the standard $\Hb_L$-module
$$N_{L,h,\chi}=\Hom_{A(h,e)}(\chi,H_*(\Bc^{\tilde s}_{L,e},\CC)).$$
By A.3.2 we have ${}^LN_{h,\chi}\subset N_{h,\chi}$.
Further ${}^LN_{h,\chi}$ is preserved by the action of the subalgebra $\Hb_L$
of $\hat\Hb$,
and the map $\Bc_{L,e}\to{}^L\Bc_e$, $\ben\to\ben\oplus\nen$
yields an isomorphism of $\Hb_L$-modules
$N_{L,h,\chi}\to {}^LN_{h,\chi}.$
Since $\len^h=\tilde\gen_{1,2}$ and $\Lie(L^h)=\tilde\gen_{1,0}$
the $(\ad L^h)$-orbit of $e$ in $\len^h$ is open, by \cite{K4, thm.~4.2}.
Thus $N_{L,h,\chi}$ is irreducible.
This implies the first claim. 

The second one follows from A.2.4, because 
$\zen_\nen(f)^h=\{0\}$ by A.2(ii) above
and, thus, $\Mc^h_e=\Bc_e^{\tilde s}$.

\qed

\vskip3mm

\proclaim{A.3.4. Proposition}
If $k<0$
the standard module $N_{h,\chi}$ has a simple top isomorphic to $L_{h,\chi}$.
\endproclaim

\noindent{\sl Proof :}
The quotient map $N_{h,\chi}\to L_{h,\chi}$ is as in 2.5.2.
The standard module is the direct sum of the subspaces
$$N_{h,\chi,j}=\Hom_{A(h,e)}(\chi,H_*(\Bc^{\tilde s}_{e,j},\CC)),$$
where the $\Bc^{\tilde s}_{e,j}$'s are 
the $Z_{\tilde G^h}(e)$-saturated sets generated
by the connected components of $\Bc^{\tilde s}_e$.
By A.3.2 we have
$$N_{h,\chi}={}^LN_{h,\chi}\oplus {}^{\not L}N_{h,\chi},$$
where 
${}^{\not L}N_{h,\chi}=\sum_jN_{h,\chi,j}$ for all $j$ such that
$\Bc^{\tilde s}_{e,j}\cap{}^L\Bc=\emptyset$.

Let $M\subset N_{h,\chi}$ be a proper $\hat\Hb$-submodule.
The spaces ${}^LN_{h,\chi}$, 
${}^{\not L}N_{h,\chi}$ are sums of weight spaces
of $N_{h,\chi}$, and any weight space of $N_{h,\chi}$ is contained either in
${}^LN_{h,\chi}$ or in ${}^{\not L}N_{h,\chi}$.
Thus any weight space of $M$ is contained either in
${}^LN_{h,\chi}$ or in ${}^{\not L}N_{h,\chi}$,
yielding a decomposition
$$M={}^LM\oplus {}^{\not L}M.$$
Since $M\neq N_{h,\chi}$ we have also 
${}^LM\neq {}^LN_{h,\chi}$ by A.3.3.
Thus ${}^LM=0$ by A.3.3 again.
Hence the sum of all proper $\hat\Hb$-submodules of
$N_{h,\chi}$ is contained into ${}^{\not L}M$,
so it is proper.
We are done.

\qed

\vskip3mm

\proclaim{A.3.5. Corollary}
If $k<0$ any module in $\Irr(O(\hat\Hb_c))$ is the top of an induced module
$\hat\Hb\otimes_{\Hb_L}{}^LN_{h,\chi}$ for some subset $I_L\subsetneq I$.
\endproclaim

\noindent{\sl Proof :}
Follows from A.2.3(c) and A.3.3, A.3.4.

\qed

\vskip3mm

Therefore, using the involution $IM$ we get the following.

\proclaim{A.3.6. Corollary}
Any module in $\Irr(O(\hat\Hb_c))$ is the top of a module induced from $\Hb_L$
for some subset $I_L\subsetneq I$.
\endproclaim

\subhead A.4. Weights of simple modules\endsubhead

We'll use the same notations as in section A.2.

\proclaim{A.4.1. Proposition}
If $M,M'\in O(\hat\Hb_c)$ are simple modules with the same weights,
then $M\simeq M'$.
\endproclaim

\noindent{\sl Proof :}
Since $M,M'$ are simple with the same weights
there is an element $h=(s,\tau,\zeta)\in \hat H$
such that $M,M'\in O_h(\hat\Hb)$.
Up to twisting both modules by $IM$ we can assume that $k<0$.
So there are $\chi,\chi'\in\Xc_h$ such that
$M=L_{h,\chi}$ and $M'=L_{h,\chi'}$.
Fix nilpotent elements $e,e'\in\Nc^h$ such that
$\chi$, $\chi'$ are supported by the $\tilde G^h$-orbits of $e,e'$ respectively.
Using respectively $e$ and $e'$,
fix group homomorphisms
$\phi,\phi'$ and Levi subgroups $L,L'\subset\tilde G$ as in $A.2.2$.
Then $M,M'$ are the top of the $\hat\Hb$-modules
$N_{h,\chi}$, $N_{h,\chi'}$ respectively by A.3.4, we have 
$$N_{h,\chi}=\hat\Hb\otimes_{\Hb_L}{}^LN_{h,\chi},
\quad
N_{h,\chi'}=\hat\Hb\otimes_{\Hb_{L'}}{}^{L'}N_{h,\chi'},$$
and ${}^LN_{h,\chi}$, ${}^{L'}N_{h,\chi'}$ are simple
modules over $\Hb_L$, $\Hb_{L'}$ by A.3.3.
Further the weights of
${}^LN_{h,\chi}$, ${}^{L'}N_{h,\chi'}$ 
are weights of $M,M'$ respectively by the proof of A.3.4.

\proclaim{A.4.2. Lemma}
(a)
We have $\Hb_L=\Hb_{L'}$.

(b)
The $\Hb_L$-modules
${}^LN_{h,\chi}$, ${}^{L'}N_{h,\chi'}$ have the same weights. 
\endproclaim

\noindent
Therefore, since $L$ is a proper subgroup of $\tilde G$ by A.2.3(c), we have 
${}^LN_{h,\chi}={}^{L'}N_{h,\chi'}$ by \cite{EM, thm.~5.5} and A.4.2(a). 
Thus $M\simeq M'$.

\qed

\vskip3mm

\noindent{\sl Proof of A.4.2 :}
For each $\l\in X_+$ the number $v_\chi(\l)=v_\l(\tilde s_\phi,\ben_\phi)$
is independent of the choice of the element 
$\ben_\phi\in{}^L\Bc^{\tilde s_\phi}$ by A.3.1(a) and is additive
with respect to $\l$.
Consider the linear forms 
$v_\chi,v_{\tilde s}\in V_\RR^*$ such that 
$$v_\chi(\l)=v_\chi(\l')-v_\chi(\l''),
\quad
v_{\tilde s}(\l)=v(\l(\tilde s)),
\quad
\forall \l\in X,
$$
with $\l=\l'-\l''$ and $\l',\l''\in X_+$.
Let $h_1=(\tilde s_1,\zeta)$ be a weight of $N_{h,\chi}$.
We have

\vskip1mm

\itemitem{(i)}
$v_\chi(-\a_i)=0$ if $i\in I_L$, and $v_\chi(-\a_i)>0$ else,

\vskip1mm

\itemitem{(ii)}
if $\l\in X_+$ 
then $v_{\tilde s_1}(\l)-v_\chi(\l)\ge 0$, 

\vskip1mm

\itemitem{(iii)}
if $\l\in X_{L,+}$ and 
$h_1$ is a weight of ${}^LN_{h,\chi}$ 
then
$v_{\tilde s_1}(\l)=v_\chi(\l)$,
and if $\l\in X_{L,++}$ and 
$v_{\tilde s_1}(\l)=v_\chi(\l)$
then
$h_1$ is a weight of ${}^LN_{h,\chi}$ .

\vskip1mm

Indeed, let $\ben_{\phi,i}$ be the parahoric Lie subalgebra of type $i$
containing $\ben_\phi$. Then
$$v_\chi(-\a_i)=v(\tr_{\ben_{\phi,i}/\ben_\phi}(\tilde s_\phi)).$$
If $i\in I_L$ then  
$\ben_{\phi,i}\subset\qen$, because $\ben_\phi\subset\qen$.
Thus
$\ben_{\phi,i}/\ben_\phi\subset(\ben_{\phi,i}\cap\len)/(\ben_\phi\cap\len)$.
So $v_\chi(-\a_i)=0$ by definition of $\len$.
If $i\notin I_L$ then $\ben_{\phi,i}\not\subset\qen$. 
Thus $v_\chi(-\a_i)>0$ by definition of $\qen$.
This proves (i).

Let $N_{h,\chi,j}$, $\Bc^{\tilde s}_{e,j}$ be as in the proof of $A.3.4$.
For any $\ben\in\Bc^{\tilde s}_{e,j}$ the element $x_\l^{-1}\in\hat\Hb$ acts on 
$N_{h,\chi,j}$ as scalar multiplication by 
$r_\l(\tilde s,\ben)$ times a unipotent operator.
Thus, since $h_1$ is a weight of the module $N_{h,\chi}$
there is a Lie algebra $\ben\in\Bc^{\tilde s}_e$ such that
$v_{\tilde s_1}(\l)=v_\l(\tilde s,\ben)$.
So (ii) follows from A.3.1(b).

Further $h_1$ is a weight of 
${}^LN_{h,\chi}$ iff it is a weight of  
$H_*({}^L\Bc_e^{\tilde s},\CC)$, i.e.,
iff there is a Lie algebra $\ben\in{}^L\Bc^{\tilde s}_e$ such that
$v_{\tilde s_1}(\l)=v_\l(\tilde s,\ben)$.
So the first part of (iii) follows from A.3.1(c), and the second one from
A.3.1(d).

Now, we define a function $v_{\chi'}$ in the same way as $v_\chi$, using
$\phi'$ instead of $\phi$.

\proclaim{A.4.3. Lemma}
We have $v_\chi=v_{\chi'}.$
\endproclaim

We have $I_L=I_{L'}$ by (i) and A.4.3.
Part (a) is proved.
We have the following inclusions of sets of weights 
$$wt({}^LN_{h,\chi})\subset wt(M)=wt(M')\supset wt ({}^LN_{h,\chi'}),
\quad
wt(M)\subset wt(N_{h,\chi})\cap wt (N_{h,\chi'}).$$
Fix $\l\in X_{L,++}$
and assume that $h_1$ is a weight of $M$.
Then (ii), (iii) imply that $h_1$ is a weight of
${}^LN_{h,\chi}$ iff we have $v_{\tilde s_1}(\l)=v_\chi(\l)$,
and that it is a weight of
${}^LN_{h,\chi'}$ iff we have $v_{\tilde s_1}(\l)=v_{\chi'}(\l)$.
So (b) follows also from A.4.3. 

\qed

\vskip3mm

To prove A.4.3 we need two geometric lemmas.
Let
%$\Hc^-=\{v\in\check V_{\RR};\delta(v)<0\}$
$\Hc=\{v\in\check V_{\RR};\delta(v)=v(\tau)\}.$
Recall that $v(\tau)<0$.
For each subset $J\subset I$ set
$$\Hc_{J,+}=\bigoplus_{i\in J}\RR_{\ge 0}\alphav_i,
\quad
\Hc_{J,-}=\Bigl(\bigoplus_{i\notin J}\RR_{<0}\omegav_i\oplus\RR\check\delta
\Bigr)\cap\Hc,
\quad
\Hc_J=\Hc_{J,-}+\Hc_{J,+}.$$
%Let $\Hc_{J,-}$, $\Hc_J$ 
%be the traces of 
%$\Hc_{J,-}^-$, $\Hc_J^-$ 
%on $\Hc$.
Notice that $\Hc_J\subset\Hc.$
Let $p:\check V_\RR\to\bigoplus_{i\in I}\RR\omegav_i$ 
be the projection along $\RR\check\delta$.

\proclaim{A.4.4. Lemma}
(a)
We have
$\Hc=\coprod_{J\subsetneq I}\Hc_J$.
%$\Hc^-=\coprod_{J\subsetneq I}\Hc^-_J$.

(b)
For each $v\in\Hc_J$ with $J\subsetneq I$ 
the decomposition $v=v_-+v_+$ with
$v_-\in\Hc_{J,-}$, $v_+\in\Hc_{J,+}$
is unique.

(c)
If $v\in\Hc$ then $v_-\in\Hc$.
The map $v\mapsto v_-$ 
factors to the projection $\pi:\check V_{0,\RR}\to\bar\Ac$
relatively to the scalar product $(\ :\ )$.
\endproclaim

\noindent{\sl Proof :}
To prove claim (a) it is enough to check that
$p(\Hc)=\coprod_{J\subsetneq I}p(\Hc_J),$ 
because $\Hc_{J,-}+\RR\check\delta=\Hc_{J,-}$ for all $J$.
A little attention shows that this follows from
the following geometric lemma (left to the reader).
Fix generators $A_i$, $i\in I$, of the vector space $\RR^{I_0}$ 
such that $\sum_iA_i=0$.
Fix linearly independent vectors
$B_i\in(\sum_{j\neq i}\RR_{\ge 0}A_j)\times\{1\}$ with $i\in I$.
The convex polyhedral cones 
$$\Bigl(\bigoplus_{i\in I\setminus J}\RR_{>0}B_i\Bigr)
+\Bigl(\bigoplus_{i\in J}\RR_{\ge 0}A_i\Bigr)$$
with $J\subsetneq I$ form a partition of the set
$\RR^{I_0}\times\RR_{>0}$.

To prove claim (b) it is enough to check that
for all $J\subsetneq I$ the vector subspaces 
$\bigoplus_{i\in J}\RR\alphav_i,$
$\bigoplus_{i\notin J}\RR\omegav_i\oplus\RR\check\delta$
intersect trivially.
Since these subspaces are orthogonal to each other
for the bilinear form $(\ :\ )$, 
any element $v$ in the intersection
is a sum
$x+y\omegav_0+z\check\delta$ with $x\in\check V_{0,\RR}$, $y,z\in\RR$, and 
$(x:x)+2yz=0.$
Since 
$v\in\bigoplus_{i\in J}\RR\alphav_i,$
we have $y=0$.
Thus $x=0$ because $(\ :\ )$ is positive definite on $\check V_{0,\RR}$.
Thus $z=0$ because $J\neq I$.

Now we prove (c).
Observe that
$p(\Hc)=\check V_{0,\RR}+v(\tau)\omegav_0.$
Thus the map
$v+v(\tau)\omegav_0\mapsto v(\tau)^{-1}v$
yields an isomorphism $\phi:p(\Hc)\to\check V_{0,\RR}$.
Set $a_0=1-\theta$ (an affine function on $\check V_{0,\RR}$),
and
$$\bar \Ac=\{x\in\check V_{0,\RR};a_i(x)\ge 0,\,\forall i\in I\},
\quad
\Ac_J=\{x\in\bar \Ac;\,a_i(x)=0\iff i\in J\}.$$
For each $J\subsetneq I$ we have 
$p(\Hc_{J,-})=v(\tau)(\Ac_J+\omegav_0).$
So 
$$\phi p(\Hc_{J,-})=\Ac_J,\quad
\phi p(\coprod_{J\subsetneq I}\Hc_{J,-})=\bar\Ac.$$
Notice that $\Ac_\emptyset$ is the alcove $\Ac$.
We have $(v+y\check\delta)_-=v_-+y\check\delta$ for each
$v\in\Hc$, $y\in\RR$,
because $\Hc_{J,-}+\RR\check\delta=\Hc_{J,-}$.
Thus the map $\Hc\to\Hc$, $v\mapsto v_-$ 
factors to a map $\check V_{0,\RR}\to\bar \Ac$.

The projection of $x\in\check V_{0,\RR}$
to the convex set $\bar \Ac$ is the unique element $\pi(x)\in\bar \Ac$
such that one of the following identities is satisfied

\vskip1mm

\itemitem{(iv)}
$\|x-\pi(x)\|=\min\{\|x-x'\|;x'\in\bar \Ac\}$,

\vskip1mm

\itemitem{(v)}
$(x-\pi(x):x'-\pi(x))\le 0$, $\forall x'\in\bar \Ac$.

\vskip1mm

Given $J\subsetneq I$ and $v\in\Hc_J$,
we must prove that
$\phi p(v_-)=\pi\phi p(v)$.
Since $\phi p(v)-\phi p(v_-)\in -p(\Hc_{J,+})$,
we must check that 
$$(\alphav_i:x'-\phi p(v_-))\ge 0,\quad \forall i\in J,\,
\forall x'\in\bar \Ac.$$
Since $\phi p(v_-)\in\Ac_J,$ we have
$$x'-\phi p(v_-)\in\{x\in\check V_{0,\RR};
\,a_i(x)\ge 0,\,\forall i\in J\}.$$
We are done.

\qed

\vskip3mm

\proclaim{A.4.5. Lemma}
If $v'\in\Hc$ and 
$v-v'\in\Hc_{I,+}$ then 
$v_--v'_-\in\Hc_{I,+}$.
\endproclaim

\noindent{\sl Proof :}
An easy induction on the cardinal of the minimal set $J$ such that
$v-v'\in\Hc_{J,+}$ shows that it is enough to assume  that
$$v-v'\in\Hc_{\{i_0\},+}$$ with $i_0\in I$.
Since $v'\in\Hc$ and $v-v'\in\Hc_{I,+}$ we have $v\in\Hc.$ 
Fix $J,J'$ such that 
$v\in\Hc_{J}$, $v'\in\Hc_{J'}$.
We'll use the same notations as in the proof of A.4.4(c).
Write also $\av_0$ for $-\check\theta$,
and $\Hc_{J,\pm}$ for $(-\Hc_{J,+})\cup\Hc_{J,+}.$

%Put $\Vc_J=\phi p(\Hc_J)$.
%We have $\Vc_J=\pi^{-1}(\Ac_J)$ by A.4.4(c).
%Set $\Vc^i=\bigcup_{i\notin J''}\Vc_{J''}$.
%Notice that $i\in J''$ iff $\Ac_{J''}\subset\bar \Ac_{\{i\}}$,
%and that $\Vc_J=\pi^{-1}(\Ac_J)$ by A.4.4(c).
%Thus the set $\Vc^i=\bigcup_{i\notin J''}\Vc_{J''}$ is equal to
%$\pi^{-1}(\bar \Ac\setminus\bar \Ac_{\{i\}})$.
Given an element $i\in I$, consider
the convex polyhedral cone
$$\Cc_i=\{x\in\check V_{0,\RR};a_j(x)\ge 0,\,\forall j\neq i\}$$ 
whose origin is the vertex $x_i\in\bar \Ac$ 
opposit to the face $\bar \Ac_{\{i\}}$.
Put $\Vc_{J''}=\phi p(\Hc_{J''})$ and $\Vc^i=\bigcup_{i\notin J''}\Vc_{J''}$.
If $x\in\Vc_{J''}$ then 
$x-\pi(x)\in -p(\Hc_{J'',+})$
and $\pi(x)\in\Ac_{J''}$
by A.4.4(c).
Hence if $i\notin J''$ then 
$$(x-\pi(x):x'-\pi(x))\le 0,\quad\forall x'\in\Cc_i,$$
because
$(\check\alpha_j:x'-\pi(x))\ge 0$ for all $j\in J''.$
Thus, since $\bar\Ac\subset\Cc_i$,
the restriction of the map $\pi$ to $\Vc^i$ is
the projection to $\Cc_i$.

%In the rest of the proof we will assume that $i_0\neq 0$ 
%(the case $i_0=0$ is similar and is left to the reader).

If $i\notin J'$ for some $i\neq i_0$ then we have also $i\notin J$,
because 
%$\phi p(v)\in V_{J'}+\RR_{\le 0}\check a_{i_0}$.
$$\phi p(v)-\phi p(v')\in\RR_{\le 0}\av_{i_0},
\quad
\Vc^i+\RR_{\le 0}\av_{i_0}\subset\Vc^i,
\quad
\phi p(v')\in\Vc^i.$$
Therefore the projections of $\phi p(v)$, $\phi p(v')$ to 
$\bar \Ac$ coincide with their projection to $\Cc_i$.
Hence
$$\phi p(v_-)-\phi p(v'_-)\in\sum_{j\neq i}\RR_{\le 0}\av_j$$
by \cite{C1, cor.~1.12}.
Thus we have
$$v_--v'_-\in
\RR\check\delta+\Hc_{I\setminus\{i\},+}.$$
Since we have also
$v-v'\in\Hc_{\{i_0\},+}$
and
$v_+,v'_+\in\Hc_{I\setminus\{i\},+},$
this yields
$$v_--v'_-\in
(\RR\check\delta+\Hc_{I\setminus\{i\},+})\cap
\Hc_{I\setminus\{i\},\pm}=
\Hc_{I\setminus\{i\},+}.$$

Now, assume that $J'=I\setminus\{i_0\}$.
So $\phi p(v'_-)=x_{i_0}$, because
$\phi p(v'_-)\in\Ac_{J'}$ and 
$\Ac_{J'}=\{x_{i_0}\}$.
The cones $\sum_{j\in J''}\RR_{<0}\av_j$, with $J''\subsetneq I$,
form a partition of $\check V_{0,\RR}$.
Since $\phi p(v_-)\in\bar \Ac$ the element,
$\phi p(v_-)-x_{i_0}$ belongs to the set
%\{x\in\check V_{0,\RR};a_i(x)\ge 0,\forall i\neq i_0\}\subset
$\check V_{0,\RR}\setminus\sum_{j\neq i_0}\RR_{<0}\av_j.$
Thus we can fix an element 
$i\neq i_0$ such that 
$$\phi p(v_-)-x_{i_0}\in\sum_{j\neq i}\RR_{\le 0}\av_j.
\leqno(A.4.6)$$
From $A.4.6$ we get
$$v_--v'_-\in\RR\check\delta+\Hc_{I\setminus\{i\},+}.$$
Since
$v-v'\in\Hc_{\{i_0\},+}$,
$v_+\in\Hc_{J,+}$,
and $v'_+\in\Hc_{J',+}$,
this yields 
$$v_--v'_-\in(\RR\check\delta+\Hc_{I\setminus\{i\},+})\cap
(\Hc_{I,+}-\Hc_{J,+}).$$
If $\phi p(v_-)\notin\bar\Ac_{\{i\}}$ then $i\notin J$,
because else we would have
$\Ac_J\subset\bar\Ac_{\{i\}}$.
So in this case we get
$$v_--v'_-\in(\RR\check\delta+\Hc_{I\setminus\{i\},+})\cap
(\Hc_{\{i\},+}+\Hc_{I\setminus\{i\},\pm})=\Hc_{I,+}.$$
If $\phi p(v_-)\in\bar\Ac_{\{i\}}$  then
$\phi p(v_-), x_{i_0}$ both belong to $\bar\Ac_{\{i\}}.$
Thus $\phi p(v_-)=x_{i_0}$ by $A.4.6$, because
$a_i(\av_j)<0$ for each $j\neq i$.
Therefore we get
$v_--v'_-\in\RR\check\delta$,
$v-v'\in\Hc_{\{i_0\},+}$,
and $v_+,v'_+\in\Hc_{J',+}$,
yielding
$$v_--v'_-
\in\RR\check\delta\cap(\Hc_{\{i_0\},+}+\Hc_{J',\pm})\subset\Hc_{I,+}.$$

\qed

\vskip3mm

\noindent{\sl Proof of A.4.3 :}
Let $h'$ be a weight of $M$.
From (i), (ii), (iii) above we get

\vskip1mm

\itemitem{(vi)}
$v_\chi(\a_i)=0$ if $i\in I_L$
and $v_\chi(\a_i)<0$ else,

\vskip1mm

\itemitem{(vii)}
$v_{\chi'}(\a_i)=0$ if $i\in I_{L'}$
and $v_{\chi'}(\a_i)<0$ else,

\vskip1mm

\itemitem{(viii)}
if $\l\in X_+$ then
$v_{\tilde s'}(\l)\ge v_\chi(\l), v_{\chi'}(\l)$,

\vskip1mm

\itemitem{(ix)}
if $\l\in X_{L,+}$ and $h'$ is a weight of ${}^LN_{h,\chi}$ then
$v_{\tilde s'}(\l)=v_\chi(\l)$,

\vskip1mm

\itemitem{(x)}
if $\l\in X_{L',+}$ and $h'$ is a weight of ${}^{L'}N_{h,\chi'}$ then 
$v_{\tilde s'}(\l)=v_{\chi'}(\l)$.

\vskip1mm

Now, assume that $h_1$, $h'_1$ are weights of $M$ (i.e., of $M'$)
such that 
$h_1$ is also a weight of ${}^LN_{h,\chi}$ and
$h'_1$ is also a weight of ${}^{L'}N_{h',\chi'}$. 
Identify $V_\RR^*$ with $\check V_\RR$ via the pairing 1.0.3. 
We have
$$v_{\tilde s_1}(\delta)=v_{\tilde s'_1}(\delta)=v(\tau).$$
Since $\pm\delta\in X_+$, from (viii) we get also
$$v_{\chi}(\delta)=v_{\chi'}(\delta)=v(\tau).$$
So $v_\chi,$ $v_{\tilde s_1}$, $v_{\chi'},$ $v_{\tilde s'_1}$ 
may be viewed as elements of $\Hc$.
We have 
$v_{\chi}\in\Hc_{I_L,-}$ by (vi), 
and
$v_{\tilde s_1}-v_\chi\in\Hc_{I_L,+}$ 
by (viii), (ix).
Thus 
$v_\chi=(v_{\tilde s_1})_-$.
In the same way we prove that
$v_{\chi'}=(v_{\tilde s'_1})_-$.
We have also
$(v_{\chi'})_-=v_{\chi'}$,
$(v_{\chi})_-=v_{\chi}$,
and
$v_{\tilde s'_1}-v_\chi,
v_{\tilde s_1}-v_{\chi'}
\in\Hc_{I,+}$
by (viii). 
Thus
$v_{\chi'}-v_\chi,
v_\chi-v_{\chi'}
\in\Hc_{I,+}$
by A.4.5.
Therefore $v_\chi=v_{\chi'}$. 

\qed

\vskip3mm

\noindent{\bf A.4.7. Remark.}
Lemmas A.4.4, A.4.5 are affine analogues of two well-known lemmas of
Langlands. See \cite{BW, lem.~6.11, 6.13} and
\cite{KL1, sec.~2.10, 7.11}. 
The latter can be interpreted in terms of orthogonal
projections of weights on the Weyl chamber as in \cite{C1, sec.~1}.
We do not know if our lemmas have a similar interpretation
(in the affine case the invariant form is not positive).

\vskip3mm

\subhead A.5. Index of notation\endsubhead

\itemitem{0 : }
$c$,
$k$,
$m$,
$K$,
$\bar K$,
$A$,
$\CC_{q,t}$,
$\CC_t$,
$\CC_\kappa$,

\itemitem{1.0 : }
$w_0$,
$N$,
$\Pi$,

\itemitem{1.1 : }
EN, RN, 
$W_0[m]$,
$I_m$,

\itemitem{1.2 : }
$\tau_\ell$,
$\tau$,
$\lambdav_\ell$,
$F$,
$F_\tau$,
$F_w$,

\itemitem{1.3 : }
TN, TNS, RS, HRS, ERS, HERS,
$\sen$,
$e_R$,
$f_R$,
$B_0$,
$B_0^-$,
$U_0$,
$U_0^-$,
$\ben_0$,
$\ben_0^-$,
$\uen_0$,
$\uen_0^-$,

\itemitem{2.1 : }
$H$,
$\hat H$,
$\hat T_0$,
$\tilde T$,
$T_d$,
$h_c$,
$h_c^\dag$,
$h^\dag$,
$h^\ddag$,
$h=(s,\tau,\zeta)$,
$\tilde s=(s,\tau)$,
$\Hb$,
$\hat\Hb$,
$\Pb$,
$\Pb^\dag$,
$\Pb_X$,
$\Pb_X^\dag$,
$x_\l$,
$t_w$,
$IM$,
$M_h$,
${}^F\CC_{q,t}$,
${}^wh$,
$O(\Hb)$,
$O_h(\Hb)$,
$\Delta_{0,(h)}$,
$\Delta_{(h)}$,
$\Delta_{\re,(h)}$,
$\Hb_{(h)}$,
$W_{(h)}$,
$I^h$,
$J^h$,
$I^E$,
$\Delta_h$,
$W_h$,
$\vartheta_\a$,
$\psi_{s_\a}$,

\itemitem{2.2 : }
$\xi_\l$,
$\CC[\xi_\l]$,
$M_\lambdav$,
$\Hb'$,
$\mod(\Hb'_c)_{uni}$,
$IM$,
$CF$,
$OP$,
$\psi'_{s_i}$,
$\Sb_c$,
$\bar\Sb_c$,
$\Sb'_c$,
$\bar\Sb'_c$,
$\HHc$,
$\HHc'$,
$\rhov_\kappa$,

\itemitem{2.3 : }
$\Hb''$,
$\mod(\Hb''_c)_{nil}$,
$\Lb'_c$,
$\Lb_c$,
$\Lb_c^\dag$,

\itemitem{2.4 : }
$\ben_i$,
$\ben_w$,
$\ben_{w,i}$,
$\CC^\times_{\omega_0}$,
$\CC_\delta^\times$,
$\hat G$,
$G_d$,
$\tilde G$,
$\Nc$,
$\dot\Nc$,
$\Bc$,
$\Bc^h_w$,
$\Bc^h_{x,w}$,
$\dot\Nc^h_w$,
$\ddot\Nc$,
$\ddot\Nc_{v,w}^h$,
$\ddot\Nc_1$,
$\ddot\Nc_{s_i}$,
$G^{h,\circ}$,
$\tilde G^h$,
$\pi^h$,
$\Bc^h_x$,
$\Delta^+_h$,
$\Den_h$,
$B^h$,
$B_w$,
$W^h$,
$n_h$,
$\Oc_\Yc(\l)$,
$\Oc^h_\Yc(\l)$,
$\widehat K(\ddot\Nc^h)$,
$\widehat H_*(\ddot\Nc^h)$,
$\Phi$,
$ch$,
$\Lambda(\Lc)$,
$e_{cox}$,
$\dot\Nc_{w,(i)},$
$\Bc_{w,(i)}$,

\itemitem{2.5 : }
$\Lc_{h,w}$,
$\Lc_h$,
$L_{h,w,\chi}$,
$L_{h,\chi}$,
$\Xc_{h,w}$,
$\Xc_h$,
$A(h,x)$,
$A(h,x)^+$,
$A(h,x)^\circ$,

\itemitem{2.6 : }
$FS$,
$\chi^\dag$,

\itemitem{2.7 : }
$I^{c,w}$,
$\widehat\Pb_{X,c}$,
$\widehat\Pb_{X,c}^\dag$,

\itemitem{3.1 : }
$w_m$,
$W_m$,
$G_0^m$,
$\gen_0^m$,
$ev$,
$A_0$,
$\aen_0$,
$\aen_\CC$,
$A(m,x_0)$,
$A(m,x_0)^\circ$,
$A_m$,
$A^\circ_m$,
$\nu$,

\itemitem{3.2 : }
$X_w$,
$A_w$,
$\een_{c,w}$,
$\Ec_{c,w}$,

\itemitem{3.3 : }
$\Den_{h,w}$,

\itemitem{3.4 : }
$V_\a$, 
$V^+_\a$, 
$V^-_\a$,
$V_{h,\reg}$,
$\Ac$,
$\Ac_w$,
$\Aen$,
$\Aen^h$,
$\Cen_c$,
$\Hb_\lambdav$,
$W_\lambdav$,
$w_c$,
$\lambdav_c$,
$\Pi_c$,

\itemitem{3.5 : }
$\underline{i}$,

\itemitem{A.1 : }
$\Res_{h,w}$,
$\Ind_{h,w}$,
$\Gamma_{\tilde G^h}^{B^h}$,
$a^h$,
$b^h$,
$d^h$,
$c^h$,
$\Extb_h$,
$\Phi'$,

\itemitem{A.2 : }
$p_v$,
$r_v$,
$W_L$,
${}^L\Bc_v$,
$\Bc_L$,
$\Bc_{L,e}$,
$\dot\Nc_L$,
$\Mc_e$,
${}^L\Mc_{e,v}$,
${}^L\widetilde\Mc_{e,v}$,
$\cen$,
$\qen$,
$\len$,
$\uen$,
$\phi$,
$\phi_\zeta$,
$s_\phi$,
$\gamma_\zeta$,
$N_{h,\chi}$,
$\Hb_L$,
$I_L$.

\vskip3mm

\subhead A.6. Remarks on \cite{V1} (by E.V.)\endsubhead

In this section we'll use the same notations as in \cite{V1}.
One source of technical problems in loc.~cit.~
was the fact that the scheme ${}^w\Xc^k$ in
Section 2.4 may be not separated.
Several people (including A. Bravermann) told us
that, since the group $G$ is of affine type, this can't occur. 
Therefore, in Appendix A (and at each place where it is used in the paper) one
may assume that the scheme $\Xc$ is separated.
The reason for that is the following (details are left to the reader).

First, since the map $\pi_\flat : \Xc\to\Xc^\flat$ factors to
a smooth morphism
${}^w\Xc^{k}\to{}^w\Xc^{\flat,k}$ with separated fibers,
it is enough to prove that the scheme 
$${}^w\Xc^{\flat,k}=U_k^-(\CC)\setminus{}^w\Xc^\flat$$ is separated. 
For each $k>0$ we'll set 
$\uen_k=\eps^k\cdot\uen_\flat$. See loc.~cit., sec.~2.10.
Then the group-$\CC$-scheme $U_k^-$ is such that
$U_k^-(\CC)=\underline G(\eps^{k+1}A).$
See loc.~ cit., sec.~ 2.2.

Now, the $\CC$-scheme $\Xc^\flat$ represents the functor
associated to the set of $\underline G$-torsors over
$\PP^1$ with a trivialisation on the formal neighborhood of $\{0\}$.
Thus ${}^w\Xc^{\flat,k}$ is an open subset of the scheme representing
the functor associated to the set of pairs $(\Pc,\phi)$, where
$\Pc$ is a $\underline G$-torsor over
$\PP^1$ and $\phi$ is a trivialisation of the restriction of
$\Pc$ to the subscheme
$$\xi=\Spec(A/\eps^{k+1}A)\subset\PP^1.$$
Here $\eps$ is 
identified with a local coordinate on $\PP^1$ centered at $0$. 
Fix a faithfull finite dimensional rational module $V$ of $\underline G$.
Consider the vector space $W=V\otimes_\CC\CC[\xi]$. 
Restricting the sections of the vector bundle
$\Pc\times_{\underline G}V$ on $\PP^1$ to the subscheme $\xi$ yields a map
$\psi$ from ${}^w\Xc^{\flat,k}$ to the Grassmannian $\Gr(W)$ of $W$. 
The fibers of $\psi$ embed in a Quot scheme.
So ${}^w\Xc^{\flat,k}$ is separated. 

Further, the following corrections should be made in loc.~ cit. :

\itemitem{p 280 :}
in line 9, replace the definition of $A_e$ by
``$A_e=\{a\in A;\ad_ae=e,\ad_af=\zeta^{-2}f\}$".

\itemitem{p 299 :}
replace the first 2 lines after CLAIM 2 by
``Taking now $\Sc_\ben=\dot\Nc\cap(\ben\times\Bc),\ben$ instead of
$\Sc_\phi,\sen_\phi$ 
the same proof as for Theorem 5.8 (with $e=0,\phi=1$), 
yields an isomorphism of $\underline\Hb_\zeta$-modules
$\Kb^A(\Sc_\ben^a)=\underline\Hb_\zeta\otimes_{\Rb_A}\CC_a$."
Then replace $\Sc_\phi$ by $\Sc_\ben$ at each place until CLAIM 3.

\vskip1mm

\itemitem{p 308 :}
in line 1, replace ``$\ZZ/k\ZZ$" by ``$\ZZ/(r+1)\ZZ$".

\enddocument